%% file: 1105.6007
\documentclass[12pt]{article}
\usepackage{amsmath,amssymb, amsthm, amsopn, amsfonts,epsfig,graphics}{}
\usepackage{hyperref}
\usepackage[french, english]{babel}
\usepackage[applemac]{inputenc}
\usepackage[all]{xy}
\newcommand{\field}[1]{\mathbb{#1}} \newcommand{\rz}{\field{R}}
\newcommand{\cz}{\field{C}} \newcommand{\nz}{\field{N}}
\newcommand{\zz}{\field{Z}} \newcommand{\tz}{\field{T}}
\newcommand{\pz}{\field{P}} \newcommand{\qz}{\field{Q}}
\newcommand{\dz}{\field{D}} \newcommand{\kz}{\field{K}}
\newcommand{\sz}{\field{S}}
\newcommand{\iw}{{\bf i}}
\newcommand{\nw}{{\bf n}}
\newcommand{\tw}{{\bf t}}
\newcommand{\Id}{{\textrm{Id}}}
\newcommand{\ima}{{\textrm{Im\,}}}
\newcommand{\rank}{{\textrm{rank\,}}}
\DeclareMathOperator\Real{\rz e\,}
\DeclareMathOperator\supp{\textrm{supp}}
\newcommand{\Hess}{{\text{Hess\,}}}
\newcommand{\Ker}{{\text{Ker\,}}}
\newcommand{\Ran}{{\text{Ran\,}}}
\newtheorem{theorem}{Theorem}[section]
\newtheorem{lemma}[theorem]{Lemma}

\newtheorem{proposition}[theorem]{Proposition}
\newtheorem{definition}[theorem]{Definition}
\newtheorem{remark}[theorem]{Remark}

\newtheorem{assumption}{Hypothesis}

\def\dim{\textrm{dim\;}}

\title{Precise Arrhenius law for $p$-forms:\\
The Witten Laplacian and\\ Morse-Barannikov complex.}
\author{D.~Le~Peutrec\thanks{DŽpartement de MathŽmatiques,
UMR-CNRS 8628, B‰t. 425, UniversitŽ Paris~11, F91405 Orsay
Cedex. \textit{email: dorian.lepeutrec@math.u-psud.fr}}\; ,
  F.~Nier\thanks{IRMAR, UMR-CNRS 6625, UniversitŽ de Rennes 1, Campus
    de Beaulieu, F35042 Rennes Cedex. \textit{email: francis.nier@univ-rennes1.fr}}\; , 
C.~Viterbo\thanks{CMLS, UMR-CNRS 7640, Ecole Polytechnique, F91128
  Palaiseau Cedex and Eilenberg Chair for Spring 2011 at Columbia
  University, New-York, USA. \textit{email: viterbo@math.polytechnique.fr}}
}
\begin{document}
\input{fig4tex}

\bibliographystyle{plain}
\selectlanguage{english}
\maketitle 
\begin{abstract}
Accurate asymptotic expressions are given for the exponentially small
eigenvalues of Witten Laplacians acting on $p$-forms. The key
ingredient, which replaces explicit formulas for global quasimodes in
the case $p=0$, is Barannikov's presentation of Morse theory in \cite{Bar}.
\end{abstract}
\section{Introduction and main statements}
\subsection{Presentation}
\label{se.present}
The Brownian motion of a particle, at position
$x(t)$ (in $\rz^{d}$ for this rapid presentation), 
and experiencing a gradient field $-2\nabla f(x)$, can be
modelled by the Smoluchowski stochastic differential equation (see
a.e. \cite{Nel}):
\begin{equation}
  \label{eq.sde}
dx=-2\nabla f(x)dt - \sqrt{2h}dW\quad,\quad  x(t)=x_{0}\in \rz^{d}\,.
\end{equation}
Local minima of the energy profile, $f$, are stable steady states when
 $h=0$ 
and become metastable states when $h>0$ and
small. If $h>0$ is thought as a temperature, the  lifetime
such a metastable states $U_{0}$
is exponentially large in term of $1/h$. Its inverse
$\tau_{U_{0}}(h)$ follows an Arrhenius law
$\tau_{U_{0}}(h)\propto e^{-\frac{E_{act}(U_{0})}{h}}$,
 where the
activation energy $E_{act}(U_{0})$ equals $2(f(U_{1})-f(U_{0}))$
and $U_{1}$ is a proper
saddle point, of the energy profile $f$, associated with $U_{0}$\,.
Those inverse lifetime are actually the exponentially small
eigenvalues of the Feller semigroup generator, associated with 
\eqref{eq.sde}, 
\begin{equation}
  \label{eq.feller}
  -2\partial_{x}f.\partial_{x}-h\Delta_{x}=\frac{1}{h}(-h\partial_{x}+2\partial_{x}f).(h\partial_{x})\quad
  \text{on}~\rz^{d},
\end{equation}
defined on $L^{2}(\rz^{d}, e^{-\frac{2f(x)}{h}}dx)$\,, while
$e^{-\frac{2f(x)}{h}}dx$ 
is
the invariant measure associated with \eqref{eq.sde}.\\
The analysis of these activation energies, or exponentially small
eigenvalues in terms of $h$, has motivated
various mathematical studies within the probabilistic approach and
simulated annealing techniques in the 80's (see for instance \cite{HKS}\cite{FrWe}).
More recently several works have been devoted to the accurate
computation of the prefactors, $P_{U_{0}}(h)$ in
$\tau_{U_{0}}(h)=P_{U_{0}}(h)e^{-\frac{E_{act}(U_{0})}{h}}$, with a
probabilistic 
and potential theory approach in \cite{BEGK}\cite{BoGaKl}, or with PDE
and spectral techniques in
\cite{HeNi2}\cite{HKN}\cite{HeNi}\cite{Lep1}\cite{Lep2}\cite{Lep3}\cite{HerHiSj2}.\\
After conjugating with $e^{-\frac{f}{h}}$ and multiplying by $h$, the operator
\eqref{eq.feller} becomes a Witten Laplacian acting of functions
($0$-forms)
\begin{equation}
  \label{eq.wittRd}
  (-h\partial_{x}+\partial_{x}f).(h\partial_{x}+\partial_{x}f)=-h^{2}\Delta +|\nabla_{x}|^{2}-h^{2}\Delta_{x}=d_{f,h}^{*}d_{f,h}=\Delta_{f,h}^{(0)}\,,
\end{equation}
with $d_{f,h}=e^{-\frac{f}{h}}(hd)e^{\frac{f}{h}}=hd+df\wedge$ and
$d_{f,h}^{*}=hd^{*}+\mathbf{i}_{\nabla f}$\,.\\
On a general configuration space, that is a manifold, the Witten
Laplacian, acting on the space of all smooth differential forms
\begin{equation}
  \label{eq.WittenLapl}
  \Delta_{f,h}=(d_{f,h}+d_{f,h}^{*})^{2}=d_{f,h}^{*}d_{f,h}+d_{f,h}d_{f,h}^{*}
\end{equation}
is decomposed as the direct sum
$\Delta_{f,h}=\oplus_{p=0}^{d}\Delta_{f,h}^{(p)}$ with
$\Delta_{f,h}^{(p)}$ acting on $p$-forms. It provides the
geometrically intrisic writing, depending on the metric $g$ and the
Morse
function $f$, and exhibits the relationship with other structures.
In his celebrated article \cite{Wit}, Witten showed that this deformation
of Hodge theory allows one to recover analytically the Morse inequalities for the function $f$. The number of eigenvalues of $\Delta_{f,h}^{(p)}$ lying in $[0,h^{3/2})$
equals, for $h>0$ small enough, the number $m_{p}$ of critical points
of $f$ with index $p$, while conjugating the differential with
$e^{\frac{f}{h}}$ provides an isomorphism in cohomology between the de~Rham chain
complex $(\Omega^{*}(M), d)$ 
and the chain complex $(\Omega^{*}_{f,h},d_{f,h})$, where $\Omega^{p}_{f,h}$ is the space generated by the eigenmodes with eigenvalue less than $h^{3/2}$ for $\Delta^{p}_{f,h}$, and $d_{f,h}$
the Witten differential.\\  
Shortly after Witten's article, it was proved in \cite{HelSj4} that the
$\mathcal{O}(h^{3/2})$ eigenvalues of $\Delta_{f,h}^{(p)}$ are
actually exponentially small,
$\lambda_{k}^{(p)}=\mathcal{O}(e^{-\frac{C_{k}^{(p)}}{h}})$ without
specifying the $C_{k}^{(p)}$'s. The values of the activation energies
$C_{k}^{(0)}$ for $p=0$ were
already known from \cite{FrWe}\cite{HKS}. The
accurate determination, in the case of functions ($p=0$), of the
prefactors, $P_{k}^{(0)}(h)$ in
$\lambda_{k}^{(0)}=P_{k}^{(0)}(h)e^{-\frac{C_{k}^{(0)}}{h}}$, 
  came later, motivated by
probabilistic questions in \cite{BEGK}\cite{BoGaKl}, or by the
analysis of the Kramers-Fokker-Planck operators in \cite{HerNi}.\\
The accurate computation of the small eigenvalues of
$\Delta_{f,h}^{(p)}$ is made difficult  by the interactions due to the  tunneling effect between the
$m_{p}$ quantum wells, with a hierarchy of weakly resonant tunneling
quantum wells, according to the terminology of
\cite{HelSj2}\cite{HeSj3}.
 This hierarchy, which orders the exponentially small quantities, can be
solved by considering the interaction via the deformed differential
$d_{f,h}$ with the eigenmodes of $\Delta_{f,h}^{(p+1)}$ and via the
deformed codifferential $d_{f,h}^{*}$ with the eigenmodes of
$\Delta_{f,h}^{(p-1)}$\,.
When $p=0$, it is simply understood within the probabilistic approach 
by ordering the exit times, following in some sense the intuition of
Arrhenius law. It actually amounts to elementary topological arguments
 by considering how the number of connected components of the sublevel set
$f^{\lambda}=\left\{x, f(x)<\lambda\right\}$ varies as $\lambda$
crosses a critical value. Moreover, the analysis carried out in
\cite{HeNi2}\cite{HKN}\cite{HeNi}\cite{Lep1}\cite{Lep2}\cite{Lep3},
relied on the important remark that the eigenvalues of
$\Delta_{f,h}^{(0)}$ are the squares of the singular values of the
differential $d_{f,h}^{(0)}$\,. The Fan inequalities
(see \cite{Sim}) for singular
values then allow to propagate \underline{relative} errors (i.e. small
errors relative to various exponentially small quantities).\\
For all these reasons, the study of exponentially small eigenvalues of
$\Delta_{f,h}^{(p)}$ for a general $p$, is a natural question which is also encountered
 in geometry (see for instance
\cite{Zha}\cite{BiLe}\cite{Bis2} and references therein)
 or in statistical physics (see \cite{TTK}).
A first attempt was done in \cite{Lep4}, extending the result for
$p=0$ to $p=0,1,2$ on surfaces, with simple duality and chain complex
arguments.
For a general $p$, the global quasimodes of the form
$\chi_{U_{0}}(x)e^{-\frac{f(x)-f(U_{0})}{h}}$ used for the case $p=0$
in \cite{HeNi2}\cite{HKN}\cite{HeNi}\cite{Lep1}\cite{Lep2}\cite{Lep3}
and which propagate the information through weakly resonant quantum
wells, are missing. 

The solution comes from the use of Barannikov's  version of the Morse complex
from  \cite{Bar} 
 which fits exactly the handling of
global quasimodes for Witten Laplacian. There are two reasons for
this: 

1) this new chain complex has nice restriction properties which
are implemented by boundary Witten Laplacians; 

2) a side result coming
from this presentation of Morse theory allows to replace the
analytical computations with
$\chi_{U_{0}}(x)e^{-\frac{f(x)-f(U_{0})}{h}}$ in the case $p=0$, 
 by a subtle repeated use of
Stokes' formula.\\

We conclude this introduction by emphasizing that the accurate
analysis of the tunnel effect required for the computation of
exponentially small eigenvalues, goes far beyond the instantonic
picture (see \cite{Bot2}), which sticks in some sense to the intuition
of classical mechanics. However, it is remarkable that discriminating
between so
small quantities (exponentially small quantities $e^{-C_{k}/h}$ as $h\to 0$) is
made possible by global topological arguments.

\subsection{Assumptions and result}
\label{se.assresult}
\begin{assumption}
\label{hyp.1}
  We shall work on an oriented compact riemannian manifold $(M,g)$ and
  $f$ will be an excellent Morse function : $f$ is smooth has  non degenerate critical points and these have
  distinct critical values. Moreover, homology and cohomology will always be with real coefficients. 
\end{assumption}
Barannikov's simple Morse complex, allows to partition the
set of critical points, $\mathcal{U}=\left\{x\in M, \nabla
  f(x)=0\right\}$ (resp. the set of critical points with index $p$, 
$\mathcal{U}^{(p)}=\left\{x\in M, \nabla f(x)=
  0\,,~\textrm{sign(Hess f)}(x)=(d-p,p)\right\}$), into upper, lower
and homological critical points:
\begin{eqnarray}
  \label{eq.partU}
  &&\mathcal{U}=\mathcal{U}_{U}\sqcup
  \mathcal{U}_{L}\sqcup\mathcal{U}_{H}\\
\text{resp.}&&
\label{eq.partUp}
\mathcal{U}^{(p)}=\mathcal{U}^{(p)}_{U}\sqcup
  \mathcal{U}^{(p)}_{L}\sqcup\mathcal{U}^{(p)}_{H}\,.
\end{eqnarray}
Homological critical points in $\mathcal{U}_{H}^{(p)}$ are associated
with the kernel $\ker(\Delta_{f,h}^{(p)})\sim \ker
\Delta_{\text{Hodge}}^{(p)}$ and their number is the $p$-th Betti number
$\beta_{p}=\dim H^{p}(M, {\mathbb R})$\,. 
The boundary operator $\partial_{\mathcal B}$ of Barannikov's chain complex,
defined on $\oplus_{U\in \mathcal{U}}\rz U$,
associates with any $U'\in \mathcal{U}_{U}^{(p)}$ an element $U\in
\mathcal{U}_{L}^{(p-1)}$  such that $f(U)<f(U')$, 
 and vanishes on all other critical points
$U'\in \mathcal{U}_{H}\cup \mathcal{U}_{L}$\,. Details are given in
Section~\ref{se.barann}.
The second assumption avoids technical (nevertheless
interesting) questions about multiplicities of non zero
 exponentially small
eigenvalues.
\begin{assumption}
\label{hyp.2}
The values $f(U')-f(U)$ obtained for $U'\in \mathcal{U}_{U}$ and
$\partial_{B}U'=U$ are all distinct.
\end{assumption}
Here is our main result
\begin{theorem}
\label{th.main}
Assume Hypotheses~\ref{hyp.1} and \ref{hyp.2} hold.
Let $\mathcal U_{H} $,
$\mathcal U_{L}$,  $\mathcal U_{U}$ respectively denote 
the sets of homological,
lower and upper critical points.\\
For $h_{0}$ small enough
and $0<h<h_{0}$, 
there exists a mapping $j$ from $\mathcal U
:=\mathcal U_{H} \cup \mathcal U_{L}\cup \mathcal U_{U}  $ onto 
$\sigma(\Delta_{f,h})\cap [0, h^{\frac32}[$
and the restriction $j_{p}:=j\big|_{\mathcal U^{(p)}}$ 
is onto $\sigma(\Delta^{(p)}_{f,h})\cap [0, h^{\frac32}[$ and one to
one provided the eigenvalues of $\sigma(\Delta^{(p)}_{f,h})$ are counted
with multiplicities.\\
Moreover, the map $j$ satisfies the following properties:
\begin{enumerate} 
\item 
For $U^{(p)}$ in $\mathcal U^{(p)}_{H}$,
$$j(U^{(p)})=0\,.$$
\item For $U^{(p)}$ in $ \mathcal U_{L}^{(p)}$,
let $U^{(p+1)}$ denote 
the element of  $\mathcal U_{U}^{(p+1)}$ s.t. $\partial_{\mathcal B}(U^{(p+1)})=U^{(p)}$. Then,
there exists a homological constant $\kappa(U^{p+1})\in\mathbb R^*$
such that 
\begin{multline*}
j(U^{(p)})
=\kappa^2(U^{(p+1)})\frac{h}{\pi}
\frac{|\lambda_{1}^{(p+1)}\cdots\lambda_{p+1}^{(p+1)}|}{
|\lambda_{1}^{(p)}\cdots\lambda_{p}^{(p)}|}
\frac{{|\Hess f(U^{(p)})|}^\frac12}
{{|\Hess f(U^{(p+1)})|}^\frac12}
\\
\times
e^{-2\frac{f(U^{(p+1)})-f(U^{(p)})}{h}}(1+\mathcal{O}(h))\,,
\end{multline*}
where $\lambda_{1}^{(\ell)},\dots,\lambda_{\ell}^{(\ell)}$
denote the negative eigenvalues of $\Hess f(U^{(\ell)})$,
for $\ell\in\{p,p+1\}$.\\
\item Finally , for $U^{(p)}$ in $\mathcal U_{U}^{(p)}$, the equality
$
j(U^{(p)})=j(\partial_{\mathcal B}(U^{(p)}))=j(U^{(p-1)})$ holds with
$U^{(p-1)}:=\partial_{\mathcal B}(U^{(p)})\in \mathcal U_{L}^{(p-1)}$,
i.e.
\begin{multline*}
j(U^{(p)})
=\kappa^2(U^{(p)})\frac{h}{\pi}
\frac{|\lambda_{1}^{(p)}\cdots\lambda_{p}^{(p)}|}{
|\lambda_{1}^{(p-1)}\cdots\lambda_{p-1}^{(p-1)}|}
\frac{{|\Hess f(U^{(p-1)})|}^\frac12}
{{|\Hess f(U^{(p)})|}^\frac12}
\\
\times
e^{-2\frac{f(U^{(p})-f(U^{(p-1)})}{h}}(1+\mathcal{O}(h))\,,
\end{multline*}
where $\lambda_{1}^{(\ell)},\dots,\lambda_{\ell}^{(\ell)}$
denote the negative eigenvalues of $\Hess f(U^{(\ell)})$,
for $\ell\in\{p-1,p\}$.
\end{enumerate} \end{theorem}
A relative version of this result, implemented with boundary Witten Laplacians, is
given  in Subsection~\ref{se.relatth} at the end of this paper.
\begin{remark}
Although we are not able to prove it in general, there is a strong
indication that the ``homological'' constant $\kappa_{p}(U_{0})$
equals $\pm 1$\,.
This is indeed the case for $p=0$ as shown in
\cite{HKN}\cite{HeNi}\cite{Lep3}.
By duality it is also true when $p=d$\,. Finally in the case of
surfaces treated in \cite{Lep4}, a combination of these results says
that it is true for $p=0,1,2$\,.\\
A general proof requires a better understanding of the topological
aspects of Morse theory and of Barannikov's construction.\\
Surely, this constant is completely
determined by the structure of the homology groups of the sublevel
sets $H_{*}(\left\{f< \lambda\right\})$, $\lambda\in
[-\infty,+\infty]$\,.
 It does not depend on $h$, on the
riemannian metric $g$ or on the Morse function $f$ (as long our generic
assumptions are fulfilled), contrary to the other factors. This is a
reason to use the attribute ``homological'' for this, up to now
unknown, constant.
\end{remark}

\section{Barannikov's simple complex and Morse theory}
\label{se.barann}

In this section, we adapt the approach of Barannikov, using notations and
definitions better suited to the treatment of Witten
Laplacians.

\subsection{Sublevel sets and bases of Morse theory}
\label{se.sublevel}
Remember that we work on a riemannian compact oriented manifold
$(M,g)$, endowed with an excellent Morse function according to
Hypothesis~\ref{hyp.1}. With such an assumption we may identify  a critical 
point $U$  with the corresponding critical value
$c=f(U)$\,. For any index $p$, $0\leq p\leq d=\dim M$, the set 
of critical points of $f$
with index $p$ is $\mathcal{U}^{(p)}=\left\{U_{k}^{(p)}, 1\leq k\leq m_{p}\right\}$\,.
It is equivalently represented  by 
a vertical line
with a $m_{p}$ points with heights
$c_{k}^{(p)}=f(U^{(p)}_{k})$, $1\leq k\leq m_{p}$\,. 
The vector space spanned by these points is denoted by
$\mathcal{C}^{(p)}(f)$ and 
we set $\mathcal{C}(f)=\oplus_{p=0}^{d}\mathcal{C}_{p}(f)$\,. We shall construct explicitly
a  differential on $\mathcal{C}(f)$ quasi-isomorphic to the
Morse chain complex associated with $f$.\\

For $\lambda\in [-\infty, +\infty]$,  $f^{\lambda}$  denotes the
sublevel set $\left\{x\in M, f(x)<\lambda\right\}$ while $f_{\lambda}$
denotes the upper level set $\left\{x\in M, f(x)>\lambda\right\}$ and
$f_{\lambda}^{\mu}=\left\{x\in M, \lambda < x < \mu\right\}$\,.
Let us recall a few elementary facts related with Morse theory known
from \cite{Mil}\cite{Bott}\cite{Lau1}\cite{Lau2}. We refer to
\cite{Hat}\cite{Mas}\cite{BoTu}\cite{Ful} for basic material in homological algebra.
When $\lambda\in (\min f, \max f)$ 
is not a critical value, $f^{\lambda}$ and $f_{\lambda}$ are boundary
manifolds, while $f^{\lambda}=M$, $f_{\lambda}=\emptyset$ for
$\lambda> \max f$ and $f_{\lambda}=\emptyset$, $f_{\lambda}=M$ for
$\lambda < \min f$\,.
When there is no critical value between $\lambda_{1}$ and
$\lambda_{2}$\,,the natural inclusion of $f^{\lambda_{1}}$ in $f^{\lambda_{2}}$
(resp. $f_{\lambda_{2}}$ in $f_{\lambda_{1}}$) induces a
homotopy equivalence and therefore induces an isomorphism of their homology groups:
$$
\forall p\in \left\{0,\ldots, d\right\}\,,\quad
H_{p}(f^{\lambda_{1}})=H_{p}(f^{\lambda_{2}})\quad
\text{and}\quad
H_{p}(f_{\lambda_{1}})=H_{p}(f_{\lambda_{2}})\,.
$$
With the help of the five lemma,
this holds also for the relative homology groups $H_{*}(f^{\mu}, f^{\lambda_{1}})$ and $H_{*}(f^{\mu}, f^{\lambda_{2}})$, and $H_{*}( f_{\lambda_{1}}, f_{\mu})$ $H_{*}( f_{\lambda_{2}}, f_{\mu})$ for $\mu > \lambda_{1}, \lambda_{2}$, as is easily seen using the long
exact sequences 
\begin{eqnarray*}
&& 
\xymatrix{
H_{*+1}(f^{\lambda})
\ar[r]^{i_{*}}&
H_{*+1}(f^{\mu})
\ar[r]^{j_{*}}&
H_{*+1}(f^{\mu},f^{\lambda})
\ar[r]^{\partial}&
H_{*}(f^{\lambda})
}
\;,\\
&&
\xymatrix{
H_{*+1}(f_{\mu})\ar[r]^{i_{*}}&
H_{*+1}(f_{\lambda})
\ar[r]^{j_{*}}&
H_{*+1}(f_{\lambda},f_{\mu})
\ar[r]^{\partial}&
H_{*}(f_{\mu})}\,.
\end{eqnarray*}
when $\mu>\lambda$ are not critical values.\\
Passing a critical point with index $p$  and critical value $c$, the pair
$(f^{c+\varepsilon},f^{c-\varepsilon})$ is homologous to the pair
$(\dz^{p},\partial \dz^{p})$, associated with the $p$-cell
$e^{p}=\dz^{p}\setminus \partial \dz^{p}$ (see \cite{Mil}).
This gives using excision
\begin{equation}
  \label{eq.exactlong}
\xymatrix@C=0.7pc{
0\ar[r]&
 H_{p}(f^{c-\varepsilon})
\ar[r]&
H_{p}(f^{c+\varepsilon})
\ar[r]&
H_{p}(f^{c+\varepsilon},f^{c-\varepsilon})
\ar[r] &
H_{p-1}(f^{c-\varepsilon})
\ar[r]&
H_{p-1}(f^{c+\varepsilon})
\ar[r]&
0}\,,
\end{equation}
with $\dim H_{p}(f^{c+\varepsilon},f^{c-\varepsilon})=1$
and
ensures that for $k\neq p, p-1$ we have the  equality of  $H_{k}(f^{c\pm \varepsilon})$\,.\\ 
This yields  two mutually exclusive cases
\begin{eqnarray}
  \label{eq.1stcasesub}
&&
\left\{
  \begin{array}[c]{l}
\xymatrix@C=1pc{
H_{p}(f^{c-\varepsilon})\ar[r]^{\sim}
&H_{p}(f^{c+\varepsilon})
}\quad\text{and}
\\
\xymatrix@C=1pc{
0\ar[r]&
 H_{p}(f^{c+\varepsilon},f^{c-\varepsilon})
\ar[r]&
 H_{p-1}(f^{c-\varepsilon})
\ar[r]&
 H_{p-1}(f^{c+\varepsilon})\to 0
}\,,
\end{array}
\right.
\\
\label{eq.2ndcasesub}
\text{or}
&&
\left\{
\begin{array}[c]{l}
\xymatrix@C=1pc{
0\ar[r]&
H_{p}(f^{c-\varepsilon})
\ar[r]&
 H_{p}(f^{c+\varepsilon})\to
H_{p}(f^{c+\varepsilon},f^{c-\varepsilon})
\ar[r]& 0
}
\\
\text{and}\quad
 \xymatrix{
H_{p-1}(f^{c-\varepsilon})\ar[r]^{\sim}&
 H_{p-1}(f^{c+\varepsilon})
}\,.
\end{array}
\right.
\end{eqnarray}
The Poincar{\'e} duality takes a nice form owing to
Theorem~3.43~\cite{Hat}  ($M$ is oriented), with the excision argument
$H_{*}(f^{\mu},f^{\lambda})=H_{*}( \overline{f_{\lambda}^{\mu}},  \left\{f=\lambda\right\})$: For two non critical values $-\infty\leq
\lambda < \mu\leq +\infty$, the cohomology group
$H^{k}(f^{\mu},f^{\lambda})$ is isomorphic to
$H_{d-k}(f_{\lambda},f_{\mu})$\,.
In \cite{Spa}-p296, this is called Alexander duality and proved,
without excision, via coverings and Mayer-Vietoris techniques.
 With $f_{\lambda}=(-f)^{-\lambda}$,
this is often summarized by changing $f$ into $-f$ and inverting indexes $p$
and $d-p$\,. Thus, the dual version of \eqref{eq.exactlong} for a critical
value with index $p$ is
\begin{equation*}
\xymatrix@C=0.7pc{
0
&
\ar[l]
H_{d-p-1}(f_{c-\varepsilon})
&\ar[l]
H_{d-p-1}(f_{c+\varepsilon})
&\ar[l]
H_{d-p}(f_{c-\varepsilon},f_{c+\varepsilon})
&\ar[l]
H_{d-p}(f_{c-\varepsilon})
\\
&&&&\ar[u]
H_{d-p}(f_{c+\varepsilon})
\\
&&&&\ar[u]
0}\,.
\end{equation*}
For this, use the excision property
$$
H_{k}(f^{c+\varepsilon},f^{c-\varepsilon})=H_{k}(
\left\{c-\varepsilon\leq  f < c+\varepsilon\right\},
\left\{f=c-\varepsilon\right\})
$$
 while noticing that according to PoincarŽ duality
$H_{k}(\left\{f=c-\varepsilon\right\})$ is isomorphic to
$H^{d-1-k}(\left\{f=c-\varepsilon\right\})$\,.\\
Hence passing a critical value with upper level sets leads to the two
exclusive cases
\begin{eqnarray}
  \label{eq.1stcaseup}
&&
\left\{
  \begin{array}[c]{l}
\xymatrix@C=1pc{
H_{d-p}(f_{c+\varepsilon})\ar[r]^{\sim}
&H_{d-p}(f_{c-\varepsilon})
}\quad\text{and}
\\
\xymatrix@C=1pc{
0\ar[r]&
 H_{d-p}(f_{c-\varepsilon},f_{c+\varepsilon})
\ar[r]&
 H_{d-p-1}(f_{c+\varepsilon})
\ar[r]&
 H_{d-p-1}(f_{c-\varepsilon})\to 0
}\,,
\end{array}
\right.
\\
\label{eq.2ndcaseup}
\text{or}
&&
\left\{
\begin{array}[c]{l}
\xymatrix@C=1pc{
0\ar[r]&
H_{d-p}(f_{c+\varepsilon})
\ar[r]&
 H_{d-p}(f_{c-\varepsilon})\to
H_{p}(f_{c-\varepsilon},f_{c+\varepsilon})
\ar[r]& 0
}
\\
\text{and}\quad
 \xymatrix@C=1pc{
H_{d-p-1}(f_{c+\varepsilon})\ar[r]^{\sim}&
 H_{d-p-1}(f_{c+\varepsilon})
}\,.
\end{array}
\right.
\end{eqnarray}

\subsection{Classification of critical points}
\label{se.class}

\subsubsection{Partition}
\label{se.partcrit}

The critical points are divided in three classes, and we  prove that these classes 
make a partition of the set of critical points, satisfying a number of  additional properties.
\begin{definition}
\label{de.class}
  \begin{enumerate}
  \item 
A critical value  (resp. point) $c$  of $f$ is called a lower critical
value (resp. point), if
  the natural mapping
$$
\xymatrix{
H_{*}(f^{c+\varepsilon},f^{c-\varepsilon})\ar[r]&H_{*}(M,f^{c-\varepsilon})
}
$$
vanishes.
\item 
A critical value  (resp. point) $c$ of $f$ is called an upper critical
value (resp. point), if
  the natural mapping
$$
\xymatrix{
H_{*}(f^{c+\varepsilon})\ar[r]& H_{*}(f^{c+\varepsilon},f^{c-\varepsilon})
}
$$
vanishes.
\item In all other cases the critical value (resp. point) $c$, is called an
  homological critical value (resp. point).
\end{enumerate}
\end{definition}
Remember the long exact sequence for the triple $(X,A,B)$ where $B\subset  A \subset X$:
$$
\xymatrix{ &\ar[r]&
H_{*}(A,B)
\ar[r]&
H_{*}(X,B)
\ar[r]&
H_{*}(X,A)
\ar[r]^{\partial}
&
H_{*-1}(A,B)\ar[r] & \\
}
$$
and the commutative diagram associated with a map $\varphi:X\to X'$ satisfying
$\varphi(A)\subset A'$ and $\varphi(B)\subset B'$:
\begin{equation}
  \label{eq.comdiag}
\hspace{-0.7cm}\xymatrix{&\ar[r]&
H_{*}(A,B)
\ar[r]\ar[d]^{\varphi_{*}}&
H_{*}(X,B)
\ar[r]
\ar[d]^{\varphi_{*}}
&
H_{*}(X,A)
\ar[r]^{\partial}
\ar[d]^{\overline{\varphi_{*}}}
&
H_{*-1}(A,B)
\ar[d]^{\varphi_{*}}\ar[r]& \\&\ar[r]&
H_{*}(A',B')
\ar[r]&
H_{*}(X',B')
\ar[r]
&
H_{*}(X',A')
\ar[r]^{\partial}
&
H_{*-1}(A',B')\ar[r]& 
}
\end{equation}
\begin{proposition}
\label{pr.lowup}
The set of lower critical values (resp. points) and upper critical
values (resp. points) are disjoint and the classification into lower,
upper and homological critical values (resp. points) is a partition.
\end{proposition}
\begin{proof}
Consider the long exact sequences correspionding to the triples
$$
(X,A,B) =(f^{c+\varepsilon},
f^{c-\varepsilon},\emptyset,)\quad \text{and}\quad(X',B',A')=(M,
f^{c-\varepsilon},\emptyset,)
$$
 with the mapping
$i^{\infty,c+\varepsilon}:f^{c+\varepsilon}\to M=f^{+\infty}$~:
$$
\xymatrix{
H_{*}(f^{c-\varepsilon})
\ar[r]\ar[d]^{\Id}&
H_{*}(f^{c+\varepsilon})
\ar[r]
\ar[d]^{i^{\infty,c+\varepsilon}_{*}}
&
H_{*}(f^{c+\varepsilon},f^{c-\varepsilon})
\ar[r]^{\partial}
\ar[d]^{\overline{i^{\infty,c+\varepsilon}_{*}}}
&
H_{*-1}(f^{c-\varepsilon})
\ar[d]^{\Id}\\
H_{*}(f^{c-\varepsilon})
\ar[r]&
H_{*}(M)
\ar[r]
&
H_{*}(M,f^{c-\varepsilon})
\ar[r]^{\partial'}
&
H_{*-1}(f^{c-\varepsilon})
}\,.
$$
Now, assume that both the mappings
  $H_{*}(f^{c+\varepsilon},f^{c-\varepsilon})\to
  H_{*}(M,f^{c-\varepsilon})$ and  $H_{*}(f^{c+\varepsilon})\to
  H_{*}(f^{c+\varepsilon},f^{c-\varepsilon})$ both vanish. This
  implies that $\partial$ is injective while
  $\overline{i^{\infty,c+\varepsilon}}=0$\,. 
This contradicts the commutativity $\partial=\partial'\circ \overline{i^{\infty,c+\varepsilon}}$\,.
\end{proof}
\subsubsection{Upper critical points}
\label{se.upper}
We first give other characterizations for upper critical points with
index $p$, which
will be used further. An additional property about the rank of
$i^{\lambda}_{*}:H_{*}(f^{\lambda})\to H_{*}(M)$ is given.
\begin{proposition}
\label{pr.upper}
A critical value $c$ with index $p$ is an upper critical value, iff one
of the following conditions is satisfied:
\begin{enumerate}
\item The critical value satisfies the condition
  \eqref{eq.1stcasesub}, namely:
$$
\left\{
  \begin{array}[c]{l}
\xymatrix@C=1pc{
H_{p}(f^{c-\varepsilon})\ar[r]^{\sim}
&H_{p}(f^{c+\varepsilon})
}\quad\text{and}
\\
\xymatrix@C=1pc{
0\ar[r]&
 H_{p}(f^{c+\varepsilon},f^{c-\varepsilon})
\ar[r]&
 H_{p-1}(f^{c-\varepsilon})
\ar[r]&
 H_{p-1}(f^{c+\varepsilon})\to 0
}
\,.
\end{array}
\right.
$$
\item  The mapping $\partial: H_{*+1}(f^{c+\varepsilon},f^{c-\varepsilon})\to
  H_{*}(f^{c-\varepsilon})$ is one to one.
\item There exists $\lambda\in ]-\infty, c)$ such that the mapping\\
$H_{*}(f^{c+\varepsilon},f^{\lambda})\to
H_{*}(f^{c+\varepsilon},
f^{c-\varepsilon})$ vanishes.
\item The exists $\lambda \in ]-\infty,c)$ such that the mapping\\
$\partial: H_{*+1}(f^{c+\varepsilon},f^{c-\varepsilon})\to
  H_{*}(f^{c+\varepsilon},f^{\lambda})$ is one to one. 
\end{enumerate}
\end{proposition}
\begin{proof}
 The condition 1) is just the explicit form of the definition of an
 upper critical value, in view of the long exact sequence
 \eqref{eq.exactlong}.\\
Condition 2) is obtained by considering  the right-hand side of the long exact sequence
$$
\xymatrix{
H_{*+1}(f^{c-\varepsilon})
\ar[r]&
H_{*+1}(f^{c+\varepsilon})
\ar[r]&
H_{*+1}(f^{c+\varepsilon},f^{c-\varepsilon})
\ar[r]^{\partial}&
H_{*}(f^{c-\varepsilon})\,.
}
$$
Similarly  condition 4) is equivalent to condition 3). \\
Consider condition 3): It is necessary (take $\lambda=-\infty$).
The middle square of the commutative diagram
\eqref{eq.comdiag}
with
the embedding $\varphi=i^{\lambda}:(X,A,B)=(f^{c+\varepsilon},f^{c-\varepsilon},f^{-\infty}=\emptyset)\to
(X',A',B')=(f^{c+\varepsilon},f^{c-\varepsilon},f^{\lambda})$:
$$
\xymatrix{
H_{*}(f^{c-\varepsilon})
\ar[r]\ar[d]^{i^{\lambda}_{*}}&
H_{*}(f^{c+\varepsilon})
\ar[r]
\ar[d]^{i^{\lambda}_{*}}
&
H_{*}(f^{c+\varepsilon},f^{c-\varepsilon})
\ar[r]^{\partial}
\ar[d]^{\Id}
&
H_{*-1}(f^{c-\varepsilon})
\ar[d]^{i^{\lambda}_{*}}\\
H_{*}(f^{c-\varepsilon},f^{\lambda})
\ar[r]&
H_{*}(f^{c+\varepsilon},f^{\lambda})
\ar[r]^{0}
&
H_{*}(f^{c+\varepsilon},f^{c-\varepsilon})
\ar[r]^{\partial}
&
H_{*-1}(f^{c-\varepsilon},f^{\lambda})\,,
}
$$
provides the sufficiency.
\end{proof}
 Consider for $-\infty\leq \lambda<\mu\leq +\infty$ 
which are not critical values,
 the embeddings $i^{\lambda}:f^{\lambda}\to M$ and $i^{\mu}:f^{\mu}\to
 M$ and $i^{\mu,\lambda}:(f^{\lambda},M)\to (f^{\mu},M)$\,. 
Then the commutative diagram
$$
\xymatrix{
H_{*}(f^{\lambda})
\ar[r]^{i^{\lambda}_{*}}\ar[d]^{i^{\mu,\lambda}_{*}}&
H_{*}(M)
\ar[r]
\ar[d]^{\Id}
&
H_{*}(M,f^{\lambda})
\ar[r]^{\partial}
\ar[d]^{\overline{i^{\mu,\lambda}_{*}}}
&
H_{*-1}(f^{\lambda})
\ar[d]^{i^{\mu,\lambda}_{*}}\\
H_{*}(f^{\mu})
\ar[r]^{i^{\mu}_{*}}&
H_{*}(M)
\ar[r]
&
H_{*}(M,f^{\mu})
\ar[r]^{\partial}
&
H_{*-1}(f^{\mu})
}
$$
implies
\begin{eqnarray}
  \label{eq.comp}
&&  i^{\lambda}_{*}=i^{\mu}_{*}\circ i^{\mu,\lambda}_{*}\quad
(\lambda<\mu)\\
\label{eq.range}
&&\ima i^{\lambda}_{*}\subset \ima i^{\mu}_{*} \quad,\quad \rank
i^{\lambda}_{*}\leq \rank i^{\mu}_{*}\,.
\end{eqnarray}
\begin{proposition}
\label{pr.rkiup}
  When $c$ is an upper critical value, the ranges of
  $i^{c+\varepsilon}_{*}:H_{*}(f^{c+\varepsilon})\to H_{*}(M)$ and
$i^{c-\varepsilon}_{*}:H_{*}(f^{c-\varepsilon})\to H_{*}(M)$ are the same.
\end{proposition}
\begin{proof}
  The condition 1) of Proposition~\ref{pr.upper} ensures that for any
  $k\in \left\{0,\ldots,d\right\}$, the mapping
$$
i^{c+\varepsilon,c-\varepsilon}_{*}:H_{k}(f^{c-\varepsilon})\to H_{k}(f^{c+\varepsilon})  
$$
is onto. Therefore $i^{c+\varepsilon}_{*}$ and
$i^{c-\varepsilon}_{*}=i^{c+\varepsilon}\circ
i^{c+\varepsilon,c-\varepsilon}_{*}$ have the same range.
\end{proof}

\subsubsection{Lower critical points}
\label{se.lower}
With $H_{*}(M,f^{c-\varepsilon})=H_{*}(f^{+\infty},f^{c+\varepsilon})$
and the duality $H^{*}(f^{\mu},f^{\lambda})\simeq
H_{d-*}(f_{\lambda},f_{\mu})$, says that the lower critical values are
the one for which the mapping
$$
\xymatrix{
H_{*}(f_{c-\varepsilon})=H_{*}(f_{c-\varepsilon},f_{+\infty})
\ar[r]&
H_{*}(f_{c-\varepsilon},f_{c+\varepsilon})}
$$
vanishes. It is therefore the dual notion to the one of upper critical
values and all the dual properties of the ones of the upper critical
values will hold for lower critical points. We shall use the following
 characterization.
\begin{proposition}
\label{pr.lower}
A critical value $c$ with index $p$ is a lower critical value, iff one
of the following conditions is satisfied:
\begin{enumerate}
\item The critical value satisfies the condition \eqref{eq.1stcaseup},
  namely:
$$
\left\{
  \begin{array}[c]{l}
\xymatrix@C=1pc{
H_{d-p}(f_{c+\varepsilon})\ar[r]^{\sim}
&H_{d-p}(f_{c-\varepsilon})
}\quad\text{and}
\\
\xymatrix@C=1pc{
0\ar[r]&
 H_{d-p}(f_{c-\varepsilon},f_{c+\varepsilon})
\ar[r]&
 H_{d-p-1}(f_{c+\varepsilon})
\ar[r]&
 H_{d-p-1}(f_{c-\varepsilon})\to 0
}\,.
\end{array}
\right.
$$
\item The mapping $\partial:H_{*+1}(M,f^{c+\varepsilon})\to
  H_{*}(f^{c+\varepsilon},f^{c-\varepsilon})$ is onto.
\item There exists $\lambda\in (c,+\infty]$ such that the mapping 
$H_{*}(f^{c+\varepsilon},f^{c-\varepsilon})\to
H_{*}(f^{\lambda},f^{c-\varepsilon})$ vanishes.
\item There exists $\lambda\in (c,+\infty]$ such that the mapping 
$\partial:H_{*+1}(f^{\lambda},f^{c+\varepsilon})\to
  H_{*}(f^{c+\varepsilon},f^{c-\varepsilon})$ is onto.
\end{enumerate}
\end{proposition}
\begin{proof}
The condition 1) is the dual statement of \eqref{eq.1stcasesub} which
is equivalent to the dual notion of upper critical values.\\
The equivalence with the condition 2) is contained  in the long exact sequence
$$
\xymatrix{
H_{*+1}(M,f^{c-\varepsilon})
\ar[r]&
H_{*+1}(M,f^{c+\varepsilon})
\ar[r]^{\partial}&
H_{*}(f^{c+\varepsilon},f^{c-\varepsilon})
\ar[r]&
H_{*}(M,f^{c-\varepsilon})
}\,.
$$
Similarly the condition 3) and 4) are equivalent.\\
Consider the condition 3): It is necessary (take
$\lambda=+\infty$). If there exists $\lambda\in (c,+\infty]$ such that
 the composition of the embeddings
 $i^{\lambda,c+\varepsilon}:(f^{c+\varepsilon},f^{c-\varepsilon})\to
 (f^{\lambda},f^{c-\varepsilon})$
and $i^{\lambda}:(f^{\lambda},f^{c-\varepsilon})\to
(M=f^{+\infty},f^{c-\varepsilon})$ implies
that the composed map 
$$
\xymatrix@C=3pc{
H_{*}(f^{c+\varepsilon},f^{c-\varepsilon})
\ar[r]^{i^{\lambda,c+\varepsilon}_{*}=0}&
H_{*}(f^{\lambda},f^{c-\varepsilon})
\ar[r]^{i^{\lambda}_{*}}&
H_{*}(M,f^{c-\varepsilon})
}\,,
$$
which equals $i^{c+\varepsilon}_{*}$, vanishes.
\end{proof}

We next prove that the lower critical values share the same property
as the upper critical values, concerning the rank of
$i^{\lambda}_{*}: H_{*}(f^{\lambda})\to H_{*}(M)$\,.

\begin{proposition}
\label{pr.rkilow}
When $c$ is a lower critical value, the ranges of
$i^{c+\varepsilon}_{*}:H_{*}(f^{c+\varepsilon})\to H_{*}(M)$ and 
$i^{c-\varepsilon}_{*}:H_{*}(f^{c-\varepsilon})\to H_{*}(M)$ are the same.
\end{proposition}
\begin{proof}
Assume that $c$ is a critical value with index $p$\,. Then for any
$k\neq p$ the map
$i^{c+\varepsilon,c-\varepsilon}_{*}:H_{k}(f^{c-\varepsilon})\to
H_{k}(f^{c+\varepsilon})$ is onto and the range of
$i^{c-\varepsilon}_{*}=i^{c+\varepsilon}_{*}\circ i^{c+\varepsilon,c-\varepsilon}_{*}$
and $i^{c+\varepsilon}_{*}$, when restricted to $H_{k}$, are equal.
For the case $k=p$, we start from the exact sequence
$$
\xymatrix@C=2pc{
0\ar[r]&
H_{p}(f^{c-\varepsilon})
\ar[r]^{i^{c+\varepsilon,c-\varepsilon}_{*}}&
H_{p}(f^{c+\varepsilon})
\ar[r]&
H_{p}(f^{c+\varepsilon},f^{c-\varepsilon})
\ar[r]^{\partial}&
H_{p-1}(f^{c-\varepsilon})
\ar[r]&
0
}\,,
$$
where the $\partial$-arrow vanishes, because $c$ cannot be an 
upper critical value. Hence the range of
$i^{c+\varepsilon,c-\varepsilon}_{*}$ is an hyperplane of
$H_{p}(f^{c+\varepsilon})$
 and the equality $i^{c-\varepsilon}_{*}=i^{c+\varepsilon}_{*}\circ
 i^{c+\varepsilon,c-\varepsilon}_{*}$ implies that
 $i^{c+\varepsilon}_{*}$ and $i^{c-\varepsilon}_{*}$ have the same
 range iff
$$
1+\dim \ker(i^{c+\varepsilon}_{*}\big|_{\ima
  i^{c+\varepsilon,c-\varepsilon}_{*}})
-\dim \ker(i^{c+\varepsilon}_{*})=0\,.
$$
We write for simplicity
$$
i=i^{c+\varepsilon,c-\varepsilon}\,.
$$
We have to prove that  there exists $\alpha\in
H_{p}(f^{c+\varepsilon})$ such that
$$
 i^{c+\varepsilon}_{*}(\alpha)=0 \quad\text{and}\quad \alpha\not\in
\ima i_{*}\,.
$$
The functoriality of the relative homology gives
\begin{equation}
  \label{eq.comdiag0}
\xymatrix{
H_{*}(f^{c-\varepsilon})
\ar[r]^{i^{c-\varepsilon}_{*}}\ar[d]^{i_{*}}&
H_{*}(M)
\ar[r]^{j^{c-\varepsilon}_{*}}
\ar[d]^{\Id}
&
H_{*}(M,f^{c-\varepsilon})
\ar[r]^{\partial^{c-\varepsilon}}
\ar[d]^{\overline{i_{*}}}
&
H_{*-1}(f^{c-\varepsilon})
\ar[d]^{i_{*}}\\
H_{*}(f^{c+\varepsilon})
\ar[r]^{i^{c+\varepsilon}_{*}}&
H_{*}(M)
\ar[r]^{j^{c+\varepsilon}_{*}}
&
H_{*}(M,f^{c+\varepsilon})
\ar[r]^{\partial^{c+\varepsilon}}
&
H_{*-1}(f^{c+\varepsilon})
}\,.
\end{equation}
The condition 2) of Proposition~\ref{pr.lower} and the long exact
sequence of relative homologies for the triple
$(M, ^{c+\varepsilon},f^{c-\varepsilon})$,  provide the exact
sequence
$$
\xymatrix{
0
\ar[r]&
H_{p+1}(M,f^{c-\varepsilon})
\ar[r]^{\overline{i_{*}}}&
H_{p+1}(M,f^{c+\varepsilon})
\ar[r]^{\partial}&
H_{p}(f^{c+\varepsilon},f^{c-\varepsilon})\to 0
}\,.
$$
Let $\alpha_{0}\in H_{p+1}(M,f^{c+\varepsilon})$ be such that
$\partial \alpha_{0}\neq 0$\,.
By the second line
 of the above commutative diagram \eqref{eq.comdiag0}, $\partial^{c+\varepsilon}\alpha_{0} \in
H_{p}(f^{c+\varepsilon})$ belongs to $\ker i^{c+\varepsilon}_{*}$\,.\\
Assume $\partial^{c+\varepsilon}\alpha_{0}\in \ima i_{*}$ and take
$\beta\in H_{p}(f^{c-\varepsilon})$ such that
$$
\partial^{c+\varepsilon}\alpha_{0}=i_{*}\beta\,.
$$
The commutative diagramm \eqref{eq.comdiag0} then implies
$$i^{c-\varepsilon}_{*}\beta=(\Id\circ
i^{c-\varepsilon}_{*})\beta=
(i^{c+\varepsilon}_{*}\circ i_{*})\beta
=i^{c+\varepsilon}_{*}(\partial^{c+\varepsilon}\alpha)=0\,.
$$
Hence $\beta\in \ker
i^{c-\varepsilon}_{*}=\ima \partial^{c-\varepsilon}$\,.
Hence there exists $\gamma\in H_{p+1}(M,f^{c-\varepsilon})$ such that
 $$
\partial^{c+\varepsilon}\alpha_{0}=
(i_{*}\circ \partial^{c-\varepsilon})\gamma
= \partial^{c+\varepsilon}\left(\overline{i_{*}}\gamma\right)\,.
$$
The cycle $\alpha_{0}-\overline{i_{*}}\gamma$ belongs to
$\ker \partial^{c+\varepsilon}=\ima j^{c+\varepsilon}_{*}$\,.
Hence there exists $\delta\in H_{p+1}(M)$ such that
$$
\alpha_{0}-\overline{i_{*}}\gamma=(j^{c+\varepsilon}_{*}\circ\Id)\delta=
(\overline{i_{*}}\circ j^{c-\varepsilon}_{*})\delta\,.
$$
We finally get
$\alpha_{0}=\overline{i_{*}}(\gamma+j^{c-\varepsilon}_{*}\delta)\in
\ima\overline{i_{*}}=\ker \partial$\,. But this contradicts the first
assumption $\partial\alpha_{0}\neq 0$\,.
Thus we have found $\alpha=\partial^{c+\varepsilon}\alpha_{0}\in
H_{p}(f^{c+\varepsilon})$ such that
$$
i^{c+\varepsilon}(\alpha)=
(i^{c+\varepsilon}\circ \partial^{c+\varepsilon})\alpha_{0}=0
\quad\text{and}\quad
\alpha \not\in \ima i_{*}\,.
$$
This ends the proof.
\end{proof}

\subsubsection{Properties of homological critical values}
\label{se.prohomcrit}

The attribute ``homological'' is justified by the following result.
Let $C_{H}^{(p)}$ be the set of homological critical values with index
$p$
and set $C_{H}=\cup_{p=0}^{d}C_{H}^{(p)}$. Remember the mappings
$i^{\lambda}_{*}:H_{*}(f^{\lambda})\to H_{*}(M)$\,.
\begin{theorem}
\label{th.homisom}
For every $p\in \left\{0,\ldots, d\right\}$, there is a one to one
mapping $\alpha^{(p)}:C_{H}^{(p)}\to H_{p}(M)$ such that:
\begin{itemize}
\item The range of $\alpha^{(p)}$ is a basis of $H_{p}(M)$;
\item For every $c\in C_{H}^{(p)}$, the quotient $\ima
  i^{c+\varepsilon}_{*}\big/ \ima i^{c-\varepsilon}_{*}$ is the
  one-dimensional space spanned by
  class of $\alpha^{(p)}(c)$\,.
\end{itemize}
Finally the cardinal of $C_{H}^{(p)}$ is the $p^{th}$ Betti number  of $M$\,.
\end{theorem}
We need the following result, where homological critical values differ
from the lower and upper critical values.
\begin{proposition}
 \label{pr.homcrit}
Assume that $c$ is an homological critical value (resp. point)
according to the Definition~\ref{de.class}.
Then the mapping
$$
H_{*}(f^{c+\varepsilon})\to H_{*}(M,f^{c-\varepsilon})
$$
is non zero.\\
Moreover the mappings $i^{c\pm\varepsilon}_{*}:H_{*}(f^{c\pm
  \varepsilon})\to H_{*}(M)$ satisfy
$$
\ima i^{c-\varepsilon}_{*}\subset \ima i^{c+\varepsilon}_{*}\quad,\quad 
\rank i^{c-\varepsilon}_{*}=\rank i^{c+\varepsilon}_{*}-1\,.
$$
\end{proposition}
\begin{proof}
  By definition, homological critical value is neither a lower
  critical value nor an upper critical value. Therefore the mappings
  \begin{eqnarray*}
&&    H_{*}(f^{c+\varepsilon},f^{c-\varepsilon})\to
H(M,f^{c-\varepsilon})\\
\text{and}
&&
H_{*}(f^{c+\varepsilon})\to H_{*}(f^{c+\varepsilon},f^{c-\varepsilon})
  \end{eqnarray*}
are non zero. Since $H_{*}(f^{c+\varepsilon},f^{c-\varepsilon})$ is
one dimensional, the second one is onto and the composed map
$$
\xymatrix{
H_{*}(f^{c+\varepsilon})
\ar[r]^{\sigma} &H_{*}(M,f^{c-\varepsilon})
}
$$
is non zero.
Consider now the second statement. We have already checked  in
\eqref{eq.comp}\eqref{eq.range} the relations
$$
i^{c-\varepsilon}_{*}=i^{c+\varepsilon}_{*}\circ
i^{c+\varepsilon,c-\varepsilon}_{*}\quad
\text{and}
\quad \ima i^{c-\varepsilon}_{*}\subset \ima i^{c+\varepsilon}_{*}\,.
$$
From the long exact sequence (when $c$ is a critical value with index $p$)
$$
\xymatrix{
0\ar[r]&
H_{p}(f^{c-\varepsilon})
\ar[r]^{i^{c+\varepsilon,c-\varepsilon}_{*}}&
H_{p}(f^{c+\varepsilon})
\ar[r]&
H_{p-1}(f^{c+\varepsilon},f^{c-\varepsilon}) 
\ar[r]&
0
}\,,
$$
we know that the codimension of $\ima
i^{c+\varepsilon,c-\varepsilon}_{*}$ is at most one. Thus
$$
\rank i^{c+\varepsilon}_{*}-1\leq \rank i^{c-\varepsilon}_{*}\leq \rank i^{c+\varepsilon}_{*}\,,
$$
and it suffices to find $\alpha \in \ima i^{c+\varepsilon}_{*}$ which
does not belong to $\ima i^{c-\varepsilon}_{*}$\,.
We use the diagram
$$
\xymatrix{
&
H_{*}(f^{c+\varepsilon})\ar[d]^{i^{c+\varepsilon}_{*}}
\ar[dr]^{\sigma}&&
\\
H_{*}(f^{c-\varepsilon})
\ar[r]^{i^{c-\varepsilon}_{*}}&
H_{*}(M)
\ar[r]^{j^{c-\varepsilon}_{*}}
&
H_{*}(M,f^{c-\varepsilon})
}\,.
$$
We know that there exists $\alpha_{0}\in H_{*}(f^{c+\varepsilon})$
such that $(j^{c-\varepsilon}_{*}\circ
i^{c+\varepsilon}_{*})(\alpha_{0})=\sigma(\alpha_{0})\neq 0$\,.
Take $\alpha= i^{c+\varepsilon}_{*}(\alpha_{0})$\,. It belongs to
$\ima i^{c+\varepsilon}_{*}$ and not in $\ker
j^{c-\varepsilon}_{*}=\ima i^{c-\varepsilon}_{*}$\,.
\end{proof}
\begin{proof}[Proof of Theorem~\ref{th.homisom}:] Fix the degree $p$,
  $0\leq p\leq d$, and  consider
  $i^{\lambda}_{*}:H_{p}(f^{\lambda})\to H_{p}(M)$\,.
Start from $\lambda=\max
  {f}+\varepsilon$ for which $\ima i^{\lambda}_{*}=\ima \Id=H_{p}(M)$
  and decrease $\lambda$ down to $\min f -\varepsilon$ for which 
$\ima  i^{\lambda}_{*}=\left\{0\right\}$\,. According to Morse theory,
  Proposition~\ref{pr.rkiup} and Proposition~\ref{pr.rkilow}, 
the range of $i^{\lambda}_{*}$
  does not change except when $\lambda$ passes an homological critical
  value with index $p$\,. For such a critical value,
  Proposition~\ref{pr.homcrit} says that the rank of
  $i^{\lambda}_{*}$ is exactly decreased by $1$\,. This yields the result.
\end{proof}

\subsection{The Morse-Barannikov chain complex}
\label{se.morbar}

\subsubsection{Definition}
\label{se.defcompl}
Remember that $\mathcal{C}(f)=\oplus_{p=0}^{d}\mathcal{C}^{(p)}(f)$ is
the vector space spanned by the critical points (identified with
the critical values and the same notation $c$ will be used for the two 
objects).
 The following definition will be proved to define a chain
complex structure on $\mathcal{C}(f)$ of which the homology groups are
isomorphic to the $H_{p}(M)$\,.
\begin{definition}
\label{de.complex}
On $\mathcal{C}(f)$ consider the linear mapping $\partial_{B}$ defined
by:
\begin{itemize}
\item When $c$ is a lower critical point or an homological critical
  point, $\partial_{B} c=0$\,.
\item When $c$ is an upper critical point, take for $c'$, according to
   the condition 3) of Proposition~\ref{pr.upper}, the supremum of
  the $\lambda'$s in $]-\infty,c)$ such that the mapping
  $H_{*}(f^{c+\varepsilon}, f^{\lambda})\to
  H_{*}(f^{c+\varepsilon},f^{c-\varepsilon})$ vanishes and
    set
$$
\partial_{B} c= c'\,.
$$
\end{itemize}
\end{definition}
\begin{theorem}
  \label{th.complex}
The mapping $\partial_{B}: \mathcal{C}(f)\to \mathcal{C}(f)$ sends
$\mathcal{C}^{(p)}(f)$ into $\mathcal{C}^{(p-1)}(f)$ and satisfies
$\partial_{B}\circ\partial_{B}=0$\,. 
Moreover the homology groups $H_{*}(\mathcal{C}(f))$ are isomorphic to
$H_{*}(M)$ and a basis of $H_{*}(M)$ is indexed by the set $C_{H}(f)$ of
homological critical points.
\end{theorem}
\begin{proof}
  It suffices to prove that when $c$ is an upper critical point with
  index $p$, the
  point $c'$ is a lower critical point with index $p-1$\,.\\
Assume that the mapping $H_{*}(f^{c+\varepsilon},
f^{c'-\varepsilon})\to H_{*}(f^{c+\varepsilon},f^{c-\varepsilon})$
vanishes while the mapping $\sigma : H_{*}(f^{c+\varepsilon},
f^{c'+\varepsilon})\to H_{*}(f^{c+\varepsilon}, f^{c-\varepsilon})$ is
non zero.
Consider the commutative diagram
$$
\xymatrix{
H_{*}(f^{c+\varepsilon},f^{c'-\varepsilon})
\ar[r]^{\varphi^{+}}\ar[d]^{0}&
H_{*}(f^{c+\varepsilon}, f^{c'+\varepsilon})
\ar[r]^{\partial^{+}}\ar[dl]^{\sigma}&
H_{*-1}(f^{c'+\varepsilon}, f^{c'-\varepsilon})\\
H_{*}(f^{c+\varepsilon}, f^{c-\varepsilon})
}
$$
where the first line is the long exact sequence for the triple
$f^{c'-\varepsilon}\subset f^{c'+\varepsilon}\subset
f^{c+\varepsilon}$\,.
Since $\sigma$ is non zero while $\sigma\circ \varphi^{+}=0$,
$\varphi^{+}$ cannot be onto and $\partial^{+}$ is non zero. We have found
$\lambda=c+\varepsilon$ such that the mapping $\partial:
H_{*}(f^{\lambda},f^{c'+\varepsilon})\to
H_{*-1}(f^{c'+\varepsilon},f^{c'-\varepsilon})$ is onto. By the
characterization 4) of Proposition~\ref{pr.lower}, $c'$ is a lower
critical point. Clearly it has the index $p-1$ when $c$ has the index $p$\,.\\
Therefore $\partial_{B}\circ \partial_{B} =0$.\\
Before we conclude, we check that if $\partial(c)=c'$, then $c$ is the
infimum of the $\lambda'$s such that $\partial : H_{*}(f^{\lambda},f^{c'+\varepsilon})\to
H_{*}(f^{c'+\varepsilon},f^{c'-\varepsilon})$ is onto.\\
We have to prove that the map $\partial^{-}:
H_{*}(f^{c-\varepsilon},f^{c'+\varepsilon})\to
H_{*}(f^{c'+\varepsilon}, f^{c'-\varepsilon})$ vanishes.
Consider the diagram
$$
\xymatrix{
H_{*}(f^{c-\varepsilon}, f^{c'-\varepsilon})
\ar[r]^{\varphi^{-}}
\ar[d]^{i^{c+\varepsilon,c-\varepsilon}_{*}}&
H_{*}(f^{c-\varepsilon},f^{c'+\varepsilon})
\ar[r]^{\partial^{-}}
\ar[d]^{i^{c+\varepsilon,c-\varepsilon}_{*}}&
H_{*-1}(f^{c'+\varepsilon},f^{c'-\varepsilon})
\ar[d]^{\Id}
\\
H_{*}(f^{c+\varepsilon}, f^{c'-\varepsilon})
\ar[r]^{\varphi^{+}}\ar[d]^{0}&
H_{*}(f^{c+\varepsilon},f^{c'+\varepsilon})
\ar[r]^{\partial^{+}}\ar[dl]^{\sigma}&
H_{*-1}(f^{c'+\varepsilon},f^{c'-\varepsilon})
\\
H_{*}(f^{c+\varepsilon},f^{c-\varepsilon})&&\,.}
$$
The maps $\sigma$ and $\partial^{+}$ have one dimensional ranges and
their kernels have the codimension $1$.
Due to $\sigma\circ \varphi^{+}=0$, we know $\ker \partial^{+}=\ima
\varphi^{+}\subset \ker \sigma$\,. With the same dimension, this
yields $\ker \partial^{+}=\ker \sigma$\,. If $\partial^{-}$ does not
vanish, there exists $u$ such that
$\partial^{+}(i^{c+\varepsilon,c-\varepsilon}_{+} u)=\partial^{-}u\neq
0$\,. Hence we get $(\sigma\circ
i^{c+\varepsilon,c-\varepsilon}_{*})u\neq 0$, which contradicts the
fact that $\sigma\circ i^{c+\varepsilon,c-\varepsilon}_{*}=0$ as a
part of the long exact sequence for the triple
$f^{c'+\varepsilon}\subset f^{c-\varepsilon}\subset
f^{c+\varepsilon}$\,.\\
Hence $c$ is the infimum of the $\lambda$'s such that 
 $\partial : H_{*}(f^{\lambda},f^{c'+\varepsilon})\to
H_{*}(f^{c'+\varepsilon},f^{c'-\varepsilon})$ is onto.\\
Now assume that $c'$ be a lower critical point.
By the characterization 4) of Proposition~\ref{pr.lower}, the infimum
of the $\lambda$'s in $(c,+\infty]$, such that  
$\partial : H_{*}(f^{\lambda},f^{c'+\varepsilon})\to
H_{*}(f^{c'+\varepsilon},f^{c'-\varepsilon})$ is onto, exists. Call it
$c$. By the dual argument of the previous one, $c$ is an upper
critical point and $c'$ is the supremum of the $\lambda$'s such that
the mapping $H_{*}(f^{c+\varepsilon},f^{\lambda})\to
H_{*}(f^{c+\varepsilon},f^{c-\varepsilon})$ vanishes. Hence $c$ is an
upper critical point such that $\partial_{B}(c)=c'$.\\
We have finally proved that the range of $\partial_{B}: \mathcal{C}(f)\to
\mathcal{C}(f)$ contains all the lower critical points.
The other statements are now straightforward consequences of 
Theorem~\ref{th.homisom}.
\end{proof}
\begin{remark}
  This result provides another proof of Morse inequalities, for
  excellent Morse functions, without making use of homotopy arguments
  to reduce the problem to self-indexed Morse functions (see
  \cite{Mil}\cite{Bott}\cite{Lau2}).
\end{remark}

\begin{center}
\newbox\surf


\figinit{0.12cm}

\figpt 1: (57,0)
\figpt 2: (20,14)

\figpt 11: (25,25)
\figpt 12: (20,45)
\figpt 13: (31,48)
\figpt 14: (57,36)
\figpt 15: (69,34)
\figpt 16: (64,12)
\figpt 17: (32,20)

\figpt 21: (15,61)
\figpt 22: (41,53)
\figpt 23: (74,63)



\figvisu{\surf}{Figure 1}{
\figinsert{Aaaaa.eps,0.6}
\figwrites 1: $H$ (0.5)
\figwrites 2: $c'_{1}$ (0.5)
\figwritew 11:  $H$ (0.1)
\figwritene 12: $H$ (0.1)
\figwrites 13:  $c'_{2}$ (0.)
\figwrites 14:  $c'_{3}$ (0.)
\figwritew 15: $H$ (0.4)
\figwritese 16:  $H$ (0.6)
\figwriten 17:  $c_{1}$ (0.)
\figwrites 21:  $c_{3}$ (0.)
\figwriten 22: $c_{2}$ (0.)
\figwritesw 23:  $H$ (0.5)
}
\hspace{2cm}\box\surf\\
Example with a compact surface with genius $2$ where
  $f$ is the height function. The homological critical points are
  labelled by $H$ while the pairing of other critical points follows
  $\partial_{\mathcal B} c_{k}=c_{k}'$\,.
\end{center}
\begin{proposition}
\label{pr.homolcoef}
  When $c$ is an upper critical point with index $p$ such that
  $\partial c=c'$, then the following commutative diagram holds:
$$
\xymatrix{
&&0\ar[d]&&\\
&&H_{p}(f^{c'+\varepsilon},f^{c'-\varepsilon})
\ar[d]^{i^{c-\varepsilon,c'+\varepsilon}_{*}}
& 0\ar[d]&\\
0
\ar[r]&
H_{p+1}(f^{c+\varepsilon}, f^{c-\varepsilon})
\ar[r]^{\partial}&
H_{p}(f^{c-\varepsilon},f^{c'-\varepsilon})
\ar[r]^{i^{c+\varepsilon,c-\varepsilon}_{*}}
\ar[d]^{j^{*}}&
 H_{p}(f^{c+\varepsilon},f^{c'-\varepsilon})
 \ar[r]\ar[d]^{j*}&0
\\
&0
\ar[r]&
H_{p}(f^{c-\varepsilon},f^{c'+\varepsilon})
\ar[r]^{i^{c+\varepsilon,c-\varepsilon}_{*}}\ar[d]&
H_{p}(f^{c+\varepsilon},f^{c'+\varepsilon})
\ar[r]\ar[d]\ar@<1ex>[u]\ar@<1ex>[l]&
0\\
&&0&0&
}\,.
$$ 
In particular, if $H_{p+1}(f^{c+\varepsilon},f^{c-\varepsilon})=\rz
[e^{p+1}]$  and
$H_{p}(f^{c'+\varepsilon},f^{c'-\varepsilon})=\rz[e^{p}]$, 
then there exists $\kappa \in \rz^{*}$ such that $\partial e^{p+1}$ and
$\kappa e^{p}$ are homologous in $f^{c-\varepsilon}$ relatively to
$f^{c'-\varepsilon}$: $[\partial e^{p+1}]=k[e^{p}]$ in $H_{p}(f^{c-\varepsilon},f^{c'-\varepsilon})$\,.
\end{proposition}
\begin{proof}
  From the definition of $\partial_{B}(c)=c'$, we know that the mapping 
$\partial^{-}: H_{p+1}(f^{c+\varepsilon},f^{c-\varepsilon}\to
H_{p}(f^{c+\varepsilon},f^{c'-\varepsilon})$ is one to one while the
mapping
$\partial^{+}: H_{p+1}(f^{c+\varepsilon},f^{c-\varepsilon}\to
H_{p}(f^{c+\varepsilon},f^{c'+\varepsilon})$ vanishes. Put in the long
exact sequences associated with the two triples
$(f^{c'-\varepsilon},f^{c-\varepsilon},f^{c+\varepsilon})$ and
$(f^{c'+\varepsilon},f^{c-\varepsilon},f^{c+\varepsilon})$, this
provides the two lines of the diagram.\\
Similarly, the relation $\partial_{B}(c)=c'$ implies that the mapping 
$\partial_{-}:H_{p+1}(f^{c+\varepsilon},f^{c'+\varepsilon})\to
H_{p}(f^{c'+\varepsilon},f^{c'-\varepsilon})$ vanishes (or
equivalently the mapping
$H_{p}(f^{c'+\varepsilon},f^{c'-\varepsilon})\to
H_{p}(f^{c+\varepsilon},f^{c'-\varepsilon})$ vanishes) and the mapping 
$\partial_{-}:H_{p+1}(f^{c-\varepsilon},f^{c'+\varepsilon})\to
H_{p}(f^{c'+\varepsilon},f^{c'-\varepsilon})$ vanishes. Inserted in
the long exact sequences associated with the triples 
$(f^{c'-\varepsilon},f^{c'+\varepsilon}, f^{c+\varepsilon})$ and
$(f^{c'-\varepsilon},f^{c'+\varepsilon}, f^{c-\varepsilon})$, this
provides the two columns of the diagram.\\
The diagram implies that the two mappings
$i^{c-\varepsilon,c'+\varepsilon}_{*}:
H_{p}(f^{c'+\varepsilon},f^{c'-\varepsilon})\to
H_{p}(f^{c-\varepsilon},f^{c'-\varepsilon})$ and
$\partial: H_{p+1}(f^{c+\varepsilon},f^{c-\varepsilon})\to
H_{p}(f^{c-\varepsilon},f^{c'-\varepsilon})$ have the same one
dimensional range. This ends the proof.
\end{proof}
\begin{center}
\newbox\surfa


\figinit{0.1cm}

\figpt 2: (46,7)

\figpt 11: (15,21)
\figpt 13: (31,35)
\figpt 14: (57,36)

\figpt 21: (15,28)
\figpt 22: (41,22)
\figpt 23: (105,20)


\figvisu{\surfa}{Figure 2}{
\figinsert{relat.eps,0.6}
\figwrites 2: $[e^{p}]$ (0.5)
\figwritew 11:  $[\partial e^{p+1}]$ (0.1)
\figwrites 13:  $c$ (0.)
\figwrites 14:  $c'_{3}$ (0.)
\figwritew 21:  $f=c-\varepsilon$ (0.)
\figwritene 22: $c'$ (0.3)
\figwritesw 23:  $f=c'-\varepsilon$ (0.5)
}
\hspace{4cm}\box\surfa\\
One example on a surface for which $[\partial e^{p+1}]=1 [e^{p}]$, $p=1$\,.
\end{center}
\subsubsection{Restriction}
\label{se.restric}
In the  previous construction the manifold $M$ equals $f^{+\infty}$
while the homology group $H_{*}(f^{\lambda})$ equals
$H_{*}(f^{\lambda},\emptyset)=H_{*}(f^{\lambda}, f^{-\infty})$\,. All
the construction can be done with sublevel sets $f^{a}$ and $f^{b}$
with $-\infty\leq a<b\leq +\infty$ which are not critical values. 
For $\lambda\in ]-\infty,+\infty[$ which is not a critical value,
consider 
$\mathcal{C}(f,\lambda)$,
the chain subcomplex of $\mathcal{C}(f)$ generated by critical values
(points) below level $\lambda$\,. Since $\partial_{B}$ preserves
$\mathcal{C}(f,\lambda)$, we can introduce the quotient
$\mathcal{C}(f,\lambda,\mu)=\mathcal{C}(f,\lambda)/\mathcal{C}(f,\mu)$
when $\mu< \lambda$\,. And there are relative homology groups 
$H_{*}(\mathcal{C}(f,\lambda),\mathcal{C}(f,\mu))$ for
$\partial_{B}$,
 which will be denoted by
$H_{*}(\mathcal{C}(f,\lambda),\mathcal{C}(f,\mu))$\,.

 All the previous definitions and proofs can be
translated to the restricted and relative homologies, after
replacing $H_{*}(M,f^{\lambda})$ by $H_{*}(f^{b},f^{\lambda})$ and
$H_{*}(f^{\lambda})$ by $H_{*}(f^{\lambda},f^{b})$, when $a<\lambda<b$ are not critical
values.
This observation gives at once.
\begin{theorem}
\label{th.relCf}
For any $a,b$, $-\infty\leq a < b\leq +\infty$, which are not critical values,
the relative homology groups
$H_{*}(\mathcal{C}(f,b),\mathcal{C}(f,a))$ are isomorphic to
$H_{*}(f^{b},f^{a})$ and the following diagram
$$
\xymatrix{
H_{*}(f^{a}) 
\ar[r]^{i^{b,a}_{*}}\ar[d]^{\simeq}&
H_{*}(f^{b})
\ar[r]^{j_{*}}\ar[d]^{\simeq}&
H_{*}(f^{b},f^{a})
\ar[r]^{\partial}\ar[d]^{\simeq}&
H_{*-1}(f^{a})
\ar[d]^{\simeq}
\\
H_{*}(\mathcal{C}(f,a))
\ar[r]&
H_{*}(\mathcal{C}(f, b))
\ar[r]&
H_{*}(\mathcal{C}(f,b), \mathcal{C}(f,a))
\ar[r]^{\partial_{B}}&
H_{*-1}(\mathcal{C}(f,a))
}
$$  
is commutative.
\end{theorem}
Since we have a good basis of the chain complex $(\mathcal{C}(f),\partial_{B})$,
where the image by $\partial_{B}$ of a generator is either $0$ or another
generator, we have a nice identification
$H_{*}(\mathcal{C}(f,b),\mathcal{C}(f,a))$\,.
\begin{proposition}
  \label{pr.restrCf}
The relative homology group $H_{*}(f^{b},f^{a})$ has a basis made of
critical values (resp. points) $c\in (a,b)$ satisfying one of the following conditions
\begin{enumerate} 
\item $c$ is an homological critical value
  (resp. point) in $M$\,;
\item or $c$ is an upper critical value (resp. point) such that
  $\partial c=c'$ is below $a$\,.
  \item or $c$ is a lower critical value
  (resp. point) in $M$, that is $\partial c'=c$ in $\mathcal{C}(f)$, but $c'$ is above $b$\,;
\end{enumerate} 
\end{proposition}

\begin{remark} 
What the theorem says is that the homological critical points for $\mathcal{C}(f,b,a)$, that should be denoted $\mathcal{C}_{H}(f,b,a)$ are {\bf not} the points of $\mathcal{C}_{H}(f)$ with critcal value in $[a,b]$, but the union of those, together with  the upper critical values in $[a,b]$ such that $\partial c$ is below $a$, and the lower critical values $c$ such that $\partial c'=c$ with $c'$ above $b$. 

\end{remark} 
\section{Relative Witten chain complex}
\label{se.wit}
The Witten Laplacian is a deformation of the Hodge Laplacian, related
with de Rham cohomology, which allows to give within a semiclassical
asymptotic framework, an analytic proof of the
Morse inequalities (see \cite{Wit, CFKS, HelSj4}). The
accurate computations of its exponentially small eigenvalues has
connections with various topics going from stochastic analysis
(\cite{FrWe, BEGK, BoGaKl})
 with kinetic theory (\cite{HerNi, HeNi2, HerSjSt, HerHiSj, HerHiSj2}),
 the computation of geometric invariants (\cite{BiLe}\cite{Bis}),
 differential topology (\cite{Mil, Bott, Lau1})
The case of manifold with boundaries has been considered in
\cite{ChLi, HeNi, KoPrSh, Lep2, Lep3} with a spectral
approach and more recently in \cite{Lau2} with a pure topological point
of view partly inspired by those previous works.
We shall consider here directly the case with boundary, which is of
interest here, and recall a few basic facts. We want to specify the
realization of the Witten Laplacian on the manifold
$\overline{f_{a}^{b}}$ with boundaries $\left\{f=a\right\}$ and
$\left\{f=b\right\}$, when 
$-\infty\leq a<b\leq +\infty$ are not critical values, which is
associated with the relative homology (after de Rham duality)
$H_{*}(f^{b},f^{a})$\,.\\

\subsection{Functional analysis}
\label{se.funcanal}

We recall that $\bigwedge T_{x}^{*}M=\oplus_{p=0}^{d}\bigwedge^{p}
T_{x}^{*}M$  is the exterior algebra on the cotangent fiber
$T_{x}^{*}M$, $\bigwedge T^{*}M$ is the corresponding fiber bundle and
$\mathcal{F}(M;\bigwedge T^{*}M)$ denotes the space of sections of
class $\mathcal{F}$ on $M$ ($\mathcal{F}$ stands for
$\mathcal{C}^{\infty}$, $L^{p}$ or $W^{m,r}$). The notation
$\mathcal{F}(\overline{f}_{a}^{b}; \bigwedge T^{*}M)$ is the set of
restrictions to $\overline{f_{a}^{b}}$ of elements in $\mathcal{F}(M;
\bigwedge T^{*}M)$\,.
 The spaces $\bigwedge
T_{x}^{*}M$ and $L^{2}(M; \bigwedge T^{*}M)$ and $L^{2}(M; \bigwedge
T^{*}M)$ are endowed with their natural scalar products inherited from
the riemannian metric $g$\,. 
A shorter notation for the Sobolev spaces will be
$$
\Lambda W^{m,r}= W^{m,r}(M; \bigwedge T^{*}M)\quad,\quad
\Lambda^{p} W^{m,r}(\overline{f_{a}^{b}})= W^{m,r}(\overline{f_{a}^{b}}; \bigwedge^{p} T^{*}M)\,.
$$
After the introduction of the
Hodge-$\star$ operator, the scalar product of two $p$-forms equals
$$
\langle\omega_{1}\,|\,\omega_{2}\rangle_{\Lambda^{p}L^{2}}
  =\int_{\Omega}\omega_{1}\wedge \star \overline{\omega_{2}}\,,
$$\\
The notation $\mathbf{t}$ and $\mathbf{n}$ are specific to the case
with boundary and useful for the analysis of boundary Hodge and Witten
Laplacians (see \cite{Sch}\cite{HeNi}\cite{Lep3}). Here is their
specific meaning: On the boundary $\partial \Omega$ of a regular
domain $\Omega$, decompose the tangent vectors $X_{i}\in
T_{\sigma}\Omega$, $\sigma\in \partial \Omega$, as 
$X_{i}=X_{i}^{T}+ x_{i}^{\perp}n_{\sigma}$ where $n_{\sigma}$ is a
normalized outgoing vector normal to $\partial \Omega$ and set
for $\omega\in \mathcal{C}^{\infty}(\Omega;\bigwedge T^{*}\Omega)$
\begin{eqnarray*}
&&   (\mathbf{t}\omega)_{\sigma}(X_{1},\ldots,
X_{p})=\omega_{\sigma}(X_{1}^{T},\ldots, X_{p}^{T})\;,\;
\forall \sigma \in \partial \Omega\;,\\
&&
\mathbf{n}\omega=\omega\Big|_{\partial\Omega}-\mathbf{t}\omega \quad
\in \mathcal{C}^{\infty}(\partial \Omega;\bigwedge^{p}T^{*}\Omega)\,.
\end{eqnarray*}
Note that the $\mathbf{t}\omega$ and $\mathbf{n}\omega$ have a natural
extension to a neighborhood of $\partial \Omega$  when the metric is
fixed. After the right choice of coordinates, with $x_{d}$
parametrizing normal curves to $\partial \Omega$, $\mathbf{t}\omega$
is the part with no $dx_{d}$ while 
$\mathbf{n}\omega$  takes the form $dx_{d}\wedge \omega'$\,.
On $\mathcal{C}^{\infty}$ differential forms, 
they satisfy various relations with the Hodge-$\star$ operator, the
differential $d$ and the codifferential $d^{*}$ (see
\cite{Sch}\cite{HeNi}\cite{Lep3} for details),
\begin{eqnarray}
\label{eq.hodge3}
&
\star d^{*,(p-1)}=(-1)^{p}d^{(d-p)} \star\;,&
\star d^{(p)}= (-1)^{p+1}d^{*,(d-p-1)}\star\;,\\
\label{eq.hodge4}
&\star \; \mathbf{n}=\mathbf{t}\; \star\;, &
\star\;\mathbf{t} =\mathbf{n}\; \star\;,\\
\label{eq.hodge5}
&\mathbf{t}\;d=d\;\mathbf{t}\;,&
\mathbf{n}\;d^{*}=d^{*}\; \mathbf{n}\;,
\end{eqnarray}
and the Stokes' formula,
$$
\forall \omega\in \mathcal{C}^{\infty}(\overline{\Omega}; 
\bigwedge^{p}T^{*}\Omega),\quad
\int_{\Omega}d\omega=\int_{\partial \Omega}\mathbf{t}\omega\;.
$$
When there is no boundary, 
the Hodge-de~Rham theory makes the relation between the spectral theory of
the Hodge~Laplacian and de~Rham duality of homology and
cohomology groups (see \cite{Ful}). For boundary manifold relative and
absolute (co-)homology groups can be considered.
By excision, remember that
$H_{*}(f^{b},f^{a})=H_{*}(\overline{f_{a}^{b}},
\left\{f=a\right\})$\,.
We briefly recall why these relative homology group are naturally
associated with specific boundary conditions for the Hodge and Witten
Laplacians. We refer the reader to \cite{Gue} for a review and
\cite{Tay} for a complete but different presentation 
relying on the isometric doubling of the boundary manifold. 
When $\gamma$ is a cycle in $\overline{f_{a}^{b}}$ relative to
$\left\{f=a\right\}$, there is a natural (i.e.
independent of the representant lying in $\left\{a\leq f\right\}$ of
$\gamma$)
integration $\int_{\gamma}\omega$ of forms
$\omega\in \mathcal{C}^{\infty}(\overline{f_{a}^{b}}; \bigwedge T^{*}M)$
such that $\mathbf{t}\omega\big|_{\left\{f=a\right\}}=0$\,. 
The dual condition along $\left\{f=b\right\}$,
$\mathbf{n}\omega\big|_{\left\{f=b\right\}}=0$, simply means 
$\mathbf{t}\omega\big|_{\left\{f=b\right\}}=\omega$ and ensures that
such form are determined by integration along chains lying in
$\left\{f\leq b\right\}$\,. 
The class $\mathcal{F}=\mathcal{C}^{\infty}_{TN}$,
$\mathcal{C}^{\infty}_{T}$ and $\mathcal{C}^{\infty}_{N}$ for
$\mathcal{C}^{\infty}$ forms fulfilling respectively both conditions,
the first one or the second one, among
\begin{equation}
  \label{eq.bc1}
\mathbf{t}\omega\big|_{\left\{f=a\right\}}=0
\quad
,
\quad
\mathbf{n}\omega\big|_{\left\{f=b\right\}}=0
\,.
\end{equation}
 Note that the differential $d$
preserves the class $\mathcal{C}^{\infty}_{T}$ while the codifferential
preserves $\mathcal{C}^{\infty}_{N}$\,. Better commutations relations appear
after considering the Hodge Laplacian (see further).
The de Rham cohomology group $H^{^{p}}(f^{b},f^{a})$ is then given by
$\ker d^{(p)}/ (\ima d^{p-1}\cap
\mathcal{C}^{\infty}_{TN}(\overline{f_{a}^{b}}, \bigwedge^{p} T^{*}M))$, when
$d$ is defined on $\mathcal{C}^{\infty}_{TN}(\overline{f_{a}^{b}},
\bigwedge T^{*}M)$\,.\\
The Witten deformation consists in introducing a small parameter $h\to
0$ and to set
$$
d_{f,h}=e^{-\frac{f}{h}}(hd)e^{\frac{f}{h}}=hd+df\wedge\quad
,\quad
d_{f,h}^{*}=e^{\frac{f}{h}}(hd^{*})e^{\frac{f}{h}}
=hd^{*}+\mathbf{i}_{\nabla f}\,.
$$
The Witten Laplacian is defined as a differential operator
in $f_{a}^{b}=\left\{a< f <b\right\}$ by
$$
\Delta_{f,h}=(d_{f,h}+d_{f,h}^{*})^{2}=d_{f,h}d_{f,h}^{*}+d_{f,h}^{*}d_{f,h}
= h^{2}(d+d^{*})^{2}+|\nabla f|^{2}+h (\mathcal{L}_{\nabla
  f}+\mathcal{L}_{\nabla f}^{*})\,.
$$
On the boundaries $\left\{f=a\right\}$ and $\left\{f=b\right\}$, the
boundary conditions  have to be completed with $f$-dependent
additional boundary
conditions in order to get a self-adjoint realization which
is elliptic up to the boundary (see \cite{Sch}).
 An additional property required here, is the commutation of the resolvent with $d_{f,h}$
and $d_{f,h}^{*}$\,. We follow the scheme of \cite{ChLi,HeNi,Lep3} where the ``Dirichlet problem'' and the ``Neumann problem'' have
been considered separately. Here the ``Dirichlet'' boundary conditions
occurs on $\left\{f=a\right\}$ while the ``Neumann'' boundary
condition appears on $\left\{f=b\right\}$\,. 
Consider in $\Lambda W^{1,2}_{TN}=\overline{\mathcal{C}^{\infty}_{TN}(\overline{f_{a}^{b}};
  \bigwedge T^{*}M)}^{W^{1,2}}$ the quadratic form given by
\begin{eqnarray*}
&&\mathcal{D}_{TN}(\omega,\eta)=
\langle d_{f,h}\omega\,\big|\,d_{f,h}\eta\rangle
+
\langle d_{f,h}^{*}\omega\,\big|\,d_{f,h}^{*}\eta\rangle\,.
\\
&&\mathcal{D}_{TN}(\omega)= \mathcal{D}_{TN}(\omega,\omega)
= \|d_{f,h}\omega\|_{L^{2}}^{2}+\|d_{f,h}^{*}\omega\|_{L^{2}}^{2}\,.
\end{eqnarray*}
Since $\left\{f=a\right\}$ and $\left\{f=b\right\}$ are disjoint,
the main  arguments are local (Sobolev trace theorem,
Lopatinski-Schapiro conditions for the ellipticity up to the boundary
and finally playing with
\eqref{eq.hodge3}\eqref{eq.hodge4}\eqref{eq.hodge5}), we can combine
without repeating the proofs the results of \cite{HeNi} and
\cite{Lep3} in order to state the following result.\\
Note that due to the boundaries of the domain $\Omega$, we avoid to consider the
closure of the differential operators
$d_{f,h}$ and $d_{f,h}^{*}$ in $\Lambda L^{2}$ which are not very
explicit.
\begin{proposition}
  \label{pr.selfadjoin}~\\
  The non-negative quadratic form $\omega\to\mathcal{D}_{TN}(\omega)$
is closed on $\Lambda W^{1,2}_{TN}$. 
The associated (self-adjoint) Friedrichs 
extension  is denoted by 
$\Delta_{f,h}^{TN}$. Its domain is
$$
D(\Delta_{f,h}^{TN})=
\left\{
\omega\in \Lambda W^{2,2}(\overline{f_{a}^{b}});\quad
\begin{array}[c]{ll}
   \mathbf{t}\omega\big|_{f=a}=0\,,&
  \mathbf{t}d_{f,h}^{*}\omega\big|_{f=a}=0
\,,\\
  \mathbf{n}\omega\big|_{f=b}=0\,,&
 \mathbf{n}d_{f,h}\omega\big|_{f=b}=0
\end{array}
\right\}\;,
$$
and acts as 
$$
\forall \omega \in D(\Delta_{f,h}^{TN}),\quad
\Delta_{f,h}^{TN}\omega=\Delta_{f,h}\omega\;.
$$
The operator $\Delta_{f,h}^{TN}$ has a compact resolvent and a
discrete spectrum.
Moreover the commutations
\begin{eqnarray*}
  &&
(z-\Delta_{f,h}^{TN})^{-1}\circ d_{f,h}\omega =
d_{f,h}\circ (z-\Delta_{f,h})^{-1}\omega\,,\\
&&
(z-\Delta_{f,h}^{TN})^{-1}\circ d_{f,h}^{*}\omega =
d_{f,h}^{*}\circ (z-\Delta_{f,h})^{-1}\omega\,,\\
&&
1_{E}(\Delta_{f,h}^{TN})\circ d_{f,h}\omega=d_{f,h}
\circ 1_{E}(\Delta_{f,h}^{TN})\omega\,,\\
\text{and}
&&
1_{E}(\Delta_{f,h}^{TN})\circ
d_{f,h}^{*}\omega=d_{f,h}^{*}\circ 1_{E}(\Delta_{f,h}^{TN})\omega\,,
\end{eqnarray*}
hold for all $z\in \cz\setminus \rz$, all Borel set $E$ in $\rz$ and all
$\omega \in \Lambda W^{1,2}_{TN}$\,.
\end{proposition}
\begin{remark}
\label{re.NT}
\begin{itemize}
\item The introduction of $\Delta_{f,h}^{TN}$, as a Friedrichs
  extenstion of a non negative closed quadratic form defined on $\Lambda
  W^{1,2}_{TN}$, ensures that it is a non negative self-adjoint
  operator. Requiring $\Delta_{f,h}^{TN}u\in \Lambda L^{2}$ for $u\in
  D(\Delta_{f,h}^{TN})$ forces the additional boundary conditions,
  after integration by part.
\item The commutation relations do not result simply of the
  commutatition of the differential operators $\Delta_{f,h}\circ
  d_{f,h}=d_{f,h}\circ
\Delta_{f,h}$ valid in the interior $f_{a}^{b}$\,. Indeed
$\Delta_{f,h}^{TN}$ can be applied only to elements of
$D(\Delta_{f,h}^{TN})$ fulfilling the boundary conditions while
$d_{f,h}$ do not preserve these boundary conditions even for ${\cal
  C}^{\infty}$-forms up to the boundaries. For details and complete
proofs, we refer again the reader to \cite{HeNi}\cite{Lep3} and \cite{ChLi}.
\item We can also define analogously the 
self-adjoint
operator $\Delta_{f,h}^{NT}$,
with domain $D(\Delta_{f,h}^{NT})$,
by switching the above conditions on 
$\mathbf n$ and $\mathbf t$.
\end{itemize}
\end{remark}

Here is the main result of this section.
\begin{theorem}
\label{th.Witt}
There are two operators $L_{\pm}\in \mathcal{L}(\Lambda L^{2};\Lambda
W^{1,2}_{TN})$, commuting for all Borel sets $E$ in $\rz$ with  $1_{E}(\Delta_{f,h}^{TN})$, such that every $u\in \Lambda L^{2}$ admits the
orthogonal decomposition 
\begin{equation}
  \label{eq.orthdecom}
u= 1_{\left\{0\right\}}(\Delta_{f,h}^{TN})u +  d_{f,h} L_{-}u + 
d_{f,h}^{*} L_{+}u\,.
\end{equation}
When $F_{M}$ denotes the finite dimensional space 
$\ima 1_{[0,M]}(\Delta_{f,h}^{TN})$ and
$\beta_{M}=d_{f,h}\big|_{F_{M}}$, its adjoint is 
 $\beta_{M}^{*}=d_{f,h}^{*}\big|_{F_{M}}$ and $F_{M}$ admits the
 orthogonal decomposition
$$
F_{M}=\ker \Delta_{f,h}^{TN}
\mathop{\oplus}^{\perp}\ima \beta_{M}
\mathop{\oplus}^{\perp}\ima \beta_{M}^{*}\,.
$$
 After setting 
$F^{(p)}_{M}=\ima  1_{[0,M]}(\Delta_{f,h}^{TN,(p)})$, the two finite
dimensional  chain complexes
\begin{equation}
  \label{eq.redcompl}
\xymatrix{
0\ar[r]&
F_{M}^{(0)}
\ldots
F_{M}^{(p-1)}
\ar[r]^-{\beta_{M}^{(p-1)}}\ar@<1ex>[l]&
F_{M}^{(p)}
\ar[r]^-{\beta_{M}^{(p)}}\ar@<1ex>[l]^-{\beta^{(p-1)*}_{M}}&
F_{M}^{(p+1)}
\ldots
F_{M}^{(d)}\ar[r]\ar@<1ex>[l]^-{\beta^{(p)*}_{M}}&0
\ar@<1ex>[l]
}
\end{equation}
are dual to each other and $\ker \beta_{M}^{(p)}/\ima \beta_{M}^{(p-1)}$ is
diffeomorphic to $H^{p}(f^{b},f^{a})$\,.\\
For any $p\in \left\{0,\ldots, d\right\}$, the spectrum of
$\sigma(\Delta_{f,h}^{TN,(p)})\cap (0,M]$ lying in $(0,M]$ and counted
with multiplicities, the set  of
$\lambda^{2}$ (counted with multiplicities) when $\lambda$ ranges over the singular values, counted
with multiplicities, of
$\beta_{M}^{(p)}\big|_{\ima \beta_{M}^{(p),*}}$ and 
$\beta_{M}^{(p-1)}\big|_{\ima \beta_{M}^{(p-1),*}}$\,.
\end{theorem}
We shall need the two following lemmas.
\begin{lemma}
\label{le.wtn}
  When $\omega$ belongs to $D(\Delta_{f,h}^{TN})$, $d_{f,h}\omega$ and
  $d_{f,h}^{*}\omega$ belong to $\Lambda W^{1,2}_{TN}$\,.
\end{lemma}
\begin{proof}
The differential operators  $d_{f,h}$ and $d_{f,h}^{*}$ are continuous
from
$D(\Delta_{f,h}^{TN})\subset \Lambda W^{2,2}$ into $\Lambda W^{1,2}$\,.
By the elliptic regularity up to the boundary of $\Delta_{f,h}^{TN}$, the set of
$\mathcal{C}^{\infty}(\overline{f_{a}^{b}}; \bigwedge T^{*}M)\cap
D(\Delta_{f,h}^{TN}))$
is dense in $D(\Delta_{f,h}^{TN})$ because $1+\Delta_{f,h}^{TN}:
D(\Delta_{f,h}^{TN})\to \Lambda L^{2}$ is an isomorphism.
For $\omega\in \mathcal{C}^{\infty}(\overline{f_{a}^{b}}; \bigwedge T^{*}M)\cap
D(\Delta_{f,h}^{TN}))$, we have
$$
\mathbf{n}(d_{f,h} \omega)\big|_{f=b}=0\quad \text{and}\quad 
\mathbf{t}(d_{f,h} \omega)\big|_{f=a}=0
$$
because $\omega\in D(\Delta_{f,h}^{TN})$\,. Moreover $\omega\in
D(\Delta_{f,h}^{TN})$ also says
$$
\mathbf{t}\omega\big|_{f=a}=0\quad \text{and}\quad \mathbf{n}\omega\big|_{f=b}=0\,.
$$
But since $\mathbf{t}e^{\pm \frac{f}{h}}\omega=e^{\pm
  \frac{f}{h}}\mathbf{t}\omega$ and $\mathbf{n}e^{\pm
  \frac{f}{h}}\omega=e^{\pm \frac{f}{h}}\mathbf{n}\omega$, the
commutations \eqref{eq.hodge5} imply
$$
\mathbf{t}(d_{f,h}\omega)\big|_{f=a}=0\quad \text{and}\quad \mathbf{n}(d_{f,h}^{*}\omega)\big|_{f=b}=0\,.
$$
This ends the proof.
\end{proof}
\begin{lemma}
\label{le.ipp} The relation
$$
\langle d_{f,h}\theta_{1}, \theta_{2}\rangle_{\Lambda L^{2}}=\langle
\theta_{1}, d_{f,h}^{*}\theta_{2} \rangle_{\Lambda L^{2}}
$$
holds for all $\theta_{1},\theta_{2}\in \Lambda W^{1,2}_{TN}$\,.
\end{lemma}
\begin{proof}
  Since both quantities are continous and
  $\mathcal{C}^{\infty}_{TN}(\overline{f_{a}^{b}}; \bigwedge T^{*}M)$
  is dense in $\Lambda W^{1,2}_{TN}$, it suffices to consider
  $\theta_{1}, \theta_{2}\in
  \mathcal{C}^{\infty}_{TN}(\overline{f_{a}^{b}}; \bigwedge
  T^{*}M)$\,.
After writing $\omega_{1}=e^{\frac{f}{h}}\theta_{1}$ and
$\omega_{2}=e^{-\frac{f}{h}}\theta_{2}$ in $\mathcal{C}^{\infty}_{TN}(\overline{f_{a}^{b}}; \bigwedge
  T^{*}M)$, our identity amounts to
$$
\langle d\omega_{1}, \omega_{2}\rangle_{\Lambda L^{2}}=\langle
\omega_{1}, d^{*}\omega_{2} \rangle_{\Lambda L^{2}}\,.
$$
But the Stokes' formula with the relations
between $d$, $\star$ and $\wedge$ gives
$$
\int_{f=a}\mathbf{t}(\omega_{1}\wedge \star \omega_{2})
+ 
\int_{f=b}\mathbf{t}(\omega_{1}\wedge \star \omega_{2})
=
\langle d\omega_{1}\,|\, \omega_{2}\rangle_{\Lambda L^{2}}
+ 
\langle \omega_{1}\,|\, d^{*}\omega_{2}\rangle\,.
$$
With the help of \eqref{eq.hodge4}, write
\begin{multline*}
\mathbf{t}(\omega_{1}\wedge \star \omega_{2})=
\mathbf{t}\left[(\mathbf{t}\omega_{1})\wedge \star
  (\mathbf{t}\omega_{2})
+ (\mathbf{t}\omega_{1})\wedge \star
  (\mathbf{n}\omega_{2})
+(\mathbf{n}\omega_{1})\wedge \star
  (\mathbf{t}\omega_{2})
\right.
\\
\left.
+(\mathbf{n}\omega_{1})\wedge \star
  (\mathbf{n}\omega_{2})
\right]
\\
=\mathbf{t}\left[(\mathbf{t}\omega_{1})\wedge
  \mathbf{n}(\star\omega_{2})
+ (\mathbf{t}\omega_{1})\wedge \mathbf{t}
  (\star \omega_{2})
+(\mathbf{n}\omega_{1})\wedge
  \mathbf{t}(\star\omega_{2})
\right.
\\
\left.
+(\mathbf{n}\omega_{1})\wedge \mathbf{n}
  (\star\omega_{2})
\right]
\end{multline*}
and notice that $\mathbf{t}((\mathbf{n}u_{1})\wedge
(\mathbf{t}u_{2}))=\mathbf{t}((\mathbf{t}u_{1})\wedge
(\mathbf{n}u_{2}))=0$\,. This leads to
$$
\mathbf{t}(\omega_{1}\wedge \star \omega_{2})=
\mathbf{t}\left[
(\mathbf{t}\omega_{1})\wedge \star
  (\mathbf{n}\omega_{2})
+(\mathbf{n}\omega_{1})\wedge \star
  (\mathbf{t}\omega_{2})
\right]\,,
$$
where both terms vanishes on $\left\{f=a\right\}\cup
\left\{f=b\right\}$ when $\omega_{1},\omega_{2}\in
\mathcal{C}^{\infty}_{TN}(\overline{f_{a}^{b}}; \bigwedge T^{*}M)$\,.
\end{proof}
\begin{proof}[Proof of Theorem~\ref{th.Witt}]
The operator $\Delta_{f,h}^{TN}$ is a self-adjoint operator  with a
compact resolvent. Therefore it is invertible when restricted to
$\ker(\Delta_{f,h}^{TN})^{\perp}= \ima
1_{(0,+\infty)}(\Delta_{f,h}^{*})$\,.
Take
$$
L_{+}=d_{f,h}(\Delta_{f,h}^{TN})^{-1}1_{(0,+\infty)}(\Delta_{f,h}^{TN})\quad
\text{and}\quad
L_{-}=d_{f,h}^{*}(\Delta_{f,h}^{TN})^{-1}1_{(0,+\infty)}(\Delta_{f,h}^{TN})\,,
$$
Since $(\Delta_{f,h}^{TN})^{-1}1_{(0,+\infty)}(\Delta_{f,h}^{TN})\in
\mathcal{L}(\Lambda L^{2}, D(\Delta_{f,h}^{TN}))$, $L_{\pm}\in
\mathcal{L}(\Lambda L^{2}, \Lambda W^{1,2}_{TN})\subset
\mathcal{L}(\Lambda L^{2}, \Lambda L^{2})$ according to
Lemma~\ref{le.wtn}. The commutation of $L_{\pm}$ with
$1_{E}(\Delta_{f,h})$ is then a consequence of the commutation stated
in Proposition~\ref{pr.selfadjoin}. These commutations of
Proposition~\ref{pr.selfadjoin} also imply
$$
L_{+}\omega=(\Delta_{f,h}^{TN})^{-1}1_{(0,+\infty)}(\Delta_{f,h}^{TN})
d_{f,h}\omega\quad
\text{and}\quad
L_{-}\omega=(\Delta_{f,h}^{TN})^{-1}1_{(0,+\infty)}(\Delta_{f,h}^{TN})d_{f,h}^{*}\omega\,,
$$
when $\omega\in \Lambda W^{1,2}_{TN}$ so that $L_{\pm}\in
\mathcal{L}(\Lambda W^{1,2}_{TN}; D(\Delta_{f,h}^{TN}))$\,.
By using again Lemma~\ref{le.wtn}, $d_{f,h}\circ L_{-}$ and
$d_{f,h}^{*}\circ L_{+}$ belong to $\mathcal{L}(\Lambda
W^{1,2}_{TN})$\,.
Consider now the decomposition
$$
\omega= 1_{\left\{0\right\}}(\Delta_{f,h}^{TN})\omega
+ d_{f,h}L_{-}\omega + d_{f,h}^{*}L_{+}\omega
$$
when $\omega\in \Lambda W^{1,2}_{TN}$\,. All the terms belong to
$\Lambda W^{1,2}_{TN}$ while 
\begin{eqnarray*}
  &&
d_{f,h}(d_{f,h}L_{-}\omega)=
d_{f,h}(1_{\left\{0\right\}}(\Delta_{f,h}^{TN})\omega)=0\\
\text{and}&&
d_{f,h}^{*}(d_{f,h}^{*}L_{+}\omega)=
d_{f,h}^{*}(1_{\left\{0\right\}}(\Delta_{f,h}^{TN})\omega)=0\,.
\end{eqnarray*}
Therefore Lemma~\ref{le.ipp} implies that the decomposition is
orthogonal when $\omega\in \Lambda W^{1,2}_{TN}$, and  this extends 
by continuity to $\omega\in \Lambda L^{2}$\,.\\
Proposition~\ref{pr.selfadjoin} ensures that $d_{f,h}$ and
$d_{f,h}^{*}$ send $F_{M}$ into itself and Lemma~\ref{le.ipp} with 
$F_{M}\subset D(\Delta_{f,h}^{TN})\subset \Lambda W^{1,2}_{TN}$
implies
$$
\beta_{M}^{*}=(1_{[0,M]}(\Delta_{f,h}^{TN})d_{f,h}1_{[0,M]}(\Delta_{f,h}^{TN}))^{*}=
(1_{[0,M]}(\Delta_{f,h}^{TN})d_{f,h}^{*}1_{[0,M]}(\Delta_{f,h}^{TN}))\,.
$$
The orthogonal decomposition follows from the result for $\omega\in
\Lambda L^{2}$\,. The chain complex structure comes from $d_{f,h}\circ
d_{f,h}=d_{f,h}^{*}\circ d_{f,h}^{*}=0$.\\
The space $\ker \beta_{M}^{(p)}/\ima \beta_{M}^{(p-1)}$ is isomorphic
to $\ker(\Delta_{f,h}^{(p)})$, which is contained like $F_{M}$ in 
$\mathcal{C}^{\infty}_{TN}(\overline{f_{a}^{b}};\bigwedge^{p}
T^{*}M)$ by the elliptic regularity up to the boundary of
$\Delta_{f,h}^{(p)}$\,. Hence $\ker \beta_{M}^{(p)}/\ima
\beta_{M}^{(p)}$ is isomorphic to $\ker d_{f,h}^{(p)}/\ima
d_{f,h}^{(p-1)}$ after considering the differential operators
$d_{f,h}$ restricted to
$\mathcal{C}^{\infty}_{TN}(\overline{f}^{a}_{b},\bigwedge T^{*}M)$\,.
Since $\omega\mapsto e^{-\frac{f}{h}}\alpha$ is an isomorphism between the two
spaces $\mathcal{C}^{\infty}_{TN}(\overline{f}^{a}_{b},\bigwedge
T^{*}M)$ respectively defined for $\Delta_{f,h}$ or $d_{f,h}$ and
$\Delta_{0,h}=h^{2}\Delta_{\text{Hodge}}$ or $d$, we obtain 
$$
ker \beta_{M}^{(p)}/\ima \beta_{M}^{(p-1)}\sim \ker d^{(p)}/\ima
d^{(p-1)}
=H^{p}(f^{a},f^{b})
$$
by Hodge-de~Rham theory (see for instance \cite{Ful} for the usual
boundaryless case and \cite{Gue}\cite{Tay} for the case with boundary).\\
The result concerned with the spectrum of $\Delta_{f,h}^{TN,(p)}$ is a
direct consequence of the orthogonal decomposition of $F_{M}$ with the
chain complex structure \eqref{eq.redcompl}. To be more specific, the
decomposition of $\Delta_{f,h}^{TN,(p)}\big|_{F_{M}^{(p)}}$ according to 
$\displaystyle F_{M}^{(p)}=\ker \Delta_{f,h}^{TN,(p)}
\mathop{\oplus}^{\perp}\ima \beta_{M}^{(p-1)}
\mathop{\oplus}^{\perp}\ima \beta_{M}^{*,(p+1)}$ writes
$$
\Delta_{f,h}^{TN,(p)}\big|_{F_{M}^{(p)}}=
0
\mathop{\oplus}^{\perp}\beta_{M}^{(p-1)}\beta^{(p-1),*}_{M}
\mathop{\oplus}^{\perp} \beta_{M}^{(p),*}\beta^{(p)}_{M}
$$
while $\beta^{(p-1),*}_{M}$ is an isomorphism from $\ima \beta^{(p-1)}_{M}$
onto $\ima \beta^{(p-1),*}$ with
$$
\beta^{(p-1),*}_{M}\left(\beta_{M}^{(p-1)}\beta^{(p-1),*}_{M}\right)=
\left(\beta^{(p-1),*}_{M}\beta_{M}^{(p-1)}\right)\beta^{(p-1),*}_{M}\,.
$$
\end{proof} 
\begin{remark}
  Note that the duality between the two chain complexes associated with
  $\beta_{M}$ and $\beta_{M}^{*}$ and their homology groups, is
  another version of the topological duality $f\to -f$\,. Actually
  changing $f$ to $-f$ and $p$-forms with $d-p$-forms with the
  Hodge-$\star$
 operator, interchanges $d_{f,h}$ and $d_{f,h}^{*}$\,.
\end{remark}

\subsection{Adapting Helffer-Sj\"ostrand analysis}
\label{se.adaptHS}

We still work in $\overline{f_{a}^{b}}$ and we introduce like  in
\cite{HelSj4} the Agmon distance $d_{Ag}$ associated with the degenerate
metric $|\nabla f|^{2}g$, where $g$ is the initial Riemannian  metric
on $M$\,. This distance satisfies
$$
d_{Ag}(x,y)\geq |f(x)-f(y)|
$$
with equality when an integral curve of $\nabla f$ joins $x$ and
$y$\,.\\
\\
Before stating the following crucial theorem, let us introduce two definitions
which will be very useful in the sequel. The first one recalls
Helffer-Sj\"ostrand notation $\tilde{\mathcal{O}}$, very convenient when
handling exponentially small quantities.
\begin{definition}
For two quantities $A(h)$,  estimated with a norm $|A(h)|$,
 and $B(h)\geq 0$, parametrized by $h\in (0,h_{0})$,
the notation $A(h)=\tilde{\mathcal{O}}(B(h))$ means:
$$
\forall \varepsilon>0,\exists C_{\varepsilon}>0,\quad
\forall h\in (0,h_{0}), |A(h)|\leq C_{\varepsilon}B(h)e^{\frac{\varepsilon}{h}}\,.
$$
\end{definition}

\begin{definition}
 Let $U\in M$ be a critical point of $f$ with index $p$
 and let $\Phi(x):=d_{Ag}(x,U)$.
 A local coordinate system $y_{1},\dots,y_{d}$ around 
 $U$ is said to be an adapted Morse coordinate system for $f$ if
 $y_{1},\dots,y_{d}$ is
 centered at $U$,
 $dy_{1},\dots,dy_{d}$ is an orthonormal positively oriented
 basis of $T^*_{U}M$, and if,
  in these coordinates, 
  the following Morse decompositions for $f$
  and $\Phi$,
$$
 f(y)=f(U)+\frac{1}{2}\sum_{j=1}^{d}\lambda_{j}y_{j}^{2}\quad,\quad
 \Phi(y)=\frac{1}{2}\sum_{j=1}^{d}|\lambda_{j}|y_{j}^{2}\,,$$
hold locally around $U$, with
$\lambda_{j}<0$ for $j\leq p$ and $\lambda_{j}>0$ for $j>p$.
\end{definition}

\noindent
Let us notice that such a coordinate system always exists,
according to  \cite{HelSj4}~pp.~272--281.
Moreover, in such a coordinate system,
the stable and unstable manifolds of $-\nabla f$, respectively
denoted by $\mathcal C_{\text{St}}$
and $\mathcal C_{\text{Unst}}$
are locally parametrized by 
\begin{equation}
\label{eq.stable}
\mathcal C_{\text{St}} =\{y\,;\ y_{1}=\cdots = y_{p}=0\}
\quad \text{and}
\quad
\mathcal C_{\text{Unst}} =\{y\,;\ y_{p+1}=\cdots = y_{d}=0\}
\,.
\end{equation}

\begin{theorem}
\label{th.HSadapt}
Let $p$ belong to $\left\{0,\ldots, d\right\}$ and denote by
$\mathcal{U}^{(p)}=\left\{U_{1}^{(p)},\ldots, U_{m_{p}}^{(p)}\right\}$ 
the set  critical points of $f$ in
$f_{a}^{b}$\,.
There exists $h_{0}>0$ such that, for all $h\in(0,h_{0}]$,  
the spectral subspace $F^{(p)}=1_{[0, Ch^{3/2}]}(\Delta_{f,h}^{TN (p)})$ is
spanned by $m_{p}$ normalized vector $v_{j}$, $1\leq j\leq m_{p}$,
which satisfy, for $x\in M$ and $\alpha\in \nz^{d}$,
\begin{eqnarray}
  \label{eq.Agmder}
&&
|\partial_{x}^{\alpha}v_{j}|=\tilde{\mathcal{O}}(e^{-\frac{d_{Ag}(x,U_{j}^{(p)})}{h}})\,,\\
&&
\label{eq.Agmderdd}
|\partial_{x}^{\alpha}d_{f,h}v_{j}|=
\tilde{\mathcal{O}}
\left(e^{-\frac{\beta_{j}^{+}(x)
}{h}}
\right)\quad,\quad
|\partial_{x}^{\alpha}d_{f,h}^{*}v_{j}|=
\tilde{\mathcal{O}}
\left(e^{-\frac{\beta_{j}^{-}(x)
}{h}}
\right)\\
\nonumber
\text{with}&&
\beta_{j}^{+}(x)=\min_{
U\in
\,\mathcal{U}^{(p+1)}\cup
    \,\mathcal{U}^{(p)}\setminus\left\{U_{j}^{(p)}\right\}
}
  d_{Ag}(U_{j},U)+d_{Ag}(U,x)
\\
\nonumber
\text{and}
&&
\beta_{j}^{-}(x)=\min_{
U\in\,
\mathcal{U}^{(p-1)}\cup\,\mathcal{U}^{(p)}\setminus\left\{U_{j}^{(p)}\right\}
}
  d_{Ag}(U_{j},U)+d_{Ag}(U,x)\,.
\end{eqnarray}
The eigenvalues of $\Delta_{f,h}^{TN (p)}$ lying in $[0,h^{3/2}]$ are
$\mathcal{O}(e^{-\frac{C}{h}})$\,.\\
When the metric $g$ is Euclidean in some adapted Morse coordinates for $f$ in
$B(U_{j}^{(p)},2\eta)$, with
$f(y)=f(U_{j})+\frac{1}{2}\sum_{j=1}^{d}\lambda_{j}y_{j}^{2}$, then the
form $v_{j}$ satisfies
$$
v_{j}=|\lambda_{1}\ldots \lambda_{d}|^{1/4}(\pi
h)^{-d/4}e^{-\frac{\sum_{j=1}^{d}|\lambda_{j}|y_{j}^{2}}{2h}}
dy_{1}\wedge\ldots \wedge dy_{p} +
{\mathcal{O}}\left(e^{-\frac{C_{\eta}}{h}}\right)
$$
in $\mathcal{C}^{\infty}(B(U_{j}^{(p)},\eta))$\,.\\
In the general case of a Riemannian metric,  there exists,
for $\eta$ small enough, some adapted Morse coordinates for $f$ in
$B(U_{j}^{(p)},2\eta)$, with
$f(y)=f(U_{j})+\frac{1}{2}\sum_{j=1}^{d}\lambda_{j}y_{j}^{2}$,
and such that the
form $v_{j}$ satisfies
$$
e^{\frac{\sum_{j=1}^{d}|\lambda_{j}|y_{j}^{2}}{2h}}v_{j}=\omega_{0}(x)
+\mathcal{O}(h^{1-d/4}) \quad 
\textrm{in}~ \mathcal{C}^{\infty}(B(U_{j}^{(p)},\eta))\,,
$$
with
$$
\omega_{0}= \frac{|\lambda_{1}\ldots \lambda_{d}|^{1/4}}{(\pi
h)^{d/4}}
dy_{1}\wedge\ldots \wedge dy_{p}\quad\textrm{along}\quad
\mathcal C_{\text{Unst}}\cap B(U_{j}^{(p)},\eta)\,,
$$
and
$$
\omega_{0}= (-1)^{p(d-p)}\frac{|\lambda_{1}\ldots \lambda_{d}|^{1/4}}{(\pi
h)^{d/4}}
\star(dy_{p+1}\wedge\ldots \wedge dy_{d})\quad\textrm{along}\quad
\mathcal C_{\text{St}}\cap B(U_{j}^{(p)},\eta)\,.
$$
\end{theorem}

\noindent We shall need an integration by part formula adapted from Lemma~4.3.3
in \cite{HeNi} and Lemma~4.3 in \cite{Lep3}.
\begin{lemma}
\label{le.ippHSadapt}
Let $\Omega$ be a regular domain of $f_{a}^{b}$ with boundary made of 
three disjoint pieces $\partial \Omega=\left\{f=a\right\}\sqcup
\left\{f=b\right\}\sqcup \Gamma$\,. Consider the self-adjoint
realization $\Delta_{f,h}^{TND}$ of $\Delta_{f,h}$ given by the form
$$
\mathcal{D}(\omega,\omega')=\langle d_{f,h}\omega\,,\,d_{f,h}\omega'\rangle
+ 
\langle
d_{f,h}^{*}\omega\,,\,d_{f,h}^{*}\omega'
\rangle
$$
with the form domain
$$
\Lambda W^{1,2}_{TND}=\left\{\omega\in \Lambda W^{1,2}(\overline{\Omega})\,;\quad
  \mathbf{t}\omega\big|_{f=a}=0\,,
\; \mathbf{n}\omega\big|_{f=b}=0\,,\;
\omega\big|_{\Gamma}=0\right\}\,,
$$
and the operator domain
$$
D(\Delta_{f,h}^{TND})=
\left\{
\omega\in \Lambda W^{2,2}(\overline{\Omega})\,;\quad
\begin{array}[c]{ll}
   \mathbf{t}\omega\big|_{f=a}=0\,,&
  \mathbf{t}d_{f,h}^{*}\omega\big|_{f=a}=0
\,,\\
  \mathbf{n}\omega\big|_{f=b}=0\,,&
 \mathbf{n}d_{f,h}\omega\big|_{f=b}=0\,,\\
 \omega \big|_{\Gamma}=0
\end{array}
\right\}\;.
$$
Let $\varphi$ be  any Lipschitz function. Then for all $\omega\in \Lambda
W^{1,2}_{TND}$ we have the integration by part formula
\begin{multline}
  \label{eq.ippphi}
  \Real \mathcal{D}(\omega, e^{\frac{2\varphi}{h}}\omega)=
h^{2}\|de^{\frac{\varphi}{h}}\omega\|^{2}+h^{2}\|d^{*}e^{\frac{\varphi}{h}}\omega\|^{2}\\
+
\langle (|\nabla f|^{2}-|\nabla\varphi|^{2}+h\mathcal{L}_{\nabla f}+h
\mathcal{L}_{\nabla f}^{*})e^{\frac{\varphi}{h}}\omega\,,\,
e^{\frac{\varphi}{h}}\omega\rangle\\
+h\left(\int_{f=b}-\int_{f=a}\right)\langle\omega\,,\,\omega\rangle_{\Lambda
  T_{\sigma}^{*}\Omega}e^{\frac{2\varphi(\sigma)}{h}}
\left(\frac{\partial f}{\partial n}\right)(\sigma)~d\sigma\,,
\end{multline}
where $\frac{\partial f}{\partial n}$ is the exterior normal
derivative.
Moreover when $\omega\in D(\Delta_{f,h}^{TND})$ then
$\mathcal{D}(\omega, e^{\frac{2\varphi}{h}}\omega)=\Real\langle
e^{\frac{2\varphi}{h}}\Delta_{f,h}^{TND}\omega, \omega\rangle$\,.
\end{lemma}
\begin{proof}[Proof of Theorem~\ref{th.HSadapt}:]
  First of all, applying the integration by part~\eqref{eq.ippphi} with
  $\varphi=0$ and the local harmonic approximation around critical
  points like in \cite{CFKS}, one obtains that the number of
  eigenvalues in $[0,Ch^{3/2}]$ is $m_{p}$, with no other eigenvalues in
  $(Ch^{3/2}, h/C]$, when $C$ is chosen large
  enough.
The boundary term in \eqref{eq.ippphi} is non
negative because $\frac{\partial f}{\partial n}$ is non positive
(resp. non negative)
on $\left\{f=a\right\}$ (resp. $\left\{f=b\right\}$). 
The IMS localization formula
$-h^{2}\Delta=-\sum_{j}\chi_{j}(h^{2}\Delta)\chi_{j}
-h^{2}\sum_{j}|\nabla \chi_{j}|$ with one $\chi_{j_{a}}$ (resp.
$\chi_{j_{b}}$) localizing around $\left\{f=a\right\}$
(resp. $\left\{f=b\right\}$) shows that eigenfunctions associated
with the
$\mathcal{O}(h^{3/2})$ eigenvalues have asymptotically no mass around
$\left\{f=a\right\}\cup \left\{f=b\right\}$\,. Contrary to \cite{HeNi}
and \cite{Lep3}, there are no generalized critical points at the
boundary and the assumption that $f$ restricted to the boundary is a
Morse function is not necessary here.\\
We construct now a global quasimode $\phi_{j}^{h}$ associated with the
critical point $U_{j}^{(p)}$. Following \cite{HelSj4}\cite{HelSj2}, consider, for a
small constant $\gamma>0$,  the
domain 
$$
\overline{\Omega_{j}}=\overline{f_{a}^{b}}\setminus \cup_{k\neq
  j}B(U_{k}^{(p)},\gamma)\,,
$$ 
with  $\partial
\Omega_{j}=\left\{f=a\right\}\cup \left\{f=b\right\}\cup
\Gamma$ and $\Gamma=\cup_{k\neq
  j}\partial B(U_{k}^{(p)},\gamma)$, and take the self-adjoint realization
$\Delta_{f,h}^{TND (p)}$ of Lemma~\ref{le.ippHSadapt} acting on
$p$-forms. 
It admits a single
eigenvalue $\mu_{j}^{h}$ which is $\mathcal{O}(h^{3/2})$ (with the rest
of the spectrum in $[h/C,+\infty)$) and take $\phi^h_{j}$ a
normalized eigenvector associated with $\mu^h_{j}$: 
$$
\|\phi_{j}^{h}\|_{\Lambda L^{2}}=1\quad,\quad
\Delta_{f,h}^{TND (p)}\phi^h_{j}=\mu_{j}^{h}\phi^h_{j}\,.
$$
Applying Lemma~\ref{le.ippHSadapt},
in the spirit of \cite{DiSj}~pp.~49--55,
with
$$\omega=\phi_{j}^{h}\ ,\ 
\varphi_{\varepsilon}(x)=(1-\varepsilon)d_{Ag}\left(x,B(U_{j}^{(p)},\varepsilon)\right)\ ,\ |\nabla \varphi_{\varepsilon}|\leq
(1-\varepsilon)|\nabla f|\,,$$
 gives  
$$
\|e^{\frac{\varphi_{\varepsilon}}{h}}\phi^h_{j}\|_{\Lambda W^{1,2}}= 
\mathcal{O}_{\varepsilon}(\frac{1}{h})\,,
$$
where the subscript $_{\varepsilon}$ recalls that the factor of
$\frac{1}{h}$ depends on the parameter $\varepsilon>0$\,.
By elliptic regularity up to the boundary of $\Delta_{f,h}^{TND (p)}$,
$\phi^h_{j}$ is  $\mathcal{C}^{\infty}$ in $\overline{\Omega_{j}}$\,.
 The differential operator 
$e^{\frac{\varphi}{h}}\Delta_{f,h}e^{-\frac{\varphi}{h}}$ equals 
$$
h^{2}(dd^{*}+d^{*}d)+|\nabla f|^{2}-|\nabla \varphi|^{2}+h
(\mathcal{L}_{\nabla \varphi}-\mathcal{L}^*_{\nabla \varphi}+
\mathcal{L}_{\nabla f}+ \mathcal{L}_{\nabla f}^{*})
$$
where the last part is a first order differential operator. 
With the boundary conditions, the form
$u_{j}^{h}=e^{\frac{\varphi_{\varepsilon}}{h}}\phi_{j}^{h}$ satisfies the system
$$
\left\{
\begin{array}[c]{l}
(dd^{*}+d^{*}d)u_{j}^{h}=r_{j}
\\
\mathbf{t}u_{j}^{h}\big|_{f=a}=0\quad\mathbf{t}d^{*}u_{j}^{h}\big|_{f=a}=\varrho_{j,a}\\
\mathbf{n}u_{j}^{h}\big|_{f=b}=0\quad\mathbf{n}d u_{j}^{h}\big|_{f=b}=\varrho_{j,b}\\
u_{j}^{h}\big|_{\Gamma}=0
\end{array}
\right.
$$
where  $\|r_{j}\|_{\Lambda L^{2}}$, 
$\|\varrho_{j,*}\|_{\Lambda W^{1/2,2}}$ are $\mathcal{O}(\frac{1}{h^3})$. This
provides a $\mathcal{O}(\frac{1}{h^3})$ estimate for $\|u_{j}^{h}\|_{\Lambda
  W^{2,2}}$ and bootstraping gives 
$$
\forall \alpha \in \nz^{d}\,, \forall x\in \overline{\Omega_{j}}\,,\quad
|\partial_{x}^{\alpha}\phi_{j}^{h}(x)|=
\tilde{\mathcal{O}}(e^{-\frac{\varphi_{\varepsilon}(x)}{h}})\,.
$$
Since this holds for all $\varepsilon>0$, the definition of
$\tilde{\mathcal{O}}$ provides the same result with $\varepsilon=0$. It means 
the following estimate holds, with $\varphi(x):=\varphi_{0}(x)=d_{Ag}(x,U_{j}^{(p)})$:
$$
\forall \alpha \in \nz^{d}\,, \forall x\in
\overline{\Omega}_{j}\,,\quad
|\partial_{x}^{\alpha}\phi^h_{j}(x)|=\tilde{\mathcal{O}}(e^{-\frac{\varphi(x)}{h}})\,.
$$
The differential  $d_{f,h}\phi^h_{j}$  solves in $\Omega_{j}$ the
differential equation
$$
\Delta_{f,h}(d_{f,h}\phi^h_{j})= d_{f,h}(\Delta_{f,h}\phi^h_{j})= \mu_{j}^{h}d_{f,h}\phi^h_{j}\,.
$$
The same argument as the one for Lemma~\ref{le.wtn} leads to the fact
that $d_{f,h}\phi^h_{j}$ satisfies the same boundary conditions as
$\phi^h_{j}$ on $\left\{f=a\right\}\cup
\left\{f=b\right\}$\,. Consider now the domain
$\overline{\Omega_{j}'}=\overline{\Omega_{j}}\setminus \cup_{U\in
  \mathcal{U}^{(p+1)}}B(U,\gamma)$, note $\mathcal{V}_{j}=\mathcal{U}^{(p+1)}\cup\,\mathcal{U}^{(p)}\setminus
  \left\{U_{j}^{(p)}\right\}$, and work with the associated
$\Delta_{f,h}^{TND (p+1)}$\,.
The form ${u'}_{j}^{h}=\chi_{\gamma}d_{f,h}\phi^h_{j}$, where $\chi_{\gamma}\in
\mathcal{C}^{\infty}(\overline{\Omega_{j}'})$
vanishes in $\cup_{U\in
  \mathcal{V}_{j}}B(U,2\gamma)$ and equals $1$
outside
$\cup_{U\in \mathcal{V}_{j}}B(U,3\gamma)$, belongs to
$D(\Delta_{f,h}^{TND (p+1)})$ and solves
$$
\Delta_{f,h}^{TDN (p+1)}({u'}_{j}^{h})=\mu_{j}^{h}{u'}_{j}^{h} + r'_{j}\,,
$$
with $\supp r'_{j}\subset \cup_{U\in
\mathcal{V}_{j}}B(U,3\gamma)$ and 
$$
|\partial_{x}^{\alpha}r'_{j}(x)|=\tilde{\mathcal{O}}(e^{\frac{-\min_{U\in
    \mathcal{V}_{j}}
d_{Ag}(U,U_{j}^{(p)})+c\gamma}{h}})\quad\ (\text{with }c>0)\,.
$$
With our choice of $\overline{\Omega_{j}'}$, $\Delta_{f,h}^{TND
  (p+1)}$ has no eigenvalue in $[0,h/C]$ and the same analysis as
above
leads to
$$
\forall \alpha\in \nz^{d}\,,
\forall x\in \overline{\Omega_{j}'}\,,\quad
|\partial_{x}^{\alpha}d_{f,h}\phi_{j}^{h}(x)|
=\tilde{\mathcal{O}}
\left(
e^{\frac{-\min_{U \in
\mathcal{V}_{j}}d_{Ag}(U_{j}^{(p)},U)-d_{Ag}(U,x)+c\gamma}{h}}
\right)\,,
$$
where the previous estimates extend the result  to  all $\overline{\Omega_{j}}$\,.
After changing $\mathcal{U}^{(p+1)}$ into
$\mathcal{U}^{(p-1)}$, a similar result holds for
$d_{f,h}^{*}\phi_{j}^{h}$\,.  A simple computation now gives 
$\mu^h_{j}=\mathcal{D}(\phi^h_{j},\phi^h_{j})=
\tilde{\mathcal{O}}(e^{-2\frac{\min_{U\in
  \mathcal{V}_{j}}d_{Ag}(U,U_{j})+c\gamma}{h}})$\,.\\
Now let us work with $\Delta_{f,h}^{TN (p)}$ on
$\overline{f_{a}^{b}}$\,. Consider the cut-off $\theta_{j,\gamma}\in
\mathcal{C}^{\infty}(\Omega_{j})$ which vanishes in $\cup_{U\in
  \mathcal{U}^{(p)}\setminus \left\{U_{j}^{(p)}\right\}}B(U,2\gamma)$ and equals $1$ in  $\cup_{U\in
  \mathcal{U}^{(p)}\setminus \left\{U_{j}^{(p)}\right\}}B(U,3\gamma)$, and set
$$
\psi^h_{j}=\theta_{j,\gamma}\phi^h_{j}\,.
$$
These $m_{p}$ vectors belong to $D(\Delta_{f,h}^{TN(p)})$ and satisfy
\begin{eqnarray*}
&&\Delta_{f,h}^{TN(p)}\psi^h_{j}=\mu^h_{j}\psi^h_{j}+r_{j}\\
\text{with}
&& \mu^h_{j}=\mathcal{O}(e^{-\frac{C}{h}})\,,\\
&&
|r_{j}(x)|=\mathcal{O}(e^{-\frac{d_{Ag}(x,U_{j})}{h}})\quad,\quad
\supp r_{j}\subset \cup_{U\in
  \mathcal{U}^{(p)}\setminus\left\{U_{j}^{(p)}\right\}}B(U,3\gamma)\,,\\
\text{and}
&&\left(\langle \psi_{j}^h\,,\, \psi_{k}^h\rangle\right)_{1\leq j,k\leq m_{p}}
=\Id+\mathcal{O}(e^{-\frac{C}{h}})\,.
\end{eqnarray*}
while $\Delta_{f,h}^{TN(p)}$ has only $m_{p}$ eigenvalues in $[0,Ch^{3/2}]$.
The Proposition~4.1 in \cite{Hel}  implies 
$$
\psi^h_{j}-1_{[0,Ch^{3/2}]}(\Delta_{f,h}^{TN (p)})\psi^h_{j}=\mathcal{O}(e^{-\frac{C}{h}})\,,
$$
and we set $u_{j}^{h}=1_{[0,Ch^{3/2}]}(\Delta_{f,h}^{TN (p)})\psi^h_{j}$\,.
The min-max principle applied with the $\psi^h_{j}$'s also implies that
the eigenvalues of $\Delta_{f,h}^{TN(p)}$ in $[0,Ch^{3/2}]$ are actually 
exponentially small. 
With the integration contour $\mathcal C_{h}=\left\{z\in\cz,
  |z|=h^{3/2}\right\}$, write
$$
u_{j}^{h}-\psi^h_{j}=\frac{1}{2i\pi}
\int_{\mathcal C_{h}}(z-\mu^h_{j})^{-1}(z-\Delta_{f,h}^{TN(p)})^{-1}r_{j}~dz\,. 
$$
The resolvent estimates of Proposition~2.2.5 in \cite{HelSj2} can be
carried over to our boundary problem thanks to
Lemma~\ref{le.ippHSadapt} and elliptic regularity up to the boundary. 
With the estimates and support condition on $r_{j}$, they
lead to
\begin{multline*}
\forall \alpha\in \nz^{d}, \forall x\in \overline{f_{a}^{b}}\,,
|\partial_{x}^{\alpha}\omega_{z}(x)|
=\tilde{\mathcal{O}}\left(e^{\frac{-\min_{U\in
        \mathcal{U}^{(p)}\setminus\left\{U_{j}^{(p)}\right\}}
      d_{Ag}(U_{j},U)+d_{Ag}(U,x)+c\gamma}{h}}\right)\\
=\tilde{\mathcal{O}}\left(e^{\frac{-d_{Ag}(U_{j},x)+c\gamma}{h}}\right)
\end{multline*}
when $\omega_{z}=[(z-\Delta_{f,h}^{TN(p)})^{-1}r_{j}]$ and $z\in
\mathcal C_{h}$\,.\\
With the estimates on $\psi_{j}=\theta_{j,h}\phi_{j}^{h}$, this leads to
$$
\forall \alpha\in \nz^{d}, \forall x\in \overline{f_{a}^{b}}\,,\quad
|\partial_{x}^{\alpha}u_{j}^{h}|=\tilde{\mathcal{O}}(e^{\frac{-d_{Ag}(x,U_{j
   })+c\gamma}{h}})\,,
$$
and we take
$$
v^h_{j}=\|u_{j}^{h}\|^{-1}u_{j}^{h}=(1+\mathcal{O}(e^{-C/h}))u_{j}^{h}\,.
$$
The estimates for $d_{f,h}v^h_{j}$ (resp. $d_{f,h}^{*}v^h_{j}$) are obtained
after writing the equation for $d_{f,h}(u_{j}^{h}-\psi^h_{j})$
(resp. $d_{f,h}^{*}(u_{j}^{h}-\psi_{j}^h)$) and using the resolvent
estimates for $\Delta_{f,h}^{TN(p+1)}$ (resp. $\Delta_{f,h}^{TN(p-1)}$)\,.\\
Finally, 
since these estimates hold for any $\gamma>0$, the definition of
$\tilde{\mathcal{O}}$ provides the same result with $\gamma=0$.\\
The rest of the proof  of Theorem~\ref{th.HSadapt} is a direct consequence of Theorem~2.5
of \cite{HelSj4} and therein related WKB construction.
\end{proof}

\subsection{An important remark}
\label{se.remark}

It is clear that the results stated for $\overline{f_{a}^{b}}$ hold
when $a=-\infty$ or $b=+\infty$, that is when one boundary is empty.\\
Another variation on it consists in deforming homotopically
$\left\{f=a\right\}$ (resp. $\left\{f=b\right\}$) while preserving the
sign conditions $\frac{\partial f}{\partial n}<0$
(resp. $\frac{\partial f}{\partial n}>0$).

\section{Barannikov-Morse chain complex and construction of accurate global quasimodes}
\label{se.barquasi}

\subsection{Properties of quasimodes associated with lower and upper critical points}
\label{se.qualowup}

Consider the operator $\Delta_{f,h}^{TN}$ defined on
$\overline{f_{a}^{b}}$ and set
 $F^{(p)}=\ima 1_{[0,h^{3/2}]}(\Delta_{f,h}^{TN (p)})$ for $p\in
 \left\{0,\ldots, d\right\}$, 
$F=\oplus_{p}F^{(p)}$,  $\beta=d_{f,h}\big|_{F}$ and
$\beta^{*}=d_{f,h}^{*}\big|_{F}$\,.
According to Section~\ref{se.wit}, $\dim F^{(p)}=m_{p}$ 
and the chain complex associated with $\beta$
$$
\xymatrix{
0\ar[r]&
F^{(0)}
\ldots
F^{(p-1)}
\ar[r]^-{\beta^{(p-1)}}\ar@<1ex>[l]&
F^{(p)}
\ar[r]^-{\beta^{(p)}}\ar@<1ex>[l]^-{\beta^{(p-1)*}}&
F^{(p+1)}
\ldots
F^{(d)}\ar[r]\ar@<1ex>[l]^-{\beta^{(p)*}}&0
\ar@<1ex>[l]
}
$$
has the homology group
$H^{*}(\overline{f_{a}^{b}},\left\{f=a\right\})$ dual to
$H_{*}(f^{b},f^{a})=H_{*}(\overline{f_{a}^{b}},\left\{f=a\right\})$\,.
Remember also  that $F^{(p)}$ admits the orthogonal decompositions
\begin{eqnarray*}
F^{(p)}&=&\ker \Delta_{f,h}^{TN(p)}\mathop{\oplus}^{\perp}\ima \beta^{(p-1)}
\mathop{\oplus}^{\perp}\ima \beta^{(p)*}\\
 &=& \ker \beta^{(p)}
\mathop{\oplus}^{\perp}\ima \beta^{(p)*}
=
\ima \beta^{(p-1)}\mathop{\oplus}^{\perp}\ker \beta^{(p-1)*}
\end{eqnarray*}
and that $F^{(p)}$ admits an almost orthonormal basis
$\left\{v_{U}, U\in \mathcal{U}^{(p)}\right\}$ fulfilling the
properties of Theorem~\ref{th.HSadapt}.
The orthogonal projection on any subspace $G$ of the above orthogonal
decomposition will be denoted by $\Pi_{G}$\,.
\begin{proposition}
\label{pr.quasinotup}
  Assume that $U\in \mathcal{U}^{(p)}$ is not an upper critical
  points in $f_{a}^{b}$, then $v_{U}$ is almost orthogonal to $\ima
  \beta$:
$$
\|\Pi_{\ima \beta^{(p-1)}}v_{U}\|=\mathcal{O}(e^{-\frac{C_{\eta}}{h}})\,,
$$
where the constant $C_{\eta}$ depends on the small radius $\eta$
fixed by the geometry in Theorem~\ref{th.HSadapt}.
\end{proposition}
\begin{proof}
  When $U\in \mathcal{U}^{(p)}$
 with $f(U)=c_{U}$ is not an upper critical point, it means that the mapping
$H_{p}(f^{c_{U}+\varepsilon})\to
H_{p}(f^{c_{U}+\varepsilon},f^{c_{U}-\varepsilon})$ does not
vanish. Since $H_{p}(f^{c_{U}+\varepsilon}, f^{c_{U}-\varepsilon})$ is
one dimensional and is generated by $e^{p}_{U}$ the unstable manifold
for $-\nabla f$ leaving $U$ and restricted to
$f_{c_{U}-\varepsilon}$,
there exists a cycle $C_{U}^{p}$ in $f^{c_{U}+\varepsilon}$, or a
cycle in $M$ supported in $f^{c_{U}+\varepsilon}$, of which the
restriction to $f_{c_{U}-\varepsilon}^{c_{U}+\varepsilon}$ is
$e_{U}^{p}$\,.
We choose $\varepsilon=\varepsilon_{\eta}>0$ so that $e_{U}^{p}$ is contained in the ball
$B(U,\eta)$ of Theorem~\ref{th.HSadapt}. In adapted Morse coordinates,
$e_{U}^{p}$ equals up to the orientation
$$
\left\{y_{p+1}=\ldots=y_{d}=0,\quad \frac{1}{2}
  \sum_{j=1}^{p}|\lambda_{j}|y_{j}^{2}< \varepsilon_{\eta} \right\}\,,
$$
while $e^{\frac{f-c_{U}}{h}}v_{U}$ equals
$$
|\lambda_{1}\ldots \lambda_{d}|^{1/4}(\pi
h)^{-d/4}e^{-\frac{\sum_{j=1}^{p}|\lambda_{j}|y_{j}^{2}}{h}}(\omega_{0}(y)+h\omega'(y,h))
 +
\mathcal{O}\left(e^{-\frac{C_{\eta}}{h}}\right)
$$
with $\omega'$ bounded in $\mathcal{C}^{\infty}(B(U,\eta))$ and
$$
\omega_{0}=
dy_{1}\wedge\ldots \wedge dy_{p} \quad\text{along}\quad
\{y_{p+1}=\cdots=y_{d}=0\}\cap B(U,\eta)\,.
$$
Decompose $v_{U}$ according to
\begin{eqnarray*}
  && v_{U}=v_{U}'+v_{U}''=\Pi_{\ima \beta^{(p-1)}}v_{U}+\Pi_{\ker
    \beta^{(p-1)*}}v_{U}\,,\\
&& v_{U}'=\sum_{U'\in \mathcal{U}^{(p)}}t_{U'}v_{U'}\quad,\quad 
v_{U}''=\sum_{U'\in \mathcal{U}^{(p)}}s_{U'}v_{U'}\,.
\end{eqnarray*}
The decomposition $v_{U}=v_{U}'+v_{U}''$ is orthogonal
$$
\|v_{U}'\|^{2}+\|v_{U}''\|^2=1\,.
$$
Meanwhile, the exponential decay estimates of $(v_{U'})_{U'\in
  \mathcal{U}^{(p)}}$ stated in Theorem~\ref{th.HSadapt} provide the
almost orthogonality
\begin{eqnarray*}
&& \sum_{U'}|t_{U'}|^{2}\leq 1+\mathcal{O}(e^{-C/h})\quad,
\sum_{U'}|s_{U'}|^{2}\leq 1+ \mathcal{O}(e^{-C/h})\,,\\
&&
t_{U}+s_{U}=1+\mathcal{O}(e^{-C/h}) \quad
\text{and}\quad t_{U'}+s_{U'}=\mathcal{O}(e^{-C/h})
\quad\text{for}\; U'\neq U\,.
\end{eqnarray*}
All the $v_{U'}$ have $\tilde{\mathcal{O}}(1)$ estimates in
$\mathcal{C}^{\infty}(\overline{f_{a}^{b}}; \bigwedge^{p} T^{*}M)$
and the support conditions on $e_{U}^{p}$ and $C_{U}^{p}$ give
$$
\int_{e_{U}^{p}}e^{\frac{f-c_{U}}{h}}v_{U}=\int_{C_{U}^{p}}
e^{\frac{f-c_{U}}{h}}v_{U}+\mathcal{O}(e^{-\frac{C_{\eta}}{h}})\,.
$$
By using $v_{U}'=d_{f,h}\omega$ we get
$$
\int_{e_{U}^{p}}e^{\frac{f-c_{U}}{h}}v_{U}
= 
h\int_{C_{U}^{p}}d\left(e^{\frac{f-c_{U}}{h}}\omega\right)
+\int_{C_{U}^{p}}e^{\frac{f-c_{U}}{h}}v_{U}''\,,
$$
and finally with $\partial C_{U}^{p}=0$,
$$
\int_{e_{U}^{p}}e^{\frac{f-c_{U}}{h}}v_{U}=\sum_{U'\in
  \mathcal{U}^{(p)}}s_{U'}
\int_{C_{U}^{p}}e^{\frac{f-c_{U}}{h}}v_{U'}
=
\sum_{U'\in
  \mathcal{U}^{(p)}}s_{U'}
\left(\int_{e^{p}_{U}}+\int_{C_{U}^{p}\setminus e^{p}_{U}}\right)e^{\frac{f-c_{U}}{h}}v_{U'}\,.
$$
With $s_{U'}=\mathcal{O}(1)$,
$e^{\frac{f-c_{U}}{h}}=\mathcal{O}(e^{-\frac{C_{\eta}}{h}})$ on
$C_{U}^{p}-e^{p}_{U}$ and $|v_{U'}|=\tilde{\mathcal{O}}(1)$ the second
integral gives an exponentially small term.
When $U'\neq U$ the exponential decay estimate of
Theorem~\ref{th.HSadapt} imply that
$s_{U'}\int_{e^{p}_{U}}e^{\frac{f-c_{U}}{h}}v_{U'}$ is
$\mathcal{O}(e^{-C/h})$\,.\\
We have proved
$$
\int_{e^{p}_{U}}e^{\frac{f-c_{U}}{h}}v_{U}=s_{U}\int_{e^{p}_{U}}e^{\frac{f-c_{U}}{h}}v_{U}+ \mathcal{O}(e^{-\frac{C_{\eta}}{h}})\,.
$$
But a direct calculation in the Morse coordinates gives
\begin{eqnarray*}
\int_{e^{p}_{U}}e^{\frac{f-c_{U}}{h}}v_{U}
&=&\left(1+\mathcal O(h)\right)\!\int_{\sum_{j=1}^{p}|\lambda_{j}|y_{j}^{2}<
2\varepsilon_{\eta}}\!\frac{|\lambda_{1}\ldots \lambda_{d}|^{1/4}}{(\pi
 h)^{d/4}}e^{-\frac{\sum_{j=1}^{p}|\lambda_{j}|y_{j}^{2}}{h}}dy_{1}\ldots dy_{p}
\\
&=&
\frac{|\lambda_{p+1}\ldots \lambda_{d}|^{1/4}}{|\lambda_{1}\ldots
  \lambda_{p}|^{1/4}}(\pi h)^{(2p-d)/4}\left(1+\mathcal O(h)\right)\,.
\end{eqnarray*}
This proves $s_{U}=1+ \mathcal{O}(e^{-\frac{C_{\eta}}{h}})$ while all
the other coefficients are $\mathcal{O}(e^{-\frac{C_{\eta}}{h}})$\,.
\end{proof}
By duality $f\to -f$, other results can be deduced.
\begin{proposition}
  \label{pr.quasiclass}
When $U\in \mathcal{U}^{(p)}$ is not a lower critical point in $f_{a}^{b}$, $v_{U}$
is almost orthogonal to  $\ima \beta^{(p)*}$:
$$
\|\Pi_{\ima \beta^{(p)*}}v_{U}\|=\mathcal{O}(e^{-\frac{C_{\eta}}{h}})\,.
$$
When $U\in \mathcal{U}^{(p)}$ is an homological critical point in
$f_{a}^{b}$, $v_{U}$ is exponentially close to $\ker \Delta_{f,h}^{(p)}$: 
$$
 \|v_{U}-\Pi_{\ker \Delta_{f,h}^{TN (p)}}v_{U}\|=\mathcal{O}(e^{-\frac{C_{\eta}}{h}})\,.
 $$
Finally, when $U\in \mathcal{U}^{(p)}$ is an upper (resp. a lower) critical point,
$v_{U}$ is exponentially close to $\ima \beta^{(p-1)}$
(resp. $\ima\beta^{(p)*}$):
$$
 \|v_{U}-\Pi_{\ima \beta^{(p-1)}}v_{U}\|=\mathcal{O}(e^{-\frac{C_{\eta}}{h}})
\quad(\text{resp.}~ \|v_{U}-\Pi_{\ima
  \beta^{(p)*}}v_{U}\|=\mathcal{O}(e^{-\frac{C_{\eta}}{h}}))\,.
$$
\end{proposition}
\begin{proof}
  The first statement is dual to the one of
  Proposition~\ref{pr.quasinotup}.
For the second one it suffices to notice that homological critical
points are neither upper nor lower critical points.
For the last one it suffices to notice that the number of homological
critical points equals the dimension of $\ker
\Delta_{f,h}^{TN(p)}$\,. Hence the set of $\Pi_{\ker \Delta_{f,h}^{TN
    (p)}}v_{U'}$ when $U'$ ranges over the homological critical points,
is an almost orthonormal basis of $\ker \Delta_{f,h}^{TN (p)}$\,. If
$U$ is an upper critical point, $\Pi_{\ker \beta^{(p)}}v_{U}=v_{U}+
\mathcal{O}(e^{-\frac{C_{\eta}}{h}})$ is almost orthogonal to the
$v_{U'}$ and therefore to $\ker \Delta_{f,h}^{TN (p)}$\,. We deduce
$$
v_{U}=\Pi_{\ima \beta^{(p-1)}}v_{U}+ \Pi_{\ker
  \Delta_{f,h}^{TN(p)}}v_{U}+ \mathcal{O}(e^{-\frac{C_{\eta}}{h}})
= \Pi_{\ima \beta^{(p-1)}}v_{U}+ \mathcal{O}(e^{-\frac{C_{\eta}}{h}})\,.
$$
\end{proof}

\subsection{Construction of accurate global quasimodes}
\label{se.globquasi}

We now define the global quasimodes for $\Delta_{f,h}$ on $M$
 which will be used in our computations.
 \begin{itemize}
 \item When $U$ is an homological critical point, take simply
$$
\omega_{U}=\Pi_{\ker \Delta_{f,h}}v_{U}\,,
$$
where $v_{U}$ is the form defined in Theorem~\ref{th.HSadapt} with 
$a=-\infty$ and $b=+\infty$\,.
\item When $U$ is an upper critical point, take
$$
\omega_{U}=\Pi_{\ima \beta}v_{U}
$$
where $v_{U}$ is the form defined in Theorem~\ref{th.HSadapt} with 
$a=-\infty$ and $b=+\infty$\,.
\item When $U\in \mathcal{U}^{(p)}$ is a lower critical point there
  exist $U_{1}\in \mathcal{U}^{(p+1)}$ with
  $f(U_{1})=c_{1}$ such that $\partial U_{1}=U$\,. In
  $f^{c_{1}-\varepsilon}$, $U$ becomes an homological critical point.
We take 
$$
\omega_{U}=1_{[0,h^{3/2}]}(\Delta_{f,h})\chi_{\varepsilon}\tilde{v}_{U}
$$
where $\tilde{v}_{U}$ is now the form defined in Theorem~\ref{th.HSadapt} with 
$a=-\infty$ and $b=c_{1}-\varepsilon$ while $\Delta_{f,h}$ is the
operator defined on all $M$\,. The function
$\chi_{\varepsilon}$ vanishes in $f_{c_{1}-\frac32\varepsilon}$ and equals
$1$ in $f^{c_{1}-2\varepsilon}$\,. The value of the parameter
$\varepsilon$ 
will be specified further according to $\eta$\,.
\end{itemize}

\subsection{Computation of the matrix of $d_{f,h}$}
\label{se.computdf}

We work with the basis $(\omega_{U})_{U\in \mathcal{U}}$ constructed
before and we will denote by $\mathcal U_{H} $,
$\mathcal U_{L}$,  $\mathcal U_{U}$ the sets of homological,
lower and upper critical points of $f$, and $\mathcal U_{H}^{(p)} $,
$\mathcal U_{L}^{(p)}$,  $\mathcal U_{U}^{(p)}$
their respective intersection with $\mathcal U^{(p)} $,
the set of critical points of $f$ with index $p$.
\\


\begin{proposition}
\label{pr.computdf}
When $U_{0}$ belongs to $\mathcal U_{U}^{(p)}\cup \mathcal U_{H}^{(p)}$,
then for any $U'\in \mathcal U^{(p+1)}$,
\begin{equation}
\label{eq.computdf1}
\langle \omega_{U'}\,|\,d_{f,h}\omega_{U_{0}}\rangle =0\,.
\end{equation}

\noindent
When $U_{0}$ belongs to $\mathcal U_{L}^{(p)}$, let $U_{1}$ denote 
 the upper critical point with index 
$p+1$
s.t. $\partial_{\mathcal B}(U_{1})
=U_{0}$. Then there exists a real constant $C>0$
and a homological constant $\kappa=\kappa(U_{1})\neq 0$ such that
for $U'\in\mathcal U^{(p+1)}$: 
\begin{eqnarray}
\label{eq.computdf2}
&\text{If } U'\neq U_{1}\,,&\  
\langle \omega_{U'}\,|\, d_{f,h}\omega_{U_{0}}\rangle =
\mathcal O(e^{-\frac{f(U_{1})-f(U_{0})+C}{h}})\,,\\
\label{eq.computdf3}
&\text{If } U'= U_{1}\,,&\ \langle \omega_{U_{1}}\,|\, d_{f,h}\omega_{U_{0}}\rangle
=\pm \kappa A(h)
e^{-\frac{f(U_{1})-f(U_{0})}{h}}
(1+\mathcal{O}(h))\,.
\end{eqnarray}
Moreover the prefactor $A(h)$ is given  by the formula
\begin{equation}
\label{eq.computdf4}
A(h)=(\frac h\pi)^\frac12\frac{|\lambda^1_{1}\cdots\lambda^1_{p+1}|^\frac14}{
|\lambda^1_{p+2}\cdots\lambda^1_{d}|^\frac14}
\frac{|\lambda^0_{p+1}\cdots\lambda^0_{d}|^\frac14}{
|\lambda^0_{1}\cdots\lambda^0_{p}|^\frac14}\,,
\end{equation}
where $\lambda_{1}^{\ell}<\cdots<\lambda_{p+\ell}^{\ell}<0<
\lambda_{p+\ell+1}^{\ell}<\cdots<\lambda_{d}^{\ell}$
denote the eigenvalues of $\Hess f(U_{\ell})$,
for $\ell\in\{0,1\}$.
\end{proposition}

\noindent 
The rest of this section is devoted to the proof of this proposition. We are first
going to prove the relations~\eqref{eq.computdf1} and \eqref{eq.computdf2}, then,
in order to prove \eqref{eq.computdf3} and \eqref{eq.computdf4}, we will in a first
time work with a metric which is locally Euclidean around the critical points of $f$
before showing that it remains valid for a general Riemannian metric.

\begin{proof}[Proof of equations \eqref{eq.computdf1} and \eqref{eq.computdf2}]
When $U_{0}$ is an upper critical point or a homological critical point,
then the definition of $\omega_{U_{0}}$ says
$$
d_{f,h}\omega_{U_{0}}=\beta \omega_{U_{0}}=0\,,
$$
which yields equation~\eqref{eq.computdf1}.\\
Let us now compute $\langle \omega_{U'}\,,\, d_{f,h}\omega_{U_{0}}\rangle$
for $U'\in \mathcal U^{(p+1)}$ when $U_{0}\in \mathcal{U}^{(p)}$ is a lower critical point with critical
value $c_{0}$. Let $U_{1}\in \mathcal{U}^{(p+1)}$ be the upper
critical point with critical value $c_{1}$ such that $\partial_{\mathcal B} U_{1}=U_{0}$\,.
The commutation
$d_{f,h}1_{[0,h^{2}]}(\Delta_{f,h})=1_{[0,h^{2}]}(\Delta_{f,h})d_{f,h}$
gives
\begin{eqnarray*}
 \langle \omega_{U'}\,,\, d_{f,h}\omega_{U_{0}}\rangle
&=& \langle
  \omega_{U'}\,,\,d_{f,h}1_{[0,h^{2}]}(\Delta_{f,h}^{(p)})\chi_{\varepsilon}\tilde{v}_{U_{0}}\rangle
=\langle
  \omega_{U'}\,,\,d_{f,h} \chi_{\varepsilon}\tilde{v}_{U_{0}}\rangle
\\
&=& 
\langle v_{U'}\,,\, d_{f,h}\chi_{\varepsilon}\tilde{v}_{U_{0}}\rangle
+
\langle\omega_{U'}- v_{U'}\,,\, d_{f,h}\chi_{\varepsilon}\tilde{v}_{U_{0}}\rangle
\end{eqnarray*} 
Since $d_{f,h}\tilde{v}_{U_{0}}=0$ in $\supp\nabla \chi_{\varepsilon}$, we have
$$
d_{f,h}\chi_{\varepsilon}\tilde{v}_{U_{0}}=hd\chi_{\varepsilon}\wedge \tilde{v}_{U_{0}}\,.
$$
Since $d\chi_{\varepsilon}$ is supported in
$f_{c_{1}-2\varepsilon}^{c_{1}-\frac32\varepsilon}$ and 
$$
|\tilde{v}_{U_{0}}(x)|=\tilde{\mathcal{O}}(e^{-\frac{d_{Ag}(x,U_{0})}{h}})=\mathcal{O}(e^{-\frac{c_{1}-c_{0}-C\varepsilon}{h}})\quad
\text{for}\; x\in \supp\nabla \chi_{\varepsilon}\,,
$$
the remainder term $\langle \omega_{U'}-v_{U'}\,,\,
d_{f,h}\chi_{\varepsilon}\tilde{v}_{U_{0}}\rangle$ is bounded by
$$
\|\omega_{U'}-v_{U'}\|\mathcal{O}(e^{-\frac{c_{1}-c_{0}-C\varepsilon}{h}})\,.
$$
When $U'$ is not a lower critical point,
the relation $\|\omega_{U'}-v_{U'}\|=\mathcal{O}(e^{-\frac{C}{h}})$ comes from
Proposition~\ref{pr.quasiclass}. When $U'$ is a lower critical point,
simply
note that both terms
of the r.h.s. in 
\begin{equation}
  \label{eq.compwvlower}
\|\omega_{U'}-v_{U'}\|\leq \|\omega_{U'}-\chi_{\varepsilon}\tilde{v}_{U'}\|+\|\chi_{\varepsilon}\tilde{v}_{U'}-v_{U'}\|
\end{equation}
are $\mathcal{O}(e^{-\frac{C}{h}})$\,. Actually, the estimate for the
second term is obtained after comparing $v_{U'}$ and
$\tilde{v}_{U'}$ with the single eigenmode of a Dirichlet realization
of $\Delta_{f,h}^{(p+1)}$ in $B(U',\eta_{0})$, again by following \cite{Hel}.\\
Hence we have proved
$$
 \langle \omega_{U'}\,,\, d_{f,h}\omega_{U_{0}}\rangle
= \langle v_{U'}\,,\, hd\chi_{\varepsilon}\wedge\tilde{v}_{U_{0}}\rangle + \mathcal{O}(e^{-\frac{c_{1}-c_{0}+C}{h}})
$$
when $\varepsilon>0$ is chosen small enough.\\
If $U'\neq U_{1}$ the exponential decay of $v_{U'}$,
$$
|v_{U'}(x)|=\tilde{\mathcal{O}}(e^{-\frac{d_{Ag}(x,U')}{h}})=
\mathcal{O}(e^{-\frac{|c_{1}-f(U')|-C\varepsilon}{h}})\quad\text{for}\;
x\in \supp \nabla \chi_{\varepsilon}\,,
$$
leads to
$$
\langle \omega_{U'}\,,\, d_{f,h}\omega_{U_{0}}\rangle =
\mathcal O(e^{-\frac{c_{1}-c_{0}+C}{h}})\,,
$$
and equation~\eqref{eq.computdf2} is proved.
\end{proof}

\subsubsection{Proof of Proposition~\ref{pr.computdf}
when the metric is Euclidean in some adapted Morse coordinates}

Let us check equations~\eqref{eq.computdf3}-\eqref{eq.computdf4},
when the metric is  Euclidean in some local  adapted Morse coordinates
for $f$ (around each critical point).\\
Note first that we have already proved, for a general metric, the following result,
$$
\forall U'\in\mathcal U^{(p+1)}\,,\quad 
 \langle \omega_{U'}\,,\, d_{f,h}\omega_{U_{0}}\rangle
= \langle v_{U'}\,,\, hd\chi_{\varepsilon}\wedge\tilde{v}_{U_{0}}\rangle +
 \mathcal{O}(e^{-\frac{c_{1}-c_{0}+C}{h}})\,,
$$
where $C$ is a positive constant. According to the choice of
$\chi_{\varepsilon}$, the first term of the right-hand side vanishes
when $\partial_{B} U_{1}\neq U_{0}$\,. Thus, we can focus on the term
$\langle v_{U_{1}}\,,\, hd\chi_{\varepsilon}\wedge
\tilde{v}_{U_{0}}\rangle$, when $U_{1}\in\mathcal U^{(p+1)}$
satisfies $\partial_{\mathcal B}(U_{1})=U_{0}$. \\
\\
In the ball $B(U_{1},\eta)$, we use the above adapted Morse coordinates $(y',y'')$
with $y'=(y_{1},\ldots, y_{p+1})$, $y''=(y_{p+2},\ldots,
y_{d})$, and
$f(y)-c_{1}=\frac{1}{2}\sum_{j=1}^{d}\lambda^1_{j}y_{j}^{2}$\,. 
The parameter $\varepsilon>0$ is chosen according to
$\eta>0$ so that
$$
f_{c_{1}-2\varepsilon}^{c_{1}-\frac32\varepsilon}\cap B(U_{1},\eta)\neq \emptyset\,.
$$
More precisely one takes $C_{1},\,C_{2}>1$  and
$\varepsilon=\varepsilon_{\eta}$
 such that
$$
f_{c_{1}-2\varepsilon}^{c_{1}-\frac32\varepsilon}\cap\left\{|y''|<
  \frac{\eta}{C_{1}}\right\}
\subset \left\{\frac{\eta}{C_{2}}<|y'|< \frac{2\eta}{C_{2}}\,,
|y''|<
  \frac{\eta}{C_{1}}\right\}
\subset B(U_{1},\eta)\,.
$$
Lemma~A.2.2 of \cite{HelSj4} says that
$d_{Ag}(x,y)=\left|f(x)-f(y)\right|$ if and only if there is a
generalized integral curve of $\nabla f$ going from $x$ to $y$\,.
Hence in $f_{c_{1}-2\varepsilon}^{c_{1}-\frac32\varepsilon}$ the only points
such that $d_{Ag}(U_{1},y)=c_{1}-f(y)$ are the points lying on the
unstable manifold for $-\nabla f$\,. Hence there exists a constant 
$C_{\eta}>0$ such that 
$$
\forall y\in f_{c_{1}-2\varepsilon}^{c_{1}-\frac32\varepsilon}
\setminus \left\{|y''|<\frac{\eta}{C_{1}}\right\}\,,\quad
d_{Ag}(U_{1},y)\geq c_{1}-f(y)+C_{\eta}\,.
$$
By combining this with the exponential decay estimates for
$\tilde{v}_{U_{0}}$ we deduce
$$
\langle \omega_{U_{1}}\,,\, hd\chi_{\varepsilon}\wedge
\tilde{v}_{U_{0}}\rangle
=\int_{|y''|\leq \frac{\eta}{C_{1}}}\langle v_{U_{1}}\,,\, hd\chi_{\varepsilon}\wedge
\tilde{v}_{U_{0}}\rangle_{\Lambda T^{*}_{y}M} \  +\  \mathcal{O}(e^{-\frac{c_{1}-c_{0}+C_{\eta}}{h}})\,.
$$
With the above inclusion and the approximation of $v_{U_{1}}$ in
$B(U_{1},\eta)$ stated in Theorem~\ref{th.HSadapt} we get
\begin{eqnarray*}
\langle \omega_{U_{1}}\,,\, hd\chi_{\varepsilon}\wedge
\tilde{v}_{U_{0}}\rangle
\!=\!K_{U_{1}}^{h}\int_{|y''|\leq \frac{\eta}{C_{1}}}\!\!\langle e^{-\frac{\sum_{j=1}^{d}|\lambda^1_{j}|y_{j}^{2}}{2h}}
dy_{1}\wedge\ldots \wedge dy_{p+1}\,,\, hd\chi_{\varepsilon}\wedge
\tilde{v}_{U_{0}}\rangle_{\Lambda T^{*}_{y}M}\\
+\quad  \mathcal{O}(e^{-\frac{c_{1}-c_{0}+C_{\eta}}{h}})\,,
\end{eqnarray*}
with $K_{U_{1}}^{h}
=\frac{|\lambda^1_{1}\ldots \lambda^1_{d}|^{1/4}}{(\pi
h)^{d/4}}$\,. With an Euclidean metric, inserting $e^{-\frac{f-c_{1}}{h}}\times
e^{\frac{f-c_{0}}{h}}$ in the bracket implies that
$\frac{e^{\frac{c_{1}-c_{0}}{h}}}{K_{U_{1}}^{h}}\langle \omega_{U_{1}}\,,\, hd\chi_{\varepsilon}\wedge
\tilde{v}_{U_{0}}\rangle$ equals
\begin{eqnarray*}
&&\int_{|y''|\leq \frac{\eta}{C_{1}}}
\langle e^{-\frac{\sum_{j=p+2}^{d}|\lambda^1_{j}|y_{j}^{2}}{h}}
dy_{1}\wedge\ldots \wedge dy_{p+1}\,,\, hd\chi_{\varepsilon}\wedge
e^{\frac{f-c_{0}}{h}}\tilde{v}_{U_{0}}\rangle_{\Lambda T^{*}_{y}M}
+\mathcal{O}(e^{-\frac{C_{\eta}}{h}})
\\
&&=
\pm\int_{|y''|\leq \frac{\eta}{C_{1}}}
e^{-\frac{\sum_{j=p+2}^{d}|\lambda^1_{j}|y_{j}^{2}}{h}}
dy_{p+2}\wedge\ldots \wedge dy_{d}\wedge(hd\chi_{\varepsilon})\wedge
e^{\frac{f-c_{0}}{h}}\tilde{v}_{U_{0}}+\mathcal{O}(e^{-\frac{C_{\eta}}{h}})\,.
\end{eqnarray*}
But our assumption says that
$d\left(e^{\frac{f}{h}}\tilde{v}_{U_{0}}\right)=0$ in $\supp \nabla
\chi_{\varepsilon}$\,.
Moreover one has clearly $d(e^{-\frac{\sum_{j=p+2}^{d}|\lambda^1_{j}|y_{j}^{2}}{h}}
dy_{p+2}\wedge\ldots \wedge dy_{d})=0$\,.
Hence the integrand is nothing but
$$
hd\left(\chi_{\varepsilon} e^{-\frac{\sum_{j=p+2}^{d}|\lambda^1_{j}|y_{j}^{2}}{h}}
dy_{p+2}\wedge\ldots \wedge dy_{d}\wedge e^{\frac{f-c_{0}}{h}}\tilde{v}_{U_{0}} \right)\,.
$$
By Stokes' formula the quantity 
$ \frac{e^{\frac{c_{1}-c_{0}}{h}}}{K_{U_{1}}^{h}}\langle \omega_{U_{1}}\,,\, hd\chi_{\varepsilon}\wedge
\tilde{v}_{U_{0}}\rangle$ equals
$$
\pm h\int_{|y''|\leq \frac{\eta}{C_{1}}}\int_{|y'|= \frac{2\eta}{C_{2}}}
e^{-\frac{\sum_{j=p+2}^{d}|\lambda^1_{j}|y_{j}^{2}}{h}}
dy_{p+2}\wedge\ldots \wedge dy_{d}\wedge
e^{\frac{f-c_{0}}{h}}\tilde{v}_{U_{0}}
+ \mathcal{O}(e^{-\frac{C_{\eta}}{h}})\,,
$$
and by introducing for every fixed  $y''$ such that $|y''|\leq
\frac{\eta}{C_{1}}$ the cycle $C_{y''}$ supported by 
$\left\{(y',y''), |y'|=\frac{2\eta}{C_{2}}\right\}$  and homotopic to 
$\partial e_{U_{1}}^{p+1}$ we get
$$
\frac{e^{\frac{c_{1}-c_{0}}{h}}}{K_{U_{1}}^{h}}\langle \omega_{U_{1}}\,,\, hd\chi_{\varepsilon}\wedge
\tilde{v}_{U_{0}}\rangle
=\pm h\int_{|y''|\leq
  \frac{\eta}{C_{1}}}e^{-\frac{\sum_{j=p+2}^{d}|\lambda^1_{j}|y_{j}^{2}}{h}}
\int_{C_{y''}}e^{\frac{f-c_{0}}{h}}\tilde{v}_{U_{0}} +\mathcal{O}(e^{-\frac{C_{\eta}}{h}})\,.
$$
For any $y''$, the cycle $C_{y''}$ is homologous to $\partial
e_{U_{1}}^{p+1}$ and according to Proposition~\ref{pr.homolcoef} to
$\kappa [e_{U_{0}}^{p}]$ in $f^{c_{1}-\varepsilon}$ relatively to
$f^{c_{0}-\gamma_{\eta}}$, with $\gamma_{\eta}>0$ small
enough. Owing to
$d\left(e^{\frac{f}{h}}\tilde{v}_{U_{0}}=0\right)$ in
$f^{c_{1}-\varepsilon}$ and the exponential decay estimate of
$\tilde{v}_{U_{0}}$ stated in Theorem~\ref{th.HSadapt} for $\tilde{v}_{U_{0}}$, we obtain
$$
\frac{e^{\frac{c_{1}-c_{0}}{h}}}{K_{U_{1}}^{h}}\langle \omega_{U_{1}}\,,\, hd\chi_{\varepsilon}\wedge
\tilde{v}_{U_{0}}\rangle
=\pm \kappa h\int_{|y''|\leq
  \frac{\eta}{C_{1}}}e^{-\frac{\sum_{j=p+2}^{d}|\lambda^1_{j}|y_{j}^{2}}{h}}
\int_{e_{U_{0}}^{p}}
e^{\frac{f-c_{0}}{h}}\tilde{v}_{U_{0}} +\mathcal{O}(e^{-\frac{C_{\eta}}{h}})\,.
$$
Using again Theorem~\ref{th.HSadapt} with
$\tilde{v}_{U_{0}}$, $f(U_{0})=c_{0}$, and decomposition
of $f$ around $U_{0}$,
$f(z)-c_{0}=\frac{1}{2}\sum_{j=1}^{d}\lambda^0_{j}z_{j}^{2}$
in some (adapted) Morse coordinates $(z',z'')
=(z_{1},\ldots, z_{p},z_{p+1},\ldots,
z_{d})$,
we get,
\begin{eqnarray*}
\frac{e^{\frac{c_{1}-c_{0}}{h}}}{K_{U_{1}}^{h}K_{U_{0}}^{h}}
\langle \omega_{U_{1}}\,,\, hd\chi_{\varepsilon}\wedge
\tilde{v}_{U_{0}}\rangle&=&
\pm \kappa h\int_{|y''|\leq
  \frac{\eta}{C_{1}}}e^{-\frac{\sum_{j=p+2}^{d}|\lambda^1_{j}|y_{j}^{2}}{h}}
\int_{e_{U_{0}}^{p}}e^{\frac{
\sum_{j=1}^{d}\lambda^0_{j}z_{j}^{2}
}{2h}}e^{-\frac{\sum_{j=1}^{d}|\lambda^0_{j}|z_{j}^{2}}{2h}}\\
&&\quad+\quad\mathcal{O}(e^{-\frac{C_{\eta}}{h}})\\
&=&
\pm \kappa h\int_{|y''|\leq
  \frac{\eta}{C_{1}}}e^{-\frac{\sum_{j=p+2}^{d}|\lambda^1_{j}|y_{j}^{2}}{h}}
\int_{e_{U_{0}}^{p}}e^{-\frac{\sum_{j=1}^{p}|\lambda^0_{j}|z_{j}^{2}}{h}}\\
&&\quad +\quad\mathcal{O}(e^{-\frac{C_{\eta}}{h}})\,,
\end{eqnarray*}
with
$K_{U_{0}}^{h}=\frac{|\lambda^0_{1}\ldots \lambda^0_{d}|^{1/4}}{(\pi
h)^{d/4}}$\,.\\
Now, writing successively two Laplace methods,
we obtain
$$
\frac{e^{\frac{c_{1}-c_{0}}{h}}}{K_{U_{1}}^{h}K_{U_{0}}^{h}}
\langle \omega_{U_{1}}\,,\, hd\chi_{\varepsilon}\wedge
\tilde{v}_{U_{0}}\rangle
=
\pm \kappa h\frac{(\pi h)^{\frac{d-1}{2}}}{|\lambda^1_{p+2}\cdots\lambda^1_{d}|^\frac12
\,|\lambda^0_{1}\cdots\lambda^0_{p}|^\frac12}(1+\mathcal{O}(h))\,,$$
which leads, finally, to the following formula: 
$$
\langle \omega_{U_{1}}\,,\, hd\chi_{\varepsilon}\wedge
\tilde{v}_{U_{0}}\rangle
=
\pm \kappa {(\frac h\pi)}^\frac12 \frac{|\lambda^1_{1}\cdots\lambda^1_{p+1}|^\frac14}{
|\lambda^1_{p+2}\cdots\lambda^1_{d}|^\frac14}
\frac{|\lambda^0_{p+1}\cdots\lambda^0_{d}|^\frac14}{
|\lambda^0_{1}\cdots\lambda^0_{p}|^\frac14}e^{-\frac{c_{1}-c_{0}}{h}}
(1+\mathcal{O}(h))\,.
$$

The picture below summarizes the scheme of the calculation and use of
Stokes' formula, for $d=3$
and $p=1$.
\begin{center}
\newbox\surfb


\figinit{0.1cm}

\figpt 2: (18,60)

\figpt 11: (32,60)
\figpt 13: (31,35)
\figpt 14: (68,36)

\figpt 21: (5,50)
\figpt 22: (44,22)


\figvisu{\surfb}{Figure 3}{
\figinsert{shema.eps,0.6}
\figwrites 2: $\partial e^{p+1}_{U_{1}}$ (0.5)
\figwrites 11:  $C_{y''}$ (0.1)
\figwritesw 13:  $\kappa\times$ (0.5)
\figwrites 14:  $f=f{(U_{0})}$ (0.)
\figwritew 21:  $f=f{(U_{1})}$ (0.)
\figwritene 22: $[e^{p}_{U_{0}}]$ (0.3)
}
\hspace{4cm}\box\surfb\\
The arrows show the use of Skokes' formula. The dotted part of
$[e_{U_{0}}^{p}]$ shows the part of $[e_{U_{0}}^{p}]$ lying below $f(U_{0})-\gamma_{\eta}$\,.
\end{center}

\subsubsection{Proof of Proposition~\ref{pr.computdf}
for a general Riemannian metric}

\noindent
As in the previous subsection, we look at the term 
$\langle v_{U_{1}}\,,\, hd\chi_{\varepsilon}\wedge
\tilde{v}_{U_{0}}\rangle$, where $U_{1}\in\mathcal U^{(p+1)}$
satisfies $\partial_{\mathcal B}(U_{1})=U_{0}$, and
we use some adapted Morse coordinates $(y',y'')=
(y_{1},\ldots, y_{p+1}, y_{p+2},\ldots,
y_{d})$
in the ball $B(U_{1},\eta)$. 
Let us recall that the function $f$ has the following decomposition
in these coordinates:
$$f(y)-c_{1}=\frac{1}{2}\sum_{j=1}^{d}\lambda^1_{j}y_{j}^{2}\,.$$ 
Again, one takes $C_{1},\,C_{2}>1$  and
$\varepsilon=\varepsilon_{\eta}$
 such that
$$
f_{c_{1}-2\varepsilon}^{c_{1}-\frac32\varepsilon}\cap\left\{|y''|<
  \frac{\eta}{C_{1}}\right\}
\subset \left\{\frac{\eta}{C_{2}}<|y'|< \frac{2\eta}{C_{2}}\,,
|y''|<
  \frac{\eta}{C_{1}}\right\}
\subset B(U_{1},\eta)\,,
$$
and we have the existence of
$C_{\eta}>0$ s.t.
$$\langle v_{U_{1}}\,,\, hd\chi_{\varepsilon}\wedge
\tilde{v}_{U_{0}}\rangle
=\int_{|y''|\leq \frac{\eta}{C_{1}}}\langle v_{U_{1}}\,,\, hd\chi_{\varepsilon}\wedge
\tilde{v}_{U_{0}}\rangle_{\Lambda T^{*}_{y}M} \  +\  \mathcal{O}(e^{-\frac{c_{1}-c_{0}+C_{\eta}}{h}})\,.
$$
Choose $\eta>0$ small enough such that 
$U_{1}$ is the only critical point of $f$
in $f^{-1}([c_{1}-2\eta,c_{1}+2\eta])$.
\\
\\
Now, let us introduce the metric $g_{1}$, 
$$   
g_{1}(y) = \chi(y) g_{e}(y) + (1-\chi(y)) g(y)\,,
$$
where $0\leq\chi\leq 1$ is a smooth cut-off function 
such that $\chi=1$ in $B(U_{1},\eta)$, 
$\chi=0$ outside $B(U_{1},\frac32\eta)$, $g$ is the usual metric on $M$, and $g_{e}$ is the Euclidean metric
$g_{e}=\sum_{i=1}^d (dy_{i})^2$.\\
\\
For $g$ and $g_{1}$, let $\tilde v_{U_{1}}$ and $\tilde v_{U_{1}}^{1}$
denote respectively the forms
defined in Theorem~\ref{th.HSadapt} with 
$a=c_{1}-2\eta$ and $b=c_{1}+2\eta$. Since $U_{1}$
is the only critical point of $f$ in $f^{-1}([c_{1}-2\eta,c_{1}+2\eta])$,
this means that $\tilde v^{i}_{U_{1}}$
is a normalized form in the one dimensional kernel of
 $\Delta^{TN,(p+1)}_{g,f,h}$ (resp. $\Delta^{TN,(p+1)}_{g_1,f,h}$),
 the Witten Laplacian corresponding to the metric $g$ (resp. $g_{1}$).\\
Note also that the boundary conditions are strictly the same for both
$\tilde v_{U_{1}}$ and $\tilde v^{1}_{U_{1}}$, since the metrics $g$
and $g_{1}$ coincide near the boundary. In particular,
$\tilde v_{U_{1}}$ and $\tilde v^{1}_{U_{1}}$ belong to the same
domain $D(\Delta^{TN,(p+1)}_{f,h}):=D(\Delta^{TN,(p+1)}_{g,f,h})=D(\Delta^{TN,(p+1)}_{g_{1},f,h})$.
Analogously,
$\star\tilde v_{U_{1}}$ and $ \star_{1}\tilde v^{1}_{U_{1}}$
belong to the same domain
$D(\Delta^{NT,(d-p-1)}_{-f,h}):=D(\Delta^{NT,(d-p-1)}_{g,-f,h})=D(\Delta^{NT,(d-p-1)}_{g_{1},-f,h})$ (we refer to Remark~\ref{re.NT}
for the meaning of $\Delta^{NT}_{-f,h})$).
\\
\noindent 
Since
$$ v_{U_{1}}=  \tilde v_{U_{1}} +\mathcal O(e^{-\frac{C_{\eta}}{h}})$$
holds in  $f^{-1}([c_{1}-2\eta,c_{1}+2\eta])\cap\supp d\chi_{\varepsilon}$\,,
it suffices to estimate
$$
\langle v_{U_{1}}\,,\, hd\chi_{\varepsilon}\wedge
\tilde{v}_{U_{0}}\rangle
=
\int_{|y''|\leq \frac{\eta}{C_{1}}}\langle \tilde v_{U_{1}}\,,\, hd\chi_{\varepsilon}\wedge
\tilde{v}_{U_{0}}\rangle_{\Lambda T^{*}_{y}M} \  +\  \mathcal{O}(e^{-\frac{c_{1}-c_{0}+C_{\eta}}{h}})\,.
$$

\noindent The following lemma gives some useful relations between
$\tilde v_{U_{1}}$ and $\tilde v^{1}_{U_{1}}$, especially the second one which will be crucial in the sequel.

\begin{lemma}
\label{le.c1}
There exist $\omega$ in $D(\Delta^{TN,(p)}_{f,h})$
and $\omega'$ in $D(\Delta^{NT,(d-p-2)}_{-f,h})$ s.t.
\begin{eqnarray}
e^{\frac{f-c_{1}}{h}}\tilde v_{U_{1}}&=& h\,d \left(e^{\frac{f-c_{1}}{h}}\omega\right)+(1+\mathcal O(h))e^{\frac{f-c_{1}}{h}}\tilde v^1_{U_{1}}\,,\\
\star\!\left(e^{-\frac{f-c_{1}}{h}}\tilde v_{U_{1}}\right)&=& h\,d \left(e^{-\frac{f-c_{1}}{h}}{\omega}'
\right)+(1+\mathcal O(h))\star_{1}\!\left(e^{-\frac{f-c_{1}}{h}}\tilde v^{1}_{U_{1}}\right)
\,.
\end{eqnarray}
\end{lemma}

\begin{proof}

The form $\tilde v_{U_{1}}$ (resp. $\tilde v^{1}_{U_{1}}$) is 
in the one dimensional kernel of $\Delta^{TN,(p+1)}_{g,f,h}$ (resp. $\Delta^{TN,(p+1)}_{g_{1},f,h}$)
and the isomorphisms
$$ 
\Ker\Delta^{TN,(p+1)}_{g,f,h}\,\sim\,\Ker d^{TN}_{f,h} /\,\Ran
d^{TN}_{f,h}\sim \Ker\Delta^{TN,(p+1)}_{g_{1},f,h}\,,
$$
with the middle set independent of the metrics, implies the existence
of a constant $\alpha_{1}\neq 0$ 
s.t.
$$
\tilde v_{U_{1}}-\alpha_{1}\,
\tilde v^{1}_{U_{1}}\ \in\ \Ran d^{TN}_{f,h}\,.
$$
This means the existence of $\omega$ in $D(\Delta^{TN,(p+1)}_{f,h})$
s.t. 
\begin{equation}
\label{eq.c1}
e^{\frac{f-c_{1}}{h}}\tilde v_{U_{1}}-\alpha_{1}\,e^{\frac{f-c_{1}}{h}}\tilde v_{U_{1}}^{1}
=h\,d\left(e^{\frac{f-c_{1}}{h}}\omega\right)\,.
\end{equation}
Moreover, the form
$\star\tilde{v}_{U_{1}}$ (resp. $\star_{1}\tilde{v}_{U_{1}}^{1}$)
 belongs to  $\ker(\Delta^{NT,(d-p-1)}_{g,-f,h})$ (resp. $\ker(\Delta^{NT,(d-p-1)}_{g_{1},-f,h})$),
and there exists another constant $\alpha'_{1}\neq 0$
s.t.
$$
\star\tilde v_{U_{1}}-\alpha'_{1}\star_{1}\tilde v_{U_{1}}^{1}\ \in\ \Ran d^{NT}_{-f,h}\,.
$$
According to the definition of $d_{-f,h}$, it means that there exists 
$\omega'$ in $D(\Delta^{NT,(d-p-2)}_{-f,h})$
s.t. 
\begin{equation}
\label{eq.alphaprim1}
\star\!\left(e^{-\frac{f-c_{1}}{h}}\tilde v_{U_{1}}\right)-\alpha'_{1}\star_{1}\!\left(e^{-\frac{f-c_{1}}{h}}\tilde v_{U_{1}}^{1}\right) = h\,d \left(e^{-\frac{f-c_{1}}{h}}\omega'\right)\,.
\end{equation}
In order to show that $\alpha_{1}=1+\mathcal O(h)$,
let us integrate equation \eqref{eq.c1} along the unstable manifold 
$\mathcal C_{\text{Unst}}$, which is 
 a $(p+1)$-cycle in $\overline{f_{c_{1}-2\eta}^{c_{1}+2\eta}}$ relatively to
$\left\{f=c_{1}-2\eta\right\}$. Using Stokes' formula,
we obtain
$$
\int_{\mathcal C_{\text{Unst}}}e^{\frac{f-c_{1}}{h}}\left(\tilde v_{U_{1}}-\alpha_{1}\tilde v_{U_{1}}^{1}
\right)
=h\int_{\mathcal C_{\text{Unst}}} d\left(e^{\frac{f-c_{1}}{h}}\omega\right)=0\,.
$$
Consequently, the constant $\alpha_{1}$ is given by
$$
\alpha_{1}=\frac{\int_{\mathcal C_{\text{Unst}}}e^{\frac{f-c_{1}}{h}}\tilde v_{U_{1}}}{
\int_{\mathcal C_{\text{Unst}}}e^{\frac{f-c_{1}}{h}}\tilde v_{U_{1}}^{1}}
=1+\mathcal O(h)\,,
$$
where the last equality comes from a Laplace method applied to each 
integral, after introducing the first order WKB approximations of $\tilde v_{U_{1}}$
and $\tilde v_{U_{1}}^{1}$ recalled from \cite{HelSj4} in the last statements of Theorem~\ref{th.HSadapt}.\\
In order to obtain the same estimate for $\alpha'_{1}$, we make the same computation
with 
equation \ref{eq.alphaprim1} along the stable manifold 
$\mathcal C_{\text{St}}$, which is 
 a $(d-p-1)$-cycle in $\overline{f_{c_{1}-2\eta}^{c_{1}+2\eta}}$ relatively to
$\left\{f=c_{1}+2\eta\right\}$. This gives, according
again to the last statements of Theorem~\ref{th.HSadapt}, 
$$
\alpha'_{1}=\frac{\int_{\mathcal C_{\text{St}}}e^{-\frac{f-c_{1}}{h}}\star\tilde v_{U_{1}}}{
\int_{\mathcal C_{\text{St}}}e^{-\frac{f-c_{1}}{h}}\star_{1}\tilde v_{U_{1}}^{1}}
=1+\mathcal O(h)\,.
$$
\end{proof}

\noindent
Consider now the quantity
\begin{eqnarray*}
A&:=&\int_{|y''|\leq \frac{\eta}{C_{1}}}\langle \tilde v_{U_{1}}\,,\, hd\chi_{\varepsilon}\wedge
\tilde{v}_{U_{0}}\rangle_{\Lambda T^{*}_{y}M}\\
&=&e^{-\frac{c_{1}-c_{0}}{h}}\int_{|y''|\leq \frac{\eta}{C_{1}}}
\star(e^{-\frac{f-c_{1}}{h}} \tilde v_{U_{1}})
\wedge(hd\chi_{\varepsilon})\wedge(e^{\frac{f-c_{0}}{h}}\tilde v_{U_{0}}) \,.
\end{eqnarray*}
By our assumption, $\tilde v_{U_{1}}$
is in $f^{-1}([c_{1}-2\eta,c_{1}+2\eta])$ solution
of  $\Delta_{g,f,h}^{TN,(p+1)}\tilde v_{U_{1}}=0$,
then we have on this domain the following equality:
$$
d\left(\star(e^{-\frac{f-c_{1}}{h}}\tilde v_{U_{1}})
\right)
=
(-1)^{d-p}e^{-\frac{f-c_{1}}{h}}
d^*_{g,f,h}\tilde v_{U_{1}}=0\,.
$$
Keeping in mind the relation
$d(e^{\frac{f-c_{0}}{h}}\tilde v_{U_{0}})=0$ in 
$\supp d\chi_{\varepsilon}$, this implies
$$
\star(e^{-\frac{f-c_{1}}{h}}\tilde v_{U_{1}})
\wedge(hd\chi)\wedge(e^{\frac{f-c_{0}}{h}}\tilde v_{U_{0}})=
hd\left( \chi_{\varepsilon}\star\!(e^{-\frac{f-c_{1}}{h}}\tilde v_{U_{1}})
\wedge(e^{\frac{f-c_{0}}{h}}\tilde v_{U_{0}})\right)\,.
$$
Then we have by Stokes' formula,
\begin{eqnarray*}
e^{\frac{c_{1}-c_{0}}{h}}A
&=& h\int_{\partial(\{|y''|\leq \frac{\eta}{C_{1}}\})}
\chi_{\varepsilon}\star(e^{-\frac{f-c_{1}}{h}}\tilde v_{U_{1}})
\wedge(e^{\frac{f-c_{0}}{h}}\tilde v_{U_{0}})
\\
&=&h\int_{|y''|\leq \frac{\eta}{C_{1}}}\int_{|y'|= \frac{2\eta}{C_{2}}}
\star(e^{-\frac{f-c_{1}}{h}}\tilde v_{U_{1}})
\wedge(e^{\frac{f-c_{0}}{h}}\tilde v_{U_{0}})+\mathcal O(e^{-\frac Ch})\,.
\end{eqnarray*}

\noindent
Using now the second  equation of Lemma~\ref{le.c1},
let us write
\begin{eqnarray}
\nonumber
e^{\frac{c_{1}-c_{0}}{h}}A
=
h(1+\mathcal O(h))\int_{|y''|\leq \frac{\eta}{C_{1}}}\int_{|y'|= \frac{2\eta}{C_{2}}}
\star_{1}(e^{-\frac{f-c_{1}}{h}}\tilde v_{U_{1}}^{1})
\wedge(e^{\frac{f-c_{0}}{h}}\tilde v_{U_{0}})\\
\label{eq.comparEuclid}
+
h^2\int_{|y''|\leq \frac{\eta}{C_{1}}}\int_{|y'|= \frac{2\eta}{C_{2}}}
d \left(e^{-\frac{f-c_{1}}{h}}\omega'\right)\wedge(e^{\frac{f-c_{0}}{h}}\tilde v_{U_{0}})
+\mathcal O(e^{-\frac Ch})\,.
\end{eqnarray}
But, looking  at the second equation of Lemma~\ref{le.c1}
and using our exponential decay estimates,
the second term of the r.h.s. is also (up to an exponentially error term),  
\begin{eqnarray*}
&h^2&\int_{\partial(\{|y''|\leq \frac{\eta}{C_{1}}\})}
d \left(e^{-\frac{f-c_{1}}{h}}\omega'\right)\wedge(e^{\frac{f-c_{0}}{h}}\tilde v_{U_{0}})
+ \mathcal O(e^{-\frac Ch})\,,
\end{eqnarray*}
where the integral term is $0$
owing to Stokes' formula.\\
Equation~\ref{eq.comparEuclid} can then be rewritten
\begin{eqnarray*}
\nonumber
e^{\frac{c_{1}-c_{0}}{h}}A
=
h(1+\mathcal O(h))\int_{|y''|\leq \frac{\eta}{C_{1}}}\int_{|y'|= \frac{2\eta}{C_{2}}}
\star_{1}(e^{-\frac{f-c_{1}}{h}}\tilde v_{U_{1}}^{1})
\wedge(e^{\frac{f-c_{0}}{h}}\tilde v_{U_{0}})
+\mathcal O(e^{-\frac Ch})\,,
\end{eqnarray*}
and we can focus on the integral part  of the r.h.s.,
$$
B:=
\int_{|y''|\leq \frac{\eta}{C_{1}}}\int_{|y'|= \frac{2\eta}{C_{2}}}
\star_{1}(e^{-\frac{f-c_{1}}{h}}\tilde v_{U_{1}}^{1})
\wedge(e^{\frac{f-c_{0}}{h}}\tilde v_{U_{0}})\,.
$$ 
Since the metric $g_{1}$ is Euclidean in the ball $B(U_{1},\eta)$,
the Morse decomposition of $f$ combined with Theorem~\ref{th.HSadapt}
gives
\begin{eqnarray*}
B
=
K_{U_{1}}^{h}\int_{|y''|\leq \frac{\eta}{C_{1}}}\int_{|y'|= \frac{2\eta}{C_{2}}}e^{-\frac{\sum_{j=p+2}^{d}|\lambda^1_{j}|y_{j}^{2}}{h}}
dy_{p+2}\wedge\ldots \wedge dy_{d}\wedge
(e^{\frac{f-c_{0}}{h}}\tilde v_{U_{0}})
+\mathcal{O}(e^{-\frac{C_{\eta}}{h}})\,,
\end{eqnarray*}
with $K_{U_{1}}^{h}
=\frac{|\lambda^1_{1}\ldots \lambda^1_{d}|^{1/4}}{(\pi
h)^{d/4}}$\,.\\
We can then follow the same proof as the one used in the locally Euclidean case
in order to obtain the wanted result. The only difference
arises in the fact that the WKB expansion of  $\tilde v_{U_{0}}$, given in the last statement 
of Theorem~\ref{th.HSadapt},
contains higher order correcting terms, because the full WKB expansion
of $\tilde{v}_{U_{0}}$ depends on the metric, and
this produces a relative $\mathcal{O}(h)$ error  term.

\subsection{End of the proof of Theorem~\ref{th.main}}
We recall firstly some notations of Section~\ref{se.qualowup}.
The operator $\Delta_{f,h}$ is defined on
$M$ (i.e. on $\overline{f_{-\infty}^{+\infty}}$)
and
$$F=\oplus_{p=0}^d F^{(p)}\quad 
\text{with}\quad 
 F^{(p)}=\ima 1_{[0,h^{2}]}(\Delta_{f,h}^{(p)})
\,.
 $$
According to Section~\ref{se.wit}, 
$F^{(p)}$ admits an almost orthonormal basis,
$\left\{v_{U}, U\in \mathcal{U}^{(p)}\right\}$, fulfilling the
properties of Theorem~\ref{th.HSadapt}.\\
Consider now the family of quasimodes,
$\left\{\omega_{U}, U\in \bigcup_{p=0}^d\mathcal{U}^{(p)}\right\}$,
constructed in Section~\ref{se.globquasi}.
For any $p$ in $\{0,\dots,d\}$
and $U\in \mathcal{U}^{(p)}$,
the $p$-form $\omega_{U}$ belongs to $F^{(p)}$
and, as already mentioned, $\omega_{U}$ satisfies
the relation $\|\omega_{U}-v_{U}\|=\mathcal O(e^{-\frac{C_{\eta}}{
    h}})$
(see Proposition~\ref{pr.quasiclass} for upper or homological critical
points and \eqref{eq.compwvlower} for lower critical points).\\
The family $\left\{\omega_{U}, U\in \mathcal{U}^{(p)}\right\}$
is then an almost orthonormal basis of
$F^{(p)}$ and, thanks to Proposition~\ref{pr.computdf},
Theorem~2.3 of \cite{Lep1} applies. For $h_{0}$ small enough and $h\in(0,h_{0}]$,
we obtain an accurate writing of the non zero eigenvalues of 
$d^*_{f,h}d_{f,h}:F\to F$.\\
More precisely, 
when restricted to $F^{(p)}$,
the non zero eigenvalues of
$d^*_{f,h}d_{f,h}$
are the quantities
$$
\kappa^2 B(h)
e^{-2\frac{f(U_{U}^{(p+1)})-f(\partial_{\mathcal B}(U_{U}^{(p+1)}))}{h}}
(1+\mathcal{O}(h))\ ,\ U_{U}^{(p+1)}\in\mathcal U_{U}^{(p+1)}\,,
$$
with
$$
B(h)=\frac h\pi\frac{|\lambda^1_{1}\cdots\lambda^1_{p+1}|^\frac12}{
|\lambda^1_{p+2}\cdots\lambda^1_{d}|^\frac12}
\frac{|\lambda^0_{p+1}\cdots\lambda^0_{d}|^\frac12}{
|\lambda^0_{1}\cdots\lambda^0_{p}|^\frac12}\,,
$$
where $\lambda_{1}^{\ell}<\cdots<\lambda_{p+\ell}^{\ell}<0<
\lambda_{p+\ell+1}^{\ell}<\cdots<\lambda_{d}^{\ell}$
denote the eigenvalues of $\Hess f(U_{\ell})$;
for $\ell\in\{0,1\}$, $U_{1}:=U_{U}^{(p+1)}$,
and $U_{0}:=\partial_{\mathcal B}(U_{U}^{(p+1)})$.\\
\\
In particular, $d^*_{f,h}d_{f,h}:F\to F$
has exactly $\text{card}\,\mathcal U_{U}=\text{card}\,\mathcal U_{L}$
non zero eigenvalues and these eigenvalues are distinct.\\
This provides all the exponentially small non zero eigenvalues of
$\Delta_{f,h}$ according to the last statement of
Theorem~\ref{th.Witt}.\\
This ends the proof of Theorem~\ref{th.main}.

\subsection{A relative version of Theorem~\ref{th.main}}
\label{se.relatth}

The analysis for Theorem~\ref{th.main} is done on
$M=f^{+\infty}_{-\infty}$\,. All the constructions and the good
restriction properties of Morse-Barannikov chain complex, when
considering $H_{*}(f^{b},f^{a})$, $-\infty\leq a<b\leq +\infty$, 
have their counterpart with the Witten Laplacians $\Delta_{f,h}^{TN}$
defined on
$\overline{f_{a}^{b}}$ in Section~\ref{se.wit}.
Hence all the proof of Theorem~\ref{th.main} is still valid for
$\Delta_{f,h}^{TN}$ on $\overline{f_{a}^{b}}$ except that some end
points of the relation $\partial_{B}U^{(p+1)}=U^{(p)}$ disappear when
they lie in $f^{a}$ or $f_{b}$\,.
We state without more detail the spectral result for
$\Delta_{f,h}^{TN}$ on $\overline{f_{a}^{b}}$\,.
\begin{theorem}
\label{th.relatmain}
With the same assumptions as in Theorem~\ref{th.main} and when $a,b$
are not critical values of $f$, $-\infty\leq a < b\leq +\infty$, the
exponentially small eigenvalues of $\Delta_{f,h}^{TN}$, defined on
$\overline{f_{a}^{b}}$ according to Proposition~\ref{pr.selfadjoin},
are given by a mapping $j_{a}^{b}:\mathcal{U}\cap f_{a}^{b}\to
\sigma(\Delta_{f,h}^{TN})$ derived from the mapping $j$ of
Theorem~\ref{th.main} by:
\begin{itemize}
\item  $j_{a}^{b}(U)=j(U)(1+\mathcal{O}(h))\neq 0$ if $U\in \mathcal{U}_{L}\cap
  f_{a}^{b}$ and  $U=\partial_{B}U'$  with
  $U'\in f_{a}^{b}$\,;
\item $j_{a}^{b}(U)=j(U)(1+\mathcal{O}(h))\neq 0$ if $U\in \mathcal{U}_{U}\cap
  f_{a}^{b}$ and  $\partial_{B}U=U'$  with
  $U'\in f_{a}^{b}$\,;
\item $j_{a}^{b}(U)=0$ else.
\end{itemize}
\end{theorem}

\noindent\textbf{Acknowledgements~:} The authors wish to thank
Franois Laudenbach from many insightful discussions concerning Morse
theory. These played an essential role in the 
elaboration of this work.

\end{document}

%% file: fig4tex.tex
%
%
%
\ifx\figfortexisloaded\relax \else\let\figfortexisloaded=\relax\fi
\message{version 1.7.7}
\newif\iftextures
\catcode`\@=11
\newdimen\epsil@n\epsil@n=0.00005pt
\newdimen\Cepsil@n\Cepsil@n=0.005pt
\newdimen\dcq@\dcq@=254pt
\newdimen\PI@\PI@=3.141592pt
\newdimen\DemiPI@deg\DemiPI@deg=90pt
\newdimen\PI@deg\PI@deg=180pt
\newdimen\DePI@deg\DePI@deg=360pt
\chardef\t@n=10
\chardef\c@nt=100
\chardef\@lxxiv=74
\chardef\@xci=91
\mathchardef\@nMnCQn=9949
\chardef\@vi=6
\chardef\@xxx=30
\chardef\@lvi=56
\chardef\@lxxi=71
\chardef\@lxxxv=85
\mathchardef\@mmmmlxviii=4068
\mathchardef\@ccclx=360
\mathchardef\@dccxx=720
\newcount\p@rtent \newcount\f@ctech \newcount\result@tent
\newdimen\v@lmin \newdimen\v@lmax \newdimen\v@leur
\newdimen\result@t\newdimen\result@@t
\newdimen\mili@u \newdimen\c@rre \newdimen\delt@
\def\degT@rd{0.017453 }  
\def\rdT@deg{57.295779 } 
{\catcode`p=12 \catcode`t=12 \gdef\v@leurseule#1pt{#1}}
\def\repdecn@mb#1{\expandafter\v@leurseule\the#1\space}
\def\arct@n#1(#2,#3){{\v@lmin=#2\v@lmax=#3%
    \maxim@m{\mili@u}{-\v@lmin}{\v@lmin}\maxim@m{\c@rre}{-\v@lmax}{\v@lmax}%
    \delt@=\mili@u\m@ech\mili@u%
    \ifdim\c@rre>\@nMnCQn\mili@u\divide\v@lmax\tw@\c@lATAN\v@leur(\z@,\v@lmax)
    \else%
    \maxim@m{\mili@u}{-\v@lmin}{\v@lmin}\maxim@m{\c@rre}{-\v@lmax}{\v@lmax}%
    \m@ech\c@rre%
    \ifdim\mili@u>\@nMnCQn\c@rre\divide\v@lmin\tw@
    \maxim@m{\mili@u}{-\v@lmin}{\v@lmin}\c@lATAN\v@leur(\mili@u,\z@)%
    \else\c@lATAN\v@leur(\delt@,\v@lmax)\fi\fi%
    \ifdim\v@lmin<\z@\v@leur=-\v@leur\ifdim\v@lmax<\z@\advance\v@leur-\PI@%
    \else\advance\v@leur\PI@\fi\fi%
    \global\result@t=\v@leur}#1=\result@t}
\def\m@ech#1{\ifdim#1>1.646pt\divide\mili@u\t@n\divide\c@rre\t@n\m@ech#1\fi}
\def\c@lATAN#1(#2,#3){{\v@lmin=#2\v@lmax=#3\v@leur=\z@\delt@=\tw@ pt%
    \un@iter{0.785398}{\v@lmax<}%
    \un@iter{0.463648}{\v@lmax<}%
    \un@iter{0.244979}{\v@lmax<}%
    \un@iter{0.124355}{\v@lmax<}%
    \un@iter{0.062419}{\v@lmax<}%
    \un@iter{0.031240}{\v@lmax<}%
    \un@iter{0.015624}{\v@lmax<}%
    \un@iter{0.007812}{\v@lmax<}%
    \un@iter{0.003906}{\v@lmax<}%
    \un@iter{0.001953}{\v@lmax<}%
    \un@iter{0.000976}{\v@lmax<}%
    \un@iter{0.000488}{\v@lmax<}%
    \un@iter{0.000244}{\v@lmax<}%
    \un@iter{0.000122}{\v@lmax<}%
    \un@iter{0.000061}{\v@lmax<}%
    \un@iter{0.000030}{\v@lmax<}%
    \un@iter{0.000015}{\v@lmax<}%
    \global\result@t=\v@leur}#1=\result@t}
\def\un@iter#1#2{%
    \divide\delt@\tw@\edef\dpmn@{\repdecn@mb{\delt@}}%
    \mili@u=\v@lmin%
    \ifdim#2\z@%
      \advance\v@lmin-\dpmn@\v@lmax\advance\v@lmax\dpmn@\mili@u%
      \advance\v@leur-#1pt%
    \else%
      \advance\v@lmin\dpmn@\v@lmax\advance\v@lmax-\dpmn@\mili@u%
      \advance\v@leur#1pt%
    \fi}
\def\c@ssin#1#2#3{\expandafter\ifx\csname COS@\number#3\endcsname\relax\c@lCS{#3pt}%
    \expandafter\xdef\csname COS@\number#3\endcsname{\repdecn@mb\result@t}%
    \expandafter\xdef\csname SIN@\number#3\endcsname{\repdecn@mb\result@@t}\fi%
    \edef#1{\csname COS@\number#3\endcsname}\edef#2{\csname SIN@\number#3\endcsname}}
\def\c@lCS#1{{\mili@u=#1\p@rtent=\@ne%
    \relax\ifdim\mili@u<\z@\red@ng<-\else\red@ng>+\fi\f@ctech=\p@rtent%
    \relax\ifdim\mili@u<\z@\mili@u=-\mili@u\f@ctech=-\f@ctech\fi\c@@lCS}}
\def\c@@lCS{\v@lmin=\mili@u\c@rre=-\mili@u\advance\c@rre\DemiPI@deg\v@lmax=\c@rre%
    \mili@u\@lxxi\mili@u\divide\mili@u\@mmmmlxviii%
    \edef\v@larg{\repdecn@mb{\mili@u}}\mili@u=-\v@larg\mili@u%
    \edef\v@lmxde{\repdecn@mb{\mili@u}}%
    \c@rre\@lxxi\c@rre\divide\c@rre\@mmmmlxviii%
    \edef\v@largC{\repdecn@mb{\c@rre}}\c@rre=-\v@largC\c@rre%
    \edef\v@lmxdeC{\repdecn@mb{\c@rre}}%
    \fctc@s\mili@u\v@lmin\global\result@t\p@rtent\v@leur%
    \let\t@mp=\v@larg\let\v@larg=\v@largC\let\v@largC=\t@mp%
    \let\t@mp=\v@lmxde\let\v@lmxde=\v@lmxdeC\let\v@lmxdeC=\t@mp%
    \fctc@s\c@rre\v@lmax\global\result@@t\f@ctech\v@leur}
\def\fctc@s#1#2{\v@leur=#1\relax\ifdim#2<\@lxxxv\p@\cosser@h\else\sinser@t\fi}
\def\cosser@h{\advance\v@leur\@lvi\p@\divide\v@leur\@lvi%
    \v@leur=\v@lmxde\v@leur\advance\v@leur\@xxx\p@%
    \v@leur=\v@lmxde\v@leur\advance\v@leur\@ccclx\p@%
    \v@leur=\v@lmxde\v@leur\advance\v@leur\@dccxx\p@\divide\v@leur\@dccxx}
\def\sinser@t{\v@leur=\v@lmxdeC\p@\advance\v@leur\@vi\p@%
    \v@leur=\v@largC\v@leur\divide\v@leur\@vi}
\def\red@ng#1#2{\relax\ifdim\mili@u#1#2\DemiPI@deg\advance\mili@u#2-\PI@deg%
    \p@rtent=-\p@rtent\red@ng#1#2\fi}
\def\invers@#1#2{{\v@leur=#2\maxim@m{\v@lmax}{-\v@leur}{\v@leur}%
    \f@ctech=\@ne\m@inv@rs%
    \multiply\v@leur\f@ctech\edef\v@lv@leur{\repdecn@mb{\v@leur}}%
    \p@rtentiere{\p@rtent}{\v@leur}\v@lmin=\p@\divide\v@lmin\p@rtent%
    \inv@rs@\multiply\v@lmax\f@ctech\global\result@t=\v@lmax}#1=\result@t}
\def\m@inv@rs{\ifdim\v@lmax<\p@\multiply\v@lmax\t@n\multiply\f@ctech\t@n\m@inv@rs\fi}
\def\inv@rs@{\v@lmax=-\v@lmin\v@lmax=\v@lv@leur\v@lmax%
    \advance\v@lmax\tw@ pt\v@lmax=\repdecn@mb{\v@lmin}\v@lmax%
    \delt@=\v@lmax\advance\delt@-\v@lmin\ifdim\delt@<\z@\delt@=-\delt@\fi%
    \ifdim\delt@>\epsil@n\v@lmin=\v@lmax\inv@rs@\fi}
\def\minim@m#1#2#3{\relax\ifdim#2<#3#1=#2\else#1=#3\fi}
\def\maxim@m#1#2#3{\relax\ifdim#2>#3#1=#2\else#1=#3\fi}
\def\p@rtentiere#1#2{#1=#2\divide#1by65536 }
\def\r@undint#1#2{{\v@leur=#2\divide\v@leur\t@n\p@rtentiere{\p@rtent}{\v@leur}%
    \v@leur=\p@rtent pt\global\result@t=\t@n\v@leur}#1=\result@t}
\def\sqrt@#1#2{{\v@leur=#2%
    \minim@m{\v@lmin}{\p@}{\v@leur}\maxim@m{\v@lmax}{\p@}{\v@leur}%
    \f@ctech=\@ne\m@sqrt@\sqrt@@%
    \mili@u=\v@lmin\advance\mili@u\v@lmax\divide\mili@u\tw@\multiply\mili@u\f@ctech%
    \global\result@t=\mili@u}#1=\result@t}
\def\m@sqrt@{\ifdim\v@leur>\dcq@\divide\v@leur\c@nt\v@lmax=\v@leur%
    \multiply\f@ctech\t@n\m@sqrt@\fi}
\def\sqrt@@{\mili@u=\v@lmin\advance\mili@u\v@lmax\divide\mili@u\tw@%
    \c@rre=\repdecn@mb{\mili@u}\mili@u%
    \ifdim\c@rre<\v@leur\v@lmin=\mili@u\else\v@lmax=\mili@u\fi%
    \delt@=\v@lmax\advance\delt@-\v@lmin\ifdim\delt@>\epsil@n\sqrt@@\fi}
\def\extrairelepremi@r#1\de#2{\expandafter\lepremi@r#2@#1#2}
\def\lepremi@r#1,#2@#3#4{\def#3{#1}\def#4{#2}\ignorespaces}
\def\@cfor#1:=#2\do#3{%
  \edef\@fortemp{#2}%
  \ifx\@fortemp\empty\else\@cforloop#2,\@nil,\@nil\@@#1{#3}\fi}
\def\@cforloop#1,#2\@@#3#4{%
  \def#3{#1}%
  \ifx#3\Fig@nnil\let\@nextwhile=\Fig@fornoop\else#4\relax\let\@nextwhile=\@cforloop\fi%
  \@nextwhile#2\@@#3{#4}}

\def\@ecfor#1:=#2\do#3{%
  \def\@@cfor{\@cfor#1:=}%
  \edef\@@@cfor{#2}%
  \expandafter\@@cfor\@@@cfor\do{#3}}
\def\Fig@nnil{\@nil}
\def\Fig@fornoop#1\@@#2#3{}
\def\trtlis@rg#1#2{\def\list@@rg{#1}%
    \@ecfor\p@rv@l:=\list@@rg\do{\expandafter#2\p@rv@l|}}
\newbox\b@xvisu
\newtoks\let@xte
\newif\ifitis@K
\newcount\s@mme
\newcount\l@mbd@un \newcount\l@mbd@de
\newcount\superc@ntr@l\superc@ntr@l=\@ne        
\newcount\typec@ntr@l\typec@ntr@l=\superc@ntr@l 
\newdimen\v@lX  \newdimen\v@lY  \newdimen\v@lZ
\newdimen\v@lXa \newdimen\v@lYa \newdimen\v@lZa
\newdimen\unit@\unit@=\p@ 
\def\unit@util{pt}
\def\ptT@ptps{0.996264 }
\def\ptpsT@pt{1.00375 }
\def\ptT@unit@{1} 
\def\setunit@#1{\def\unit@util{#1}\setunit@@#1:\invers@{\result@t}{\unit@}%
    \edef\ptT@unit@{\repdecn@mb\result@t}}
\def\setunit@@#1#2:{\ifcat#1a\unit@=\@ne#1#2\else\unit@=#1#2\fi}
\def\d@fm@cdim#1#2{{\v@leur=#2\v@leur=\ptT@unit@\v@leur\xdef#1{\repdecn@mb\v@leur}}}
\newif\ifBdingB@x\BdingB@xtrue
\newdimen\c@@rdXmin \newdimen\c@@rdYmin  
\newdimen\c@@rdXmax \newdimen\c@@rdYmax
\def\b@undb@x#1#2{\ifBdingB@x%
    \relax\ifdim#1<\c@@rdXmin\global\c@@rdXmin=#1\fi%
    \relax\ifdim#2<\c@@rdYmin\global\c@@rdYmin=#2\fi%
    \relax\ifdim#1>\c@@rdXmax\global\c@@rdXmax=#1\fi%
    \relax\ifdim#2>\c@@rdYmax\global\c@@rdYmax=#2\fi\fi}
\def\b@undb@xP#1{{\Figg@tXY{#1}\b@undb@x{\v@lX}{\v@lY}}}
\def\ellBB@x#1;#2,#3(#4,#5,#6){{\s@uvc@ntr@l\et@tellBB@x%
    \setc@ntr@l{2}\figptell-2::#1;#2,#3(#4,#6)\b@undb@xP{-2}%
    \figptell-2::#1;#2,#3(#5,#6)\b@undb@xP{-2}%
    \c@ssin{\C@}{\S@}{#6}\v@lmin=\C@ pt\v@lmax=\S@ pt%
    \mili@u=#3\v@lmin\delt@=#2\v@lmax\arct@n\v@leur(\delt@,\mili@u)%
    \mili@u=-#3\v@lmax\delt@=#2\v@lmin\arct@n\c@rre(\delt@,\mili@u)%
    \v@leur=\rdT@deg\v@leur\advance\v@leur-\DePI@deg%
    \c@rre=\rdT@deg\c@rre\advance\c@rre-\DePI@deg%
    \v@lmin=#4pt\v@lmax=#5pt%
    \loop\ifdim\v@leur<\v@lmax\ifdim\v@leur>\v@lmin%
    \edef\@ngle{\repdecn@mb\v@leur}\figptell-2::#1;#2,#3(\@ngle,#6)%
    \b@undb@xP{-2}\fi\advance\v@leur\PI@deg\repeat%
    \loop\ifdim\c@rre<\v@lmax\ifdim\c@rre>\v@lmin%
    \edef\@ngle{\repdecn@mb\c@rre}\figptell-2::#1;#2,#3(\@ngle,#6)%
    \b@undb@xP{-2}\fi\advance\c@rre\PI@deg\repeat%
    \resetc@ntr@l\et@tellBB@x}\ignorespaces}
\def\initb@undb@x{\c@@rdXmin=\maxdimen\c@@rdYmin=\maxdimen%
    \c@@rdXmax=-\maxdimen\c@@rdYmax=-\maxdimen}
\def\c@ntr@lnum#1{%
    \relax\ifnum\typec@ntr@l=\@ne%
    \ifnum#1<\z@%
    \immediate\write16{*** Forbidden point number (#1). Abort.}\end\fi\fi%
    \set@bjc@de{#1}}
\def\set@bjc@de#1{\edef\objc@de{@BJ\ifnum#1<\z@ M\romannumeral-#1\else\romannumeral#1\fi}}
\def\setc@ntr@l#1{\ifnum\superc@ntr@l>#1\typec@ntr@l=\superc@ntr@l%
    \else\typec@ntr@l=#1\fi}
\def\resetc@ntr@l#1{\global\superc@ntr@l=#1\setc@ntr@l{#1}}
\def\s@uvc@ntr@l#1{\edef#1{\the\superc@ntr@l}}
\def\c@lproscalDD#1[#2,#3]{{\Figg@tXY{#2}%
    \edef\Xu@{\repdecn@mb{\v@lX}}\edef\Yu@{\repdecn@mb{\v@lY}}\Figg@tXY{#3}%
    \global\result@t=\Xu@\v@lX\global\advance\result@t\Yu@\v@lY}#1=\result@t}
\def\c@lproscalTD#1[#2,#3]{{\Figg@tXY{#2}\edef\Xu@{\repdecn@mb{\v@lX}}%
    \edef\Yu@{\repdecn@mb{\v@lY}}\edef\Zu@{\repdecn@mb{\v@lZ}}%
    \Figg@tXY{#3}\global\result@t=\Xu@\v@lX\global\advance\result@t\Yu@\v@lY%
    \global\advance\result@t\Zu@\v@lZ}#1=\result@t}
\def\c@lprovec#1{%
    \det@rmC\v@lZa(\v@lX,\v@lY,\v@lmin,\v@lmax)%
    \det@rmC\v@lXa(\v@lY,\v@lZ,\v@lmax,\v@leur)%
    \det@rmC\v@lYa(\v@lZ,\v@lX,\v@leur,\v@lmin)%
    \Figv@ctCreg#1(\v@lXa,\v@lYa,\v@lZa)}
\def\det@rm#1[#2,#3]{{\Figg@tXY{#2}\Figg@tXYa{#3}%
    \delt@=\repdecn@mb{\v@lX}\v@lYa\advance\delt@-\repdecn@mb{\v@lY}\v@lXa%
    \global\result@t=\delt@}#1=\result@t}
\def\det@rmC#1(#2,#3,#4,#5){{\global\result@t=\repdecn@mb{#2}#5%
    \global\advance\result@t-\repdecn@mb{#3}#4}#1=\result@t}
\def\getredf@ctDD#1(#2,#3){{\maxim@m{\v@lXa}{-#2}{#2}\maxim@m{\v@lYa}{-#3}{#3}%
    \maxim@m{\v@lXa}{\v@lXa}{\v@lYa}
    \ifdim\v@lXa>\@xci pt\divide\v@lXa\@xci%
    \p@rtentiere{\p@rtent}{\v@lXa}\advance\p@rtent\@ne\else\p@rtent=\@ne\fi%
    \global\result@tent=\p@rtent}#1=\result@tent\ignorespaces}
\def\getredf@ctTD#1(#2,#3,#4){{\maxim@m{\v@lXa}{-#2}{#2}\maxim@m{\v@lYa}{-#3}{#3}%
    \maxim@m{\v@lZa}{-#4}{#4}\maxim@m{\v@lXa}{\v@lXa}{\v@lYa}%
    \maxim@m{\v@lXa}{\v@lXa}{\v@lZa}
    \ifdim\v@lXa>\@lxxiv pt\divide\v@lXa\@lxxiv%
    \p@rtentiere{\p@rtent}{\v@lXa}\advance\p@rtent\@ne\else\p@rtent=\@ne\fi%
    \global\result@tent=\p@rtent}#1=\result@tent\ignorespaces}
\def\FigptintercircB@zDD#1:#2:#3,#4[#5,#6,#7,#8]{{\s@uvc@ntr@l\et@tfigptintercircB@zDD%
    \setc@ntr@l{2}\figvectPDD-1[#5,#8]\Figg@tXY{-1}\getredf@ctDD\f@ctech(\v@lX,\v@lY)%
    \mili@u=#4\unit@\divide\mili@u\f@ctech\c@rre=\repdecn@mb{\mili@u}\mili@u%
    \figptBezierDD-5::#3[#5,#6,#7,#8]%
    \v@lmin=#3\p@\v@lmax=\v@lmin\advance\v@lmax0.1\p@%
    \loop\edef\T@{\repdecn@mb{\v@lmax}}\figptBezierDD-2::\T@[#5,#6,#7,#8]%
    \figvectPDD-1[-5,-2]\n@rmeucCDD{\delt@}{-1}\ifdim\delt@<\c@rre\v@lmin=\v@lmax%
    \advance\v@lmax0.1\p@\repeat%
    \loop\mili@u=\v@lmin\advance\mili@u\v@lmax%
    \divide\mili@u\tw@\edef\T@{\repdecn@mb{\mili@u}}\figptBezierDD-2::\T@[#5,#6,#7,#8]%
    \figvectPDD-1[-5,-2]\n@rmeucCDD{\delt@}{-1}\ifdim\delt@>\c@rre\v@lmax=\mili@u%
    \else\v@lmin=\mili@u\fi\v@leur=\v@lmax\advance\v@leur-\v@lmin%
    \ifdim\v@leur>\epsil@n\repeat\figptcopyDD#1:#2/-2/%
    \resetc@ntr@l\et@tfigptintercircB@zDD}\ignorespaces}
\def\inters@cDD#1:#2[#3,#4;#5,#6]{{\s@uvc@ntr@l\et@tinters@cDD%
    \setc@ntr@l{2}\vecunit@{-1}{#4}\vecunit@{-2}{#6}%
    \Figg@tXY{-1}\setc@ntr@l{1}\Figg@tXYa{#3}%
    \edef\A@{\repdecn@mb{\v@lX}}\edef\B@{\repdecn@mb{\v@lY}}%
    \v@lmin=\B@\v@lXa\advance\v@lmin-\A@\v@lYa%
    \Figg@tXYa{#5}\setc@ntr@l{2}\Figg@tXY{-2}%
    \edef\C@{\repdecn@mb{\v@lX}}\edef\D@{\repdecn@mb{\v@lY}}%
    \v@lmax=\D@\v@lXa\advance\v@lmax-\C@\v@lYa%
    \delt@=\A@\v@lY\advance\delt@-\B@\v@lX%
    \invers@{\v@leur}{\delt@}\edef\v@ldelta{\repdecn@mb{\v@leur}}%
    \v@lXa=\A@\v@lmax\advance\v@lXa-\C@\v@lmin%
    \v@lYa=\B@\v@lmax\advance\v@lYa-\D@\v@lmin%
    \v@lXa=\v@ldelta\v@lXa\v@lYa=\v@ldelta\v@lYa%
    \setc@ntr@l{1}\Figp@intregDD#1:{#2}(\v@lXa,\v@lYa)%
    \resetc@ntr@l\et@tinters@cDD}\ignorespaces}
\def\inters@cTD#1:#2[#3,#4;#5,#6]{{\s@uvc@ntr@l\et@tinters@cTD%
    \setc@ntr@l{2}\figvectNVTD-1[#4,#6]\figvectNVTD-2[#6,-1]\figvectPTD-1[#3,#5]%
    \r@pPSTD\v@leur[-2,-1,#4]\edef\v@lcoef{\repdecn@mb{\v@leur}}%
    \figpttraTD#1:{#2}=#3/\v@lcoef,#4/\resetc@ntr@l\et@tinters@cTD}\ignorespaces}
\def\r@pPSTD#1[#2,#3,#4]{{\Figg@tXY{#2}\edef\Xu@{\repdecn@mb{\v@lX}}%
    \edef\Yu@{\repdecn@mb{\v@lY}}\edef\Zu@{\repdecn@mb{\v@lZ}}%
    \Figg@tXY{#3}\v@lmin=\Xu@\v@lX\advance\v@lmin\Yu@\v@lY\advance\v@lmin\Zu@\v@lZ%
    \Figg@tXY{#4}\v@lmax=\Xu@\v@lX\advance\v@lmax\Yu@\v@lY\advance\v@lmax\Zu@\v@lZ%
    \invers@{\v@leur}{\v@lmax}\global\result@t=\repdecn@mb{\v@leur}\v@lmin}%
    #1=\result@t}
\def\n@rminfDD#1#2{{\Figg@tXY{#2}\maxim@m{\v@lX}{\v@lX}{-\v@lX}%
    \maxim@m{\v@lY}{\v@lY}{-\v@lY}\maxim@m{\global\result@t}{\v@lX}{\v@lY}}%
    #1=\result@t}
\def\n@rminfTD#1#2{{\Figg@tXY{#2}\maxim@m{\v@lX}{\v@lX}{-\v@lX}%
    \maxim@m{\v@lY}{\v@lY}{-\v@lY}\maxim@m{\v@lZ}{\v@lZ}{-\v@lZ}%
    \maxim@m{\v@lX}{\v@lX}{\v@lY}\maxim@m{\global\result@t}{\v@lX}{\v@lZ}}%
    #1=\result@t}
\def\n@rmeucCDD#1#2{\Figg@tXY{#2}\divide\v@lX\f@ctech\divide\v@lY\f@ctech%
    #1=\repdecn@mb{\v@lX}\v@lX\v@lX=\repdecn@mb{\v@lY}\v@lY\advance#1\v@lX}
\def\n@rmeucCTD#1#2{\Figg@tXY{#2}%
    \divide\v@lX\f@ctech\divide\v@lY\f@ctech\divide\v@lZ\f@ctech%
    #1=\repdecn@mb{\v@lX}\v@lX\v@lX=\repdecn@mb{\v@lY}\v@lY\advance#1\v@lX%
    \v@lX=\repdecn@mb{\v@lZ}\v@lZ\advance#1\v@lX}
\def\n@rmeucSVDD#1#2{{\Figg@tXY{#2}%
    \v@lXa=\repdecn@mb{\v@lX}\v@lX\v@lYa=\repdecn@mb{\v@lY}\v@lY%
    \advance\v@lXa\v@lYa\sqrt@{\global\result@t}{\v@lXa}}#1=\result@t}
\def\n@rmeucSVTD#1#2{{\Figg@tXY{#2}\v@lXa=\repdecn@mb{\v@lX}\v@lX%
    \v@lYa=\repdecn@mb{\v@lY}\v@lY\v@lZa=\repdecn@mb{\v@lZ}\v@lZ%
    \advance\v@lXa\v@lYa\advance\v@lXa\v@lZa\sqrt@{\global\result@t}{\v@lXa}}#1=\result@t}
\def\n@rmeucDD#1#2{{\Figg@tXY{#2}\getredf@ctDD\f@ctech(\v@lX,\v@lY)%
    \divide\v@lX\f@ctech\divide\v@lY\f@ctech%
    \v@lXa=\repdecn@mb{\v@lX}\v@lX\v@lYa=\repdecn@mb{\v@lY}\v@lY%
    \advance\v@lXa\v@lYa\sqrt@{\global\result@t}{\v@lXa}%
    \global\multiply\result@t\f@ctech}#1=\result@t}
\def\n@rmeucTD#1#2{{\Figg@tXY{#2}\getredf@ctTD\f@ctech(\v@lX,\v@lY,\v@lZ)%
    \divide\v@lX\f@ctech\divide\v@lY\f@ctech\divide\v@lZ\f@ctech%
    \v@lXa=\repdecn@mb{\v@lX}\v@lX%
    \v@lYa=\repdecn@mb{\v@lY}\v@lY\v@lZa=\repdecn@mb{\v@lZ}\v@lZ%
    \advance\v@lXa\v@lYa\advance\v@lXa\v@lZa\sqrt@{\global\result@t}{\v@lXa}%
    \global\multiply\result@t\f@ctech}#1=\result@t}
\def\vecunit@DD#1#2{{\Figg@tXY{#2}\getredf@ctDD\f@ctech(\v@lX,\v@lY)%
    \divide\v@lX\f@ctech\divide\v@lY\f@ctech%
    \Figv@ctCreg#1(\v@lX,\v@lY)\n@rmeucSV{\v@lYa}{#1}%
    \invers@{\v@lXa}{\v@lYa}\edef\v@lv@lXa{\repdecn@mb{\v@lXa}}%
    \v@lX=\v@lv@lXa\v@lX\v@lY=\v@lv@lXa\v@lY%
    \Figv@ctCreg#1(\v@lX,\v@lY)\multiply\v@lYa\f@ctech\global\result@t=\v@lYa}}
\def\vecunit@TD#1#2{{\Figg@tXY{#2}\getredf@ctTD\f@ctech(\v@lX,\v@lY,\v@lZ)%
    \divide\v@lX\f@ctech\divide\v@lY\f@ctech\divide\v@lZ\f@ctech%
    \Figv@ctCreg#1(\v@lX,\v@lY,\v@lZ)\n@rmeucSV{\v@lYa}{#1}%
    \invers@{\v@lXa}{\v@lYa}\edef\v@lv@lXa{\repdecn@mb{\v@lXa}}%
    \v@lX=\v@lv@lXa\v@lX\v@lY=\v@lv@lXa\v@lY\v@lZ=\v@lv@lXa\v@lZ%
    \Figv@ctCreg#1(\v@lX,\v@lY,\v@lZ)\multiply\v@lYa\f@ctech\global\result@t=\v@lYa}}
\def\vecunitC@TD[#1,#2]{\Figg@tXYa{#1}\Figg@tXY{#2}%
    \advance\v@lX-\v@lXa\advance\v@lY-\v@lYa\advance\v@lZ-\v@lZa\c@lvecunitTD}
\def\vecunitCV@TD#1{\Figg@tXY{#1}\c@lvecunitTD}
\def\c@lvecunitTD{\getredf@ctTD\f@ctech(\v@lX,\v@lY,\v@lZ)%
    \divide\v@lX\f@ctech\divide\v@lY\f@ctech\divide\v@lZ\f@ctech%
    \v@lXa=\repdecn@mb{\v@lX}\v@lX%
    \v@lYa=\repdecn@mb{\v@lY}\v@lY\v@lZa=\repdecn@mb{\v@lZ}\v@lZ%
    \advance\v@lXa\v@lYa\advance\v@lXa\v@lZa\sqrt@{\v@lYa}{\v@lXa}%
    \invers@{\v@lXa}{\v@lYa}\edef\v@lv@lXa{\repdecn@mb{\v@lXa}}%
    \v@lX=\v@lv@lXa\v@lX\v@lY=\v@lv@lXa\v@lY\v@lZ=\v@lv@lXa\v@lZ}
\def\figgetangleDD#1[#2,#3,#4]{\ifps@cri{\s@uvc@ntr@l\et@tfiggetangleDD\setc@ntr@l{2}%
    \figvectPDD-1[#2,#3]\figvectPDD-2[#2,#4]\vecunit@{-1}{-1}%
    \c@lproscalDD\delt@[-2,-1]\figvectNVDD-1[-1]\c@lproscalDD\v@leur[-2,-1]%
    \arct@n\v@lmax(\delt@,\v@leur)\v@lmax=\rdT@deg\v@lmax%
    \ifdim\v@lmax<\z@\advance\v@lmax\DePI@deg\fi\xdef#1{\repdecn@mb{\v@lmax}}%
    \resetc@ntr@l\et@tfiggetangleDD}\ignorespaces\fi}
\def\figgetangleTD#1[#2,#3,#4,#5]{\ifps@cri{\s@uvc@ntr@l\et@tfiggetangleTD\setc@ntr@l{2}%
    \figvectPTD-1[#2,#3]\figvectPTD-2[#2,#5]\figvectNVTD-3[-1,-2]%
    \figvectPTD-2[#2,#4]\figvectNVTD-4[-3,-1]%
    \vecunit@{-1}{-1}\c@lproscalTD\delt@[-2,-1]\c@lproscalTD\v@leur[-2,-4]%
    \arct@n\v@lmax(\delt@,\v@leur)\v@lmax=\rdT@deg\v@lmax%
    \ifdim\v@lmax<\z@\advance\v@lmax\DePI@deg\fi\xdef#1{\repdecn@mb{\v@lmax}}%
    \resetc@ntr@l\et@tfiggetangleTD}\ignorespaces\fi}    
\def\figgetdist#1[#2,#3]{\ifps@cri{\s@uvc@ntr@l\et@tfiggetdist\setc@ntr@l{2}%
    \figvectP-1[#2,#3]\n@rmeuc{\v@lX}{-1}\v@lX=\ptT@unit@\v@lX\xdef#1{\repdecn@mb{\v@lX}}%
    \resetc@ntr@l\et@tfiggetdist}\ignorespaces\fi}
\def\Figg@tT#1{\c@ntr@lnum{#1}%
    {\expandafter\expandafter\expandafter\extr@ctT\csname\objc@de\endcsname:%
     \ifnum\B@@ltxt=\z@\ptn@me{#1}\else\csname\objc@de T\endcsname\fi}}
\def\extr@ctT#1,#2,#3/#4:{\def\B@@ltxt{#3}}
\def\Figg@tXY#1{\c@ntr@lnum{#1}%
    \expandafter\expandafter\expandafter\extr@ctC\csname\objc@de\endcsname:}
\def\extr@ctCDD#1/#2,#3,#4:{\v@lX=#2\v@lY=#3}
\def\extr@ctCTD#1/#2,#3,#4:{\v@lX=#2\v@lY=#3\v@lZ=#4}
\def\Figg@tXYa#1{\c@ntr@lnum{#1}%
    \expandafter\expandafter\expandafter\extr@ctCa\csname\objc@de\endcsname:}
\def\extr@ctCaDD#1/#2,#3,#4:{\v@lXa=#2\v@lYa=#3}
\def\extr@ctCaTD#1/#2,#3,#4:{\v@lXa=#2\v@lYa=#3\v@lZa=#4}
\def\figinit#1{\initpr@lim\Figinit@#1,:\initpss@ttings\ignorespaces}
\def\Figinit@#1,#2:{\setunit@{#1}\def\t@xt@{#2}\ifx\t@xt@\empty\else\Figinit@@#2:\fi}
\def\Figinit@@#1#2:{\if#12 \else\Figs@tproj{#1}\initTD@\fi}
\newif\ifTr@isDim
\def\UnD@fined{UNDEFINED}
\def\ifundefined#1{\expandafter\ifx\csname#1\endcsname\relax}
\def\initpr@lim{\initb@undb@x\figsetmark{}\figsetptname{$A_{##1}$}\def\Sc@leFact{1}%
    \initDD@\figsetroundcoord{yes}\ps@critrue\expandafter\setupd@te\defaultupdate:%
    \edef\disob@unit{\UnD@fined}\edef\t@rgetpt{\UnD@fined}}
\def\initDD@{\Tr@isDimfalse%
    \let\c@lDCUn=\c@lDCUnDD%
    \let\c@lDCDeux=\c@lDCDeuxDD%
    \let\c@ldefproj=\relax%
    \let\c@lproscal=\c@lproscalDD%
    \let\c@lprojSP=\relax%
    \let\extr@ctC=\extr@ctCDD%
    \let\extr@ctCa=\extr@ctCaDD%
    \let\extr@ctCF=\extr@ctCFDD%
    \let\Figp@intreg=\Figp@intregDD%
    \let\Figpts@xes=\Figpts@xesDD%
    \let\n@rmeucSV=\n@rmeucSVDD\let\n@rmeuc=\n@rmeucDD\let\n@rminf=\n@rminfDD%
    \let\pr@dMatV=\pr@dMatVDD%
    \let\ps@xes=\ps@xesDD%
    \let\vecunit@=\vecunit@DD%
    \let\figcoord=\figcoordDD%
    \let\figgetangle=\figgetangleDD%
    \let\figpt=\figptDD%
    \let\figptBezier=\figptBezierDD%
    \let\figptbary=\figptbaryDD%
    \let\figptcirc=\figptcircDD%
    \let\figptcircumcenter=\figptcircumcenterDD%
    \let\figptcopy=\figptcopyDD%
    \let\figptcurvcenter=\figptcurvcenterDD%
    \let\figptell=\figptellDD%
    \let\figptendnormal=\figptendnormalDD%
    \def\figptinterlineplane{\un@v@ilable{figptinterlineplane}}%
    \let\figptinterlines=\inters@cDD%
    \let\figptorthocenter=\figptorthocenterDD%
    \let\figptorthoprojline=\figptorthoprojlineDD%
    \def\figptorthoprojplane{\un@v@ilable{figptorthoprojplane}}%
    \let\figptrot=\figptrotDD%
    \let\figptscontrol=\figptscontrolDD%
    \let\figptsintercirc=\figptsintercircDD%
    \let\figptsinterlinell=\figptsinterlinellDD%
    \let\figptsorthoprojline=\figptsorthoprojlineDD%
    \let\figptsrot=\figptsrotDD%
    \let\figptssym=\figptssymDD%
    \let\figptstra=\figptstraDD%
    \let\figptsym=\figptsymDD%
    \let\figpttraC=\figpttraCDD%
    \let\figpttra=\figpttraDD%
    \def\figptvisilimSL{\un@v@ilable{figptvisilimSL}}%
    \def\figsetobdist{\un@v@ilable{figsetobdist}}%
    \def\figsettarget{\un@v@ilable{figsettarget}}%
    \def\figsetview{\un@v@ilable{figsetview}}%
    \let\figvectDBezier=\figvectDBezierDD%
    \let\figvectN=\figvectNDD%
    \let\figvectNV=\figvectNVDD%
    \let\figvectP=\figvectPDD%
    \let\figvectU=\figvectUDD%
    \let\psarccircP=\psarccircPDD%
    \let\psarccirc=\psarccircDD%
    \let\psarcell=\psarcellDD%
    \let\psarcellPA=\psarcellPADD%
    \let\psarrowBezier=\psarrowBezierDD%
    \let\psarrowcircP=\psarrowcircPDD%
    \let\psarrowcirc=\psarrowcircDD%
    \let\psarrowhead=\psarrowheadDD%
    \let\psarrow=\psarrowDD%
    \let\psBezier=\psBezierDD%
    \let\pscirc=\pscircDD%
    \let\pscurve=\pscurveDD%
    \let\psnormal=\psnormalDD%
    }
\def\initTD@{\Tr@isDimtrue\initb@undb@xTD\newt@rgetptfalse\newdis@bfalse%
    \let\c@lDCUn=\c@lDCUnTD%
    \let\c@lDCDeux=\c@lDCDeuxTD%
    \let\c@ldefproj=\c@ldefprojTD%
    \let\c@lproscal=\c@lproscalTD%
    \let\extr@ctC=\extr@ctCTD%
    \let\extr@ctCa=\extr@ctCaTD%
    \let\extr@ctCF=\extr@ctCFTD%
    \let\Figp@intreg=\Figp@intregTD%
    \let\Figpts@xes=\Figpts@xesTD%
    \let\n@rmeucSV=\n@rmeucSVTD\let\n@rmeuc=\n@rmeucTD\let\n@rminf=\n@rminfTD%
    \let\pr@dMatV=\pr@dMatVTD%
    \let\ps@xes=\ps@xesTD%
    \let\vecunit@=\vecunit@TD%
    \let\figcoord=\figcoordTD%
    \let\figgetangle=\figgetangleTD%
    \let\figpt=\figptTD%
    \let\figptBezier=\figptBezierTD%
    \let\figptbary=\figptbaryTD%
    \let\figptcirc=\figptcircTD%
    \let\figptcircumcenter=\figptcircumcenterTD%
    \let\figptcopy=\figptcopyTD%
    \let\figptcurvcenter=\figptcurvcenterTD%
    \let\figptinterlineplane=\figptinterlineplaneTD%
    \let\figptinterlines=\inters@cTD%
    \let\figptorthocenter=\figptorthocenterTD%
    \let\figptorthoprojline=\figptorthoprojlineTD%
    \let\figptorthoprojplane=\figptorthoprojplaneTD%
    \let\figptrot=\figptrotTD%
    \let\figptscontrol=\figptscontrolTD%
    \let\figptsintercirc=\figptsintercircTD%
    \let\figptsorthoprojline=\figptsorthoprojlineTD%
    \let\figptsorthoprojplane=\figptsorthoprojplaneTD%
    \let\figptsrot=\figptsrotTD%
    \let\figptssym=\figptssymTD%
    \let\figptstra=\figptstraTD%
    \let\figptsym=\figptsymTD%
    \let\figpttraC=\figpttraCTD%
    \let\figpttra=\figpttraTD%
    \let\figptvisilimSL=\figptvisilimSLTD%
    \let\figsetobdist=\figsetobdistTD%
    \let\figsettarget=\figsettargetTD%
    \let\figsetview=\figsetviewTD%
    \let\figvectDBezier=\figvectDBezierTD%
    \let\figvectN=\figvectNTD%
    \let\figvectNV=\figvectNVTD%
    \let\figvectP=\figvectPTD%
    \let\figvectU=\figvectUTD%
    \let\psarccircP=\psarccircPTD%
    \let\psarccirc=\psarccircTD%
    \let\psarcell=\psarcellTD%
    \let\psarcellPA=\psarcellPATD%
    \let\psarrowBezier=\psarrowBezierTD%
    \let\psarrowcircP=\psarrowcircPTD%
    \let\psarrowcirc=\psarrowcircTD%
    \let\psarrowhead=\psarrowheadTD%
    \let\psarrow=\psarrowTD%
    \let\psBezier=\psBezierTD%
    \let\pscirc=\pscircTD%
    \let\pscurve=\pscurveTD%
    }
\def\un@v@ilable#1{\immediate\write16{*** The macro #1 is not available in the current context.}}
\def\figinsert#1{\Figinsert@#1,:\ignorespaces}
\def\Figinsert@#1,#2:{{\def\t@xt@{#2}\ifx\t@xt@\empty\xdef\Sc@leFact{1}\else%
    \def\Xarg@##1,{\def\@rgdeux{##1}}\Xarg@#2\xdef\Sc@leFact{\@rgdeux}\fi\@psfgetbb{#1}%
    \v@lX=\@psfllx\p@\v@lX=\ptpsT@pt\v@lX\v@lX=\Sc@leFact\v@lX%
    \v@lY=\@psflly\p@\v@lY=\ptpsT@pt\v@lY\v@lY=\Sc@leFact\v@lY%
    \b@undb@x{\v@lX}{\v@lY}%
    \v@lX=\@psfurx\p@\v@lX=\ptpsT@pt\v@lX\v@lX=\Sc@leFact\v@lX%
    \v@lY=\@psfury\p@\v@lY=\ptpsT@pt\v@lY\v@lY=\Sc@leFact\v@lY%
    \b@undb@x{\v@lX}{\v@lY}%
    \v@lX=\c@nt pt\v@lX=\Sc@leFact\v@lX\edef\F@ct{\repdecn@mb{\v@lX}}%
    \iftextures\special{postscriptfile #1 vscale=\F@ct\space hscale=\F@ct}%
    \else\includegraphics{#1}\fi%
    \message{[#1]}}\ignorespaces}
\def\figptDD#1:#2(#3,#4){\ifps@cri\c@ntr@lnum{#1}%
    {\v@lX=#3\unit@\v@lY=#4\unit@\Fig@dmpt{#2}{\z@}}\ignorespaces\fi}
\def\Fig@dmpt#1#2{\def\t@xt@{#1}\ifx\t@xt@\empty\def\B@@ltxt{\z@}%
    \else\expandafter\gdef\csname\objc@de T\endcsname{#1}\def\B@@ltxt{\@ne}\fi%
    \expandafter\xdef\csname\objc@de\endcsname{\ifitis@vect@r\C@dCl@svect%
    \else\C@dCl@spt\fi,\z@,\B@@ltxt/\the\v@lX,\the\v@lY,#2}}
\def\C@dCl@spt{P}
\def\C@dCl@svect{V}
\def\figptTD#1:#2(#3,#4){\ifps@cri\c@ntr@lnum{#1}%
    \def\c@@rdYZ{#4,0,0}\extrairelepremi@r\c@@rdY\de\c@@rdYZ%
    \extrairelepremi@r\c@@rdZ\de\c@@rdYZ%
    {\v@lX=#3\unit@\v@lY=\c@@rdY\unit@\v@lZ=\c@@rdZ\unit@\Fig@dmpt{#2}{\the\v@lZ}%
    \b@undb@xTD{\v@lX}{\v@lY}{\v@lZ}}\ignorespaces\fi}
\def\Figp@intregDD#1:#2(#3,#4){\c@ntr@lnum{#1}%
    {\result@t=#4\v@lX=#3\v@lY=\result@t\Fig@dmpt{#2}{\z@}}\ignorespaces}
\def\Figp@intregTD#1:#2(#3,#4){\c@ntr@lnum{#1}%
    \def\c@@rdYZ{#4,\z@,\z@}\extrairelepremi@r\c@@rdY\de\c@@rdYZ%
    \extrairelepremi@r\c@@rdZ\de\c@@rdYZ%
    {\v@lX=#3\v@lY=\c@@rdY\v@lZ=\c@@rdZ\Fig@dmpt{#2}{\the\v@lZ}%
    \b@undb@xTD{\v@lX}{\v@lY}{\v@lZ}}\ignorespaces}
\def\figptBezierDD#1:#2:#3[#4,#5,#6,#7]{\ifps@cri{\s@uvc@ntr@l\et@tfigptBezierDD%
    \FigptBezier@#3[#4,#5,#6,#7]\Figp@intregDD#1:{#2}(\v@lX,\v@lY)%
    \resetc@ntr@l\et@tfigptBezierDD}\ignorespaces\fi}
\def\figptBezierTD#1:#2:#3[#4,#5,#6,#7]{\ifps@cri{\s@uvc@ntr@l\et@tfigptBezierTD%
    \FigptBezier@#3[#4,#5,#6,#7]\Figp@intregTD#1:{#2}(\v@lX,\v@lY,\v@lZ)%
    \resetc@ntr@l\et@tfigptBezierTD}\ignorespaces\fi}
\def\FigptBezier@#1[#2,#3,#4,#5]{\setc@ntr@l{2}%
    \edef\T@{#1}\v@leur=\p@\advance\v@leur-#1pt\edef\UNmT@{\repdecn@mb{\v@leur}}%
    \figptcopy-4:/#2/\figptcopy-3:/#3/\figptcopy-2:/#4/\figptcopy-1:/#5/%
    \l@mbd@un=-4 \l@mbd@de=-\thr@@\p@rtent=\m@ne\c@lDecast%
    \l@mbd@un=-4 \l@mbd@de=-\thr@@\p@rtent=-\tw@\c@lDecast%
    \l@mbd@un=-4 \l@mbd@de=-\thr@@\p@rtent=-\thr@@\c@lDecast\Figg@tXY{-4}}
\def\c@lDCUnDD#1#2{\Figg@tXY{#1}\v@lX=\UNmT@\v@lX\v@lY=\UNmT@\v@lY%
    \Figg@tXYa{#2}\advance\v@lX\T@\v@lXa\advance\v@lY\T@\v@lYa%
    \Figp@intregDD#1:(\v@lX,\v@lY)}
\def\c@lDCUnTD#1#2{\Figg@tXY{#1}\v@lX=\UNmT@\v@lX\v@lY=\UNmT@\v@lY\v@lZ=\UNmT@\v@lZ%
    \Figg@tXYa{#2}\advance\v@lX\T@\v@lXa\advance\v@lY\T@\v@lYa\advance\v@lZ\T@\v@lZa%
    \Figp@intregTD#1:(\v@lX,\v@lY,\v@lZ)}
\def\c@lDecast{\relax\ifnum\l@mbd@un<\p@rtent\c@lDCUn{\l@mbd@un}{\l@mbd@de}%
    \advance\l@mbd@un\@ne\advance\l@mbd@de\@ne\c@lDecast\fi}
\def\figptmap#1:#2=#3/#4/#5/{\ifps@cri{\s@uvc@ntr@l\et@tfigptmap%
    \setc@ntr@l{2}\figvectP-1[#4,#3]\Figg@tXY{-1}%
    \pr@dMatV/#5/\figpttra#1:{#2}=#4/1,-1/%
    \resetc@ntr@l\et@tfigptmap}\ignorespaces\fi}
\def\pr@dMatVDD/#1,#2;#3,#4/{\v@lXa=#1\v@lX\advance\v@lXa#2\v@lY%
    \v@lYa=#3\v@lX\advance\v@lYa#4\v@lY\Figv@ctCreg-1(\v@lXa,\v@lYa)}
\def\pr@dMatVTD/#1,#2,#3;#4,#5,#6;#7,#8,#9/{%
    \v@lXa=#1\v@lX\advance\v@lXa#2\v@lY\advance\v@lXa#3\v@lZ%
    \v@lYa=#4\v@lX\advance\v@lYa#5\v@lY\advance\v@lYa#6\v@lZ%
    \v@lZa=#7\v@lX\advance\v@lZa#8\v@lY\advance\v@lZa#9\v@lZ%
    \Figv@ctCreg-1(\v@lXa,\v@lYa,\v@lZa)}
\def\figptbaryDD#1:#2[#3;#4]{\ifps@cri{\edef\list@num{#3}\extrairelepremi@r\p@int\de\list@num%
    \s@mme=\z@\@ecfor\c@ef:=#4\do{\advance\s@mme\c@ef}%
    \edef\listec@ef{#4,0}\extrairelepremi@r\c@ef\de\listec@ef%
    \Figg@tXY{\p@int}\divide\v@lX\s@mme\divide\v@lY\s@mme%
    \multiply\v@lX\c@ef\multiply\v@lY\c@ef%
    \@ecfor\p@int:=\list@num\do{\extrairelepremi@r\c@ef\de\listec@ef%
           \Figg@tXYa{\p@int}\divide\v@lXa\s@mme\divide\v@lYa\s@mme%
           \multiply\v@lXa\c@ef\multiply\v@lYa\c@ef%
           \advance\v@lX\v@lXa\advance\v@lY\v@lYa}%
    \Figp@intregDD#1:{#2}(\v@lX,\v@lY)}\ignorespaces\fi}
\def\figptbaryTD#1:#2[#3;#4]{\ifps@cri{\edef\list@num{#3}\extrairelepremi@r\p@int\de\list@num%
    \s@mme=\z@\@ecfor\c@ef:=#4\do{\advance\s@mme\c@ef}%
    \edef\listec@ef{#4,0}\extrairelepremi@r\c@ef\de\listec@ef%
    \Figg@tXY{\p@int}\divide\v@lX\s@mme\divide\v@lY\s@mme\divide\v@lZ\s@mme%
    \multiply\v@lX\c@ef\multiply\v@lY\c@ef\multiply\v@lZ\c@ef%
    \@ecfor\p@int:=\list@num\do{\extrairelepremi@r\c@ef\de\listec@ef%
           \Figg@tXYa{\p@int}\divide\v@lXa\s@mme\divide\v@lYa\s@mme\divide\v@lZa\s@mme%
           \multiply\v@lXa\c@ef\multiply\v@lYa\c@ef\multiply\v@lZa\c@ef%
           \advance\v@lX\v@lXa\advance\v@lY\v@lYa\advance\v@lZ\v@lZa}%
    \Figp@intregTD#1:{#2}(\v@lX,\v@lY,\v@lZ)}\ignorespaces\fi}
\def\figptbaryR#1:#2[#3;#4]{\ifps@cri{%
    \v@leur=\z@\@ecfor\c@ef:=#4\do{\maxim@m{\v@lmax}{\c@ef pt}{-\c@ef pt}%
    \ifdim\v@lmax>\v@leur\v@leur=\v@lmax\fi}%
    \ifdim\v@leur<\p@\f@ctech=\@M\else\ifdim\v@leur<\t@n\p@\f@ctech=\@m\else%
    \ifdim\v@leur<\c@nt\p@\f@ctech=\c@nt\else\ifdim\v@leur<\@m\p@\f@ctech=\t@n\else%
    \f@ctech=\@ne\fi\fi\fi\fi%
    \def\listec@ef{0}%
    \@ecfor\c@ef:=#4\do{\sc@lec@nvRI{\c@ef pt}\edef\listec@ef{\listec@ef,\the\s@mme}}%
    \extrairelepremi@r\c@ef\de\listec@ef\figptbary#1:#2[#3;\listec@ef]}\ignorespaces\fi}
\def\sc@lec@nvRI#1{\v@leur=#1\p@rtentiere{\s@mme}{\v@leur}\advance\v@leur-\s@mme\p@%
    \multiply\v@leur\f@ctech\p@rtentiere{\p@rtent}{\v@leur}%
    \multiply\s@mme\f@ctech\advance\s@mme\p@rtent}
\def\figptcircDD#1:#2:#3;#4(#5){\ifps@cri{\s@uvc@ntr@l\et@tfigptcircDD%
    \c@lptellDD#1:{#2}:#3;#4,#4(#5)\resetc@ntr@l\et@tfigptcircDD}\ignorespaces\fi}
\def\figptcircTD#1:#2:#3,#4,#5;#6(#7){\ifps@cri{\s@uvc@ntr@l\et@tfigptcircTD%
    \setc@ntr@l{2}\c@lExtAxes#3,#4,#5(#6)\figptellP#1:{#2}:#3,-4,-5(#7)%
    \resetc@ntr@l\et@tfigptcircTD}\ignorespaces\fi}
\def\figptcircumcenterDD#1:#2[#3,#4,#5]{\ifps@cri{\s@uvc@ntr@l\et@tfigptcircumcenterDD%
    \setc@ntr@l{2}\figvectNDD-5[#3,#4]\figptbaryDD-3:[#3,#4;1,1]%
                  \figvectNDD-6[#4,#5]\figptbaryDD-4:[#4,#5;1,1]%
    \resetc@ntr@l{2}\inters@cDD#1:{#2}[-3,-5;-4,-6]%
    \resetc@ntr@l\et@tfigptcircumcenterDD}\ignorespaces\fi}
\def\figptcircumcenterTD#1:#2[#3,#4,#5]{\ifps@cri{\s@uvc@ntr@l\et@tfigptcircumcenterTD%
    \setc@ntr@l{2}\figvectNTD-1[#3,#4,#5]%
    \figvectPTD-3[#3,#4]\figvectNVTD-5[-1,-3]\figptbaryTD-3:[#3,#4;1,1]%
    \figvectPTD-4[#4,#5]\figvectNVTD-6[-1,-4]\figptbaryTD-4:[#4,#5;1,1]%
    \resetc@ntr@l{2}\inters@cTD#1:{#2}[-3,-5;-4,-6]%
    \resetc@ntr@l\et@tfigptcircumcenterTD}\ignorespaces\fi}
\def\figptcopyDD#1:#2/#3/{\ifps@cri{\Figg@tXY{#3}%
    \Figp@intregDD#1:{#2}(\v@lX,\v@lY)}\ignorespaces\fi}
\def\figptcopyTD#1:#2/#3/{\ifps@cri{\Figg@tXY{#3}%
    \Figp@intregTD#1:{#2}(\v@lX,\v@lY,\v@lZ)}\ignorespaces\fi}
\def\figptcurvcenterDD#1:#2:#3[#4,#5,#6,#7]{\ifps@cri{\s@uvc@ntr@l\et@tfigptcurvcenterDD%
    \setc@ntr@l{2}\c@lcurvradDD#3[#4,#5,#6,#7]\edef\Sprim@{\repdecn@mb{\result@t}}%
    \figptBezierDD-1::#3[#4,#5,#6,#7]\figpttraDD#1:{#2}=-1/\Sprim@,-5/%
    \resetc@ntr@l\et@tfigptcurvcenterDD}\ignorespaces\fi}
\def\figptcurvcenterTD#1:#2:#3[#4,#5,#6,#7]{\ifps@cri{\s@uvc@ntr@l\et@tfigptcurvcenterTD%
    \setc@ntr@l{2}\figvectDBezierTD -5:1,#3[#4,#5,#6,#7]%
    \figvectDBezierTD -6:2,#3[#4,#5,#6,#7]\vecunit@TD{-5}{-5}%
    \edef\Sprim@{\repdecn@mb{\result@t}}\figvectNVTD-1[-6,-5]%
    \figvectNVTD-5[-5,-1]\c@lproscalTD\v@leur[-6,-5]%
    \invers@{\v@leur}{\v@leur}\v@leur=\Sprim@\v@leur\v@leur=\Sprim@\v@leur%
    \figptBezierTD-1::#3[#4,#5,#6,#7]\edef\Sprim@{\repdecn@mb{\v@leur}}%
    \figpttraTD#1:{#2}=-1/\Sprim@,-5/\resetc@ntr@l\et@tfigptcurvcenterTD}\ignorespaces\fi}
\def\c@lcurvradDD#1[#2,#3,#4,#5]{{\figvectDBezierDD -5:1,#1[#2,#3,#4,#5]%
    \figvectDBezierDD -6:2,#1[#2,#3,#4,#5]\vecunit@DD{-5}{-5}%
    \edef\Sprim@{\repdecn@mb{\result@t}}\figvectNVDD-5[-5]\c@lproscalDD\v@leur[-6,-5]%
    \invers@{\v@leur}{\v@leur}\v@leur=\Sprim@\v@leur\v@leur=\Sprim@\v@leur%
    \global\result@t=\v@leur}}
\def\figptellDD#1:#2:#3;#4,#5(#6,#7){\ifps@cri{\s@uvc@ntr@l\et@tfigptell%
    \c@lptellDD#1::#3;#4,#5(#6)\figptrotDD#1:{#2}=#1/#3,#7/%
    \resetc@ntr@l\et@tfigptell}\ignorespaces\fi}
\def\c@lptellDD#1:#2:#3;#4,#5(#6){\c@ssin{\C@}{\S@}{#6}\v@lmin=\C@ pt\v@lmax=\S@ pt%
    \v@lmin=#4\v@lmin\v@lmax=#5\v@lmax%
    \edef\Xc@mp{\repdecn@mb{\v@lmin}}\edef\Yc@mp{\repdecn@mb{\v@lmax}}%
    \setc@ntr@l{2}\figvectC-1(\Xc@mp,\Yc@mp)\figpttraDD#1:{#2}=#3/1,-1/}
\def\figptellP#1:#2:#3,#4,#5(#6){\ifps@cri{\s@uvc@ntr@l\et@tfigptellP%
    \setc@ntr@l{2}\figvectP-1[#3,#4]\figvectP-2[#3,#5]%
    \v@leur=#6pt\c@lptellP{#3}{-1}{-2}\figptcopy#1:{#2}/-3/%
    \resetc@ntr@l\et@tfigptellP}\ignorespaces\fi}
\def\c@lptellP#1#2#3{\edef\@ngle{\repdecn@mb\v@leur}\c@ssin{\C@}{\S@}{\@ngle}%
    \figpttra-3:=#1/\C@,#2/\figpttra-3:=-3/\S@,#3/}
\def\figptendnormalDD#1:#2:#3,#4[#5,#6]{\ifps@cri{\s@uvc@ntr@l\et@tfigptendnormal%
    \Figg@tXYa{#5}\Figg@tXY{#6}%
    \advance\v@lX-\v@lXa\advance\v@lY-\v@lYa%
    \setc@ntr@l{2}\Figv@ctCreg-1(\v@lX,\v@lY)\vecunit@{-1}{-1}\Figg@tXY{-1}%
    \delt@=#3\unit@\maxim@m{\delt@}{\delt@}{-\delt@}\edef\l@ngueur{\repdecn@mb{\delt@}}%
    \v@lX=\l@ngueur\v@lX\v@lY=\l@ngueur\v@lY%
    \delt@=\p@\advance\delt@-#4pt\edef\l@ngueur{\repdecn@mb{\delt@}}%
    \figptbaryR-1:[#5,#6;#4,\l@ngueur]\Figg@tXYa{-1}%
    \advance\v@lXa\v@lY\advance\v@lYa-\v@lX%
    \setc@ntr@l{1}\Figp@intregDD#1:{#2}(\v@lXa,\v@lYa)\resetc@ntr@l\et@tfigptendnormal}%
    \ignorespaces\fi}
\def\figptexcenter#1:#2[#3,#4,#5]{\ifps@cri{\let@xte={-}%
    \Figptexinsc@nter#1:#2[#3,#4,#5]}\ignorespaces\fi}
\def\figptincenter#1:#2[#3,#4,#5]{\ifps@cri{\let@xte={}%
    \Figptexinsc@nter#1:#2[#3,#4,#5]}\ignorespaces\fi}
\let\figptinscribedcenter=\figptincenter
\def\Figptexinsc@nter#1:#2[#3,#4,#5]{%
    \figgetdist\LA@[#4,#5]\figgetdist\LB@[#3,#5]\figgetdist\LC@[#3,#4]%
    \figptbaryR#1:{#2}[#3,#4,#5;\the\let@xte\LA@,\LB@,\LC@]}
\def\figptinterlineplaneTD#1:#2[#3,#4;#5,#6]{\ifps@cri{\s@uvc@ntr@l\et@tfigptinterlineplane%
    \setc@ntr@l{2}\figvectPTD-1[#3,#5]\vecunit@TD{-2}{#6}%
    \r@pPSTD\v@leur[-2,-1,#4]\edef\v@lcoef{\repdecn@mb{\v@leur}}%
    \figpttraTD#1:{#2}=#3/\v@lcoef,#4/\resetc@ntr@l\et@tfigptinterlineplane}\ignorespaces\fi}
\def\figptorthocenterDD#1:#2[#3,#4,#5]{\ifps@cri{\s@uvc@ntr@l\et@tfigptorthocenterDD%
    \setc@ntr@l{2}\figvectNDD-3[#3,#4]\figvectNDD-4[#4,#5]%
    \resetc@ntr@l{2}\inters@cDD#1:{#2}[#5,-3;#3,-4]%
    \resetc@ntr@l\et@tfigptorthocenterDD}\ignorespaces\fi}
\def\figptorthocenterTD#1:#2[#3,#4,#5]{\ifps@cri{\s@uvc@ntr@l\et@tfigptorthocenterTD%
    \setc@ntr@l{2}\figvectNTD-1[#3,#4,#5]%
    \figvectPTD-2[#3,#4]\figvectNVTD-3[-1,-2]%
    \figvectPTD-2[#4,#5]\figvectNVTD-4[-1,-2]%
    \resetc@ntr@l{2}\inters@cTD#1:{#2}[#5,-3;#3,-4]%
    \resetc@ntr@l\et@tfigptorthocenterTD}\ignorespaces\fi}
\def\figptorthoprojlineDD#1:#2=#3/#4,#5/{\ifps@cri{\s@uvc@ntr@l\et@tfigptorthoprojlineDD%
    \setc@ntr@l{2}\figvectPDD-3[#4,#5]\figvectNVDD-4[-3]\resetc@ntr@l{2}%
    \inters@cDD#1:{#2}[#3,-4;#4,-3]\resetc@ntr@l\et@tfigptorthoprojlineDD}\ignorespaces\fi}
\def\figptorthoprojlineTD#1:#2=#3/#4,#5/{\ifps@cri{\s@uvc@ntr@l\et@tfigptorthoprojlineTD%
    \setc@ntr@l{2}\figvectPTD-1[#4,#3]\figvectPTD-2[#4,#5]\vecunit@TD{-2}{-2}%
    \c@lproscalTD\v@leur[-1,-2]\edef\v@lcoef{\repdecn@mb{\v@leur}}%
    \figpttraTD#1:{#2}=#4/\v@lcoef,-2/\resetc@ntr@l\et@tfigptorthoprojlineTD}\ignorespaces\fi}
\def\figptorthoprojplaneTD#1:#2=#3/#4,#5/{\ifps@cri{\s@uvc@ntr@l\et@tfigptorthoprojplane%
    \setc@ntr@l{2}\figvectPTD-1[#3,#4]\vecunit@TD{-2}{#5}%
    \c@lproscalTD\v@leur[-1,-2]\edef\v@lcoef{\repdecn@mb{\v@leur}}%
    \figpttraTD#1:{#2}=#3/\v@lcoef,-2/\resetc@ntr@l\et@tfigptorthoprojplane}\ignorespaces\fi}
\def\figpthom#1:#2=#3/#4,#5/{\ifps@cri{\s@uvc@ntr@l\et@tfigpthom%
    \setc@ntr@l{2}\figvectP-1[#4,#3]\figpttra#1:{#2}=#4/#5,-1/%
    \resetc@ntr@l\et@tfigpthom}\ignorespaces\fi}
\def\figptrotDD#1:#2=#3/#4,#5/{\ifps@cri{\s@uvc@ntr@l\et@tfigptrotDD%
    \c@ssin{\C@}{\S@}{#5}\setc@ntr@l{2}\figvectPDD-1[#4,#3]\Figg@tXY{-1}%
    \v@lXa=\C@\v@lX\advance\v@lXa-\S@\v@lY%
    \v@lYa=\S@\v@lX\advance\v@lYa\C@\v@lY%
    \Figv@ctCreg-1(\v@lXa,\v@lYa)\figpttraDD#1:{#2}=#4/1,-1/%
    \resetc@ntr@l\et@tfigptrotDD}\ignorespaces\fi}
\def\figptrotTD#1:#2=#3/#4,#5,#6/{\ifps@cri{\s@uvc@ntr@l\et@tfigptrotTD%
    \c@ssin{\C@}{\S@}{#5}%
    \setc@ntr@l{2}\figptorthoprojplaneTD-3:=#4/#3,#6/\figvectPTD-2[-3,#3]%
    \n@rmeucTD\v@leur{-2}\ifdim\v@leur<\Cepsil@n\Figg@tXYa{#3}\else%
    \edef\v@lcoef{\repdecn@mb{\v@leur}}\figvectNVTD-1[#6,-2]%
    \Figg@tXYa{-1}\v@lXa=\v@lcoef\v@lXa\v@lYa=\v@lcoef\v@lYa\v@lZa=\v@lcoef\v@lZa%
    \v@lXa=\S@\v@lXa\v@lYa=\S@\v@lYa\v@lZa=\S@\v@lZa\Figg@tXY{-2}%
    \advance\v@lXa\C@\v@lX\advance\v@lYa\C@\v@lY\advance\v@lZa\C@\v@lZ%
    \Figg@tXY{-3}\advance\v@lXa\v@lX\advance\v@lYa\v@lY\advance\v@lZa\v@lZ\fi%
    \Figp@intregTD#1:{#2}(\v@lXa,\v@lYa,\v@lZa)\resetc@ntr@l\et@tfigptrotTD}\ignorespaces\fi}
\def\figptsymDD#1:#2=#3/#4,#5/{\ifps@cri{\s@uvc@ntr@l\et@tfigptsymDD%
    \resetc@ntr@l{2}\figptorthoprojlineDD-5:=#3/#4,#5/\figvectPDD-2[#3,-5]%
    \figpttraDD#1:{#2}=#3/2,-2/\resetc@ntr@l\et@tfigptsymDD}\ignorespaces\fi}
\def\figptsymTD#1:#2=#3/#4,#5/{\ifps@cri{\s@uvc@ntr@l\et@tfigptsymTD%
    \resetc@ntr@l{2}\figptorthoprojplaneTD-3:=#3/#4,#5/\figvectPTD-2[#3,-3]%
    \figpttraTD#1:{#2}=#3/2,-2/\resetc@ntr@l\et@tfigptsymTD}\ignorespaces\fi}
\def\figpttraDD#1:#2=#3/#4,#5/{\ifps@cri{\Figg@tXYa{#5}\v@lXa=#4\v@lXa\v@lYa=#4\v@lYa%
    \Figg@tXY{#3}\advance\v@lX\v@lXa\advance\v@lY\v@lYa%
    \Figp@intregDD#1:{#2}(\v@lX,\v@lY)}\ignorespaces\fi}
\def\figpttraTD#1:#2=#3/#4,#5/{\ifps@cri{\Figg@tXYa{#5}\v@lXa=#4\v@lXa\v@lYa=#4\v@lYa%
    \v@lZa=#4\v@lZa\Figg@tXY{#3}\advance\v@lX\v@lXa\advance\v@lY\v@lYa%
    \advance\v@lZ\v@lZa\Figp@intregTD#1:{#2}(\v@lX,\v@lY,\v@lZ)}\ignorespaces\fi}
\def\figpttraCDD#1:#2=#3/#4,#5/{\ifps@cri{\v@lXa=#4\unit@\v@lYa=#5\unit@%
    \Figg@tXY{#3}\advance\v@lX\v@lXa\advance\v@lY\v@lYa%
    \Figp@intregDD#1:{#2}(\v@lX,\v@lY)}\ignorespaces\fi}
\def\figpttraCTD#1:#2=#3/#4,#5,#6/{\ifps@cri{\v@lXa=#4\unit@\v@lYa=#5\unit@\v@lZa=#6\unit@%
    \Figg@tXY{#3}\advance\v@lX\v@lXa\advance\v@lY\v@lYa\advance\v@lZ\v@lZa%
    \Figp@intregTD#1:{#2}(\v@lX,\v@lY,\v@lZ)}\ignorespaces\fi}
\def\figptsaxes#1:#2(#3){\ifps@cri{\an@lys@xes#3,:\ifx\t@xt@\empty%
    \ifTr@isDim\Figpts@xes#1:#2(0,#3,0,#3,0,#3)\else\Figpts@xes#1:#2(0,#3,0,#3)\fi%
    \else\Figpts@xes#1:#2(#3)\fi}\ignorespaces\fi}
\def\Figpts@xesDD#1:#2(#3,#4,#5,#6){%
    \s@mme=#1\figpttraC\the\s@mme:$x$=#2/#4,0/%
    \advance\s@mme\@ne\figpttraC\the\s@mme:$y$=#2/0,#6/}
\def\Figpts@xesTD#1:#2(#3,#4,#5,#6,#7,#8){%
    \s@mme=#1\figpttraC\the\s@mme:$x$=#2/#4,0,0/%
    \advance\s@mme\@ne\figpttraC\the\s@mme:$y$=#2/0,#6,0/%
    \advance\s@mme\@ne\figpttraC\the\s@mme:$z$=#2/0,0,#8/}
\def\figptsmap#1=#2/#3/#4/{\ifps@cri{\s@uvc@ntr@l\et@tfigptsmap%
    \setc@ntr@l{2}\def\list@num{#2}\s@mme=#1%
    \@ecfor\p@int:=\list@num\do{\figvectP-1[#3,\p@int]\Figg@tXY{-1}%
    \pr@dMatV/#4/\figpttra\the\s@mme:=#3/1,-1/\advance\s@mme\@ne}%
    \resetc@ntr@l\et@tfigptsmap}\ignorespaces\fi}
\def\figptscontrolDD#1[#2,#3,#4,#5]{\ifps@cri{\s@uvc@ntr@l\et@tfigptscontrolDD%
    \setc@ntr@l{2}\figptcopyDD-2:/#2/\s@mme=#1\advance\s@mme\@ne%
    \Figg@tXY{-2}\v@lX=8\v@lX\v@lY=8\v@lY\Figp@intregDD-1:(\v@lX,\v@lY)%
    \figpttraDD-1:=-1/-27,#3/\figpttraDD-1:=-1/1,#5/%
    \figpttraDD-2:=-2/-27,#4/\figpttraDD-2:=-2/8,#5/%
    \c@lptCtrlDD#1[-2,-1]\c@lptCtrlDD\the\s@mme[-1,-2]%
    \resetc@ntr@l\et@tfigptscontrolDD}\ignorespaces\fi}
\def\figptscontrolTD#1[#2,#3,#4,#5]{\ifps@cri{\s@uvc@ntr@l\et@tfigptscontrolTD%
    \setc@ntr@l{2}\figptcopyTD-2:/#2/\s@mme=#1\advance\s@mme\@ne%
    \Figg@tXY{-2}\v@lX=8\v@lX\v@lY=8\v@lY\v@lZ=8\v@lZ\Figp@intregTD-1:(\v@lX,\v@lY,\v@lZ)%
    \figpttraTD-1:=-1/-27,#3/\figpttraTD-1:=-1/1,#5/%
    \figpttraTD-2:=-2/-27,#4/\figpttraTD-2:=-2/8,#5/%
    \c@lptCtrlTD#1[-2,-1]\c@lptCtrlTD\the\s@mme[-1,-2]%
    \resetc@ntr@l\et@tfigptscontrolTD}\ignorespaces\fi}
\def\figptscontrolcurve#1,#2[#3]{\ifps@cri{\s@uvc@ntr@l\et@tfigptscontrolcurve%
    \def\list@num{#3}\extrairelepremi@r\Ak@\de\list@num%
    \extrairelepremi@r\Ai@\de\list@num\extrairelepremi@r\Aj@\de\list@num%
    \s@mme=#1\figptcopy\the\s@mme:/\Ai@/%
    \setc@ntr@l{2}\figvectP -1[\Ak@,\Aj@]%
    \@ecfor\Ak@:=\list@num\do{\advance\s@mme\@ne\figpttra\the\s@mme:=\Ai@/\curv@roundness,-1/%
       \figvectP -1[\Ai@,\Ak@]\advance\s@mme\@ne\figpttra\the\s@mme:=\Aj@/-\curv@roundness,-1/%
       \advance\s@mme\@ne\figptcopy\the\s@mme:/\Aj@/%
       \edef\Ai@{\Aj@}\edef\Aj@{\Ak@}}\advance\s@mme-#1\divide\s@mme\thr@@%
       \xdef#2{\the\s@mme}%
    \resetc@ntr@l\et@tfigptscontrolcurve}\ignorespaces\fi}
\def\figptsintercircDD#1[#2,#3;#4,#5]{\ifps@cri{\s@uvc@ntr@l\et@tfigptsintercircDD%
    \setc@ntr@l{2}\let\c@lNVintc=\c@lNVintcDD\Figptsintercirc@#1[#2,#3;#4,#5]%
    \resetc@ntr@l\et@tfigptsintercircDD}\ignorespaces\fi}
\def\figptsintercircTD#1[#2,#3;#4,#5;#6]{\ifps@cri{\s@uvc@ntr@l\et@tfigptsintercircTD%
    \setc@ntr@l{2}\let\c@lNVintc=\c@lNVintcTD\vecunitC@TD[#2,#6]%
    \Figv@ctCreg-3(\v@lX,\v@lY,\v@lZ)\Figptsintercirc@#1[#2,#3;#4,#5]%
    \resetc@ntr@l\et@tfigptsintercircTD}\ignorespaces\fi}
\def\Figptsintercirc@#1[#2,#3;#4,#5]{\figvectP-1[#2,#4]%
    \vecunit@{-1}{-1}\delt@=\result@t\f@ctech=\result@tent%
    \s@mme=#1\advance\s@mme\@ne\figptcopy#1:/#2/\figptcopy\the\s@mme:/#4/%
    \ifdim\delt@=\z@\else%
    \v@lmin=#3\unit@\v@lmax=#5\unit@\v@leur=\v@lmin\advance\v@leur\v@lmax%
    \ifdim\v@leur>\delt@%
    \v@leur=\v@lmin\advance\v@leur-\v@lmax\maxim@m{\v@leur}{\v@leur}{-\v@leur}%
    \ifdim\v@leur<\delt@%
    \divide\v@lmin\f@ctech\divide\v@lmax\f@ctech\divide\delt@\f@ctech%
    \v@lmin=\repdecn@mb{\v@lmin}\v@lmin\v@lmax=\repdecn@mb{\v@lmax}\v@lmax%
    \invers@{\v@leur}{\delt@}\advance\v@lmax-\v@lmin%
    \v@lmax=-\repdecn@mb{\v@leur}\v@lmax\advance\delt@\v@lmax\delt@=.5\delt@%
    \v@lmax=\delt@\multiply\v@lmax\f@ctech%
    \edef\t@ille{\repdecn@mb{\v@lmax}}\figpttra-2:=#2/\t@ille,-1/%
    \delt@=\repdecn@mb{\delt@}\delt@\advance\v@lmin-\delt@%
    \sqrt@{\v@leur}{\v@lmin}\multiply\v@leur\f@ctech\edef\t@ille{\repdecn@mb{\v@leur}}%
    \c@lNVintc\figpttra#1:=-2/-\t@ille,-1/\figpttra\the\s@mme:=-2/\t@ille,-1/\fi\fi\fi}
\def\c@lNVintcDD{\Figg@tXY{-1}\Figv@ctCreg-1(-\v@lY,\v@lX)} 
\def\c@lNVintcTD{{\Figg@tXY{-3}\v@lmin=\v@lX\v@lmax=\v@lY\v@leur=\v@lZ%
    \Figg@tXY{-1}\c@lprovec{-3}\vecunit@{-3}{-3}
    \Figg@tXY{-1}\v@lmin=\v@lX\v@lmax=\v@lY%
    \v@leur=\v@lZ\Figg@tXY{-3}\c@lprovec{-1}}} 
\def\figptsinterlinellDD#1[#2,#3,#4,#5;#6,#7]{\ifps@cri{\s@uvc@ntr@l\et@tfigptsinterlinellDD%
    \figptcopy#1:/#6/\s@mme=#1\advance\s@mme\@ne\figptcopy\the\s@mme:/#7/%
    \v@lmin=#3\unit@\v@lmax=#4\unit@
    \setc@ntr@l{2}\figptbaryDD-4:[#6,#7;1,1]\figptsrotDD-3=-4,#7/#2,-#5/
    \Figg@tXY{-3}\Figg@tXYa{#2}\advance\v@lX-\v@lXa\advance\v@lY-\v@lYa
    \figvectP-1[-3,-2]\Figg@tXYa{-1}\figvectP-3[-4,#7]\Figptsint@rLE{#1}
    \resetc@ntr@l\et@tfigptsinterlinellDD}\ignorespaces\fi}
\def\figptsinterlinellP#1[#2,#3,#4;#5,#6]{\ifps@cri{\s@uvc@ntr@l\et@tfigptsinterlinellP%
    \figptcopy#1:/#5/\s@mme=#1\advance\s@mme\@ne\figptcopy\the\s@mme:/#6/\setc@ntr@l{2}%
    \figvectP-1[#2,#3]\vecunit@{-1}{-1}\v@lmin=\result@t
    \figvectP-2[#2,#4]\vecunit@{-2}{-2}\v@lmax=\result@t
    \figptbary-4:[#5,#6;1,1]
    \figvectP-3[#2,-4]\c@lproscal\v@lX[-3,-1]\c@lproscal\v@lY[-3,-2]
    \figvectP-3[-4,#6]\c@lproscal\v@lXa[-3,-1]\c@lproscal\v@lYa[-3,-2]
    \Figptsint@rLE{#1}\resetc@ntr@l\et@tfigptsinterlinellP}\ignorespaces\fi}
\def\Figptsint@rLE#1{%
    \getredf@ctDD\f@ctech(\v@lmin,\v@lmax)%
    \getredf@ctDD\p@rtent(\v@lX,\v@lY)\ifnum\p@rtent>\f@ctech\f@ctech=\p@rtent\fi%
    \getredf@ctDD\p@rtent(\v@lXa,\v@lYa)\ifnum\p@rtent>\f@ctech\f@ctech=\p@rtent\fi%
    \divide\v@lmin\f@ctech\divide\v@lmax\f@ctech\divide\v@lX\f@ctech\divide\v@lY\f@ctech%
    \divide\v@lXa\f@ctech\divide\v@lYa\f@ctech%
    \c@rre=\repdecn@mb\v@lXa\v@lmax\mili@u=\repdecn@mb\v@lYa\v@lmin%
    \getredf@ctDD\f@ctech(\c@rre,\mili@u)%
    \c@rre=\repdecn@mb\v@lX\v@lmax\mili@u=\repdecn@mb\v@lY\v@lmin%
    \getredf@ctDD\p@rtent(\c@rre,\mili@u)\ifnum\p@rtent>\f@ctech\f@ctech=\p@rtent\fi%
    \divide\v@lmin\f@ctech\divide\v@lmax\f@ctech\divide\v@lX\f@ctech\divide\v@lY\f@ctech%
    \divide\v@lXa\f@ctech\divide\v@lYa\f@ctech%
    \v@lmin=\repdecn@mb{\v@lmin}\v@lmin\v@lmax=\repdecn@mb{\v@lmax}\v@lmax%
    \edef\G@xde{\repdecn@mb\v@lmin}\edef\P@xde{\repdecn@mb\v@lmax}%
    \c@rre=-\v@lmax\v@leur=\repdecn@mb\v@lY\v@lY\advance\c@rre\v@leur\c@rre=\G@xde\c@rre%
    \v@leur=\repdecn@mb\v@lX\v@lX\v@leur=\P@xde\v@leur\advance\c@rre\v@leur
    \v@lmin=\repdecn@mb\v@lYa\v@lmin\v@lmax=\repdecn@mb\v@lXa\v@lmax%
    \mili@u=\repdecn@mb\v@lX\v@lmax\advance\mili@u\repdecn@mb\v@lY\v@lmin
    \v@lmax=\repdecn@mb\v@lXa\v@lmax\advance\v@lmax\repdecn@mb\v@lYa\v@lmin
    \ifdim\v@lmax>\epsil@n%
    \maxim@m{\v@leur}{\c@rre}{-\c@rre}\maxim@m{\v@lmin}{\mili@u}{-\mili@u}%
    \maxim@m{\v@leur}{\v@leur}{\v@lmin}\maxim@m{\v@lmin}{\v@lmax}{-\v@lmax}%
    \maxim@m{\v@leur}{\v@leur}{\v@lmin}\p@rtentiere{\p@rtent}{\v@leur}\advance\p@rtent\@ne%
    \divide\c@rre\p@rtent\divide\mili@u\p@rtent\divide\v@lmax\p@rtent%
    \delt@=\repdecn@mb{\mili@u}\mili@u\v@leur=\repdecn@mb{\v@lmax}\c@rre%
    \advance\delt@-\v@leur\ifdim\delt@<\z@\else\sqrt@\delt@\delt@%
    \invers@\v@lmax\v@lmax\edef\Uns@rAp{\repdecn@mb\v@lmax}%
    \v@leur=-\mili@u\advance\v@leur-\delt@\v@leur=\Uns@rAp\v@leur%
    \edef\t@ille{\repdecn@mb\v@leur}\figpttra#1:=-4/\t@ille,-3/\s@mme=#1\advance\s@mme\@ne%
    \v@leur=-\mili@u\advance\v@leur\delt@\v@leur=\Uns@rAp\v@leur%
    \edef\t@ille{\repdecn@mb\v@leur}\figpttra\the\s@mme:=-4/\t@ille,-3/\fi\fi}
\def\figptsorthoprojlineDD#1=#2/#3,#4/{\ifps@cri{\s@uvc@ntr@l\et@tfigptsorthoprojlineDD%
    \setc@ntr@l{2}\figvectPDD-3[#3,#4]\figvectNVDD-4[-3]\resetc@ntr@l{2}%
    \def\list@num{#2}\s@mme=#1\@ecfor\p@int:=\list@num\do{%
    \inters@cDD\the\s@mme:[\p@int,-4;#3,-3]\advance\s@mme\@ne}%
    \resetc@ntr@l\et@tfigptsorthoprojlineDD}\ignorespaces\fi}
\def\figptsorthoprojlineTD#1=#2/#3,#4/{\ifps@cri{\s@uvc@ntr@l\et@tfigptsorthoprojlineTD%
    \setc@ntr@l{2}\figvectPTD-2[#3,#4]\vecunit@TD{-2}{-2}%
    \def\list@num{#2}\s@mme=#1\@ecfor\p@int:=\list@num\do{%
    \figvectPTD-1[#3,\p@int]\c@lproscalTD\v@leur[-1,-2]%
    \edef\v@lcoef{\repdecn@mb{\v@leur}}\figpttraTD\the\s@mme:=#3/\v@lcoef,-2/%
    \advance\s@mme\@ne}\resetc@ntr@l\et@tfigptsorthoprojlineTD}\ignorespaces\fi}
\def\figptsorthoprojplaneTD#1=#2/#3,#4/{\ifps@cri{\s@uvc@ntr@l\et@tfigptsorthoprojplane%
    \setc@ntr@l{2}\vecunit@TD{-2}{#4}%
    \def\list@num{#2}\s@mme=#1\@ecfor\p@int:=\list@num\do{\figvectPTD-1[\p@int,#3]%
    \c@lproscalTD\v@leur[-1,-2]\edef\v@lcoef{\repdecn@mb{\v@leur}}%
    \figpttraTD\the\s@mme:=\p@int/\v@lcoef,-2/\advance\s@mme\@ne}%
    \resetc@ntr@l\et@tfigptsorthoprojplane}\ignorespaces\fi}
\def\figptshom#1=#2/#3,#4/{\ifps@cri{\s@uvc@ntr@l\et@tfigptshom%
    \setc@ntr@l{2}\def\list@num{#2}\s@mme=#1%
    \@ecfor\p@int:=\list@num\do{\figvectP-1[#3,\p@int]%
    \figpttra\the\s@mme:=#3/#4,-1/\advance\s@mme\@ne}%
    \resetc@ntr@l\et@tfigptshom}\ignorespaces\fi}
\def\figptsrotDD#1=#2/#3,#4/{\ifps@cri{\s@uvc@ntr@l\et@tfigptsrotDD%
    \c@ssin{\C@}{\S@}{#4}\setc@ntr@l{2}\def\list@num{#2}\s@mme=#1%
    \@ecfor\p@int:=\list@num\do{\figvectPDD-1[#3,\p@int]\Figg@tXY{-1}%
    \v@lXa=\C@\v@lX\advance\v@lXa-\S@\v@lY%
    \v@lYa=\S@\v@lX\advance\v@lYa\C@\v@lY%
    \Figv@ctCreg-1(\v@lXa,\v@lYa)\figpttraDD\the\s@mme:=#3/1,-1/\advance\s@mme\@ne}%
    \resetc@ntr@l\et@tfigptsrotDD}\ignorespaces\fi}
\def\figptsrotTD#1=#2/#3,#4,#5/{\ifps@cri{\s@uvc@ntr@l\et@tfigptsrotTD%
    \c@ssin{\C@}{\S@}{#4}%
    \setc@ntr@l{2}\def\list@num{#2}\s@mme=#1%
    \@ecfor\p@int:=\list@num\do{\figptorthoprojplaneTD-3:=#3/\p@int,#5/%
    \figvectPTD-2[-3,\p@int]%
    \figvectNVTD-1[#5,-2]\n@rmeucTD\v@leur{-2}\edef\v@lcoef{\repdecn@mb{\v@leur}}%
    \Figg@tXYa{-1}\v@lXa=\v@lcoef\v@lXa\v@lYa=\v@lcoef\v@lYa\v@lZa=\v@lcoef\v@lZa%
    \v@lXa=\S@\v@lXa\v@lYa=\S@\v@lYa\v@lZa=\S@\v@lZa\Figg@tXY{-2}%
    \advance\v@lXa\C@\v@lX\advance\v@lYa\C@\v@lY\advance\v@lZa\C@\v@lZ%
    \Figg@tXY{-3}\advance\v@lXa\v@lX\advance\v@lYa\v@lY\advance\v@lZa\v@lZ%
    \Figp@intregTD\the\s@mme:(\v@lXa,\v@lYa,\v@lZa)\advance\s@mme\@ne}%
    \resetc@ntr@l\et@tfigptsrotTD}\ignorespaces\fi}
\def\figptssymDD#1=#2/#3,#4/{\ifps@cri{\s@uvc@ntr@l\et@tfigptssymDD%
    \setc@ntr@l{2}\figvectPDD-3[#3,#4]\Figg@tXY{-3}\Figv@ctCreg-4(-\v@lY,\v@lX)%
    \resetc@ntr@l{2}\def\list@num{#2}\s@mme=#1%
    \@ecfor\p@int:=\list@num\do{\inters@cDD-5:[#3,-3;\p@int,-4]\figvectPDD-2[\p@int,-5]%
    \figpttraDD\the\s@mme:=\p@int/2,-2/\advance\s@mme\@ne}%
    \resetc@ntr@l\et@tfigptssymDD}\ignorespaces\fi}
\def\figptssymTD#1=#2/#3,#4/{\ifps@cri{\s@uvc@ntr@l\et@tfigptssymTD%
    \setc@ntr@l{2}\vecunit@TD{-2}{#4}\def\list@num{#2}\s@mme=#1%
    \@ecfor\p@int:=\list@num\do{\figvectPTD-1[\p@int,#3]%
    \c@lproscalTD\v@leur[-1,-2]\v@leur=2\v@leur\edef\v@lcoef{\repdecn@mb{\v@leur}}%
    \figpttraTD\the\s@mme:=\p@int/\v@lcoef,-2/\advance\s@mme\@ne}%
    \resetc@ntr@l\et@tfigptssymTD}\ignorespaces\fi}
\def\figptstraDD#1=#2/#3,#4/{\ifps@cri{\Figg@tXYa{#4}\v@lXa=#3\v@lXa\v@lYa=#3\v@lYa%
    \def\list@num{#2}\s@mme=#1\@ecfor\p@int:=\list@num\do{\Figg@tXY{\p@int}%
    \advance\v@lX\v@lXa\advance\v@lY\v@lYa%
    \Figp@intregDD\the\s@mme:(\v@lX,\v@lY)\advance\s@mme\@ne}}\ignorespaces\fi}
\def\figptstraTD#1=#2/#3,#4/{\ifps@cri{\Figg@tXYa{#4}\v@lXa=#3\v@lXa\v@lYa=#3\v@lYa%
    \v@lZa=#3\v@lZa\def\list@num{#2}\s@mme=#1\@ecfor\p@int:=\list@num\do{\Figg@tXY{\p@int}%
    \advance\v@lX\v@lXa\advance\v@lY\v@lYa\advance\v@lZ\v@lZa%
    \Figp@intregTD\the\s@mme:(\v@lX,\v@lY,\v@lZ)\advance\s@mme\@ne}}\ignorespaces\fi}
\def\figptvisilimSLTD#1:#2[#3,#4;#5,#6]{\ifps@cri{\s@uvc@ntr@l\et@tfigptvisilimSLTD%
    \setc@ntr@l{2}\figvectP-1[#3,#4]\n@rminf{\delt@}{-1}%
    \ifcase\curr@ntproj\v@lX=\cxa@\p@\v@lY=-\p@\v@lZ=\cxb@\p@
    \Figv@ctCreg-2(\v@lX,\v@lY,\v@lZ)\figvectP-3[#5,#6]\figvectNV-1[-2,-3]%
    \or\figvectP-1[#5,#6]\vecunitCV@TD{-1}\v@lmin=\v@lX\v@lmax=\v@lY
    \v@leur=\v@lZ\v@lX=\cza@\p@\v@lY=\czb@\p@\v@lZ=\czc@\p@\c@lprovec{-1}%
    \or\c@ley@pt{-2}\figvectN-1[#5,#6,-2]\fi
    \edef\Ai@{#3}\edef\Aj@{#4}\figvectP-2[#5,\Ai@]\c@lproscal\v@leur[-1,-2]%
    \ifdim\v@leur>\z@\p@rtent=\@ne\else\p@rtent=\m@ne\fi%
    \figvectP-2[#5,\Aj@]\c@lproscal\v@leur[-1,-2]%
    \ifdim\p@rtent\v@leur>\z@\figptcopy#1:#2/#3/%
    \message{*** \BS@ figptvisilimSL: points are on the same side.}\else%
    \figptcopy-3:/#3/\figptcopy-4:/#4/%
    \loop\figptbary-5:[-3,-4;1,1]\figvectP-2[#5,-5]\c@lproscal\v@leur[-1,-2]%
    \ifdim\p@rtent\v@leur>\z@\figptcopy-3:/-5/\else\figptcopy-4:/-5/\fi%
    \divide\delt@\tw@\ifdim\delt@>\epsil@n\repeat%
    \figptbary#1:#2[-3,-4;1,1]\fi\resetc@ntr@l\et@tfigptvisilimSLTD}\ignorespaces\fi}
\def\c@ley@pt#1{\t@stp@r\ifitis@K\v@lX=\cza@\p@\v@lY=\czb@\p@\v@lZ=\czc@\p@%
    \Figv@ctCreg-1(\v@lX,\v@lY,\v@lZ)\Figp@intreg-2:(\wd\Bt@rget,\ht\Bt@rget,\dp\Bt@rget)%
    \figpttra#1:=-2/-\disob@intern,-1/\else\end\fi}
\def\t@stp@r{\itis@Ktrue\ifnewt@rgetpt\else\itis@Kfalse%
    \message{*** \BS@ figptvisilimXX: target point undefined.}\fi\ifnewdis@b\else%
    \itis@Kfalse\message{*** \BS@ figptvisilimXX: observation distance undefined.}\fi%
    \ifitis@K\else\message{*** This macro must be called after \BS@ psbeginfig or after
    having set the missing parameter(s) with \BS@ figset proj()}\fi}
\def\figscan#1(#2,#3){{\s@uvc@ntr@l\et@tfigscan\@psfgetbb{#1}\if@psfbbfound\else%
    \def\@psfllx{0}\def\@psflly{20}\def\@psfurx{540}\def\@psfury{640}\fi%
    \unit@=\@ne bp\setc@ntr@l{2}\figsetmark{}%
    \def\minst@p{20pt}%
    \v@lX=\@psfllx\p@\v@lX=\Sc@leFact\v@lX\r@undint\v@lX\v@lX%
    \v@lY=\@psflly\p@\v@lY=\Sc@leFact\v@lY\ifdim\v@lY>\z@\r@undint\v@lY\v@lY\fi%
    \delt@=\@psfury\p@\delt@=\Sc@leFact\delt@%
    \advance\delt@-\v@lY\v@lXa=\@psfurx\p@\v@lXa=\Sc@leFact\v@lXa\v@leur=\minst@p%
    \edef\valv@lY{\repdecn@mb{\v@lY}}\edef\LgTr@it{\the\delt@}%
    \loop\ifdim\v@lX<\v@lXa\edef\valv@lX{\repdecn@mb{\v@lX}}%
    \figptDD -1:(\valv@lX,\valv@lY)\figwriten -1:\hbox{\vrule height\LgTr@it}(0)%
    \ifdim\v@leur<\minst@p\else\figsetmark{\raise-8bp\hbox{$\scriptscriptstyle\triangle$}}%
    \figwrites -1:\@ffichnb{0}{\valv@lX}(6)\v@leur=\z@\figsetmark{}\fi%
    \advance\v@leur#2pt\advance\v@lX#2pt\repeat%
    \def\minst@p{10pt}%
    \v@lX=\@psfllx\p@\v@lX=\Sc@leFact\v@lX\ifdim\v@lX>\z@\r@undint\v@lX\v@lX\fi%
    \v@lY=\@psflly\p@\v@lY=\Sc@leFact\v@lY\r@undint\v@lY\v@lY%
    \delt@=\@psfurx\p@\delt@=\Sc@leFact\delt@%
    \advance\delt@-\v@lX\v@lYa=\@psfury\p@\v@lYa=\Sc@leFact\v@lYa\v@leur=\minst@p%
    \edef\valv@lX{\repdecn@mb{\v@lX}}\edef\LgTr@it{\the\delt@}%
    \loop\ifdim\v@lY<\v@lYa\edef\valv@lY{\repdecn@mb{\v@lY}}%
    \figptDD -1:(\valv@lX,\valv@lY)\figwritee -1:\vbox{\hrule width\LgTr@it}(0)%
    \ifdim\v@leur<\minst@p\else\figsetmark{$\triangleright$\kern4bp}%
    \figwritew -1:\@ffichnb{0}{\valv@lY}(6)\v@leur=\z@\figsetmark{}\fi%
    \advance\v@leur#3pt\advance\v@lY#3pt\repeat%
    \resetc@ntr@l\et@tfigscan}\ignorespaces}
\def\figshowpts[#1,#2]{{\figsetmark{$\bullet$}\figsetptname{\bf ##1}%
    \p@rtent=#2\relax\ifnum\p@rtent<\z@\p@rtent=\z@\fi%
    \s@mme=#1\relax\ifnum\s@mme<\z@\s@mme=\z@\fi%
    \loop\ifnum\s@mme<\p@rtent\pt@rvect{\s@mme}%
    \ifitis@K\figwriten{\the\s@mme}:(4pt)\fi\advance\s@mme\@ne\repeat%
    \pt@rvect{\s@mme}\ifitis@K\figwriten{\the\s@mme}:(4pt)\fi}\ignorespaces}
\def\pt@rvect#1{\set@bjc@de{#1}%
    \expandafter\expandafter\expandafter\inqpt@rvec\csname\objc@de\endcsname:}
\def\inqpt@rvec#1#2:{\if#1\C@dCl@spt\itis@Ktrue\else\itis@Kfalse\fi}
\def\figshowsettings{{%
    \immediate\write16{====================================================================}%
    \immediate\write16{ Current settings about:}%
    \immediate\write16{ --- GENERAL ---}%
    \immediate\write16{Scale factor and Unit = \unit@util\space (\the\unit@)
     \space -> \BS@ figinit{ScaleFactorUnit}}%
    \immediate\write16{Update mode = \ifpstestm@de yes\else no\fi
     \space-> \BS@ psset(update=yes/no) or \BS@ pssetdefault(update=yes/no)}%
    \immediate\write16{ --- PRINTING ---}%
    \immediate\write16{Implicit point name = \ptn@me{i} \space-> \BS@ figsetptname{Name}}%
    \immediate\write16{Point marker = \the\c@nsymb \space -> \BS@ figsetmark{Mark}}%
    \immediate\write16{Print rounded coordinates = \ifr@undcoord yes\else no\fi
     \space-> \BS@ figsetroundcoord{yes/no}}%
    \immediate\write16{ --- POSTSCRIPT (general) ---}%
    \immediate\write16{First-level (or primary) settings:}%
    \immediate\write16{ Color = \curr@ntcolor \space-> \BS@ psset(color=ColorDefinition)}%
    \immediate\write16{ Filling mode = \iffillm@de yes\else no\fi
     \space-> \BS@ psset(fillmode=yes/no)}%
    \immediate\write16{ Line join = \curr@ntjoin \space-> \BS@ psset(join=miter/round/bevel)}%
    \immediate\write16{ Line style = \curr@ntdash \space-> \BS@ psset(dash=Index/Pattern)}%
    \immediate\write16{ Line width = \curr@ntwidth
     \space-> \BS@ psset(width=real in PostScript units)}%
    \immediate\write16{Second-level (or secondary) settings:}%
    \immediate\write16{ Color = \sec@ndcolor \space-> \BS@ psset second(color=ColorDefinition)}%
    \immediate\write16{ Line style = \curr@ntseconddash
     \space-> \BS@ psset second(dash=Index/Pattern)}%
    \immediate\write16{ Line width = \curr@ntsecondwidth
     \space-> \BS@ psset second(width=real in PostScript units)}%
    \immediate\write16{Third-level (or ternary) settings:}%
    \immediate\write16{ Color = \th@rdcolor \space-> \BS@ psset third(color=ColorDefinition)}%
    \immediate\write16{ --- POSTSCRIPT (specific) ---}%
    \immediate\write16{Arrow-head:}%
    \immediate\write16{ (half-)Angle = \@rrowheadangle
     \space-> \BS@ psset arrowhead(angle=real in degrees)}%
    \immediate\write16{ Filling mode = \if@rrowhfill yes\else no\fi
     \space-> \BS@ psset arrowhead(fillmode=yes/no)}%
    \immediate\write16{ "Outside" = \if@rrowhout yes\else no\fi
     \space-> \BS@ psset arrowhead(out=yes/no)}%
    \immediate\write16{ Length = \@rrowheadlength
     \if@rrowratio\space(not active)\else\space(active)\fi
     \space-> \BS@ psset arrowhead(length=real in user coord.)}%
    \immediate\write16{ Ratio = \@rrowheadratio
     \if@rrowratio\space(active)\else\space(not active)\fi
     \space-> \BS@ psset arrowhead(ratio=real in [0,1])}%
    \immediate\write16{Curve: Roundness = \curv@roundness
     \space-> \BS@ psset curve(roundness=real in [0,0.5])}%
    \immediate\write16{Mesh: Diagonal = \the\c@ntrolmesh
     \space-> \BS@ psset mesh(diag=integer in {-1,0,1})}%
    \immediate\write16{Flow chart:}%
    \immediate\write16{ Arrow position = \@rrowp@s
     \space-> \BS@ psset flowchart(arrowposition=real in [0,1])}%
    \immediate\write16{ Arrow reference point = \ifcase\@rrowr@fpt start\else end\fi
     \space-> \BS@ psset flowchart(arrowrefpt = start/end)}%
    \immediate\write16{ Line type = \ifcase\fclin@typ@ curve\else polygon\fi
     \space-> \BS@ psset flowchart(line=polygon/curve)}%
    \immediate\write16{ Padding = (\Xp@dd, \Yp@dd)
     \space-> \BS@ psset flowchart(padding = real in user coord.)}%
    \immediate\write16{\space\space\space\space(or
     \BS@ psset flowchart(xpadding=real, ypadding=real) )}%
    \immediate\write16{ Radius = \fclin@r@d
     \space-> \BS@ psset flowchart(radius=positive real in user coord.)}%
    \immediate\write16{ Shape = \fcsh@pe
     \space-> \BS@ psset flowchart(shape = rectangle, ellipse or lozenge)}%
    \immediate\write16{ Thickness = \thickn@ss
     \space-> \BS@ psset flowchart(thickness = real in user coord.)}%
    \ifTr@isDim%
    \immediate\write16{ --- 3D to 2D PROJECTION ---}%
    \immediate\write16{Projection : \typ@proj \space-> \BS@ figinit{ScaleFactorUnit, ProjType}}%
    \immediate\write16{Longitude (psi) = \v@lPsi \space-> \BS@ figset proj(psi=real in degrees)}%
    \ifcase\curr@ntproj\immediate\write16{Depth coeff. (Lambda)
     \space = \v@lTheta \space-> \BS@ figset proj(lambda=real in [0,1])}%
    \else\immediate\write16{Latitude (theta)
     \space = \v@lTheta \space-> \BS@ figset proj(theta=real in degrees)}%
    \fi%
    \ifnum\curr@ntproj=\tw@%
    \immediate\write16{Observation distance = \disob@unit
     \space-> \BS@ figset proj(dist=real in user coord.)}%
    \immediate\write16{Target point = \t@rgetpt \space-> \BS@ figset proj(targetpt=pt number)}%
     \v@lX=\ptT@unit@\wd\Bt@rget\v@lY=\ptT@unit@\ht\Bt@rget\v@lZ=\ptT@unit@\dp\Bt@rget%
    \immediate\write16{ Its coordinates are
     (\repdecn@mb{\v@lX}, \repdecn@mb{\v@lY}, \repdecn@mb{\v@lZ})}%
    \fi%
    \fi%
    \immediate\write16{====================================================================}%
    \ignorespaces}}
{\catcode`\/=0 \catcode`/\=12 /gdef/BS@{\}}
\newif\ifitis@vect@r
\def\figvectC#1(#2,#3){{\itis@vect@rtrue\figpt#1:(#2,#3)}\ignorespaces}
\def\Figv@ctCreg#1(#2,#3){{\itis@vect@rtrue\Figp@intreg#1:(#2,#3)}\ignorespaces}
\def\figvectDBezierDD#1:#2,#3[#4,#5,#6,#7]{\ifps@cri{\s@uvc@ntr@l\et@tfigvectDBezierDD%
    \FigvectDBezier@#2,#3[#4,#5,#6,#7]\v@lX=\c@ef\v@lX\v@lY=\c@ef\v@lY%
    \Figv@ctCreg#1(\v@lX,\v@lY)\resetc@ntr@l\et@tfigvectDBezierDD}\ignorespaces\fi}
\def\figvectDBezierTD#1:#2,#3[#4,#5,#6,#7]{\ifps@cri{\s@uvc@ntr@l\et@tfigvectDBezierTD%
    \FigvectDBezier@#2,#3[#4,#5,#6,#7]\v@lX=\c@ef\v@lX\v@lY=\c@ef\v@lY\v@lZ=\c@ef\v@lZ%
    \Figv@ctCreg#1(\v@lX,\v@lY,\v@lZ)\resetc@ntr@l\et@tfigvectDBezierTD}\ignorespaces\fi}
\def\FigvectDBezier@#1,#2[#3,#4,#5,#6]{\setc@ntr@l{2}%
    \edef\T@{#2}\v@leur=\p@\advance\v@leur-#2pt\edef\UNmT@{\repdecn@mb{\v@leur}}%
    \ifnum#1=\tw@\def\c@ef{6}\else\def\c@ef{3}\fi%
    \figptcopy-4:/#3/\figptcopy-3:/#4/\figptcopy-2:/#5/\figptcopy-1:/#6/%
    \l@mbd@un=-4 \l@mbd@de=-\thr@@\p@rtent=\m@ne\c@lDecast%
    \ifnum#1=\tw@\c@lDCDeux{-4}{-3}\c@lDCDeux{-3}{-2}\c@lDCDeux{-4}{-3}\else%
    \l@mbd@un=-4 \l@mbd@de=-\thr@@\p@rtent=-\tw@\c@lDecast%
    \c@lDCDeux{-4}{-3}\fi\Figg@tXY{-4}}
\def\c@lDCDeuxDD#1#2{\Figg@tXY{#2}\Figg@tXYa{#1}%
    \advance\v@lX-\v@lXa\advance\v@lY-\v@lYa\Figp@intregDD#1:(\v@lX,\v@lY)}
\def\c@lDCDeuxTD#1#2{\Figg@tXY{#2}\Figg@tXYa{#1}\advance\v@lX-\v@lXa%
    \advance\v@lY-\v@lYa\advance\v@lZ-\v@lZa\Figp@intregTD#1:(\v@lX,\v@lY,\v@lZ)}
\def\figvectNDD#1[#2,#3]{\ifps@cri{\Figg@tXYa{#2}\Figg@tXY{#3}%
    \advance\v@lX-\v@lXa\advance\v@lY-\v@lYa%
    \Figv@ctCreg#1(-\v@lY,\v@lX)}\ignorespaces\fi}
\def\figvectNTD#1[#2,#3,#4]{\ifps@cri{\vecunitC@TD[#2,#4]\v@lmin=\v@lX\v@lmax=\v@lY%
    \v@leur=\v@lZ\vecunitC@TD[#2,#3]\c@lprovec{#1}}\ignorespaces\fi}
\def\figvectNVDD#1[#2]{\ifps@cri{\Figg@tXY{#2}\Figv@ctCreg#1(-\v@lY,\v@lX)}\ignorespaces\fi}
\def\figvectNVTD#1[#2,#3]{\ifps@cri{\vecunitCV@TD{#3}\v@lmin=\v@lX\v@lmax=\v@lY%
    \v@leur=\v@lZ\vecunitCV@TD{#2}\c@lprovec{#1}}\ignorespaces\fi}
\def\figvectPDD#1[#2,#3]{\ifps@cri{\Figg@tXYa{#2}\Figg@tXY{#3}%
    \advance\v@lX-\v@lXa\advance\v@lY-\v@lYa%
    \Figv@ctCreg#1(\v@lX,\v@lY)}\ignorespaces\fi}
\def\figvectPTD#1[#2,#3]{\ifps@cri{\Figg@tXYa{#2}\Figg@tXY{#3}%
    \advance\v@lX-\v@lXa\advance\v@lY-\v@lYa\advance\v@lZ-\v@lZa%
    \Figv@ctCreg#1(\v@lX,\v@lY,\v@lZ)}\ignorespaces\fi}
\def\figvectUDD#1[#2]{\ifps@cri{\n@rmeuc\v@leur{#2}\invers@\v@leur\v@leur%
    \delt@=\repdecn@mb{\v@leur}\unit@\edef\v@ldelt@{\repdecn@mb{\delt@}}%
    \Figg@tXY{#2}\v@lX=\v@ldelt@\v@lX\v@lY=\v@ldelt@\v@lY%
    \Figv@ctCreg#1(\v@lX,\v@lY)}\ignorespaces\fi}
\def\figvectUTD#1[#2]{\ifps@cri{\n@rmeuc\v@leur{#2}\invers@\v@leur\v@leur%
    \delt@=\repdecn@mb{\v@leur}\unit@\edef\v@ldelt@{\repdecn@mb{\delt@}}%
    \Figg@tXY{#2}\v@lX=\v@ldelt@\v@lX\v@lY=\v@ldelt@\v@lY\v@lZ=\v@ldelt@\v@lZ%
    \Figv@ctCreg#1(\v@lX,\v@lY,\v@lZ)}\ignorespaces\fi}
\def\figvisu#1#2#3{\c@ldefproj\initb@undb@x\setbox\b@xvisu=\hbox{\ignorespaces#3}%
    \v@lXa=-\c@@rdYmin\v@lYa=\c@@rdYmax\advance\v@lYa-\c@@rdYmin%
    \v@lX=\c@@rdXmax\advance\v@lX-\c@@rdXmin%
    \setbox#1=\hbox{#2}\v@lY=-\v@lX\maxim@m{\v@lX}{\v@lX}{\wd#1}%
    \advance\v@lY\v@lX\divide\v@lY\tw@\advance\v@lY-\c@@rdXmin%
    \setbox#1=\vbox{\parindent0mm\hsize=\v@lX\vskip\v@lYa%
    \rlap{\hskip\v@lY\smash{\raise\v@lXa\box\b@xvisu}}%
    \def\t@xt@{#2}\ifx\t@xt@\empty\else\medskip\centerline{#2}\fi}\wd#1=\v@lX}
\newbox\Bt@rget\setbox\Bt@rget=\null
\newbox\BminTD@\setbox\BminTD@=\null
\newbox\BmaxTD@\setbox\BmaxTD@=\null
\newif\ifnewt@rgetpt\newif\ifnewdis@b
\def\b@undb@xTD#1#2#3{%
    \relax\ifdim#1<\wd\BminTD@\global\wd\BminTD@=#1\fi%
    \relax\ifdim#2<\ht\BminTD@\global\ht\BminTD@=#2\fi%
    \relax\ifdim#3<\dp\BminTD@\global\dp\BminTD@=#3\fi%
    \relax\ifdim#1>\wd\BmaxTD@\global\wd\BmaxTD@=#1\fi%
    \relax\ifdim#2>\ht\BmaxTD@\global\ht\BmaxTD@=#2\fi%
    \relax\ifdim#3>\dp\BmaxTD@\global\dp\BmaxTD@=#3\fi}
\def\c@ldefdisob{{\ifdim\wd\BminTD@<\maxdimen\v@leur=\wd\BmaxTD@\advance\v@leur-\wd\BminTD@%
    \delt@=\ht\BmaxTD@\advance\delt@-\ht\BminTD@\maxim@m{\v@leur}{\v@leur}{\delt@}%
    \delt@=\dp\BmaxTD@\advance\delt@-\dp\BminTD@\maxim@m{\v@leur}{\v@leur}{\delt@}%
    \v@leur=5\v@leur\else\v@leur=800pt\fi\c@ldefdisob@{\v@leur}}}
\def\c@ldefdisob@#1{{\v@leur=#1\ifdim\v@leur<\p@\v@leur=800pt\fi%
    \xdef\disob@intern{\repdecn@mb{\v@leur}}%
    \delt@=\ptT@unit@\v@leur\xdef\disob@unit{\repdecn@mb{\delt@}}%
    \f@ctech=\@ne\loop\ifdim\v@leur>\t@n pt\divide\v@leur\t@n\multiply\f@ctech\t@n\repeat%
    \xdef\disob@{\repdecn@mb{\v@leur}}\xdef\divf@ctproj{\the\f@ctech}}%
    \global\newdis@btrue}
\def\c@ldeft@rgetpt{\newt@rgetpttrue\def\t@rgetpt{CenterBoundBox}{%
    \delt@=\wd\BmaxTD@\advance\delt@-\wd\BminTD@\divide\delt@\tw@%
    \v@leur=\wd\BminTD@\advance\v@leur\delt@\global\wd\Bt@rget=\v@leur%
    \delt@=\ht\BmaxTD@\advance\delt@-\ht\BminTD@\divide\delt@\tw@%
    \v@leur=\ht\BminTD@\advance\v@leur\delt@\global\ht\Bt@rget=\v@leur%
    \delt@=\dp\BmaxTD@\advance\delt@-\dp\BminTD@\divide\delt@\tw@%
    \v@leur=\dp\BminTD@\advance\v@leur\delt@\global\dp\Bt@rget=\v@leur}}
\def\c@ldefprojTD{\ifnewt@rgetpt\else\c@ldeft@rgetpt\fi\ifnewdis@b\else\c@ldefdisob\fi}
\def\c@lprojcav{
    \v@lZa=\cxa@\v@lY\advance\v@lX\v@lZa%
    \v@lZa=\cxb@\v@lY\v@lY=\v@lZ\advance\v@lY\v@lZa\ignorespaces}
\def\c@lprojrea{
    \advance\v@lX-\wd\Bt@rget\advance\v@lY-\ht\Bt@rget\advance\v@lZ-\dp\Bt@rget%
    \v@lZa=\cza@\v@lX\advance\v@lZa\czb@\v@lY\advance\v@lZa\czc@\v@lZ%
    \divide\v@lZa\divf@ctproj\advance\v@lZa\disob@ pt\invers@{\v@lZa}{\v@lZa}%
    \v@lZa=\disob@\v@lZa\edef\v@lcoef{\repdecn@mb{\v@lZa}}%
    \v@lXa=\cxa@\v@lX\advance\v@lXa\cxb@\v@lY\v@lXa=\v@lcoef\v@lXa%
    \v@lY=\cyb@\v@lY\advance\v@lY\cya@\v@lX\advance\v@lY\cyc@\v@lZ%
    \v@lY=\v@lcoef\v@lY\v@lX=\v@lXa\ignorespaces}
\def\c@lprojort{
    \v@lXa=\cxa@\v@lX\advance\v@lXa\cxb@\v@lY%
    \v@lY=\cyb@\v@lY\advance\v@lY\cya@\v@lX\advance\v@lY\cyc@\v@lZ%
    \v@lX=\v@lXa\ignorespaces}
\def\Figptpr@j#1:#2/#3/{{\Figg@tXY{#3}\superc@lprojSP%
    \Figp@intregDD#1:{#2}(\v@lX,\v@lY)}\ignorespaces}
\def\figsetobdistTD(#1){{\ifcurr@ntPS%
    \immediate\write16{*** \BS@ figsetobdist is ignored inside a
     \BS@ psbeginfig-\BS@ psendfig block.}%
    \else\v@leur=#1\unit@\c@ldefdisob@{\v@leur}\fi}\ignorespaces}
\def\Figs@tproj#1{%
    \if#13 \d@faultproj\else\if#1c\d@faultproj%
    \else\if#1o\xdef\curr@ntproj{1}\xdef\typ@proj{orthogonal}%
         \figsetviewTD(\def@ultpsi,\def@ulttheta)%
         \global\let\c@lprojSP=\c@lprojort\global\let\superc@lprojSP=\c@lprojort%
    \else\if#1r\xdef\curr@ntproj{2}\xdef\typ@proj{realistic}%
         \figsetviewTD(\def@ultpsi,\def@ulttheta)%
         \global\let\c@lprojSP=\c@lprojrea\global\let\superc@lprojSP=\c@lprojrea%
    \else\d@faultproj\message{*** Unknown projection. Cavalier projection assumed.}%
    \fi\fi\fi\fi}
\def\d@faultproj{\xdef\curr@ntproj{0}\xdef\typ@proj{cavalier}\figsetviewTD(\def@ultpsi,0.5)%
         \global\let\c@lprojSP=\c@lprojcav\global\let\superc@lprojSP=\c@lprojcav}
\def\figsettargetTD[#1]{{\ifcurr@ntPS%
    \immediate\write16{*** \BS@ figsettarget is ignored inside a
     \BS@ psbeginfig-\BS@ psendfig block.}%
    \else\global\newt@rgetpttrue\xdef\t@rgetpt{#1}\Figg@tXY{#1}\global\wd\Bt@rget=\v@lX%
    \global\ht\Bt@rget=\v@lY\global\dp\Bt@rget=\v@lZ\fi}\ignorespaces}
\def\figsetviewTD(#1){\ifcurr@ntPS%
     \immediate\write16{*** \BS@ figsetview is ignored inside a
     \BS@ psbeginfig-\BS@ psendfig block.}\else\Figsetview@#1,:\fi\ignorespaces}
\def\Figsetview@#1,#2:{{\xdef\v@lPsi{#1}\def\t@xt@{#2}%
    \ifx\t@xt@\empty\def\@rgdeux{\v@lTheta}\else%
    \def\Xarg@##1,{\edef\@rgdeux{##1}}\Xarg@#2\fi%
    \c@ssin{\costhet@}{\sinthet@}{#1}\v@lmin=\costhet@ pt\v@lmax=\sinthet@ pt%
    \ifcase\curr@ntproj%
    \v@leur=\@rgdeux\v@lmin\xdef\cxa@{\repdecn@mb{\v@leur}}%
    \v@leur=\@rgdeux\v@lmax\xdef\cxb@{\repdecn@mb{\v@leur}}\v@leur=\@rgdeux pt%
    \relax\ifdim\v@leur>\p@\message{*** Lambda too large ! See \BS@ figset proj() !}\fi%
    \else%
    \v@lmax=-\v@lmax\xdef\cxa@{\repdecn@mb{\v@lmax}}\xdef\cxb@{\costhet@}%
    \ifx\t@xt@\empty\edef\@rgdeux{\def@ulttheta}\fi\c@ssin{\C@}{\S@}{\@rgdeux}%
    \v@lmax=-\S@ pt%
    \v@leur=\v@lmax\v@leur=\costhet@\v@leur\xdef\cya@{\repdecn@mb{\v@leur}}%
    \v@leur=\v@lmax\v@leur=\sinthet@\v@leur\xdef\cyb@{\repdecn@mb{\v@leur}}%
    \xdef\cyc@{\C@}\v@lmin=-\C@ pt%
    \v@leur=\v@lmin\v@leur=\costhet@\v@leur\xdef\cza@{\repdecn@mb{\v@leur}}%
    \v@leur=\v@lmin\v@leur=\sinthet@\v@leur\xdef\czb@{\repdecn@mb{\v@leur}}%
    \xdef\czc@{\repdecn@mb{\v@lmax}}\fi%
    \xdef\v@lTheta{\@rgdeux}}}
\def\def@ultpsi{40}
\def\def@ulttheta{25}
\def\figset#1(#2){\def\t@xt@{#1}\ifx\t@xt@\empty\trtlis@rg{#2}{\Figsetwr@te}
    \else\keln@mde#1|%
    \def\n@mref{pr}\ifx\l@debut\n@mref\ifcurr@ntPS
     \immediate\write16{*** \BS@ figset proj(...) is ignored inside a
     \BS@ psbeginfig-\BS@ psendfig block.}\else\trtlis@rg{#2}{\Figsetpr@j}\fi\else%
    \def\n@mref{wr}\ifx\l@debut\n@mref\trtlis@rg{#2}{\Figsetwr@te}\else
    \immediate\write16{*** Unknown keyword: \BS@ figset #1(...)}%
    \fi\fi\fi\ignorespaces}
\def\Figsetpr@j#1=#2|{\keln@mtr#1|%
    \def\n@mref{dep}\ifx\l@debut\n@mref\Figsetd@p{#2}\else
    \def\n@mref{dis}\ifx\l@debut\n@mref%
     \ifnum\curr@ntproj=\tw@\figsetobdist(#2)\else\Figset@rr\fi\else
    \def\n@mref{lam}\ifx\l@debut\n@mref\Figsetd@p{#2}\else
    \def\n@mref{lat}\ifx\l@debut\n@mref\Figsetth@{#2}\else
    \def\n@mref{lon}\ifx\l@debut\n@mref\figsetview(#2)\else
    \def\n@mref{psi}\ifx\l@debut\n@mref\figsetview(#2)\else
    \def\n@mref{tar}\ifx\l@debut\n@mref%
     \ifnum\curr@ntproj=\tw@\figsettarget[#2]\else\Figset@rr\fi\else
    \def\n@mref{the}\ifx\l@debut\n@mref\Figsetth@{#2}\else
    \immediate\write16{*** Unknown attribute: \BS@ figset proj(..., #1=...).}%
    \fi\fi\fi\fi\fi\fi\fi\fi}
\def\Figsetd@p#1{\ifnum\curr@ntproj=\z@\figsetview(\v@lPsi,#1)\else\Figset@rr\fi}
\def\Figsetth@#1{\ifnum\curr@ntproj=\z@\Figset@rr\else\figsetview(\v@lPsi,#1)\fi}
\def\Figset@rr{\message{*** \BS@ figset proj(): Attribute "\n@mref" ignored, incompatible
    with current projection}}
\def\initb@undb@xTD{\wd\BminTD@=\maxdimen\ht\BminTD@=\maxdimen\dp\BminTD@=\maxdimen%
    \wd\BmaxTD@=-\maxdimen\ht\BmaxTD@=-\maxdimen\dp\BmaxTD@=-\maxdimen}
\newbox\Gb@x      
\newbox\Gb@xSC    
\newtoks\c@nsymb  
\newif\ifr@undcoord\newif\ifunitpr@sent
\def\unssqrttw@{0.707106 }
\def\figAst{\raise-1.15ex\hbox{$\ast$}}
\def\figBullet{\raise-1.15ex\hbox{$\bullet$}}
\def\figCirc{\raise-1.15ex\hbox{$\circ$}}
\def\figDiamond{\raise-1.15ex\hbox{$\diamond$}}%
\def\boxit#1#2{\leavevmode\hbox{\vrule\vbox{\hrule\vglue#1%
    \vtop{\hbox{\kern#1{#2}\kern#1}\vglue#1\hrule}}\vrule}}
\def\centertext#1#2{\vbox{\hsize#1\parindent0cm%
  \leftskip=0pt plus 1fil\rightskip=0pt plus 1fil\parfillskip=0pt{#2}}}
\def\lefttext#1#2{\vbox{\hsize#1\parindent0cm\rightskip=0pt plus 1fil#2}}
\def\c@nterpt{\ignorespaces%
    \kern-.5\wd\Gb@xSC%
    \raise-.5\ht\Gb@xSC\rlap{\hbox{\raise.5\dp\Gb@xSC\hbox{\copy\Gb@xSC}}}%
    \kern .5\wd\Gb@xSC\ignorespaces}
\def\b@undb@xSC#1#2{{\v@lXa=#1\v@lYa=#2%
    \v@leur=\ht\Gb@xSC\advance\v@leur\dp\Gb@xSC%
    \advance\v@lXa-.5\wd\Gb@xSC\advance\v@lYa-.5\v@leur\b@undb@x{\v@lXa}{\v@lYa}%
    \advance\v@lXa\wd\Gb@xSC\advance\v@lYa\v@leur\b@undb@x{\v@lXa}{\v@lYa}}}
\def\@keldist#1#2{\edef\Dist@n{#2}\y@tiunit{\Dist@n}%
    \ifunitpr@sent#1=\Dist@n\else#1=\Dist@n\unit@\fi}
\def\y@tiunit#1{\unitpr@sentfalse\expandafter\y@tiunit@#1:}
\def\y@tiunit@#1#2:{\ifcat#1a\unitpr@senttrue\else\def\l@suite{#2}%
    \ifx\l@suite\empty\else\y@tiunit@#2:\fi\fi}
\def\figcoordDD#1{{\v@lX=\ptT@unit@\v@lX\v@lY=\ptT@unit@\v@lY%
    \ifr@undcoord\ifcase#1\v@leur=0.5pt\or\v@leur=0.05pt\or\v@leur=0.005pt%
    \or\v@leur=0.0005pt\else\v@leur=\z@\fi%
    \ifdim\v@lX<\z@\advance\v@lX-\v@leur\else\advance\v@lX\v@leur\fi%
    \ifdim\v@lY<\z@\advance\v@lY-\v@leur\else\advance\v@lY\v@leur\fi\fi%
    (\@ffichnb{#1}{\repdecn@mb{\v@lX}},\ifmmode\else\thinspace\fi%
    \@ffichnb{#1}{\repdecn@mb{\v@lY}})}}
\def\@ffichnb#1#2{\def\@@ffich{\@ffich#1(}\edef\n@mbre{#2}%
    \expandafter\@@ffich\n@mbre)}
\def\@ffich#1(#2.#3){{#2\ifnum#1>\z@.\fi\def\dig@ts{#3}\s@mme=\z@%
    \loop\ifnum\s@mme<#1\expandafter\@ffichdec\dig@ts:\advance\s@mme\@ne\repeat}}
\def\@ffichdec#1#2:{\relax#1\def\dig@ts{#20}}
\def\figcoordTD#1{{\v@lX=\ptT@unit@\v@lX\v@lY=\ptT@unit@\v@lY\v@lZ=\ptT@unit@\v@lZ%
    \ifr@undcoord\ifcase#1\v@leur=0.5pt\or\v@leur=0.05pt\or\v@leur=0.005pt%
    \or\v@leur=0.0005pt\else\v@leur=\z@\fi%
    \ifdim\v@lX<\z@\advance\v@lX-\v@leur\else\advance\v@lX\v@leur\fi%
    \ifdim\v@lY<\z@\advance\v@lY-\v@leur\else\advance\v@lY\v@leur\fi%
    \ifdim\v@lZ<\z@\advance\v@lZ-\v@leur\else\advance\v@lZ\v@leur\fi\fi%
    (\@ffichnb{#1}{\repdecn@mb{\v@lX}},\ifmmode\else\thinspace\fi%
     \@ffichnb{#1}{\repdecn@mb{\v@lY}},\ifmmode\else\thinspace\fi%
     \@ffichnb{#1}{\repdecn@mb{\v@lZ}})}}
\def\figsetroundcoord#1{\expandafter\Figsetr@undcoord#1:\ignorespaces}
\def\Figsetr@undcoord#1#2:{\if#1n\r@undcoordfalse\else\r@undcoordtrue\fi}
\def\Figsetwr@te#1=#2|{\keln@mun#1|%
    \def\n@mref{m}\ifx\l@debut\n@mref\figsetmark{#2}\else
    \immediate\write16{*** Unknown attribute: \BS@ figset (..., #1=...)}%
    \fi}
\def\figsetmark#1{\c@nsymb={#1}\setbox\Gb@xSC=\hbox{\the\c@nsymb}\ignorespaces}
\def\figsetptname#1{\def\ptn@me##1{#1}\ignorespaces}
\def\FigWrit@L#1:#2(#3,#4){\ignorespaces\@keldist\v@leur{#3}\@keldist\delt@{#4}%
    \C@rp@r@m\def\list@num{#1}\@ecfor\p@int:=\list@num\do{\FigWrit@pt{\p@int}{#2}}}
\def\FigWrit@pt#1#2{\FigWp@r@m{#1}{#2}\Vc@rrect\figWp@si%
    \ifdim\wd\Gb@xSC>\z@\b@undb@xSC{\v@lX}{\v@lY}\fi\figWBB@x}
\def\FigWp@r@m#1#2{\Figg@tXY{#1}%
    \setbox\Gb@x=\hbox{\def\t@xt@{#2}\ifx\t@xt@\empty\Figg@tT{#1}\else#2\fi}\c@lprojSP}
\let\Vc@rrect=\relax
\let\C@rp@r@m=\relax
\def\figwrite[#1]#2{{\ignorespaces\def\list@num{#1}\@ecfor\p@int:=\list@num\do{%
    \setbox\Gb@x=\hbox{\def\t@xt@{#2}\ifx\t@xt@\empty\Figg@tT{\p@int}\else#2\fi}%
    \Figwrit@{\p@int}}}\ignorespaces}
\def\Figwrit@#1{\Figg@tXY{#1}\c@lprojSP%
    \rlap{\kern\v@lX\raise\v@lY\hbox{\unhcopy\Gb@x}}\v@leur=\v@lY%
    \advance\v@lY\ht\Gb@x\b@undb@x{\v@lX}{\v@lY}\advance\v@lX\wd\Gb@x%
    \v@lY=\v@leur\advance\v@lY-\dp\Gb@x\b@undb@x{\v@lX}{\v@lY}}
\def\figwritec[#1]#2{{\ignorespaces\def\list@num{#1}%
    \@ecfor\p@int:=\list@num\do{\Figwrit@c{\p@int}{#2}}}\ignorespaces}
\def\Figwrit@c#1#2{\FigWp@r@m{#1}{#2}%
    \rlap{\kern\v@lX\raise\v@lY\hbox{\rlap{\kern-.5\wd\Gb@x%
    \raise-.5\ht\Gb@x\hbox{\raise.5\dp\Gb@x\hbox{\unhcopy\Gb@x}}}}}%
    \v@leur=\ht\Gb@x\advance\v@leur\dp\Gb@x%
    \advance\v@lX-.5\wd\Gb@x\advance\v@lY-.5\v@leur\b@undb@x{\v@lX}{\v@lY}%
    \advance\v@lX\wd\Gb@x\advance\v@lY\v@leur\b@undb@x{\v@lX}{\v@lY}}
\def\figwritep[#1]{{\ignorespaces\def\list@num{#1}\setbox\Gb@x=\hbox{\c@nterpt}%
    \@ecfor\p@int:=\list@num\do{\Figwrit@{\p@int}}}\ignorespaces}
\def\figwritew#1:#2(#3){\figwritegcw#1:{#2}(#3,0pt)}
\def\figwritee#1:#2(#3){\figwritegce#1:{#2}(#3,0pt)}
\def\figwriten#1:#2(#3){{\def\Vc@rrect{\v@lZ=\v@leur\advance\v@lZ\dp\Gb@x}%
    \Figwrit@NS#1:{#2}(#3)}\ignorespaces}
\def\figwrites#1:#2(#3){{\def\Vc@rrect{\v@lZ=-\v@leur\advance\v@lZ-\ht\Gb@x}%
    \Figwrit@NS#1:{#2}(#3)}\ignorespaces}
\def\Figwrit@NS#1:#2(#3){\let\figWp@si=\FigWp@siNS\let\figWBB@x=\FigWBB@xNS%
    \FigWrit@L#1:{#2}(#3,0pt)}
\def\FigWp@siNS{\rlap{\kern\v@lX\raise\v@lY\hbox{\rlap{\kern-.5\wd\Gb@x%
    \raise\v@lZ\hbox{\unhcopy\Gb@x}}\c@nterpt}}}
\def\FigWBB@xNS{\advance\v@lY\v@lZ%
    \advance\v@lY-\dp\Gb@x\advance\v@lX-.5\wd\Gb@x\b@undb@x{\v@lX}{\v@lY}%
    \advance\v@lY\ht\Gb@x\advance\v@lY\dp\Gb@x%
    \advance\v@lX\wd\Gb@x\b@undb@x{\v@lX}{\v@lY}}
\def\figwritenw#1:#2(#3){{\let\figWp@si=\FigWp@sigW\let\figWBB@x=\FigWBB@xgWE%
    \def\C@rp@r@m{\v@leur=\unssqrttw@\v@leur\delt@=\v@leur%
    \ifdim\delt@=\z@\delt@=\epsil@n\fi}\let@xte={-}\FigWrit@L#1:{#2}(#3,0pt)}\ignorespaces}
\def\figwritesw#1:#2(#3){{\let\figWp@si=\FigWp@sigW\let\figWBB@x=\FigWBB@xgWE%
    \def\C@rp@r@m{\v@leur=\unssqrttw@\v@leur\delt@=-\v@leur%
    \ifdim\delt@=\z@\delt@=-\epsil@n\fi}\let@xte={-}\FigWrit@L#1:{#2}(#3,0pt)}\ignorespaces}
\def\figwritene#1:#2(#3){{\let\figWp@si=\FigWp@sigE\let\figWBB@x=\FigWBB@xgWE%
    \def\C@rp@r@m{\v@leur=\unssqrttw@\v@leur\delt@=\v@leur%
    \ifdim\delt@=\z@\delt@=\epsil@n\fi}\let@xte={}\FigWrit@L#1:{#2}(#3,0pt)}\ignorespaces}
\def\figwritese#1:#2(#3){{\let\figWp@si=\FigWp@sigE\let\figWBB@x=\FigWBB@xgWE%
    \def\C@rp@r@m{\v@leur=\unssqrttw@\v@leur\delt@=-\v@leur%
    \ifdim\delt@=\z@\delt@=-\epsil@n\fi}\let@xte={}\FigWrit@L#1:{#2}(#3,0pt)}\ignorespaces}
\def\figwritegw#1:#2(#3,#4){{\let\figWp@si=\FigWp@sigW\let\figWBB@x=\FigWBB@xgWE%
    \let@xte={-}\FigWrit@L#1:{#2}(#3,#4)}\ignorespaces}
\def\figwritege#1:#2(#3,#4){{\let\figWp@si=\FigWp@sigE\let\figWBB@x=\FigWBB@xgWE%
    \let@xte={}\FigWrit@L#1:{#2}(#3,#4)}\ignorespaces}
\def\FigWp@sigW{\v@lXa=\z@\v@lYa=\ht\Gb@x\advance\v@lYa\dp\Gb@x%
    \ifdim\delt@>\z@\relax%
    \rlap{\kern\v@lX\raise\v@lY\hbox{\rlap{\kern-\wd\Gb@x\kern-\v@leur%
          \raise\delt@\hbox{\raise\dp\Gb@x\hbox{\unhcopy\Gb@x}}}\c@nterpt}}%
    \else\ifdim\delt@<\z@\relax\v@lYa=-\v@lYa%
    \rlap{\kern\v@lX\raise\v@lY\hbox{\rlap{\kern-\wd\Gb@x\kern-\v@leur%
          \raise\delt@\hbox{\raise-\ht\Gb@x\hbox{\unhcopy\Gb@x}}}\c@nterpt}}%
    \else\v@lXa=-.5\v@lYa%
    \rlap{\kern\v@lX\raise\v@lY\hbox{\rlap{\kern-\wd\Gb@x\kern-\v@leur%
          \raise-.5\ht\Gb@x\hbox{\raise.5\dp\Gb@x\hbox{\unhcopy\Gb@x}}}\c@nterpt}}%
    \fi\fi}
\def\FigWp@sigE{\v@lXa=\z@\v@lYa=\ht\Gb@x\advance\v@lYa\dp\Gb@x%
    \ifdim\delt@>\z@\relax%
    \rlap{\kern\v@lX\raise\v@lY\hbox{\c@nterpt\kern\v@leur%
          \raise\delt@\hbox{\raise\dp\Gb@x\hbox{\unhcopy\Gb@x}}}}%
    \else\ifdim\delt@<\z@\relax\v@lYa=-\v@lYa%
    \rlap{\kern\v@lX\raise\v@lY\hbox{\c@nterpt\kern\v@leur%
          \raise\delt@\hbox{\raise-\ht\Gb@x\hbox{\unhcopy\Gb@x}}}}%
    \else\v@lXa=-.5\v@lYa%
    \rlap{\kern\v@lX\raise\v@lY\hbox{\c@nterpt\kern\v@leur%
          \raise-.5\ht\Gb@x\hbox{\raise.5\dp\Gb@x\hbox{\unhcopy\Gb@x}}}}%
    \fi\fi}
\def\FigWBB@xgWE{\advance\v@lY\delt@%
    \advance\v@lX\the\let@xte\v@leur\advance\v@lY\v@lXa\b@undb@x{\v@lX}{\v@lY}%
    \advance\v@lX\the\let@xte\wd\Gb@x\advance\v@lY\v@lYa\b@undb@x{\v@lX}{\v@lY}}
\def\figwritegcw#1:#2(#3,#4){{\let\figWp@si=\FigWp@sigcW\let\figWBB@x=\FigWBB@xgcWE%
    \let@xte={-}\FigWrit@L#1:{#2}(#3,#4)}\ignorespaces}
\def\figwritegce#1:#2(#3,#4){{\let\figWp@si=\FigWp@sigcE\let\figWBB@x=\FigWBB@xgcWE%
    \let@xte={}\FigWrit@L#1:{#2}(#3,#4)}\ignorespaces}
\def\FigWp@sigcW{\rlap{\kern\v@lX\raise\v@lY\hbox{\rlap{\kern-\wd\Gb@x\kern-\v@leur%
     \raise-.5\ht\Gb@x\hbox{\raise\delt@\hbox{\raise.5\dp\Gb@x\hbox{\unhcopy\Gb@x}}}}%
     \c@nterpt}}}
\def\FigWp@sigcE{\rlap{\kern\v@lX\raise\v@lY\hbox{\c@nterpt\kern\v@leur%
    \raise-.5\ht\Gb@x\hbox{\raise\delt@\hbox{\raise.5\dp\Gb@x\hbox{\unhcopy\Gb@x}}}}}}
\def\FigWBB@xgcWE{\v@lZ=\ht\Gb@x\advance\v@lZ\dp\Gb@x%
    \advance\v@lX\the\let@xte\v@leur\advance\v@lY\delt@\advance\v@lY.5\v@lZ%
    \b@undb@x{\v@lX}{\v@lY}%
    \advance\v@lX\the\let@xte\wd\Gb@x\advance\v@lY-\v@lZ\b@undb@x{\v@lX}{\v@lY}}
\def\figwritebn#1:#2(#3){{\def\Vc@rrect{\v@lZ=\v@leur}\Figwrit@NS#1:{#2}(#3)}\ignorespaces}
\def\figwritebs#1:#2(#3){{\def\Vc@rrect{\v@lZ=-\v@leur}\Figwrit@NS#1:{#2}(#3)}\ignorespaces}
\def\figwritebw#1:#2(#3){{\let\figWp@si=\FigWp@sibW\let\figWBB@x=\FigWBB@xbWE%
    \let@xte={-}\FigWrit@L#1:{#2}(#3,0pt)}\ignorespaces}
\def\figwritebe#1:#2(#3){{\let\figWp@si=\FigWp@sibE\let\figWBB@x=\FigWBB@xbWE%
    \let@xte={}\FigWrit@L#1:{#2}(#3,0pt)}\ignorespaces}
\def\FigWp@sibW{\rlap{\kern\v@lX\raise\v@lY\hbox{\rlap{\kern-\wd\Gb@x\kern-\v@leur%
          \hbox{\unhcopy\Gb@x}}\c@nterpt}}}
\def\FigWp@sibE{\rlap{\kern\v@lX\raise\v@lY\hbox{\c@nterpt\kern\v@leur%
          \hbox{\unhcopy\Gb@x}}}}
\def\FigWBB@xbWE{\v@lZ=\ht\Gb@x\advance\v@lZ\dp\Gb@x%
    \advance\v@lX\the\let@xte\v@leur\advance\v@lY\ht\Gb@x\b@undb@x{\v@lX}{\v@lY}%
    \advance\v@lX\the\let@xte\wd\Gb@x\advance\v@lY-\v@lZ\b@undb@x{\v@lX}{\v@lY}}
\newread\frf@g  \newwrite\fwf@g
\newif\ifcurr@ntPS
\newif\ifps@cri
\newif\ifUse@llipse
\newif\ifpsdebugmode \psdebugmodefalse 
\newif\ifmored@ta
\def\c@pypsfile#1#2{\def\blankline{\par}{\catcode`\%=12
    \loop\ifeof#2\mored@tafalse\else\mored@tatrue\immediate\read#2 to\tr@c
    \ifx\tr@c\blankline\else\immediate\write#1{\tr@c}\fi\fi\ifmored@ta\repeat}}
\def\keln@mun#1#2|{\def\l@debut{#1}\def\l@suite{#2}}
\def\keln@mde#1#2#3|{\def\l@debut{#1#2}\def\l@suite{#3}}
\def\keln@mtr#1#2#3#4|{\def\l@debut{#1#2#3}\def\l@suite{#4}}
\def\keln@mqu#1#2#3#4#5|{\def\l@debut{#1#2#3#4}\def\l@suite{#5}}
\let\@psffilein=\frf@g 
\newif\if@psffileok    
\newif\if@psfbbfound   
\newif\if@psfverbose   
\@psfverbosetrue
\def\@psfgetbb#1{\global\@psfbbfoundfalse%
\global\def\@psfllx{0}\global\def\@psflly{0}%
\global\def\@psfurx{30}\global\def\@psfury{30}%
\openin\@psffilein=#1
\ifeof\@psffilein\errmessage{I couldn't open #1, will ignore it}\else
   \edef\setcolonc@tcode{\catcode`\noexpand\:\the\catcode`\:\relax}%
   {\@psffileoktrue \chardef\other=12
    \def\do##1{\catcode`##1=\other}\dospecials \catcode`\ =10 \setcolonc@tcode
    \loop
       \read\@psffilein to \@psffileline
       \ifeof\@psffilein\@psffileokfalse\else
          \expandafter\@psfaux\@psffileline:. \\%
       \fi
   \if@psffileok\repeat
   \if@psfbbfound\else
    \if@psfverbose\message{No bounding box comment in #1; using defaults}\fi\fi
   }\closein\@psffilein\fi}%
{\catcode`\%=12 \global\let\@psfpercent=
\long\def\@psfaux#1#2:#3\\{\ifx#1\@psfpercent
   \def\testit{#2}\ifx\testit\@psfbblit
      \@psfgrab #3 . . . \\%
      \@psffileokfalse
      \global\@psfbbfoundtrue
   \fi\else\ifx#1\par\else\@psffileokfalse\fi\fi}%
\def\@psfempty{}%
\def\@psfgrab #1 #2 #3 #4 #5\\{%
\global\def\@psfllx{#1}\ifx\@psfllx\@psfempty
      \@psfgrab #2 #3 #4 #5 .\\\else
   \global\def\@psflly{#2}%
   \global\def\@psfurx{#3}\global\def\@psfury{#4}\fi}%
\def\PSwrit@cmd#1#2#3{{\Figg@tXY{#1}\c@lprojSP\b@undb@x{\v@lX}{\v@lY}%
    \v@lX=\ptT@ptps\v@lX\v@lY=\ptT@ptps\v@lY%
    \immediate\write#3{\repdecn@mb{\v@lX}\space\repdecn@mb{\v@lY}\space#2}}}
\def\PSwrit@cmdS#1#2#3#4#5{{\Figg@tXY{#1}\c@lprojSP\b@undb@x{\v@lX}{\v@lY}%
    \global\result@t=\v@lX\global\result@@t=\v@lY%
    \v@lX=\ptT@ptps\v@lX\v@lY=\ptT@ptps\v@lY%
    \immediate\write#3{\repdecn@mb{\v@lX}\space\repdecn@mb{\v@lY}\space#2}}%
    \edef#4{\the\result@t}\edef#5{\the\result@@t}}
\def\psaltitude#1[#2,#3,#4]{{\ifcurr@ntPS\ifps@cri%
    \PSc@mment{psaltitude Square Dim=#1, Triangle=[#2 / #3,#4]}%
    \s@uvc@ntr@l\et@tpsaltitude\resetc@ntr@l{2}\figptorthoprojline-5:=#2/#3,#4/%
    \figvectP -1[#3,#4]\n@rminf{\v@leur}{-1}\vecunit@{-3}{-1}%
    \figvectP -1[-5,#3]\n@rminf{\v@lmin}{-1}\figvectP -2[-5,#4]\n@rminf{\v@lmax}{-2}%
    \ifdim\v@lmin<\v@lmax\s@mme=#3\else\v@lmax=\v@lmin\s@mme=#4\fi%
    \figvectP -4[-5,#2]\vecunit@{-4}{-4}\delt@=#1\unit@%
    \edef\t@ille{\repdecn@mb{\delt@}}\figpttra-1:=-5/\t@ille,-3/%
    \figptstra-3=-5,-1/\t@ille,-4/\psline[#2,-5]\psline[-1,-2,-3]%
    \ifdim\v@leur<\v@lmax\Pss@tsecondSt\psline[-5,\the\s@mme]\Psrest@reSt\fi%
    \PSc@mment{End psaltitude}\resetc@ntr@l\et@tpsaltitude\fi\fi}}
\def\Ps@rcerc#1;#2(#3,#4){\ellBB@x#1;#2,#2(#3,#4,0)%
    \immediate\write\fwf@g{newpath}{\delt@=#2\unit@\delt@=\ptT@ptps\delt@%
    \BdingB@xfalse%
    \PSwrit@cmd{#1}{\repdecn@mb{\delt@}\space #3\space #4\space arc}{\fwf@g}}}
\def\psarccircDD#1;#2(#3,#4){\ifcurr@ntPS\ifps@cri%
    \PSc@mment{psarccircDD Center=#1 ; Radius=#2 (Ang1=#3, Ang2=#4)}%
    \iffillm@de\Ps@rcerc#1;#2(#3,#4)%
    \immediate\write\fwf@g{\fillc@md}%
    \else\Ps@rcerc#1;#2(#3,#4)\immediate\write\fwf@g{stroke}\fi%
    \PSc@mment{End psarccircDD}\fi\fi}
\def\psarccircTD#1,#2,#3;#4(#5,#6){{\ifcurr@ntPS\ifps@cri\s@uvc@ntr@l\et@tpsarccircTD%
    \PSc@mment{psarccircTD Center=#1,P1=#2,P2=#3 ; Radius=#4 (Ang1=#5, Ang2=#6)}%
    \setc@ntr@l{2}\c@lExtAxes#1,#2,#3(#4)\psarcellPATD#1,-4,-5(#5,#6)%
    \PSc@mment{End psarccircTD}\resetc@ntr@l\et@tpsarccircTD\fi\fi}}
\def\c@lExtAxes#1,#2,#3(#4){%
    \figvectPTD-5[#1,#2]\vecunit@{-5}{-5}\figvectNTD-4[#1,#2,#3]\vecunit@{-4}{-4}%
    \figvectNVTD-3[-4,-5]\delt@=#4\unit@\edef\r@yon{\repdecn@mb{\delt@}}%
    \figpttra-4:=#1/\r@yon,-5/\figpttra-5:=#1/\r@yon,-3/}
\def\psarccircPDD#1;#2[#3,#4]{{\ifcurr@ntPS\ifps@cri\s@uvc@ntr@l\et@tpsarccircPDD%
    \PSc@mment{psarccircPDD Center=#1; Radius=#2, [P1=#3, P2=#4]}%
    \Ps@ngleparam#1;#2[#3,#4]\ifdim\v@lmin>\v@lmax\advance\v@lmax\DePI@deg\fi%
    \edef\@ngdeb{\repdecn@mb{\v@lmin}}\edef\@ngfin{\repdecn@mb{\v@lmax}}%
    \psarccirc#1;\r@dius(\@ngdeb,\@ngfin)%
    \PSc@mment{End psarccircPDD}\resetc@ntr@l\et@tpsarccircPDD\fi\fi}}
\def\psarccircPTD#1;#2[#3,#4,#5]{{\ifcurr@ntPS\ifps@cri\s@uvc@ntr@l\et@tpsarccircPTD%
    \PSc@mment{psarccircPTD Center=#1; Radius=#2, [P1=#3, P2=#4, P3=#5]}%
    \setc@ntr@l{2}\c@lExtAxes#1,#3,#5(#2)\psarcellPP#1,-4,-5[#3,#4]%
    \PSc@mment{End psarccircPTD}\resetc@ntr@l\et@tpsarccircPTD\fi\fi}}
\def\Ps@ngleparam#1;#2[#3,#4]{\setc@ntr@l{2}%
    \figvectPDD-1[#1,#3]\vecunit@{-1}{-1}\Figg@tXY{-1}\arct@n\v@lmin(\v@lX,\v@lY)%
    \figvectPDD-2[#1,#4]\vecunit@{-2}{-2}\Figg@tXY{-2}\arct@n\v@lmax(\v@lX,\v@lY)%
    \v@lmin=\rdT@deg\v@lmin\v@lmax=\rdT@deg\v@lmax%
    \v@leur=#2pt\maxim@m{\mili@u}{-\v@leur}{\v@leur}%
    \edef\r@dius{\repdecn@mb{\mili@u}}}
\def\Ps@rell#1;#2,#3(#4,#5,#6){\ellBB@x#1;#2,#3(#4,#5,#6)%
    \immediate\write\fwf@g{newpath}{\v@lmin=#2\unit@\v@lmin=\ptT@ptps\v@lmin%
    \v@lmax=#3\unit@\v@lmax=\ptT@ptps\v@lmax\BdingB@xfalse%
    \PSwrit@cmd{#1}%
    {#6\space\repdecn@mb{\v@lmin}\space\repdecn@mb{\v@lmax}\space #4\space #5\space ellipse}{\fwf@g}}%
    \global\Use@llipsetrue}
\def\psarcellDD#1;#2,#3(#4,#5,#6){{\ifcurr@ntPS\ifps@cri%
    \PSc@mment{psarcellDD Center=#1 ; XRad=#2, YRad=#3 (Ang1=#4, Ang2=#5, Inclination=#6)}%
    \iffillm@de\Ps@rell#1;#2,#3(#4,#5,#6)%
    \immediate\write\fwf@g{\fillc@md}%
    \else\Ps@rell#1;#2,#3(#4,#5,#6)\immediate\write\fwf@g{stroke}\fi%
    \PSc@mment{End psarcellDD}\fi\fi}}
\def\psarcellTD#1;#2,#3(#4,#5,#6){{\ifcurr@ntPS\ifps@cri\s@uvc@ntr@l\et@tpsarcellTD%
    \PSc@mment{psarcellTD Center=#1 ; XRad=#2, YRad=#3 (Ang1=#4, Ang2=#5, Inclination=#6)}%
    \setc@ntr@l{2}\figpttraC -8:=#1/#2,0,0/\figpttraC -7:=#1/0,#3,0/%
    \figvectC -4(0,0,1)\figptsrot -8=-8,-7/#1,#6,-4/\psarcellPATD#1,-8,-7(#4,#5)%
    \PSc@mment{End psarcellTD}\resetc@ntr@l\et@tpsarcellTD\fi\fi}}
\def\psarcellPADD#1,#2,#3(#4,#5){{\ifcurr@ntPS\ifps@cri\s@uvc@ntr@l\et@tpsarcellPADD%
    \PSc@mment{psarcellPADD Center=#1,PtAxis1=#2,PtAxis2=#3 (Ang1=#4, Ang2=#5)}%
    \setc@ntr@l{2}\figvectPDD-1[#1,#2]\vecunit@DD{-1}{-1}\v@lX=\ptT@unit@\result@t%
    \edef\XR@d{\repdecn@mb{\v@lX}}\Figg@tXY{-1}\arct@n\v@lmin(\v@lX,\v@lY)%
    \v@lmin=\rdT@deg\v@lmin\edef\Inclin@{\repdecn@mb{\v@lmin}}%
    \figgetdist\YR@d[#1,#3]\psarcellDD#1;\XR@d,\YR@d(#4,#5,\Inclin@)%
    \PSc@mment{End psarcellPADD}\resetc@ntr@l\et@tpsarcellPADD\fi\fi}}
\def\psarcellPATD#1,#2,#3(#4,#5){{\ifcurr@ntPS\ifps@cri\s@uvc@ntr@l\et@tpsarcellPATD%
    \PSc@mment{psarcellPATD Center=#1,PtAxis1=#2,PtAxis2=#3 (Ang1=#4, Ang2=#5)}%
    \iffillm@de\Ps@rellPATD#1,#2,#3(#4,#5)%
    \immediate\write\fwf@g{\fillc@md}%
    \else\Ps@rellPATD#1,#2,#3(#4,#5)\immediate\write\fwf@g{stroke}\fi%
    \PSc@mment{End psarcellPATD}\resetc@ntr@l\et@tpsarcellPATD\fi\fi}}
\def\Ps@rellPATD#1,#2,#3(#4,#5){\let\c@lprojSP=\relax%
    \setc@ntr@l{2}\figvectPTD-1[#1,#2]\figvectPTD-2[#1,#3]\c@lNbarcs{#4}{#5}%
    \v@leur=#4pt\c@lptellP{#1}{-1}{-2}\Figptpr@j-5:/-3/%
    \immediate\write\fwf@g{newpath}\PSwrit@cmdS{-5}{moveto}{\fwf@g}{\X@un}{\Y@un}%
    \edef\C@nt@r{#1}\s@mme=\z@\bcl@rellPATD}
\def\bcl@rellPATD{\relax%
    \ifnum\s@mme<\p@rtent\advance\s@mme\@ne%
    \Figg@tXY{-5}\v@lX=8\v@lX\v@lY=8\v@lY\Figp@intregDD-4:(\v@lX,\v@lY)%
    \advance\v@leur\delt@\c@lptellP{\C@nt@r}{-1}{-2}\Figptpr@j-6:/-3/\figpttraDD-4:=-4/-27,-6/%
    \advance\v@leur\delt@\c@lptellP{\C@nt@r}{-1}{-2}\Figptpr@j-6:/-3/\figpttraDD-5:=-5/-27,-6/%
    \advance\v@leur\delt@\c@lptellP{\C@nt@r}{-1}{-2}\Figptpr@j-6:/-3/%
    \figpttraDD-4:=-4/1,-6/\figpttraDD-5:=-5/8,-6/\BdingB@xfalse%
    \c@lptCtrlDD-3[-5,-4]\PSwrit@cmdS{-3}{}{\fwf@g}{\X@de}{\Y@de}%
    \c@lptCtrlDD-3[-4,-5]\PSwrit@cmdS{-3}{}{\fwf@g}{\X@tr}{\Y@tr}%
    \BdingB@xtrue\PSwrit@cmdS{-6}{curveto}{\fwf@g}{\X@qu}{\Y@qu}%
    \B@zierBB@x{1}{\Y@un}(\X@un,\X@de,\X@tr,\X@qu)%
    \B@zierBB@x{2}{\X@un}(\Y@un,\Y@de,\Y@tr,\Y@qu)%
    \edef\X@un{\X@qu}\edef\Y@un{\Y@qu}\figptcopyDD-5:/-6/\bcl@rellPATD\fi}
\def\c@lNbarcs#1#2{%
    \delt@=#2pt\advance\delt@-#1pt\maxim@m{\v@lmax}{\delt@}{-\delt@}%
    \v@leur=\v@lmax\divide\v@leur45 \p@rtentiere{\p@rtent}{\v@leur}\advance\p@rtent\@ne%
    \s@mme=\p@rtent\multiply\s@mme\thr@@\divide\delt@\s@mme}
\def\c@lptCtrlDD#1[#2,#3]{\figptcopyDD#1:/#2/\figpttraDD#1:=#1/-2,#3/%
    \Figg@tXY{#1}\divide\v@lX18 \divide\v@lY18 \Figp@intregDD#1:(\v@lX,\v@lY)}
\def\c@lptCtrlTD#1[#2,#3]{\figptcopyTD#1:/#2/\figpttraTD#1:=#1/-2,#3/%
    \Figg@tXY{#1}\divide\v@lX18 \divide\v@lY18 \divide\v@lZ18 %
    \Figp@intregTD#1:(\v@lX,\v@lY,\v@lZ)}
\def\psarcellPP#1,#2,#3[#4,#5]{{\ifcurr@ntPS\ifps@cri\s@uvc@ntr@l\et@tpsarcellPP%
    \PSc@mment{psarcellPP Center=#1,PtAxis1=#2,PtAxis2=#3 [Point1=#4, Point2=#5]}%
    \setc@ntr@l{2}\figvectP-2[#1,#3]\vecunit@{-2}{-2}\v@lmin=\result@t%
    \invers@{\v@lmax}{\v@lmin}%
    \figvectP-1[#1,#2]\vecunit@{-1}{-1}\v@leur=\result@t%
    \v@leur=\repdecn@mb{\v@lmax}\v@leur\edef\AsB@{\repdecn@mb{\v@leur}}
    \c@lAngle{#1}{#4}{\v@lmin}\edef\@ngdeb{\repdecn@mb{\v@lmin}}%
    \c@lAngle{#1}{#5}{\v@lmax}\ifdim\v@lmin>\v@lmax\advance\v@lmax\DePI@deg\fi%
    \edef\@ngfin{\repdecn@mb{\v@lmax}}\psarcellPA#1,#2,#3(\@ngdeb,\@ngfin)%
    \PSc@mment{End psarcellPP}\resetc@ntr@l\et@tpsarcellPP\fi\fi}}
\def\c@lAngle#1#2#3{\figvectP-3[#1,#2]%
    \c@lproscal\delt@[-3,-1]\c@lproscal\v@leur[-3,-2]%
    \v@leur=\AsB@\v@leur\arct@n#3(\delt@,\v@leur)#3=\rdT@deg#3}
\newif\if@rrowratio\@rrowratiotrue
\newif\if@rrowhfill
\newif\if@rrowhout
\def\Psset@rrowhe@d#1=#2|{\keln@mun#1|%
    \def\n@mref{a}\ifx\l@debut\n@mref\pssetarrowheadangle{#2}\else
    \def\n@mref{f}\ifx\l@debut\n@mref\pssetarrowheadfill{#2}\else
    \def\n@mref{l}\ifx\l@debut\n@mref\pssetarrowheadlength{#2}\else
    \def\n@mref{o}\ifx\l@debut\n@mref\pssetarrowheadout{#2}\else
    \def\n@mref{r}\ifx\l@debut\n@mref\pssetarrowheadratio{#2}\else
    \immediate\write16{*** Unknown attribute: \BS@ psset arrowhead(..., #1=...)}%
    \fi\fi\fi\fi\fi}
\def\pssetarrowheadangle#1{\edef\@rrowheadangle{#1}{\c@ssin{\C@}{\S@}{#1}%
    \xdef\C@AHANG{\C@}\xdef\S@AHANG{\S@}\v@lmax=\S@ pt%
    \invers@{\v@leur}{\v@lmax}\maxim@m{\v@leur}{\v@leur}{-\v@leur}%
    \xdef\UNSS@N{\the\v@leur}}}
\def\pssetarrowheadfill#1{\expandafter\set@rrowhfill#1:}
\def\set@rrowhfill#1#2:{\if#1n\@rrowhfillfalse\else\@rrowhfilltrue\fi}
\def\pssetarrowheadout#1{\expandafter\set@rrowhout#1:}
\def\set@rrowhout#1#2:{\if#1n\@rrowhoutfalse\else\@rrowhouttrue\fi}
\def\pssetarrowheadlength#1{\edef\@rrowheadlength{#1}\@rrowratiofalse}
\def\pssetarrowheadratio#1{\edef\@rrowheadratio{#1}\@rrowratiotrue}
\def\psresetarrowhead{%
    \pssetarrowheadangle{\defaultarrowheadangle}%
    \pssetarrowheadfill{\defaultarrowheadfill}%
    \pssetarrowheadout{\defaultarrowheadout}%
    \pssetarrowheadratio{\defaultarrowheadratio}%
    \d@fm@cdim\defaultarrowheadlength{\defaulth@rdahlength}
    \pssetarrowheadlength{\defaultarrowheadlength}}
\def\defaultarrowheadratio{0.1}
\def\defaultarrowheadangle{20}
\def\defaultarrowheadfill{no}
\def\defaultarrowheadout{no}
\def\defaulth@rdahlength{8pt}
\def\psarrowDD[#1,#2]{{\ifcurr@ntPS\ifps@cri\s@uvc@ntr@l\et@tpsarrow%
    \PSc@mment{psarrowDD [Pt1,Pt2]=[#1,#2]}\pssetfillmode{no}%
    \psarrowheadDD[#1,#2]\setc@ntr@l{2}\psline[#1,-3]%
    \PSc@mment{End psarrowDD}\resetc@ntr@l\et@tpsarrow\fi\fi}}
\def\psarrowTD[#1,#2]{{\ifcurr@ntPS\ifps@cri\s@uvc@ntr@l\et@tpsarrowTD%
    \PSc@mment{psarrowTD [Pt1,Pt2]=[#1,#2]}\resetc@ntr@l{2}%
    \Figptpr@j-5:/#1/\Figptpr@j-6:/#2/\let\c@lprojSP=\relax\psarrowDD[-5,-6]%
    \PSc@mment{End psarrowTD}\resetc@ntr@l\et@tpsarrowTD\fi\fi}}
\def\psarrowheadDD[#1,#2]{{\ifcurr@ntPS\ifps@cri\s@uvc@ntr@l\et@tpsarrowheadDD%
    \if@rrowhfill\def\@hangle{-\@rrowheadangle}\else\def\@hangle{\@rrowheadangle}\fi%
    \if@rrowratio%
    \if@rrowhout\def\@hratio{-\@rrowheadratio}\else\def\@hratio{\@rrowheadratio}\fi%
    \PSc@mment{psarrowheadDD Ratio=\@hratio, Angle=\@hangle, [Pt1,Pt2]=[#1,#2]}%
    \Ps@rrowhead\@hratio,\@hangle[#1,#2]%
    \else%
    \if@rrowhout\def\@hlength{-\@rrowheadlength}\else\def\@hlength{\@rrowheadlength}\fi%
    \PSc@mment{psarrowheadDD Length=\@hlength, Angle=\@hangle, [Pt1,Pt2]=[#1,#2]}%
    \Ps@rrowheadfd\@hlength,\@hangle[#1,#2]%
    \fi%
    \PSc@mment{End psarrowheadDD}\resetc@ntr@l\et@tpsarrowheadDD\fi\fi}}
\def\psarrowheadTD[#1,#2]{{\ifcurr@ntPS\ifps@cri\s@uvc@ntr@l\et@tpsarrowheadTD%
    \PSc@mment{psarrowheadTD [Pt1,Pt2]=[#1,#2]}\resetc@ntr@l{2}%
    \Figptpr@j-5:/#1/\Figptpr@j-6:/#2/\let\c@lprojSP=\relax\psarrowheadDD[-5,-6]%
    \PSc@mment{End psarrowheadTD}\resetc@ntr@l\et@tpsarrowheadTD\fi\fi}}
\def\Ps@rrowhead#1,#2[#3,#4]{\v@leur=#1\p@\maxim@m{\v@leur}{\v@leur}{-\v@leur}%
    \ifdim\v@leur>\Cepsil@n{
    \PSc@mment{ps@rrowhead Ratio=#1, Angle=#2, [Pt1,Pt2]=[#3,#4]}\v@leur=\UNSS@N%
    \v@leur=\curr@ntwidth\v@leur\v@leur=\ptpsT@pt\v@leur\delt@=.5\v@leur
    \setc@ntr@l{2}\figvectPDD-3[#4,#3]%
    \Figg@tXY{-3}\v@lX=#1\v@lX\v@lY=#1\v@lY\Figv@ctCreg-3(\v@lX,\v@lY)%
    \vecunit@{-4}{-3}\mili@u=\result@t%
    \ifdim#2pt>\z@\v@lXa=-\C@AHANG\delt@%
     \edef\c@ef{\repdecn@mb{\v@lXa}}\figpttraDD-3:=-3/\c@ef,-4/\fi%
    \edef\c@ef{\repdecn@mb{\delt@}}%
    \v@lXa=\mili@u\v@lXa=\C@AHANG\v@lXa%
    \v@lYa=\ptpsT@pt\p@\v@lYa=\curr@ntwidth\v@lYa\v@lYa=\sDcc@ngle\v@lYa%
    \advance\v@lXa-\v@lYa\gdef\sDcc@ngle{0}%
    \ifdim\v@lXa>\v@leur\edef\c@efendpt{\repdecn@mb{\v@leur}}%
    \else\edef\c@efendpt{\repdecn@mb{\v@lXa}}\fi%
    \Figg@tXY{-3}\v@lmin=\v@lX\v@lmax=\v@lY%
    \v@lXa=\C@AHANG\v@lmin\v@lYa=\S@AHANG\v@lmax\advance\v@lXa\v@lYa%
    \v@lYa=-\S@AHANG\v@lmin\v@lX=\C@AHANG\v@lmax\advance\v@lYa\v@lX%
    \setc@ntr@l{1}\Figg@tXY{#4}\advance\v@lX\v@lXa\advance\v@lY\v@lYa%
    \setc@ntr@l{2}\Figp@intregDD-2:(\v@lX,\v@lY)%
    \v@lXa=\C@AHANG\v@lmin\v@lYa=-\S@AHANG\v@lmax\advance\v@lXa\v@lYa%
    \v@lYa=\S@AHANG\v@lmin\v@lX=\C@AHANG\v@lmax\advance\v@lYa\v@lX%
    \setc@ntr@l{1}\Figg@tXY{#4}\advance\v@lX\v@lXa\advance\v@lY\v@lYa%
    \setc@ntr@l{2}\Figp@intregDD-1:(\v@lX,\v@lY)%
    \ifdim#2pt<\z@\fillm@detrue\psline[-2,#4,-1]
    \else\figptstraDD-3=#4,-2,-1/\c@ef,-4/\psline[-2,-3,-1]\fi
    \ifdim#1pt>\z@\figpttraDD-3:=#4/\c@efendpt,-4/\else\figptcopyDD-3:/#4/\fi%
    \PSc@mment{End ps@rrowhead}}\fi}
\def\sDcc@ngle{0}
\def\Ps@rrowheadfd#1,#2[#3,#4]{{%
    \PSc@mment{ps@rrowheadfd Length=#1, Angle=#2, [Pt1,Pt2]=[#3,#4]}%
    \setc@ntr@l{2}\figvectPDD-1[#3,#4]\n@rmeucDD{\v@leur}{-1}\v@leur=\ptT@unit@\v@leur%
    \invers@{\v@leur}{\v@leur}\v@leur=#1\v@leur\edef\R@tio{\repdecn@mb{\v@leur}}%
    \Ps@rrowhead\R@tio,#2[#3,#4]\PSc@mment{End ps@rrowheadfd}}}
\def\psarrowBezierDD[#1,#2,#3,#4]{{\ifcurr@ntPS\ifps@cri\s@uvc@ntr@l\et@tpsarrowBezierDD%
    \PSc@mment{psarrowBezierDD Control points=#1,#2,#3,#4}\setc@ntr@l{2}%
    \if@rrowratio\c@larclengthDD\v@leur,10[#1,#2,#3,#4]\else\v@leur=\z@\fi%
    \Ps@rrowB@zDD\v@leur[#1,#2,#3,#4]%
    \PSc@mment{End psarrowBezierDD}\resetc@ntr@l\et@tpsarrowBezierDD\fi\fi}}
\def\psarrowBezierTD[#1,#2,#3,#4]{{\ifcurr@ntPS\ifps@cri\s@uvc@ntr@l\et@tpsarrowBezierTD%
    \PSc@mment{psarrowBezierTD Control points=#1,#2,#3,#4}\resetc@ntr@l{2}%
    \Figptpr@j-7:/#1/\Figptpr@j-8:/#2/\Figptpr@j-9:/#3/\Figptpr@j-10:/#4/%
    \let\c@lprojSP=\relax\ifnum\curr@ntproj<\tw@\psarrowBezierDD[-7,-8,-9,-10]%
    \else\immediate\write\fwf@g{newpath}\PSwrit@cmd{-7}{moveto}{\fwf@g}%
    \if@rrowratio\c@larclengthDD\mili@u,10[-7,-8,-9,-10]\else\mili@u=\z@\fi%
    \p@rtent=\NBz@rcs\advance\p@rtent\m@ne\subB@zierTD\p@rtent[#1,#2,#3,#4]%
    \immediate\write\fwf@g{stroke}%
    \advance\v@lmin\p@rtent\delt@
    \v@leur=\v@lmin\advance\v@leur0.33333 \delt@\edef\unti@rs{\repdecn@mb{\v@leur}}%
    \v@leur=\v@lmin\advance\v@leur0.66666 \delt@\edef\deti@rs{\repdecn@mb{\v@leur}}%
    \figptcopyDD-8:/-10/\c@lsubBzarc\unti@rs,\deti@rs[#1,#2,#3,#4]%
    \figptcopyDD-8:/-4/\figptcopyDD-9:/-3/\Ps@rrowB@zDD\mili@u[-7,-8,-9,-10]\fi%
    \PSc@mment{End psarrowBezierTD}\resetc@ntr@l\et@tpsarrowBezierTD\fi\fi}}
\def\c@larclengthDD#1,#2[#3,#4,#5,#6]{{\p@rtent=#2\figptcopyDD-5:/#3/%
    \delt@=\p@\divide\delt@\p@rtent\c@rre=\z@\v@leur=\z@\s@mme=\z@%
    \loop\ifnum\s@mme<\p@rtent\advance\s@mme\@ne\advance\v@leur\delt@%
    \edef\T@{\repdecn@mb{\v@leur}}\figptBezierDD-6::\T@[#3,#4,#5,#6]%
    \figvectPDD-1[-5,-6]\n@rmeucDD{\mili@u}{-1}\advance\c@rre\mili@u%
    \figptcopyDD-5:/-6/\repeat\global\result@t=\ptT@unit@\c@rre}#1=\result@t}
\def\Ps@rrowB@zDD#1[#2,#3,#4,#5]{{\pssetfillmode{no}%
    \if@rrowratio\delt@=\@rrowheadratio#1\else\delt@=\@rrowheadlength pt\fi%
    \v@leur=\C@AHANG\delt@\edef\R@dius{\repdecn@mb{\v@leur}}%
    \FigptintercircB@zDD-5::0,\R@dius[#5,#4,#3,#2]%
    \pssetarrowheadlength{\repdecn@mb{\delt@}}\psarrowheadDD[-5,#5]%
    \let\n@rmeuc=\n@rmeucDD\figgetdist\R@dius[#5,-3]%
    \FigptintercircB@zDD-6::0,\R@dius[#5,#4,#3,#2]%
    \figptBezierDD-5::0.33333[#5,#4,#3,#2]\figptBezierDD-3::0.66666[#5,#4,#3,#2]%
    \figptscontrolDD-5[-6,-5,-3,#2]\psBezierDD1[-6,-5,-4,#2]}}
\def\psarrowcircDD#1;#2(#3,#4){{\ifcurr@ntPS\ifps@cri\s@uvc@ntr@l\et@tpsarrowcircDD%
    \PSc@mment{psarrowcircDD Center=#1 ; Radius=#2 (Ang1=#3,Ang2=#4)}%
    \pssetfillmode{no}\Pscirc@rrowhead#1;#2(#3,#4)%
    \setc@ntr@l{2}\figvectPDD -4[#1,-3]\vecunit@{-4}{-4}%
    \Figg@tXY{-4}\arct@n\v@lmin(\v@lX,\v@lY)%
    \v@lmin=\rdT@deg\v@lmin\v@leur=#4pt\advance\v@leur-\v@lmin%
    \maxim@m{\v@leur}{\v@leur}{-\v@leur}%
    \ifdim\v@leur>\DemiPI@deg\relax\ifdim\v@lmin<#4pt\advance\v@lmin\DePI@deg%
    \else\advance\v@lmin-\DePI@deg\fi\fi\edef\ar@ngle{\repdecn@mb{\v@lmin}}%
    \ifdim#3pt<#4pt\psarccirc#1;#2(#3,\ar@ngle)\else\psarccirc#1;#2(\ar@ngle,#3)\fi%
    \PSc@mment{End psarrowcircDD}\resetc@ntr@l\et@tpsarrowcircDD\fi\fi}}
\def\psarrowcircTD#1,#2,#3;#4(#5,#6){{\ifcurr@ntPS\ifps@cri\s@uvc@ntr@l\et@tpsarrowcircTD%
    \PSc@mment{psarrowcircTD Center=#1,P1=#2,P2=#3 ; Radius=#4 (Ang1=#5, Ang2=#6)}%
    \resetc@ntr@l{2}\c@lExtAxes#1,#2,#3(#4)\let\c@lprojSP=\relax%
    \figvectPTD-11[#1,-4]\figvectPTD-12[#1,-5]\c@lNbarcs{#5}{#6}%
    \if@rrowratio\v@lmax=\degT@rd\v@lmax\edef\D@lpha{\repdecn@mb{\v@lmax}}\fi%
    \advance\p@rtent\m@ne\mili@u=\z@%
    \v@leur=#5pt\c@lptellP{#1}{-11}{-12}\Figptpr@j-9:/-3/%
    \immediate\write\fwf@g{newpath}\PSwrit@cmdS{-9}{moveto}{\fwf@g}{\X@un}{\Y@un}%
    \edef\C@nt@r{#1}\s@mme=\z@\bcl@rcircTD\immediate\write\fwf@g{stroke}%
    \advance\v@leur\delt@\c@lptellP{#1}{-11}{-12}\Figptpr@j-5:/-3/%
    \advance\v@leur\delt@\c@lptellP{#1}{-11}{-12}\Figptpr@j-6:/-3/%
    \advance\v@leur\delt@\c@lptellP{#1}{-11}{-12}\Figptpr@j-10:/-3/%
    \figptscontrolDD-8[-9,-5,-6,-10]%
    \if@rrowratio\c@lcurvradDD0.5[-9,-8,-7,-10]\advance\mili@u\result@t%
    \maxim@m{\mili@u}{\mili@u}{-\mili@u}\mili@u=\ptT@unit@\mili@u%
    \mili@u=\D@lpha\mili@u\advance\p@rtent\@ne\divide\mili@u\p@rtent\fi%
    \Ps@rrowB@zDD\mili@u[-9,-8,-7,-10]%
    \PSc@mment{End psarrowcircTD}\resetc@ntr@l\et@tpsarrowcircTD\fi\fi}}
\def\bcl@rcircTD{\relax%
    \ifnum\s@mme<\p@rtent\advance\s@mme\@ne%
    \advance\v@leur\delt@\c@lptellP{\C@nt@r}{-11}{-12}\Figptpr@j-5:/-3/%
    \advance\v@leur\delt@\c@lptellP{\C@nt@r}{-11}{-12}\Figptpr@j-6:/-3/%
    \advance\v@leur\delt@\c@lptellP{\C@nt@r}{-11}{-12}\Figptpr@j-10:/-3/%
    \figptscontrolDD-8[-9,-5,-6,-10]\BdingB@xfalse%
    \PSwrit@cmdS{-8}{}{\fwf@g}{\X@de}{\Y@de}\PSwrit@cmdS{-7}{}{\fwf@g}{\X@tr}{\Y@tr}%
    \BdingB@xtrue\PSwrit@cmdS{-10}{curveto}{\fwf@g}{\X@qu}{\Y@qu}%
    \if@rrowratio\c@lcurvradDD0.5[-9,-8,-7,-10]\advance\mili@u\result@t\fi%
    \B@zierBB@x{1}{\Y@un}(\X@un,\X@de,\X@tr,\X@qu)%
    \B@zierBB@x{2}{\X@un}(\Y@un,\Y@de,\Y@tr,\Y@qu)%
    \edef\X@un{\X@qu}\edef\Y@un{\Y@qu}\figptcopyDD-9:/-10/\bcl@rcircTD\fi}
\def\Pscirc@rrowhead#1;#2(#3,#4){{%
    \PSc@mment{pscirc@rrowhead Center=#1 ; Radius=#2 (Ang1=#3,Ang2=#4)}%
    \v@leur=#2\unit@\edef\s@glen{\repdecn@mb{\v@leur}}\v@lY=\z@\v@lX=\v@leur%
    \resetc@ntr@l{2}\Figv@ctCreg-3(\v@lX,\v@lY)\figpttraDD-5:=#1/1,-3/%
    \figptrotDD-5:=-5/#1,#4/%
    \figvectPDD-3[#1,-5]\Figg@tXY{-3}\v@leur=\v@lX%
    \ifdim#3pt<#4pt\v@lX=\v@lY\v@lY=-\v@leur\else\v@lX=-\v@lY\v@lY=\v@leur\fi%
    \Figv@ctCreg-3(\v@lX,\v@lY)\vecunit@{-3}{-3}%
    \if@rrowratio\v@leur=#4pt\advance\v@leur-#3pt\maxim@m{\mili@u}{-\v@leur}{\v@leur}%
    \mili@u=\degT@rd\mili@u\v@leur=\s@glen\mili@u\edef\s@glen{\repdecn@mb{\v@leur}}%
    \mili@u=#2\mili@u\mili@u=\@rrowheadratio\mili@u\else\mili@u=\@rrowheadlength pt\fi%
    \figpttraDD-6:=-5/\s@glen,-3/\v@leur=#2pt\v@leur=2\v@leur%
    \invers@{\v@leur}{\v@leur}\c@rre=\repdecn@mb{\v@leur}\mili@u
    \mili@u=\c@rre\mili@u=\repdecn@mb{\c@rre}\mili@u%
    \v@leur=\p@\advance\v@leur-\mili@u
    \invers@{\mili@u}{2\v@leur}\delt@=\c@rre\delt@=\repdecn@mb{\mili@u}\delt@%
    \xdef\sDcc@ngle{\repdecn@mb{\delt@}}
    \sqrt@{\mili@u}{\v@leur}\arct@n\v@leur(\mili@u,\c@rre)%
    \v@leur=\rdT@deg\v@leur
    \ifdim#3pt<#4pt\v@leur=-\v@leur\fi%
    \if@rrowhout\v@leur=-\v@leur\fi\edef\cor@ngle{\repdecn@mb{\v@leur}}%
    \figptrotDD-6:=-6/-5,\cor@ngle/\psarrowheadDD[-6,-5]%
    \PSc@mment{End pscirc@rrowhead}}}
\def\psarrowcircPDD#1;#2[#3,#4]{{\ifcurr@ntPS\ifps@cri%
    \PSc@mment{psarrowcircPDD Center=#1; Radius=#2, [P1=#3,P2=#4]}%
    \s@uvc@ntr@l\et@tpsarrowcircPDD\Ps@ngleparam#1;#2[#3,#4]%
    \ifdim\v@leur>\z@\ifdim\v@lmin>\v@lmax\advance\v@lmax\DePI@deg\fi%
    \else\ifdim\v@lmin<\v@lmax\advance\v@lmin\DePI@deg\fi\fi%
    \edef\@ngdeb{\repdecn@mb{\v@lmin}}\edef\@ngfin{\repdecn@mb{\v@lmax}}%
    \psarrowcirc#1;\r@dius(\@ngdeb,\@ngfin)%
    \PSc@mment{End psarrowcircPDD}\resetc@ntr@l\et@tpsarrowcircPDD\fi\fi}}
\def\psarrowcircPTD#1;#2[#3,#4,#5]{{\ifcurr@ntPS\ifps@cri\s@uvc@ntr@l\et@tpsarrowcircPTD%
    \PSc@mment{psarrowcircPTD Center=#1; Radius=#2, [P1=#3,P2=#4,P3=#5]}%
    \figgetangleTD\@ngfin[#1,#3,#4,#5]\v@leur=#2pt%
    \maxim@m{\mili@u}{-\v@leur}{\v@leur}\edef\r@dius{\repdecn@mb{\mili@u}}%
    \ifdim\v@leur<\z@\v@lmax=\@ngfin pt\advance\v@lmax-\DePI@deg%
    \edef\@ngfin{\repdecn@mb{\v@lmax}}\fi\psarrowcircTD#1,#3,#5;\r@dius(0,\@ngfin)%
    \PSc@mment{End psarrowcircPTD}\resetc@ntr@l\et@tpsarrowcircPTD\fi\fi}}
\def\psaxes#1(#2){{\ifcurr@ntPS\ifps@cri\s@uvc@ntr@l\et@tpsaxes%
    \PSc@mment{psaxes Origin=#1 Range=(#2)}\an@lys@xes#2,:\resetc@ntr@l{2}%
    \ifx\t@xt@\empty\ifTr@isDim\ps@xes#1(0,#2,0,#2,0,#2)\else\ps@xes#1(0,#2,0,#2)\fi%
    \else\ps@xes#1(#2)\fi\PSc@mment{End psaxes}\resetc@ntr@l\et@tpsaxes\fi\fi}}
\def\an@lys@xes#1,#2:{\def\t@xt@{#2}}
\def\ps@xesDD#1(#2,#3,#4,#5){%
    \figpttraC-5:=#1/#2,0/\figpttraC-6:=#1/#3,0/\psarrowDD[-5,-6]%
    \figpttraC-5:=#1/0,#4/\figpttraC-6:=#1/0,#5/\psarrowDD[-5,-6]}
\def\ps@xesTD#1(#2,#3,#4,#5,#6,#7){%
    \figpttraC-7:=#1/#2,0,0/\figpttraC-8:=#1/#3,0,0/\psarrowTD[-7,-8]%
    \figpttraC-7:=#1/0,#4,0/\figpttraC-8:=#1/0,#5,0/\psarrowTD[-7,-8]%
    \figpttraC-7:=#1/0,0,#6/\figpttraC-8:=#1/0,0,#7/\psarrowTD[-7,-8]}
\def\psbeginfig#1{\ifcurr@ntPS\else%
    \edef\PSfilen@me{#1}\edef\auxfilen@me{\jobname.anx}%
    \ifpstestm@de\ps@critrue\else\openin\frf@g=\PSfilen@me\relax%
    \ifeof\frf@g\ps@critrue\else\ps@crifalse\fi\closein\frf@g\fi%
    \curr@ntPStrue\c@ldefproj\expandafter\setupd@te\defaultupdate:%
    \ifps@cri\initb@undb@x%
    \immediate\openout\fwf@g=\auxfilen@me\initpss@ttings\fi%
    \fi}
\def\initpss@ttings{\psreset{arrowhead,curve,first,flowchart,mesh,second,third}%
    \Use@llipsefalse}
\def\B@zierBB@x#1#2(#3,#4,#5,#6){{\c@rre=\t@n\epsil@n
    \v@lmax=#4\advance\v@lmax-#5\v@lmax=\thr@@\v@lmax\advance\v@lmax#6\advance\v@lmax-#3%
    \mili@u=#4\mili@u=-\tw@\mili@u\advance\mili@u#3\advance\mili@u#5%
    \v@lmin=#4\advance\v@lmin-#3\maxim@m{\v@leur}{-\v@lmax}{\v@lmax}%
    \maxim@m{\delt@}{-\mili@u}{\mili@u}\maxim@m{\v@leur}{\v@leur}{\delt@}%
    \maxim@m{\delt@}{-\v@lmin}{\v@lmin}\maxim@m{\v@leur}{\v@leur}{\delt@}%
    \ifdim\v@leur>\c@rre\invers@{\v@leur}{\v@leur}\edef\Uns@rM@x{\repdecn@mb{\v@leur}}%
    \v@lmax=\Uns@rM@x\v@lmax\mili@u=\Uns@rM@x\mili@u\v@lmin=\Uns@rM@x\v@lmin%
    \maxim@m{\v@leur}{-\v@lmax}{\v@lmax}\ifdim\v@leur<\c@rre%
    \maxim@m{\v@leur}{-\mili@u}{\mili@u}\ifdim\v@leur<\c@rre\else%
    \invers@{\mili@u}{\mili@u}\v@leur=-0.5\v@lmin%
    \v@leur=\repdecn@mb{\mili@u}\v@leur\m@jBBB@x{\v@leur}{#1}{#2}(#3,#4,#5,#6)\fi%
    \else\delt@=\repdecn@mb{\mili@u}\mili@u\v@leur=\repdecn@mb{\v@lmax}\v@lmin%
    \advance\delt@-\v@leur\ifdim\delt@<\z@\else\invers@{\v@lmax}{\v@lmax}%
    \edef\Uns@rAp{\repdecn@mb{\v@lmax}}\sqrt@{\delt@}{\delt@}%
    \v@leur=-\mili@u\advance\v@leur\delt@\v@leur=\Uns@rAp\v@leur%
    \m@jBBB@x{\v@leur}{#1}{#2}(#3,#4,#5,#6)%
    \v@leur=-\mili@u\advance\v@leur-\delt@\v@leur=\Uns@rAp\v@leur%
    \m@jBBB@x{\v@leur}{#1}{#2}(#3,#4,#5,#6)\fi\fi\fi}}
\def\m@jBBB@x#1#2#3(#4,#5,#6,#7){{\relax\ifdim#1>\z@\ifdim#1<\p@%
    \edef\T@{\repdecn@mb{#1}}\v@lX=\p@\advance\v@lX-#1\edef\UNmT@{\repdecn@mb{\v@lX}}%
    \v@lX=#4\v@lY=#5\v@lZ=#6\v@lXa=#7\v@lX=\UNmT@\v@lX\advance\v@lX\T@\v@lY%
    \v@lY=\UNmT@\v@lY\advance\v@lY\T@\v@lZ\v@lZ=\UNmT@\v@lZ\advance\v@lZ\T@\v@lXa%
    \v@lX=\UNmT@\v@lX\advance\v@lX\T@\v@lY\v@lY=\UNmT@\v@lY\advance\v@lY\T@\v@lZ%
    \v@lX=\UNmT@\v@lX\advance\v@lX\T@\v@lY%
    \ifcase#2\or\v@lY=#3\or\v@lY=\v@lX\v@lX=#3\fi\b@undb@x{\v@lX}{\v@lY}\fi\fi}}
\def\PsB@zier#1[#2]{{\immediate\write\fwf@g{newpath}%
    \s@mme=\z@\def\list@num{#2,0}\extrairelepremi@r\p@int\de\list@num%
    \PSwrit@cmdS{\p@int}{moveto}{\fwf@g}{\X@un}{\Y@un}\p@rtent=#1\bclB@zier}}
\def\bclB@zier{\relax%
    \ifnum\s@mme<\p@rtent\advance\s@mme\@ne\BdingB@xfalse%
    \extrairelepremi@r\p@int\de\list@num\PSwrit@cmdS{\p@int}{}{\fwf@g}{\X@de}{\Y@de}%
    \extrairelepremi@r\p@int\de\list@num\PSwrit@cmdS{\p@int}{}{\fwf@g}{\X@tr}{\Y@tr}%
    \BdingB@xtrue%
    \extrairelepremi@r\p@int\de\list@num\PSwrit@cmdS{\p@int}{curveto}{\fwf@g}{\X@qu}{\Y@qu}%
    \B@zierBB@x{1}{\Y@un}(\X@un,\X@de,\X@tr,\X@qu)%
    \B@zierBB@x{2}{\X@un}(\Y@un,\Y@de,\Y@tr,\Y@qu)%
    \edef\X@un{\X@qu}\edef\Y@un{\Y@qu}\bclB@zier\fi}
\def\psBezierDD#1[#2]{\ifcurr@ntPS\ifps@cri%
    \PSc@mment{psBezierDD N arcs=#1, Control points=#2}%
    \iffillm@de\PsB@zier#1[#2]%
    \immediate\write\fwf@g{\fillc@md}%
    \else\PsB@zier#1[#2]\immediate\write\fwf@g{stroke}\fi%
    \PSc@mment{End psBezierDD}\fi\fi}
\def\psBezierTD#1[#2]{\ifcurr@ntPS\ifps@cri\s@uvc@ntr@l\et@tpsBezierTD%
    \PSc@mment{psBezierTD N arcs=#1, Control points=#2}%
    \iffillm@de\PsB@zierTD#1[#2]%
    \immediate\write\fwf@g{\fillc@md}%
    \else\PsB@zierTD#1[#2]\immediate\write\fwf@g{stroke}\fi%
    \PSc@mment{End psBezierTD}\resetc@ntr@l\et@tpsBezierTD\fi\fi}
\def\PsB@zierTD#1[#2]{\ifnum\curr@ntproj<\tw@\PsB@zier#1[#2]\else\PsB@zier@TD#1[#2]\fi}
\def\PsB@zier@TD#1[#2]{{\immediate\write\fwf@g{newpath}%
    \s@mme=\z@\def\list@num{#2,0}\extrairelepremi@r\p@int\de\list@num%
    \let\c@lprojSP=\relax\setc@ntr@l{2}\Figptpr@j-7:/\p@int/%
    \PSwrit@cmd{-7}{moveto}{\fwf@g}%
    \loop\ifnum\s@mme<#1\advance\s@mme\@ne\extrairelepremi@r\p@intun\de\list@num%
    \extrairelepremi@r\p@intde\de\list@num\extrairelepremi@r\p@inttr\de\list@num%
    \subB@zierTD\NBz@rcs[\p@int,\p@intun,\p@intde,\p@inttr]\edef\p@int{\p@inttr}\repeat}}
\def\subB@zierTD#1[#2,#3,#4,#5]{\delt@=\p@\divide\delt@\NBz@rcs\v@lmin=\z@%
    {\Figg@tXY{-7}\edef\X@un{\the\v@lX}\edef\Y@un{\the\v@lY}%
    \s@mme=\z@\loop\ifnum\s@mme<#1\advance\s@mme\@ne%
    \v@leur=\v@lmin\advance\v@leur0.33333 \delt@\edef\unti@rs{\repdecn@mb{\v@leur}}%
    \v@leur=\v@lmin\advance\v@leur0.66666 \delt@\edef\deti@rs{\repdecn@mb{\v@leur}}%
    \advance\v@lmin\delt@\edef\trti@rs{\repdecn@mb{\v@lmin}}%
    \figptBezierTD-8::\trti@rs[#2,#3,#4,#5]\Figptpr@j-8:/-8/%
    \c@lsubBzarc\unti@rs,\deti@rs[#2,#3,#4,#5]\BdingB@xfalse%
    \PSwrit@cmdS{-4}{}{\fwf@g}{\X@de}{\Y@de}\PSwrit@cmdS{-3}{}{\fwf@g}{\X@tr}{\Y@tr}%
    \BdingB@xtrue\PSwrit@cmdS{-8}{curveto}{\fwf@g}{\X@qu}{\Y@qu}%
    \B@zierBB@x{1}{\Y@un}(\X@un,\X@de,\X@tr,\X@qu)%
    \B@zierBB@x{2}{\X@un}(\Y@un,\Y@de,\Y@tr,\Y@qu)%
    \edef\X@un{\X@qu}\edef\Y@un{\Y@qu}\figptcopyDD-7:/-8/\repeat}}
\def\NBz@rcs{2}
\def\c@lsubBzarc#1,#2[#3,#4,#5,#6]{\figptBezierTD-5::#1[#3,#4,#5,#6]%
    \figptBezierTD-6::#2[#3,#4,#5,#6]\Figptpr@j-4:/-5/\Figptpr@j-5:/-6/%
    \figptscontrolDD-4[-7,-4,-5,-8]}
\def\pscircDD#1(#2){\ifcurr@ntPS\ifps@cri\PSc@mment{pscircDD Center=#1 (Radius=#2)}%
    \psarccircDD#1;#2(0,360)\PSc@mment{End pscircDD}\fi\fi}
\def\pscircTD#1,#2,#3(#4){\ifcurr@ntPS\ifps@cri%
    \PSc@mment{pscircTD Center=#1,P1=#2,P2=#3 (Radius=#4)}%
    \psarccircTD#1,#2,#3;#4(0,360)\PSc@mment{End pscircTD}\fi\fi}
{\catcode`\%=12\gdef\p@urcent{
\def\PSc@mment#1{\ifpsdebugmode\immediate\write\fwf@g{\p@urcent\space#1}\fi}
{\catcode`\[=1\catcode`\{=12\gdef\acc@louv[{}}
{\catcode`\]=2\catcode`\}=12\gdef\acc@lfer{}]]
\def\PSdict@{\ifUse@llipse%
    \immediate\write\fwf@g{/ellipsedict 9 dict def ellipsedict /mtrx matrix put}%
    \immediate\write\fwf@g{/ellipse \acc@louv ellipsedict begin}%
    \immediate\write\fwf@g{ /endangle exch def /startangle exch def}%
    \immediate\write\fwf@g{ /yrad exch def /xrad exch def}%
    \immediate\write\fwf@g{ /rotangle exch def /y exch def /x exch def}%
    \immediate\write\fwf@g{ /savematrix mtrx currentmatrix def}%
    \immediate\write\fwf@g{ x y translate rotangle rotate xrad yrad scale}%
    \immediate\write\fwf@g{ 0 0 1 startangle endangle arc}%
    \immediate\write\fwf@g{ savematrix setmatrix end\acc@lfer def}%
    \fi\PShe@der{EndProlog}}
\def\Pssetc@rve#1=#2|{\keln@mun#1|%
    \def\n@mref{r}\ifx\l@debut\n@mref\pssetroundness{#2}\else
    \immediate\write16{*** Unknown attribute: \BS@ psset curve(..., #1=...)}%
    \fi}
\def\pssetroundness#1{\edef\curv@roundness{#1}}
\def\defaultroundness{0.2} 
\def\pscurveDD[#1]{{\ifcurr@ntPS\ifps@cri\PSc@mment{pscurveDD Points=#1}%
    \s@uvc@ntr@l\et@tpscurveDD%
    \iffillm@de\Psc@rveDD\curv@roundness[#1]%
    \immediate\write\fwf@g{\fillc@md}%
    \else\Psc@rveDD\curv@roundness[#1]\immediate\write\fwf@g{stroke}\fi%
    \PSc@mment{End pscurveDD}\resetc@ntr@l\et@tpscurveDD\fi\fi}}
\def\pscurveTD[#1]{{\ifcurr@ntPS\ifps@cri%
    \PSc@mment{pscurveTD Points=#1}\s@uvc@ntr@l\et@tpscurveTD\let\c@lprojSP=\relax%
    \iffillm@de\Psc@rveTD\curv@roundness[#1]%
    \immediate\write\fwf@g{\fillc@md}%
    \else\Psc@rveTD\curv@roundness[#1]\immediate\write\fwf@g{stroke}\fi%
    \PSc@mment{End pscurveTD}\resetc@ntr@l\et@tpscurveTD\fi\fi}}
\def\Psc@rveDD#1[#2]{%
    \def\list@num{#2}\extrairelepremi@r\Ak@\de\list@num%
    \extrairelepremi@r\Ai@\de\list@num\extrairelepremi@r\Aj@\de\list@num%
    \immediate\write\fwf@g{newpath}\PSwrit@cmdS{\Ai@}{moveto}{\fwf@g}{\X@un}{\Y@un}%
    \setc@ntr@l{2}\figvectPDD -1[\Ak@,\Aj@]%
    \@ecfor\Ak@:=\list@num\do{\figpttraDD-2:=\Ai@/#1,-1/\BdingB@xfalse%
       \PSwrit@cmdS{-2}{}{\fwf@g}{\X@de}{\Y@de}%
       \figvectPDD -1[\Ai@,\Ak@]\figpttraDD-2:=\Aj@/-#1,-1/%
       \PSwrit@cmdS{-2}{}{\fwf@g}{\X@tr}{\Y@tr}\BdingB@xtrue%
       \PSwrit@cmdS{\Aj@}{curveto}{\fwf@g}{\X@qu}{\Y@qu}%
       \B@zierBB@x{1}{\Y@un}(\X@un,\X@de,\X@tr,\X@qu)%
       \B@zierBB@x{2}{\X@un}(\Y@un,\Y@de,\Y@tr,\Y@qu)%
       \edef\X@un{\X@qu}\edef\Y@un{\Y@qu}\edef\Ai@{\Aj@}\edef\Aj@{\Ak@}}}
\def\Psc@rveTD#1[#2]{\ifnum\curr@ntproj<\tw@\Psc@rvePPTD#1[#2]\else\Psc@rveCPTD#1[#2]\fi}
\def\Psc@rvePPTD#1[#2]{\setc@ntr@l{2}%
    \def\list@num{#2}\extrairelepremi@r\Ak@\de\list@num\Figptpr@j-5:/\Ak@/%
    \extrairelepremi@r\Ai@\de\list@num\Figptpr@j-3:/\Ai@/%
    \extrairelepremi@r\Aj@\de\list@num\Figptpr@j-4:/\Aj@/%
    \immediate\write\fwf@g{newpath}\PSwrit@cmdS{-3}{moveto}{\fwf@g}{\X@un}{\Y@un}%
    \figvectPDD -1[-5,-4]%
    \@ecfor\Ak@:=\list@num\do{\Figptpr@j-5:/\Ak@/\figpttraDD-2:=-3/#1,-1/%
       \BdingB@xfalse\PSwrit@cmdS{-2}{}{\fwf@g}{\X@de}{\Y@de}%
       \figvectPDD -1[-3,-5]\figpttraDD-2:=-4/-#1,-1/%
       \PSwrit@cmdS{-2}{}{\fwf@g}{\X@tr}{\Y@tr}\BdingB@xtrue%
       \PSwrit@cmdS{-4}{curveto}{\fwf@g}{\X@qu}{\Y@qu}%
       \B@zierBB@x{1}{\Y@un}(\X@un,\X@de,\X@tr,\X@qu)%
       \B@zierBB@x{2}{\X@un}(\Y@un,\Y@de,\Y@tr,\Y@qu)%
       \edef\X@un{\X@qu}\edef\Y@un{\Y@qu}\figptcopyDD-3:/-4/\figptcopyDD-4:/-5/}}
\def\Psc@rveCPTD#1[#2]{\setc@ntr@l{2}%
    \def\list@num{#2}\extrairelepremi@r\Ak@\de\list@num%
    \extrairelepremi@r\Ai@\de\list@num\extrairelepremi@r\Aj@\de\list@num%
    \Figptpr@j-7:/\Ai@/%
    \immediate\write\fwf@g{newpath}\PSwrit@cmd{-7}{moveto}{\fwf@g}%
    \figvectPTD -9[\Ak@,\Aj@]%
    \@ecfor\Ak@:=\list@num\do{\figpttraTD-10:=\Ai@/#1,-9/%
       \figvectPTD -9[\Ai@,\Ak@]\figpttraTD-11:=\Aj@/-#1,-9/%
       \subB@zierTD\NBz@rcs[\Ai@,-10,-11,\Aj@]\edef\Ai@{\Aj@}\edef\Aj@{\Ak@}}}
\def\psendfig{\ifcurr@ntPS\ifps@cri\immediate\closeout\fwf@g%
    \immediate\openout\fwf@g=\PSfilen@me%
    \immediate\write\fwf@g{\p@urcent\string!PS-Adobe-2.0 EPSF-2.0}%
    \PShe@der{Creator\string: TeX (fig4tex.tex)}%
    \PShe@der{Title\string: \PSfilen@me}%
    \PShe@der{CreationDate\string: \the\day/\the\month/\the\year}%
    {\v@lX=\ptT@ptps\c@@rdXmin\v@lY=\ptT@ptps\c@@rdYmin%
     \v@lXa=\ptT@ptps\c@@rdXmax\v@lYa=\ptT@ptps\c@@rdYmax%
    \PShe@der{BoundingBox\string: \repdecn@mb{\v@lX}\space\repdecn@mb{\v@lY}%
              \space\repdecn@mb{\v@lXa}\space\repdecn@mb{\v@lYa}}}%
    \PShe@der{EndComments}\PSdict@%
    \openin\frf@g=\auxfilen@me\c@pypsfile\fwf@g\frf@g\closein\frf@g%
    \PSc@mment{End of file.}\immediate\closeout\fwf@g%
    \immediate\openout\fwf@g=\auxfilen@me\immediate\closeout\fwf@g%
    \immediate\write16{File \PSfilen@me\space created.}\fi\fi\curr@ntPSfalse\ps@critrue}
\def\PShe@der#1{\immediate\write\fwf@g{\p@urcent\p@urcent#1}}
\def\psfcconnect[#1]{{\ifcurr@ntPS\ifps@cri\PSc@mment{psfcconnect Points=#1}%
    \pssetfillmode{no}\s@uvc@ntr@l\et@tpsfcconnect\resetc@ntr@l{2}%
    \fcc@nnect@[#1]\resetc@ntr@l\et@tpsfcconnect\PSc@mment{End psfcconnect}\fi\fi}}
\def\fcc@nnect@[#1]{\let\N@rm=\n@rmeucDD\def\list@num{#1}%
    \extrairelepremi@r\Ai@\de\list@num\edef\pr@m{\Ai@}\v@leur=\z@\p@rtent=\@ne\c@llgtot%
    \ifcase\fclin@typ@\edef\list@num{[\pr@m,#1,\Ai@}\expandafter\pscurve\list@num]%
    \else\ifdim\fclin@r@d\p@>\z@\Pslin@conge[#1]\else\psline[#1]\fi\fi%
    \v@leur=\@rrowp@s\v@leur\edef\list@num{#1,\Ai@,0}%
    \extrairelepremi@r\Ai@\de\list@num\mili@u=\epsil@n\c@llgpart%
    \advance\mili@u-\epsil@n\advance\mili@u-\delt@\advance\v@leur-\mili@u%
    \ifcase\fclin@typ@\invers@\mili@u\delt@%
    \ifnum\@rrowr@fpt>\z@\advance\delt@-\v@leur\v@leur=\delt@\fi%
    \v@leur=\repdecn@mb\v@leur\mili@u\edef\v@lt{\repdecn@mb\v@leur}%
    \extrairelepremi@r\Ak@\de\list@num%
    \figvectPDD-1[\pr@m,\Aj@]\figpttraDD-6:=\Ai@/\curv@roundness,-1/%
    \figvectPDD-1[\Ak@,\Ai@]\figpttraDD-7:=\Aj@/\curv@roundness,-1/%
    \delt@=\@rrowheadlength\p@\delt@=\C@AHANG\delt@\edef\R@dius{\repdecn@mb{\delt@}}%
    \ifcase\@rrowr@fpt%
    \FigptintercircB@zDD-8::\v@lt,\R@dius[\Ai@,-6,-7,\Aj@]\psarrowheadDD[-5,-8]\else%
    \FigptintercircB@zDD-8::\v@lt,\R@dius[\Aj@,-7,-6,\Ai@]\psarrowheadDD[-8,-5]\fi%
    \else\advance\delt@-\v@leur%
    \p@rtentiere{\p@rtent}{\delt@}\edef\C@efun{\the\p@rtent}%
    \p@rtentiere{\p@rtent}{\v@leur}\edef\C@efde{\the\p@rtent}%
    \figptbaryDD-5:[\Ai@,\Aj@;\C@efun,\C@efde]\ifcase\@rrowr@fpt%
    \delt@=\@rrowheadlength\unit@\delt@=\C@AHANG\delt@\edef\t@ille{\repdecn@mb{\delt@}}%
    \figvectPDD-2[\Ai@,\Aj@]\vecunit@{-2}{-2}\figpttraDD-5:=-5/\t@ille,-2/\fi%
    \psarrowheadDD[\Ai@,-5]\fi}
\def\c@llgtot{\@ecfor\Aj@:=\list@num\do{\figvectP-1[\Ai@,\Aj@]\N@rm\delt@{-1}%
    \advance\v@leur\delt@\advance\p@rtent\@ne\edef\Ai@{\Aj@}}}
\def\c@llgpart{\extrairelepremi@r\Aj@\de\list@num\figvectP-1[\Ai@,\Aj@]\N@rm\delt@{-1}%
    \advance\mili@u\delt@\ifdim\mili@u<\v@leur\edef\pr@m{\Ai@}\edef\Ai@{\Aj@}\c@llgpart\fi}
\def\Pslin@conge[#1]{\ifnum\p@rtent>\tw@{\def\list@num{#1}%
    \extrairelepremi@r\Ai@\de\list@num\extrairelepremi@r\Aj@\de\list@num%
    \figptcopy-6:/\Ai@/\figvectP-3[\Ai@,\Aj@]\vecunit@{-3}{-3}\v@lmax=\result@t%
    \@ecfor\Ak@:=\list@num\do{\figvectP-4[\Aj@,\Ak@]\vecunit@{-4}{-4}%
    \minim@m\v@lmin\v@lmax\result@t\v@lmax=\result@t%
    \det@rm\delt@[-3,-4]\maxim@m\mili@u{\delt@}{-\delt@}\ifdim\mili@u>\Cepsil@n%
    \ifdim\delt@>\z@\figgetangleDD\Angl@[\Aj@,\Ak@,\Ai@]\else%
    \figgetangleDD\Angl@[\Aj@,\Ai@,\Ak@]\fi%
    \v@leur=\PI@deg\advance\v@leur-\Angl@\p@\divide\v@leur\tw@%
    \edef\Angl@{\repdecn@mb\v@leur}\c@ssin{\C@}{\S@}{\Angl@}\v@leur=\fclin@r@d\unit@%
    \v@leur=\S@\v@leur\mili@u=\C@\p@\invers@\mili@u\mili@u%
    \v@leur=\repdecn@mb{\mili@u}\v@leur%
    \minim@m\v@leur\v@leur\v@lmin\edef\t@ille{\repdecn@mb{\v@leur}}%
    \figpttra-5:=\Aj@/-\t@ille,-3/\psline[-6,-5]\figpttra-6:=\Aj@/\t@ille,-4/%
    \figvectNVDD-3[-3]\figvectNVDD-8[-4]\inters@cDD-7:[-5,-3;-6,-8]%
    \ifdim\delt@>\z@\psarccircP-7;\fclin@r@d[-5,-6]\else\psarccircP-7;\fclin@r@d[-6,-5]\fi%
    \else\psline[-6,\Aj@]\figptcopy-6:/\Aj@/\fi
    \edef\Ai@{\Aj@}\edef\Aj@{\Ak@}\figptcopy-3:/-4/}\psline[-6,\Aj@]}\else\psline[#1]\fi}
\def\psfcnode[#1]#2{{\ifcurr@ntPS\ifps@cri\PSc@mment{psfcnode Points=#1}%
    \s@uvc@ntr@l\et@tpsfcnode\resetc@ntr@l{2}%
    \def\t@xt@{#2}\ifx\t@xt@\empty\def\g@tt@xt{\setbox\Gb@x=\hbox{\Figg@tT{\p@int}}}%
    \else\def\g@tt@xt{\setbox\Gb@x=\hbox{#2}}\fi%
    \v@lmin=\h@rdfcXp@dd\advance\v@lmin\Xp@dd\unit@\multiply\v@lmin\tw@%
    \v@lmax=\h@rdfcYp@dd\advance\v@lmax\Yp@dd\unit@\multiply\v@lmax\tw@%
    \Figv@ctCreg-8(\unit@,-\unit@)\def\list@num{#1}%
    \delt@=\curr@ntwidth bp\divide\delt@\tw@%
    \fcn@de\PSc@mment{End psfcnode}\resetc@ntr@l\et@tpsfcnode\fi\fi}}
\def\d@butn@de{\g@tt@xt\v@lX=\wd\Gb@x%
    \v@lY=\ht\Gb@x\advance\v@lY\dp\Gb@x\advance\v@lX\v@lmin\advance\v@lY\v@lmax}
\def\fcn@deE{%
    \@ecfor\p@int:=\list@num\do{\d@butn@de\v@lX=\unssqrttw@\v@lX\v@lY=\unssqrttw@\v@lY%
    \ifdim\thickn@ss\p@>\z@
    \v@lXa=\v@lX\advance\v@lXa\delt@\v@lXa=\ptT@unit@\v@lXa\edef\XR@d{\repdecn@mb\v@lXa}%
    \v@lYa=\v@lY\advance\v@lYa\delt@\v@lYa=\ptT@unit@\v@lYa\edef\YR@d{\repdecn@mb\v@lYa}%
    \arct@n\v@leur(\v@lXa,\v@lYa)\v@leur=\rdT@deg\v@leur\edef\@nglde{\repdecn@mb\v@leur}%
    {\c@lptellDD-2::\p@int;\XR@d,\YR@d(\@nglde)}
    \advance\v@leur-\PI@deg\edef\@nglun{\repdecn@mb\v@leur}%
    {\c@lptellDD-3::\p@int;\XR@d,\YR@d(\@nglun)}%
    \figptstra-6=-3,-2,\p@int/\thickn@ss,-8/\pssetfillmode{yes}\us@secondC@lor%
    \psline[-2,-3,-6,-5]\psarcell-4;\XR@d,\YR@d(\@nglun,\@nglde,0)\fi
    \v@lX=\ptT@unit@\v@lX\v@lY=\ptT@unit@\v@lY%
    \edef\XR@d{\repdecn@mb\v@lX}\edef\YR@d{\repdecn@mb\v@lY}%
    \pssetfillmode{yes}\us@thirdC@lor\psarcell\p@int;\XR@d,\YR@d(0,360,0)%
    \pssetfillmode{no}\us@primarC@lor\psarcell\p@int;\XR@d,\YR@d(0,360,0)}}
\def\fcn@deL{\delt@=\ptT@unit@\delt@\edef\t@ille{\repdecn@mb\delt@}%
    \@ecfor\p@int:=\list@num\do{\Figg@tXYa{\p@int}\d@butn@de%
    \ifdim\v@lX>\v@lY\itis@Ktrue\else\itis@Kfalse\fi%
    \advance\v@lXa-\v@lX\Figp@intreg-1:(\v@lXa,\v@lYa)%
    \advance\v@lXa\v@lX\advance\v@lYa-\v@lY\Figp@intreg-2:(\v@lXa,\v@lYa)%
    \advance\v@lXa\v@lX\advance\v@lYa\v@lY\Figp@intreg-3:(\v@lXa,\v@lYa)%
    \advance\v@lXa-\v@lX\advance\v@lYa\v@lY\Figp@intreg-4:(\v@lXa,\v@lYa)%
    \ifdim\thickn@ss\p@>\z@\Figg@tXYa{\p@int}\pssetfillmode{yes}\us@secondC@lor
    \c@lpt@xt{-1}{-4}\c@lpt@xt@\v@lXa\v@lYa\v@lX\v@lY\c@rre\delt@%
    \Figp@intregDD-9:(\v@lZ,\v@lYa)\Figp@intregDD-11:(\v@lZa,\v@lYa)%
    \c@lpt@xt{-4}{-3}\c@lpt@xt@\v@lYa\v@lXa\v@lY\v@lX\delt@\c@rre%
    \Figp@intregDD-12:(\v@lXa,\v@lZ)\Figp@intregDD-10:(\v@lXa,\v@lZa)%
    \ifitis@K\figptstra-7=-9,-10,-11/\thickn@ss,-8/\psline[-9,-11,-5,-6,-7]\else%
    \figptstra-7=-10,-11,-12/\thickn@ss,-8/\psline[-10,-12,-5,-6,-7]\fi\fi
    \pssetfillmode{yes}\us@thirdC@lor\psline[-1,-2,-3,-4]%
    \pssetfillmode{no}\us@primarC@lor\psline[-1,-2,-3,-4,-1]}}
\def\c@lpt@xt#1#2{\figvectN-7[#1,#2]\vecunit@{-7}{-7}\figpttra-5:=#1/\t@ille,-7/%
    \figvectP-7[#1,#2]\Figg@tXY{-7}\c@rre=\v@lX\delt@=\v@lY\Figg@tXY{-5}}
\def\c@lpt@xt@#1#2#3#4#5#6{\v@lZ=#6\invers@{\v@lZ}{\v@lZ}\v@leur=\repdecn@mb{#5}\v@lZ%
    \v@lZ=#2\advance\v@lZ-#4\mili@u=\repdecn@mb{\v@leur}\v@lZ%
    \v@lZ=#3\advance\v@lZ\mili@u\v@lZa=-\v@lZ\advance\v@lZa\tw@#1}
\def\fcn@deR{\@ecfor\p@int:=\list@num\do{\Figg@tXYa{\p@int}\d@butn@de%
    \advance\v@lXa-0.5\v@lX\advance\v@lYa-0.5\v@lY\Figp@intreg-1:(\v@lXa,\v@lYa)%
    \advance\v@lXa\v@lX\Figp@intreg-2:(\v@lXa,\v@lYa)%
    \advance\v@lYa\v@lY\Figp@intreg-3:(\v@lXa,\v@lYa)%
    \advance\v@lXa-\v@lX\Figp@intreg-4:(\v@lXa,\v@lYa)%
    \ifdim\thickn@ss\p@>\z@\pssetfillmode{yes}\us@secondC@lor
    \Figv@ctCreg-5(-\delt@,-\delt@)\figpttra-9:=-1/1,-5/%
    \Figv@ctCreg-5(\delt@,-\delt@)\figpttra-10:=-2/1,-5/%
    \Figv@ctCreg-5(\delt@,\delt@)\figpttra-11:=-3/1,-5/%
    \figptstra-7=-9,-10,-11/\thickn@ss,-8/\psline[-9,-11,-5,-6,-7]\fi
    \pssetfillmode{yes}\us@thirdC@lor\psline[-1,-2,-3,-4]%
    \pssetfillmode{no}\us@primarC@lor\psline[-1,-2,-3,-4,-1]}}
\def\Pssetfl@wchart#1=#2|{\keln@mtr#1|%
    \def\n@mref{arr}\ifx\l@debut\n@mref\expandafter\keln@mtr\l@suite|%
     \def\n@mref{owp}\ifx\l@debut\n@mref\edef\@rrowp@s{#2}\else
     \def\n@mref{owr}\ifx\l@debut\n@mref\setfcr@fpt#2|\else
     \immediate\write16{*** Unknown attribute: \BS@ psset flowchart(..., #1=...)}%
     \fi\fi\else%
    \def\n@mref{lin}\ifx\l@debut\n@mref\setfccurv@#2|\else
    \def\n@mref{pad}\ifx\l@debut\n@mref\edef\Xp@dd{#2}\edef\Yp@dd{#2}\else
    \def\n@mref{rad}\ifx\l@debut\n@mref\edef\fclin@r@d{#2}\else
    \def\n@mref{sha}\ifx\l@debut\n@mref\setfcshap@#2|\else
    \def\n@mref{thi}\ifx\l@debut\n@mref\edef\thickn@ss{#2}\else
    \def\n@mref{xpa}\ifx\l@debut\n@mref\edef\Xp@dd{#2}\else
    \def\n@mref{ypa}\ifx\l@debut\n@mref\edef\Yp@dd{#2}\else
    \immediate\write16{*** Unknown attribute: \BS@ psset flowchart(..., #1=...)}%
    \fi\fi\fi\fi\fi\fi\fi\fi}
\def\setfcr@fpt#1#2|{\if#1e\def\@rrowr@fpt{1}\else\def\@rrowr@fpt{0}\fi}
\def\setfccurv@#1#2|{\if#1c\def\fclin@typ@{0}\else\def\fclin@typ@{1}\fi}
\def\setfcshap@#1#2|{%
    \if#1e\let\fcn@de=\fcn@deE\def\h@rdfcXp@dd{4pt}\def\h@rdfcYp@dd{4pt}%
     \edef\fcsh@pe{ellipse}\else%
    \if#1l\let\fcn@de=\fcn@deL\def\h@rdfcXp@dd{4pt}\def\h@rdfcYp@dd{4pt}%
     \edef\fcsh@pe{lozenge}\else%
          \let\fcn@de=\fcn@deR\def\h@rdfcXp@dd{6pt}\def\h@rdfcYp@dd{6pt}%
     \edef\fcsh@pe{rectangle}\fi\fi}
\def\psline[#1]{{\ifcurr@ntPS\ifps@cri\PSc@mment{psline Points=#1}%
    \let\pslign@=\Pslign@P\Pslin@{#1}\PSc@mment{End psline}\fi\fi}}
\def\pslineF#1{{\ifcurr@ntPS\ifps@cri\PSc@mment{pslineF Filename=#1}%
    \let\pslign@=\Pslign@F\Pslin@{#1}\PSc@mment{End pslineF}\fi\fi}}
\def\pslineC(#1){{\ifcurr@ntPS\ifps@cri\PSc@mment{pslineC}%
    \let\pslign@=\Pslign@C\Pslin@{#1}\PSc@mment{End pslineC}\fi\fi}}
\def\Pslin@#1{\iffillm@de\pslign@{#1}%
    \immediate\write\fwf@g{\fillc@md}%
    \else\pslign@{#1}\ifx\derp@int\premp@int%
    \immediate\write\fwf@g{closepath\space stroke}%
    \else\immediate\write\fwf@g{stroke}\fi\fi}
\def\Pslign@P#1{\def\list@num{#1}\extrairelepremi@r\p@int\de\list@num%
    \edef\premp@int{\p@int}\immediate\write\fwf@g{newpath}%
    \PSwrit@cmd{\p@int}{moveto}{\fwf@g}%
    \@ecfor\p@int:=\list@num\do{\PSwrit@cmd{\p@int}{lineto}{\fwf@g}%
    \edef\derp@int{\p@int}}}
\def\Pslign@F#1{\s@uvc@ntr@l\et@tPslign@F\setc@ntr@l{2}\openin\frf@g=#1\relax%
    \ifeof\frf@g\message{*** File #1 not found !}\end\else%
    \read\frf@g to\tr@c\edef\premp@int{\tr@c}\expandafter\extr@ctCF\tr@c:%
    \immediate\write\fwf@g{newpath}\PSwrit@cmd{-1}{moveto}{\fwf@g}%
    \loop\read\frf@g to\tr@c\ifeof\frf@g\mored@tafalse\else\mored@tatrue\fi%
    \ifmored@ta\expandafter\extr@ctCF\tr@c:\PSwrit@cmd{-1}{lineto}{\fwf@g}%
    \edef\derp@int{\tr@c}\repeat\fi\closein\frf@g\resetc@ntr@l\et@tPslign@F}
\def\extr@ctCFDD#1 #2:{\v@lX=#1\unit@\v@lY=#2\unit@\Figp@intregDD-1:(\v@lX,\v@lY)}
\def\extr@ctCFTD#1 #2 #3:{\v@lX=#1\unit@\v@lY=#2\unit@\v@lZ=#3\unit@%
    \Figp@intregTD-1:(\v@lX,\v@lY,\v@lZ)}
\def\Pslign@C#1{\s@uvc@ntr@l\et@tPslign@C\setc@ntr@l{2}%
    \def\list@num{#1}\extrairelepremi@r\p@int\de\list@num%
    \edef\premp@int{\p@int}\immediate\write\fwf@g{newpath}%
    \expandafter\Pslign@C@\p@int:\PSwrit@cmd{-1}{moveto}{\fwf@g}%
    \@ecfor\p@int:=\list@num\do{\expandafter\Pslign@C@\p@int:%
    \PSwrit@cmd{-1}{lineto}{\fwf@g}\edef\derp@int{\p@int}}%
    \resetc@ntr@l\et@tPslign@C}
\def\Pslign@C@#1 #2:{{\def\t@xt@{#1}\ifx\t@xt@\empty\Pslign@C@#2:
    \else\extr@ctCF#1 #2:\fi}}
\newcount\c@ntrolmesh
\def\Pssetm@sh#1=#2|{\keln@mun#1|%
    \def\n@mref{d}\ifx\l@debut\n@mref\pssetmeshdiag{#2}\else
    \immediate\write16{*** Unknown attribute: \BS@ psset mesh(..., #1=...)}%
    \fi}
\def\pssetmeshdiag#1{\c@ntrolmesh=#1}
\def\defaultmeshdiag{0}    
\def\psmesh#1,#2[#3,#4,#5,#6]{{\ifcurr@ntPS\ifps@cri%
    \PSc@mment{psmesh N1=#1, N2=#2, Quadrangle=[#3,#4,#5,#6]}%
    \s@uvc@ntr@l\et@tpsmesh\Pss@tsecondSt\setc@ntr@l{2}%
    \ifnum#1>\@ne\Psmeshp@rt#1[#3,#4,#5,#6]\fi%
    \ifnum#2>\@ne\Psmeshp@rt#2[#4,#5,#6,#3]\fi%
    \ifnum\c@ntrolmesh>\z@\Psmeshdi@g#1,#2[#3,#4,#5,#6]\fi%
    \ifnum\c@ntrolmesh<\z@\Psmeshdi@g#2,#1[#4,#5,#6,#3]\fi\Psrest@reSt%
    \psline[#3,#4,#5,#6,#3]\PSc@mment{End psmesh}\resetc@ntr@l\et@tpsmesh\fi\fi}}
\def\Psmeshp@rt#1[#2,#3,#4,#5]{{\l@mbd@un=\@ne\l@mbd@de=#1\loop%
    \ifnum\l@mbd@un<#1\advance\l@mbd@de\m@ne\figptbary-1:[#2,#3;\l@mbd@de,\l@mbd@un]%
    \figptbary-2:[#5,#4;\l@mbd@de,\l@mbd@un]\psline[-1,-2]\advance\l@mbd@un\@ne\repeat}}
\def\Psmeshdi@g#1,#2[#3,#4,#5,#6]{\figptcopy-2:/#3/\figptcopy-3:/#6/%
    \l@mbd@un=\z@\l@mbd@de=#1\loop\ifnum\l@mbd@un<#1%
    \advance\l@mbd@un\@ne\advance\l@mbd@de\m@ne\figptcopy-1:/-2/\figptcopy-4:/-3/%
    \figptbary-2:[#3,#4;\l@mbd@de,\l@mbd@un]%
    \figptbary-3:[#6,#5;\l@mbd@de,\l@mbd@un]\Psmeshdi@gp@rt#2[-1,-2,-3,-4]\repeat}
\def\Psmeshdi@gp@rt#1[#2,#3,#4,#5]{{\l@mbd@un=\z@\l@mbd@de=#1\loop%
    \ifnum\l@mbd@un<#1\figptbary-5:[#2,#5;\l@mbd@de,\l@mbd@un]%
    \advance\l@mbd@de\m@ne\advance\l@mbd@un\@ne%
    \figptbary-6:[#3,#4;\l@mbd@de,\l@mbd@un]\psline[-5,-6]\repeat}}
\def\psnormalDD#1,#2[#3,#4]{{\ifcurr@ntPS\ifps@cri%
    \PSc@mment{psnormal Length=#1, Lambda=#2 [Pt1,Pt2]=[#3,#4]}%
    \s@uvc@ntr@l\et@tpsnormal\resetc@ntr@l{2}\figptendnormal-6::#1,#2[#3,#4]%
    \figptcopyDD-5:/-1/\psarrow[-5,-6]%
    \PSc@mment{End psnormal}\resetc@ntr@l\et@tpsnormal\fi\fi}}
\def\psreset#1{\trtlis@rg{#1}{\Psreset@}}
\def\Psreset@#1|{\keln@mde#1|%
    \def\n@mref{ar}\ifx\l@debut\n@mref\psresetarrowhead\else
    \def\n@mref{cu}\ifx\l@debut\n@mref\psset curve(roundness=\defaultroundness)\else
    \def\n@mref{fi}\ifx\l@debut\n@mref\psset (color=\defaultcolor,dash=\defaultdash,%
         fill=\defaultfill,join=\defaultjoin,width=\defaultwidth)\else
    \def\n@mref{fl}\ifx\l@debut\n@mref\psset flowchart(arrowp=\defaultfcarrowposition,%
	arrowr=\defaultfcarrowrefpt,line=\defaultfcline,xpadd=\defaultfcxpadding,%
	ypadd=\defaultfcypadding,radius=\defaultfcradius,shape=\defaultfcshape,%
	thick=\defaultfcthickness)\else
    \def\n@mref{me}\ifx\l@debut\n@mref\psset mesh(diag=\defaultmeshdiag)\else
    \def\n@mref{se}\ifx\l@debut\n@mref\psresetsecondsettings\else
    \def\n@mref{th}\ifx\l@debut\n@mref\psset third(color=\defaultthirdcolor)\else
    \immediate\write16{*** Unknown keyword #1 (\BS@ psreset).}%
    \fi\fi\fi\fi\fi\fi\fi}
\def\psset#1(#2){\def\t@xt@{#1}\ifx\t@xt@\empty\trtlis@rg{#2}{\Pssetf@rst}
    \else\keln@mde#1|%
    \def\n@mref{ar}\ifx\l@debut\n@mref\trtlis@rg{#2}{\Psset@rrowhe@d}\else
    \def\n@mref{cu}\ifx\l@debut\n@mref\trtlis@rg{#2}{\Pssetc@rve}\else
    \def\n@mref{fi}\ifx\l@debut\n@mref\trtlis@rg{#2}{\Pssetf@rst}\else
    \def\n@mref{fl}\ifx\l@debut\n@mref\trtlis@rg{#2}{\Pssetfl@wchart}\else
    \def\n@mref{me}\ifx\l@debut\n@mref\trtlis@rg{#2}{\Pssetm@sh}\else
    \def\n@mref{se}\ifx\l@debut\n@mref\trtlis@rg{#2}{\Pssets@cond}\else
    \def\n@mref{th}\ifx\l@debut\n@mref\trtlis@rg{#2}{\Pssetth@rd}\else
    \immediate\write16{*** Unknown keyword: \BS@ psset #1(...)}%
    \fi\fi\fi\fi\fi\fi\fi\fi}
\def\pssetdefault#1(#2){\ifcurr@ntPS\immediate\write16{*** \BS@ pssetdefault is ignored
    inside a \BS@ psbeginfig-\BS@ psendfig block.}%
    \immediate\write16{*** It must be called before \BS@ psbeginfig.}\else%
    \def\t@xt@{#1}\ifx\t@xt@\empty\trtlis@rg{#2}{\Pssd@f@rst}\else\keln@mde#1|%
    \def\n@mref{ar}\ifx\l@debut\n@mref\trtlis@rg{#2}{\Pssd@@rrowhe@d}\else
    \def\n@mref{cu}\ifx\l@debut\n@mref\trtlis@rg{#2}{\Pssd@c@rve}\else
    \def\n@mref{fi}\ifx\l@debut\n@mref\trtlis@rg{#2}{\Pssd@f@rst}\else
    \def\n@mref{fl}\ifx\l@debut\n@mref\trtlis@rg{#2}{\Pssd@fl@wchart}\else
    \def\n@mref{me}\ifx\l@debut\n@mref\trtlis@rg{#2}{\Pssd@m@sh}\else
    \def\n@mref{se}\ifx\l@debut\n@mref\trtlis@rg{#2}{\Pssd@s@cond}\else
    \def\n@mref{th}\ifx\l@debut\n@mref\trtlis@rg{#2}{\Pssd@th@rd}\else
    \immediate\write16{*** Unknown keyword: \BS@ pssetdefault #1(...)}%
    \fi\fi\fi\fi\fi\fi\fi\fi\initpss@ttings\fi}
\def\Pssd@f@rst#1=#2|{\keln@mun#1|%
    \def\n@mref{c}\ifx\l@debut\n@mref\edef\defaultcolor{#2}\else
    \def\n@mref{d}\ifx\l@debut\n@mref\edef\defaultdash{#2}\else
    \def\n@mref{f}\ifx\l@debut\n@mref\edef\defaultfill{#2}\else
    \def\n@mref{j}\ifx\l@debut\n@mref\edef\defaultjoin{#2}\else
    \def\n@mref{u}\ifx\l@debut\n@mref\edef\defaultupdate{#2}\pssetupdate{#2}\else
    \def\n@mref{w}\ifx\l@debut\n@mref\edef\defaultwidth{#2}\else
    \immediate\write16{*** Unknown attribute: \BS@ pssetdefault (..., #1=...)}%
    \fi\fi\fi\fi\fi\fi}
\def\Pssd@@rrowhe@d#1=#2|{\keln@mun#1|%
    \def\n@mref{a}\ifx\l@debut\n@mref\edef\defaultarrowheadangle{#2}\else
    \def\n@mref{f}\ifx\l@debut\n@mref\edef\defaultarrowheadangle{#2}\else
    \def\n@mref{l}\ifx\l@debut\n@mref\y@tiunit{#2}\ifunitpr@sent%
     \edef\defaulth@rdahlength{#2}\else\edef\defaulth@rdahlength{#2pt}%
     \message{*** \BS@ pssetdefault (..., #1=#2, ...) : unit is missing, pt is assumed.}%
     \fi\else
    \def\n@mref{o}\ifx\l@debut\n@mref\edef\defaultarrowheadout{#2}\else
    \def\n@mref{r}\ifx\l@debut\n@mref\edef\defaultarrowheadratio{#2}\else
    \immediate\write16{*** Unknown attribute: \BS@ pssetdefault arrowhead(..., #1=...)}%
    \fi\fi\fi\fi\fi}
\def\Pssd@c@rve#1=#2|{\keln@mun#1|%
    \def\n@mref{r}\ifx\l@debut\n@mref\edef\defaultroundness{#2}\else%
    \immediate\write16{*** Unknown attribute: \BS@ pssetdefault curve(..., #1=...)}%
    \fi}
\def\Pssd@fl@wchart#1=#2|{\keln@mtr#1|%
    \def\n@mref{arr}\ifx\l@debut\n@mref\expandafter\keln@mtr\l@suite|%
     \def\n@mref{owp}\ifx\l@debut\n@mref\edef\defaultfcarrowposition{#2}\else
     \def\n@mref{owr}\ifx\l@debut\n@mref\edef\defaultfcarrowrefpt{#2}\else
     \immediate\write16{*** Unknown attribute: \BS@ pssetdefault flowchart(..., #1=...)}%
     \fi\fi\else%
    \def\n@mref{lin}\ifx\l@debut\n@mref\edef\defaultfcline{#2}\else
    \def\n@mref{pad}\ifx\l@debut\n@mref\edef\defaultfcxpadding{#2}%
                    \edef\defaultfcypadding{#2}\else
    \def\n@mref{rad}\ifx\l@debut\n@mref\edef\defaultfcradius{#2}\else
    \def\n@mref{sha}\ifx\l@debut\n@mref\edef\defaultfcshape{#2}\else
    \def\n@mref{thi}\ifx\l@debut\n@mref\edef\defaultfcthickness{#2}\else
    \def\n@mref{xpa}\ifx\l@debut\n@mref\edef\defaultfcxpadding{#2}\else
    \def\n@mref{ypa}\ifx\l@debut\n@mref\edef\defaultfcypadding{#2}\else
    \immediate\write16{*** Unknown attribute: \BS@ pssetdefault flowchart(..., #1=...)}%
    \fi\fi\fi\fi\fi\fi\fi\fi}
\def\defaultfcarrowposition{0.5}\let\defaultfcarrowpos=\defaultfcarrowposition
\def\defaultfcarrowrefpt{start}
\def\defaultfcline{polygon}
\def\defaultfcradius{0}
\def\defaultfcshape{rectangle}
\def\defaultfcthickness{0}\let\defaultfcthick=\defaultfcthickness
\def\defaultfcxpadding{0}\let\defaultfcxpad=\defaultfcxpadding
\def\defaultfcypadding{0}\let\defaultfcypad=\defaultfcypadding
\def\Pssd@m@sh#1=#2|{\keln@mun#1|%
    \def\n@mref{d}\ifx\l@debut\n@mref\edef\defaultmeshdiag{#2}\else%
    \immediate\write16{*** Unknown attribute: \BS@ pssetdefault mesh(..., #1=...)}%
    \fi}
\def\Pssd@s@cond#1=#2|{\keln@mun#1|%
    \def\n@mref{c}\ifx\l@debut\n@mref\edef\defaultsecondcolor{#2}\else%
    \def\n@mref{d}\ifx\l@debut\n@mref\edef\defaultseconddash{#2}\else%
    \def\n@mref{w}\ifx\l@debut\n@mref\edef\defaultsecondwidth{#2}\else%
    \immediate\write16{*** Unknown attribute: \BS@ pssetdefault second(..., #1=...)}%
    \fi\fi\fi}
\def\Pssd@th@rd#1=#2|{\keln@mun#1|%
    \def\n@mref{c}\ifx\l@debut\n@mref\edef\defaultthirdcolor{#2}\else%
    \immediate\write16{*** Unknown attribute: \BS@ pssetdefault third(..., #1=...)}%
    \fi}
\newif\iffillm@de
\def\pssetfillmode#1{\expandafter\setfillm@de#1:}
\def\setfillm@de#1#2:{\if#1n\fillm@defalse\else\fillm@detrue\fi}
\def\defaultfill{no}     
\newif\ifpstestm@de
\def\pssetupdate#1{\ifcurr@ntPS\immediate\write16{*** \BS@ pssetupdate is ignored inside a
     \BS@ psbeginfig-\BS@ psendfig block.}%
    \immediate\write16{*** It must be called before \BS@ psbeginfig.}%
    \else\expandafter\setupd@te#1:\fi}
\def\setupd@te#1#2:{\if#1n\pstestm@defalse\else\pstestm@detrue\fi}
\def\defaultupdate{no}     
\def\Pssetc@lor#1{\ifps@cri\result@tent=\@ne\expandafter\c@lnbV@l#1 :%
    \def\curr@ntcolor{}\def\curr@ntcolorc@md{}%
    \ifcase\result@tent\or\pssetgray{#1}\or\or\pssetrgb{#1}\or\pssetcmyk{#1}\fi\fi}
\def\pssetcmyk#1{\ifps@cri\def\curr@ntcolor{#1}\def\curr@ntcolorc@md{setcmykcolor}%
    \ifcurr@ntPS\PSc@mment{pssetcmyk Color=#1}\us@primarC@lor\fi\fi}
\def\pssetrgb#1{\ifps@cri\def\curr@ntcolor{#1}\def\curr@ntcolorc@md{setrgbcolor}%
    \ifcurr@ntPS\PSc@mment{pssetrgb Color=#1}\us@primarC@lor\fi\fi}
\def\pssetgray#1{\ifps@cri\def\curr@ntcolor{#1}\def\curr@ntcolorc@md{setgray}%
    \ifcurr@ntPS\PSc@mment{pssetgray Gray level=#1}\us@primarC@lor\fi\fi}
\def\us@primarC@lor{\immediate\write\fwf@g{\curr@ntcolor\space\curr@ntcolorc@md}%
    \let\fillc@md=\prfillc@md}
\def\prfillc@md{\curr@ntcolor\space\curr@ntcolorc@md\space fill}
\def\defaultcolor{0}       
\def\c@lnbV@l#1 #2:{\def\t@xt@{#1}\relax\ifx\t@xt@\empty\c@lnbV@l#2:
    \else\c@lnbV@l@#1 #2:\fi}
\def\c@lnbV@l@#1 #2:{\def\t@xt@{#2}\ifx\t@xt@\empty%
    \def\t@xt@{#1}\ifx\t@xt@\empty\advance\result@tent\m@ne\fi
    \else\advance\result@tent\@ne\c@lnbV@l@#2:\fi}
\def\Blackcmyk{0 0 0 1}
\def\Whitecmyk{0 0 0 0}
\def\Cyancmyk{1 0 0 0}
\def\Magentacmyk{0 1 0 0}
\def\Yellowcmyk{0 0 1 0}
\def\Redcmyk{0 1 1 0}
\def\Greencmyk{1 0 1 0}
\def\Bluecmyk{1 1 0 0}
\def\Graycmyk{0 0 0 0.50}
\def\BrickRedcmyk{0 0.89 0.94 0.28} 
\def\Browncmyk{0 0.81 1 0.60} 
\def\ForestGreencmyk{0.91 0 0.88 0.12} 
\def\Goldenrodcmyk{ 0 0.10 0.84 0} 
\def\Marooncmyk{0 0.87 0.68 0.32} 
\def\Orangecmyk{0 0.61 0.87 0} 
\def\Purplecmyk{0.45 0.86 0 0} 
\def\RoyalBluecmyk{1. 0.50 0 0} 
\def\Violetcmyk{0.79 0.88 0 0} 
\def\Blackrgb{0 0 0}
\def\Whitergb{1 1 1}
\def\Redrgb{1 0 0}
\def\Greenrgb{0 1 0}
\def\Bluergb{0 0 1}
\def\Cyanrgb{0 1 1}
\def\Magentargb{1 0 1}
\def\Yellowrgb{1 1 0}
\def\Grayrgb{0.5 0.5 0.5}
\def\Chocolatergb{0.824 0.412 0.118}
\def\DarkGoldenrodrgb{0.722 0.525 0.043}
\def\DarkOrangergb{1 0.549 0}
\def\Firebrickrgb{0.698 0.133 0.133}
\def\ForestGreenrgb{0.133 0.545 0.133}
\def\Goldrgb{1 0.843 0}
\def\HotPinkrgb{1 0.412 0.706}
\def\Maroonrgb{0.690 0.188 0.376}
\def\Pinkrgb{1 0.753 0.796}
\def\RoyalBluergb{0.255 0.412 0.882}
\def\Pssetf@rst#1=#2|{\keln@mun#1|%
    \def\n@mref{c}\ifx\l@debut\n@mref\Pssetc@lor{#2}\else
    \def\n@mref{d}\ifx\l@debut\n@mref\pssetdash{#2}\else
    \def\n@mref{f}\ifx\l@debut\n@mref\pssetfillmode{#2}\else
    \def\n@mref{j}\ifx\l@debut\n@mref\pssetjoin{#2}\else
    \def\n@mref{u}\ifx\l@debut\n@mref\pssetupdate{#2}\else
    \def\n@mref{w}\ifx\l@debut\n@mref\pssetwidth{#2}\else
    \immediate\write16{*** Unknown attribute: \BS@ psset (..., #1=...)}%
    \fi\fi\fi\fi\fi\fi}
\def\s@uvdash#1{\edef#1{\curr@ntdash}}
\def\defaultdash{1}        
\def\pssetdash#1{\ifps@cri\edef\curr@ntdash{#1}\ifcurr@ntPS\expandafter\Pssetd@sh#1 :\fi\fi}
\def\Pssetd@shI#1{\PSc@mment{pssetdash Index=#1}\ifcase#1%
    \or\immediate\write\fwf@g{[] 0 setdash}
    \or\immediate\write\fwf@g{[6 2] 0 setdash}
    \or\immediate\write\fwf@g{[4 2] 0 setdash}
    \or\immediate\write\fwf@g{[2 2] 0 setdash}
    \or\immediate\write\fwf@g{[1 2] 0 setdash}
    \or\immediate\write\fwf@g{[2 4] 0 setdash}
    \or\immediate\write\fwf@g{[3 5] 0 setdash}
    \or\immediate\write\fwf@g{[3 3] 0 setdash}
    \or\immediate\write\fwf@g{[3 5 1 5] 0 setdash}
    \or\immediate\write\fwf@g{[6 4 2 4] 0 setdash}
    \fi}
\def\Pssetd@sh#1 #2:{{\def\t@xt@{#1}\ifx\t@xt@\empty\Pssetd@sh#2:
    \else\def\t@xt@{#2}\ifx\t@xt@\empty\Pssetd@shI{#1}\else\s@mme=\@ne\def\debutp@t{#1}%
    \an@lysd@sh#2:\ifodd\s@mme\edef\debutp@t{\debutp@t\space\finp@t}\def\finp@t{0}\fi%
    \PSc@mment{pssetdash Pattern=#1 #2}%
    \immediate\write\fwf@g{[\debutp@t] \finp@t\space setdash}\fi\fi}}
\def\an@lysd@sh#1 #2:{\def\t@xt@{#2}\ifx\t@xt@\empty\def\finp@t{#1}\else%
    \edef\debutp@t{\debutp@t\space#1}\advance\s@mme\@ne\an@lysd@sh#2:\fi}
\def\s@uvwidth#1{\edef#1{\curr@ntwidth}}
\def\defaultwidth{0.4}     
\def\pssetwidth#1{\ifps@cri\edef\curr@ntwidth{#1}\ifcurr@ntPS%
    \PSc@mment{pssetwidth Width=#1}\immediate\write\fwf@g{#1 setlinewidth}\fi\fi}
\def\pssetjoin#1{\ifps@cri\edef\curr@ntjoin{#1}\ifcurr@ntPS\expandafter\Pssetj@in#1:\fi\fi}
\def\Pssetj@in#1#2:{\PSc@mment{pssetjoin join=#1}%
    \if#1r\def\t@xt@{1}\else\if#1b\def\t@xt@{2}\else\def\t@xt@{0}\fi\fi%
    \immediate\write\fwf@g{\t@xt@\space setlinejoin}}
\def\defaultjoin{miter}   
\def\Pssets@cond#1=#2|{\keln@mun#1|%
    \def\n@mref{c}\ifx\l@debut\n@mref\Pssets@condcolor{#2}\else%
    \def\n@mref{d}\ifx\l@debut\n@mref\pssetseconddash{#2}\else%
    \def\n@mref{w}\ifx\l@debut\n@mref\pssetsecondwidth{#2}\else%
    \immediate\write16{*** Unknown attribute: \BS@ psset second(..., #1=...)}%
    \fi\fi\fi}
\def\pssetseconddash#1{\edef\curr@ntseconddash{#1}}
\def\defaultseconddash{4}  
\def\pssetsecondwidth#1{\edef\curr@ntsecondwidth{#1}}
\edef\defaultsecondwidth{\defaultwidth} 
\def\psresetsecondsettings{%
    \pssetseconddash{\defaultseconddash}\pssetsecondwidth{\defaultsecondwidth}%
    \Pssets@condcolor{\defaultsecondcolor}}
\def\Pssets@condcolor#1{\ifps@cri\result@tent=\@ne\expandafter\c@lnbV@l#1 :%
    \def\sec@ndcolor{}\def\sec@ndcolorc@md{}%
    \ifcase\result@tent\or\pssetsecondgray{#1}\or\or\pssetsecondrgb{#1}%
    \or\pssetsecondcmyk{#1}\fi\fi}
\def\pssetsecondcmyk#1{\def\sec@ndcolor{#1}\def\sec@ndcolorc@md{setcmykcolor}}
\def\pssetsecondrgb#1{\def\sec@ndcolor{#1}\def\sec@ndcolorc@md{setrgbcolor}}
\def\pssetsecondgray#1{\def\sec@ndcolor{#1}\def\sec@ndcolorc@md{setgray}}
\def\us@secondC@lor{\immediate\write\fwf@g{\sec@ndcolor\space\sec@ndcolorc@md}%
    \let\fillc@md=\sdfillc@md}
\def\sdfillc@md{\sec@ndcolor\space\sec@ndcolorc@md\space fill}
\edef\defaultsecondcolor{\defaultcolor} 
\def\Pss@tsecondSt{%
    \s@uvdash{\typ@dash}\pssetdash{\curr@ntseconddash}%
    \s@uvwidth{\typ@width}\pssetwidth{\curr@ntsecondwidth}\us@secondC@lor}
\def\Psrest@reSt{\pssetwidth{\typ@width}\pssetdash{\typ@dash}\us@primarC@lor}
\def\Pssetth@rd#1=#2|{\keln@mun#1|%
    \def\n@mref{c}\ifx\l@debut\n@mref\Pssetth@rdcolor{#2}\else%
    \immediate\write16{*** Unknown attribute: \BS@ psset third(..., #1=...)}%
    \fi}
\def\Pssetth@rdcolor#1{\ifps@cri\result@tent=\@ne\expandafter\c@lnbV@l#1 :%
    \def\th@rdcolor{}\def\th@rdcolorc@md{}%
    \ifcase\result@tent\or\Pssetth@rdgray{#1}\or\or\Pssetth@rdrgb{#1}%
    \or\Pssetth@rdcmyk{#1}\fi\fi}
\def\Pssetth@rdcmyk#1{\def\th@rdcolor{#1}\def\th@rdcolorc@md{setcmykcolor}}
\def\Pssetth@rdrgb#1{\def\th@rdcolor{#1}\def\th@rdcolorc@md{setrgbcolor}}
\def\Pssetth@rdgray#1{\def\th@rdcolor{#1}\def\th@rdcolorc@md{setgray}}
\def\us@thirdC@lor{\immediate\write\fwf@g{\th@rdcolor\space\th@rdcolorc@md}%
    \let\fillc@md=\thfillc@md}
\def\thfillc@md{\th@rdcolor\space\th@rdcolorc@md\space fill}
\def\defaultthirdcolor{1}  
\def\pstrimesh#1[#2,#3,#4]{{\ifcurr@ntPS\ifps@cri%
    \PSc@mment{pstrimesh Type=#1, Triangle=[#2,#3,#4]}%
    \s@uvc@ntr@l\et@tpstrimesh\ifnum#1>\@ne\Pss@tsecondSt\setc@ntr@l{2}%
    \Pstrimeshp@rt#1[#2,#3,#4]\Pstrimeshp@rt#1[#3,#4,#2]%
    \Pstrimeshp@rt#1[#4,#2,#3]\Psrest@reSt\fi\psline[#2,#3,#4,#2]%
    \PSc@mment{End pstrimesh}\resetc@ntr@l\et@tpstrimesh\fi\fi}}
\def\Pstrimeshp@rt#1[#2,#3,#4]{{\l@mbd@un=\@ne\l@mbd@de=#1\loop\ifnum\l@mbd@de>\@ne%
    \advance\l@mbd@de\m@ne\figptbary-1:[#2,#3;\l@mbd@de,\l@mbd@un]%
    \figptbary-2:[#2,#4;\l@mbd@de,\l@mbd@un]\psline[-1,-2]%
    \advance\l@mbd@un\@ne\repeat}}
\initpr@lim\initpss@ttings
\catcode`\@=12

%% file: 1105.bbl
\begin{thebibliography}{99}

 \bibitem[{Bar}]{Bar} S. A.~Barannikov. 
  \newblock The framed Morse complex and its invariants. 
  \newblock Singularities and bifurcations, pp.~93--115,
Adv. Soviet Math., 21, Amer. Math. Soc. (1994). 


 \bibitem[{Bis}]{Bis} J.~M.~Bismut.
  \newblock The Witten complex and the degenerate Morse inequalities. 
  \newblock J.~Differ.~Geom.~23, pp.~207-240 (1986). 
  
  \bibitem[{BiLe}]{BiLe}
J.~M.~Bismut and G.~Lebeau.
\newblock{\sl The hypoelliptic Laplacian and Ray-Singer metrics.}
\newblock Annals of Mathematics Studies, 167, Princeton University Press (2008).

\bibitem[{Bis2}]{Bis2} J.~M.~Bismut.
\newblock Hypoelliptic Laplacian and Bott-Chern cohomology,
\newblock Comptes Rendus Mathematique, \textbf{349}-2 (2011), pp. 75--80 

\bibitem[{BEGK}]{BEGK} A.~Bovier, M.~Eckhoff, V.~Gayrard, and  M.~Klein.
  \newblock Metastability in reversible diffusion processes  I:
 Sharp asymptotics for capacities and exit times.
 \newblock  JEMS 6 (4), pp.~399-424 (2004).


 \bibitem[{BGK}]{BoGaKl} A.~Bovier,  V.~Gayrard,   and  M.~Klein.
  \newblock Metastability in reversible diffusion processes  II:
 Precise asymptotics for small eigenvalues.
 \newblock JEMS 7 (1), pp.~69-99 (2004).

\bibitem[{Bot}]{Bott} R.~Bott.
\newblock Lectures on Morse theory, old and new.
\newblock Bull. Amer. Math. Soc. (N.S.) \textbf{7}2 (1982)
pp. 331--358.

\bibitem[Bot2]{Bot2} R.~Bott.
\newblock Morse theory indomitable.
\newblock Publications IHES, tome 68 (1988) pp 99-114.

\bibitem[{BoTu}]{BoTu}
\newblock R.~Bott, L.~Tu.
\newblock\textit{Differential forms in algebraic topology.}
\newblock  Graduate Texts in Mathematics, 82. Springer-Verlag, (1982).

 \bibitem[{ChLi}]{ChLi} K.~C.~Chang and  J.~Liu.
 \newblock A cohomology complex for manifolds with boundary.
 \newblock Topological Methods in Non Linear Analysis,
 Vol.~5, pp.~325-340 (1995).

\bibitem[{CFKS}]{CFKS}
H.~L.~Cycon, R.~G.~Froese, W.~Kirsch, and B.~Simon.
\newblock{\it Schr{\"o}dinger operators with application to quantum
  mechanics and global geometry}.
\newblock Text and Monographs in Physics. Springer Verlag (1987).

\bibitem[{DiSj}]{DiSj} M.~Dimassi and  J.~Sj{\"o}strand.
 \newblock{\it  Spectral Asymptotics in the semi-classical limit}.
 \newblock London Mathematical Society, Lecture Note Series 268, 
 Cambridge University Press (1999).

\bibitem[{FrWe}]{FrWe} M.~I.~Freidlin and  A.~D.~Wentzell.
 \newblock{\it Random perturbations of dynamical systems. Transl. from the Russian by Joseph Szuecs. 2nd ed.}
 \newblock Grundlehren der Mathematischen Wissenschaften, 260, Springer-Verlag (1998). 

\bibitem[{Ful}]{Ful} W.~Fulton.
\newblock \textit{Algebraic topology. A first course.}
\newblock Graduate Texts in Mathematics, 153. Springer-Verlag, (1995)

\bibitem[{Gue}]{Gue} P.~Guerini.
\newblock Prescription du spectre du laplacien de Hodge-de Rham.
\newblock  Ann. Sci. \'Ecole Norm. Sup. (4), \textbf{37}-2 (2004)
pp. 270–-303.

\bibitem[{Hat}]{Hat} A.~Hatcher.
 \newblock{\it  Algebraic topology}.
 \newblock Cambridge University Press (2002).

\bibitem[{Hel}]{Hel} B.~Helffer.
 \newblock {\it Introduction to the semi-classical Analysis for the
 Schr{\"o}dinger operator and applications.}
\newblock Springer-Verlag.  Lecture Notes in Mathematics 1336
 (1988).  
 
\bibitem[{HKN}]{HKN} B.~Helffer, M.~Klein, and F.~Nier.
 \newblock Quantitative analysis of metastability in reversible
 diffusion processes via  a Witten complex approach.
 \newblock Matematica Contemporanea, 26, pp.~41-85 (2004).

 \bibitem[{HeNi1}]{HeNi} B.~Helffer and F.~Nier.
 \newblock \textit{Quantitative analysis of metastability in 
reversible diffusion processes via a Witten complex approach: the case with boundary}.
 \newblock M\'emoire 105, Soci\'et\'e Math\'ematique de France (2006).

 \bibitem[{HeNi2}]{HeNi2} B.~Helffer and F.~Nier.
\newblock \textit{Hypoelliptic estimates and spectral theory for Fokker-Planck operators and Witten Laplacians}.
Springer-Verlag. Lecture Notes in Mathematics 1862 (2005).

\bibitem[{HeSj2}]{HelSj2} B.~Helffer and J.~Sj{\"o}strand.
\newblock Puits multiples en limite semi-classique II
-Interaction mol{\'e}culaire-Sym{\'e}tries-Perturbations.
\newblock Ann. Inst. H. Poincar\'e Phys. Th\'eor. 42 (2),
pp.~127--212 (1985). 

\bibitem[HeSj3]{HeSj3} B.~Helffer and J.~Sj{\"o}strand.
\newblock Multiple wells in the semi-classical limit III,
\newblock Math.Nachrichte 124 (1985) p. 263-313.

\bibitem[{HeSj4}]{HelSj4} B.~Helffer and J.~Sj{\"o}strand.
\newblock   Puits multiples en limite semi-classique IV
-Etude du complexe de Witten -.
\newblock  Comm. Partial Differential Equations 10 (3), pp.~245-340
(1985).


 \bibitem[{HerNi}]{HerNi}
F.~H\'erau and F.~Nier.
\newblock Isotropic hypoellipticity and trend to the equilibrium
  for the Fokker-Planck equation with high degree potential. 
\newblock Archive for Rational Mechanics and Analysis 171 (2), pp.~151-218
(2004).
 
 \bibitem[{HSS}]{HerSjSt}
F.~H\'erau, J.~Sj{\"o}strand, and C.~Stolk.
\newblock Semiclassical analysis for the Kramers-Fokker-Planck equation.
\newblock Comm. Partial Differential Equations Vol. 30
no. 4-6, pp.~689-760 (2005).

\bibitem[{HHS1}]{HerHiSj}
F.~H\'erau, M. Hitrik, and J.~Sj{\"o}strand. 
\newblock Tunnel effect for Kramers-Fokker-Planck type operators.
\newblock Ann. Henri Poincar\'e 9, no. 2, pp.~209--274 (2008).

\bibitem[{HHS2}]{HerHiSj2}
F.~H\'erau, M. Hitrik, and J.~Sj{\"o}strand. 
\newblock Tunnel effect and symmetries for Kramers Fokker-Planck type operators.
\newblock Preprint,~69 pages, available at http://arxiv.org/abs/1007.0838 (2010).

\bibitem[HKS]{HKS} R.~Holley, R.~Kusuoka, D.~Stroock.
\newblock Asymptotic of the spectral gap with applications to the
theory of simulated annealing, 
\newblock J.~Funct.~Anal. \textbf{83}-2 (1989), pp 333--347.

\bibitem[{KPS}]{KoPrSh} N.~Koldan, I.~Prokhorenkov, and M.~Shubin.
\newblock Semiclassical Asymptotics on Manifolds with Boundary.
 \newblock Spectral analysis in geometry and number theory, pp.~239--266,
Contemp. Math., 484, Amer. Math. Soc. (2009).

\bibitem[{Lau1}]{Lau1} F.~Laudenbach.
 \newblock On the Thom-Smale complex.
 \newblock Ast\'erisque 205, pp.~219--233 (1992).

\bibitem[{Lau2}]{Lau2} F.~Laudenbach.
 \newblock A Morse complex on manifolds with boundary.
 \newblock Preprint,~11 pages, available at http://arxiv.org/abs/1003.5077 (2010).


\bibitem[{Lep1}]{Lep1} D.~Le~Peutrec.
 \newblock Small singular values of an extracted matrix of a Witten
  complex.
 \newblock Cubo, A Mathematical Journal, Vol.~11 (4), pp.~49-57 (2009).


\bibitem[{Lep2}]{Lep2} D.~Le~Peutrec.
 \newblock Local WKB construction for Witten
Laplacians on manifolds with boundary.
 \newblock Analysis \& PDE, Vol. 3, No. 3, pp. 227--260 (2010).

\bibitem[{Lep3}]{Lep3} D.~Le~Peutrec.
 \newblock Small eigenvalues of the Neumann realization of the semiclassical Witten Laplacian.
 \newblock Annales de la Facult\'e des Sciences de Toulouse,
 Vol.~19, no 3--4, pp.~735--809 (2010).
 
 \bibitem[{Lep4}]{Lep4} D.~Le~Peutrec.
 \newblock Small eigenvalues of the Witten Laplacian
acting on $p$-forms on a surface.
 \newblock To appear in Asymptotic Analysis.
Preprint,~16 pages, available at http://hal.archives-ouvertes.fr/hal-00473114/fr/
(2010).

\bibitem[{Mas}]{Mas} W.~S.~Massey.
\newblock \textit{A basic course in algebraic topology.}
\newblock Graduate Texts in Mathematics, 127. Springer-Verlag, (1991).

\bibitem[{Mil}]{Mil} J.~Milnor. 
\newblock \textit{Morse theory.}
\newblock Annals of Mathematics Studies, No. 51 Princeton University
Press, (1963). 

\bibitem[{Nel}]{Nel} E.~Nelson.
\newblock {\it Dynamical Theories of Brownian Motion}.
\newblock 2nd Ed., Princeton University Press, (2002).

 \bibitem[{Sch}]{Sch} G.~Schwarz.
 \newblock{\it Hodge decomposition. A method for Solving Boundary Value
 Problems}.
 \newblock Lecture Notes in Mathematics 1607, Springer-Verlag (1995).
 
\bibitem[{Sim}]{Sim} B.~Simon.
\newblock {\it Trace ideals and their applications}.
\newblock Cambridge University Press IX, Lecture Notes Series vol. 35 (1979).

\bibitem[{Spa}]{Spa} E.H.~Spanier.
\newblock{\it Algebraic Topology.}
\newblock McGraw-Hill Series in Higher Mathematics (1966).
\bibitem[{Tay}]{Tay} M.E.~Taylor.
\newblock{\it Partial Differential Equations~1, Basic Theory}.
\newblock Applied Mathematical Sciences 115, Springer (1997).

\bibitem[{TTK}]{TTK} J.~Tailleur, S.~Tanase-Nicola, J.~Kurchan.
\newblock Kramers equation and supersymmetry.
\newblock J. Stat. Phys. \textbf{122}-4 (2006) pp. 557-–595.

\bibitem[{Wit}]{Wit} E.~Witten.
\newblock Supersymmetry and Morse inequalities.
\newblock J. Diff. Geom. 17, p.~661-692 (1982).

\bibitem[{Zha}]{Zha} W.~Zhang.
\newblock \textit{Lectures on Chern-Weil theory and Witten deformations.}
\newblock Nankai Tracts in Mathematics, 4. World Scientific Publishing Co. (2001).
\end{thebibliography}
